\pdfoutput=1
\RequirePackage{ifpdf}
\ifpdf % We are running pdfTeX in pdf mode
\documentclass[pdftex]{sigma}
\else
\documentclass{sigma}
\fi

\newcommand{\aki}{{}}    
\newcommand{\mn}{\bar{n}}
\newcommand{\mone}{\bar{1}}
\newcommand{\mtwo}{\bar{2}}
\newcommand{\mthree}{\bar{3}}
\newcommand{\mfour}{\bar{4}}
\newcommand{\mfive}{\bar{5}}
\newcommand{\minusi}{\bar{i}}
\newcommand{\tone}{t_1}
\newcommand{\ttwo}{t_2}
\newcommand{\tN}{t_N}

\usepackage[enableskew]{youngtab}
\usepackage[all]{xy}

\numberwithin{equation}{section}

\newtheorem{Theorem}{Theorem}[section]
\newtheorem{Corollary}[Theorem]{Corollary}
\newtheorem{Lemma}[Theorem]{Lemma}
\newtheorem{Proposition}[Theorem]{Proposition}
\newtheorem{Conjecture}[Theorem]{Conjecture}
 { \theoremstyle{definition}
\newtheorem{Definition}[Theorem]{Definition}
\newtheorem{Example}[Theorem]{Example}
\newtheorem{Remark}[Theorem]{Remark} }

\begin{document}

\allowdisplaybreaks

\renewcommand{\thefootnote}{$\star$}

\renewcommand{\PaperNumber}{028}

\FirstPageHeading

\ShortArticleName{Rigged Conf\/igurations and Kashiwara Operators}

\ArticleName{Rigged Conf\/igurations and Kashiwara Operators\footnote{This paper is a~contribution to the Special Issue in
honor of Anatol Kirillov and Tetsuji Miwa.
The full collection is available at \href{http://www.emis.de/journals/SIGMA/InfiniteAnalysis2013.html}
{http://www.emis.de/journals/SIGMA/InfiniteAnalysis2013.html}}}

\Author{Reiho SAKAMOTO}

\AuthorNameForHeading{R.~Sakamoto}

\Address{Department of Physics, Tokyo University of Science, Kagurazaka, Shinjuku, Tokyo, Japan}
\Email{\href{mailto:reiho@rs.tus.ac.jp}{reiho@rs.tus.ac.jp}, \href{mailto:reihosan@08.alumni.u-tokyo.ac.jp}{reihosan@08.alumni.u-tokyo.ac.jp}}
\URLaddress{\url{https://sites.google.com/site/affinecrystal}}

\ArticleDates{Received April 02, 2013, in f\/inal form February 28, 2014; Published online March 23, 2014}

\Abstract{For types $A^{(1)}_n$ and $D^{(1)}_n$ we prove that the rigged conf\/iguration bijection intertwines the
classical Kashiwara operators on tensor products of the arbitrary Kirillov--Reshetikhin crystals and the set of the
rigged conf\/igurations.}


\Classification{17B37; 05E10; 82B23}

\begin{flushright}
\it Dedicated to Professor Anatol N.~Kirillov with admiration
\end{flushright}

\renewcommand{\thefootnote}{\arabic{footnote}}
\setcounter{footnote}{0}

{\small \tableofcontents}

\section{Introduction}

A search of a~natural presentation of the basis for representations of inf\/inite dimensional algebras sometimes reveals
an unexpected and intriguing connection with other models of mathematical physics.
Take the Feigin--Fuchs representation of the Virasoro algebra~\cite{FeFu} as an example.
In this case, a~cerebrated integrable quantum many body problem called the Calogero--Sutherland
model~\cite{Sutherland1,Sutherland2} plays the role.
Indeed, it is conjectured in~\cite{SSAFR} and partially proved in~\cite{CJ} that the excited states of the
Calogero--Sutherland model (usually called the Jack symmetric functions) behave particularly nicely as a~basis for the
Feigin--Fuchs module.
For example, matrix elements of the Virasoro generators with respect to the Jack basis become factorized rational
functions whose major part consisting of degree one polynomials of positive integer coef\/f\/icients which admit
a~transparent combinatorial interpretation.
This phenomenon seems not accidental since the theory includes Mimachi and Yamada's f\/inding~\cite{MY} which claims that
the singular vectors of the Virasoro algebra coincide with special cases of the Jack symmetric functions.
See, for example,~\cite{AFLT,BB,DLM,EPSS,FHV} for recent developments on the subject.

The purpose of the present paper is to pursue a~parallel investigation for the settings of the crystal
bases~\cite{Kashiwara:1991} of an important class of f\/inite dimensional representations (called the
Kirillov--Reshetikhin modules) of the quantum af\/f\/ine algebras.
We call the corresponding crystal the Kirillov--Reshetikhin (KR) crystal.
In this case, the corresponding physical model is the box-ball system.
The box-ball system is a~prototypical example of the ultradiscrete (or tropical) soliton system whose simplest example
was discovered by Takahashi and Satsuma~\cite{TS}.
For our purpose, an important aspect of the box-ball system is that the dynamics of the model is determined by the
combinatorial $R$-matrices for the KR crystal~\cite{FOY, HHIKTT}.
This formalism allows us to def\/ine the box-ball system for any tensor products of the KR crystals for arbitrary quantum
af\/f\/ine algebras (see, for example,~\cite{HKT1,HKT2}).
See~\cite{S:review} for an introductory review on studies on the box-ball systems.

Recently, a~complete set of the action and angle variables of the type $A^{(1)}_n$ box-ball systems is found to be
identical to the set of the rigged conf\/igurations~\cite{KOSTY}.
The rigged conf\/igurations are combinatorial objects discovered by an insightful analysis of the Bethe ansatz for quantum
integrable models~\cite{KKR,KR}.
Originally the rigged conf\/igurations are used to prove a~combinatorial identity for the Kostka--Foulkes polynomials
called the Fermionic formula as well as to give a~description of the branching coef\/f\/icients of f\/inite dimensional
representations of the quantum af\/f\/ine algebras with respect to the corresponding f\/inite dimensional subalgebras.
In such context the rigged conf\/igurations are labels of representations and thus naturally correspond to the classically
highest weight tensor products of the KR crystals.
On the other hand, since the box-ball system is def\/ined for not necessarily highest weight tensor products, we need to
consider generalized rigged conf\/igurations corresponding to any tensor products.

Therefore it is natural to introduce the Kashiwara operators on the set of the rigged con\-f\/i\-gu\-ra\-tions.
Let $I$ be the set of the Dynkin nodes of a~quantum af\/f\/ine algebra $\mathfrak{g}$ following Kac's convention~\cite{Kac}
and set $I_0=I
\setminus\{0\}
$.
In~\cite{Sch:2006} Schilling gave a~def\/inition of the classical Kashiwara operators $\widetilde{e}_i$ and
$\widetilde{f}_i$ ($i\in I_0$) for simply laced algebras and proved that the crystal structure on the rigged
conf\/igurations is actually isomorphic to the one on the tensor products of the KR crystals by using Stembridge's local
characterization of crystals.
In Section~\ref{sec:crystal_axiom_RC} of the present paper, we provide a~direct proof of the fact that the classical
Kashiwara operators on the set of the rigged conf\/igurations of types $A^{(1)}_n$ and $D^{(1)}_n$ indeed satisfy the
axiom of the crystals.

The crux of the rigged conf\/iguration theory is the so called rigged conf\/iguration bijection which gives an one to one
correspondence between the elements of tensor products of the KR crystals and the rigged conf\/igurations:
\begin{gather*}
\Phi:\text{ rigged conf\/igurations }\longmapsto\text{ tensor products.}
\end{gather*}
The algorithm of~$\Phi$ is a~rather complicated combinatorial procedure which is explained in
Section~\ref{sec:def_RCbijection}\footnote{For the most general tensor products $B^{r_1,s_1}\otimes
B^{r_2,s_2}\otimes\dots\otimes B^{r_L,s_L}$ of type $D^{(1)}_n$, a~proof of the well-def\/inedness of the bijection~$\Phi$
is the subject of~\cite{OSSS:2012}.
However we will not assume this fact in the present paper since the paper~\cite{OSSS:2012} is yet not available for the
public.
After Conjecture~\ref{conj:well_def}, we summarize the known cases of this result.}. Nevertheless it is known that the
map~$\Phi$ has a~lot of nice properties and seems to have a~deep mathematical origin.
For example, one of the most important properties of the rigged conf\/iguration bijection is that we can regard the
algorithm~$\Phi$ as a~general form of the combinatorial $R$-matrices (\cite{KSS:2002} for type $A^{(1)}_n$ case
and~\cite{S:2005,SS} for type $D^{(1)}_n$ case).
See Remark~\ref{rem:combR} for more precise meanings.
This property is the basis for providing the action and angle va\-riab\-les of the box-ball systems~\cite{KOSTY,KSY2}.
Moreover, as an application of the relationship with the box-ball system we have a~representation theoretical
interpretation of the algorithm of $\Phi^{-1}$ for type $A^{(1)}_n$~\cite{Sak2}, though the full understanding of the
rigged conf\/igurations themselves is yet to be seen.

The main result of the present paper is the proof of the following compatibility of the rigged conf\/iguration bijection
and the classical Kashiwara operators
\begin{gather}
\label{eq:intro}
[\widetilde{e}_i,\Phi]=0,
\qquad
[\widetilde{f}_i,\Phi]=0,
\qquad
i\in I_0
\end{gather}
for $\mathfrak{g}=A^{(1)}_n$ and $D^{(1)}_n$ (see Theorem~\ref{th:main}) under the assumption that the bijection~$\Phi$
is well-def\/ined.
Here some remarks are in order.
\begin{itemize}\itemsep=0pt
\item
In Appendix C of the preprint version of~\cite{DS:2006}, Deka and Schilling asserted the relation~\eqref{eq:intro} for
type $A^{(1)}_n$.
However, according to the author's opinion, they did not f\/inish their proof (see Section~\ref{sec:DekaSchilling}).
Given the fundamental importance of the problem, it seems appropriate to prove the result completely for type
$A^{(1)}_n$ also.
Indeed, the proof for type $A^{(1)}_n$ is a~part of the proof for type $D^{(1)}_n$ and easy to extract from
it\footnote{Ignore all arguments related with the parameter $\ell_{(a)}$ by formally setting $\ell_{(a)}=\infty$.
See Section~\ref{sec:def_RCbijection} for the notation.}.
\item
The proof of~\eqref{eq:intro} is direct and does not depend on other results.
\end{itemize}

Since our approach is direct, we expect that our result paves the way to deal with the rigged conf\/iguration bijection
for arbitrary quantum af\/f\/ine algebras.
Currently the following cases are established:
\begin{itemize}\itemsep=0pt
\item
All cases of type $A^{(1)}_n$~\cite{KKR,KR,KSS:2002}.
\item
Special cases of type $D^{(1)}_n$~\cite{OSS:2012, S:2005,SS} (see the comments after Conjecture~\ref{conj:well_def}).
\item
Tensor products of the vector representations of arbitrary non-exceptional quantum af\/f\/ine algebras~\cite{OSS:2003} as
well as of type $E^{(1)}_6$~\cite{OSano}.
\end{itemize}
In all cases the algorithms for the rigged conf\/iguration bijection share many common features.
We expect that the rigged conf\/iguration bijection as well as all major properties will extend to the arbitrary quantum
af\/f\/ine algebras.

We expect that such an extension of the rigged conf\/iguration theory will provide meaningful insights for both
mathematical physics side and mathematical side which are not easily seen by other methods.
One of the supporting evidences for our expectation is a~bijection constructed in~\cite{OS:2011}.
In that paper, we constructed a~combinatorial bijection between the rigged conf\/igurations for arbitrary non-exceptional
quantum af\/f\/ine algebras (under certain restrictions on ranks) and the set of pairs of a~type $A^{(1)}_n$ rigged
conf\/iguration and a~Littlewood--Richardson tableau.
Remarkably, the construction is uniform for all types of algebras and the Littlewood--Richardson tableaux naturally
appear as the recording tableaux of the algorithm.
We expect that the bijection coincides with (and generalizes) a~canonical Dynkin diagram automorphism (see Remark 3.2
of~\cite{OS:2011}).

\looseness=-1
Another supporting evidence is a~numerical comparison of the rigged conf\/iguration bijection~\cite{OSS:2003} for the
vector representations of arbitrary non-exceptional quantum af\/f\/ine algebras and the corresponding combinatorial
$R$-matrices (see Section~\ref{sec:RC_foudations} and subsequent sections of~\cite{KOSTY}).
The conjecture provided there claims that seemingly complicated algorithms presented in~\cite{OSS:2003} with many case
by case def\/inition possess uniform description in terms of the corresponding combinatorial $R$-matrices.
In fact, we expect that the rigged conf\/iguration bijection gives a~ge\-ne\-ral form of the combinatorial $R$-matrices for
general quantum af\/f\/ine algebras and thus provides action and angle variables for the corresponding box-ball systems.
Finally, it should be mentioned that in~\cite{OSS:2012} we see that the af\/f\/ine crystal structure of type $D^{(1)}_n$ is
essentially go\-ver\-ned by the rigged conf\/igurations.
Therefore we expect that the construction of the rigged conf\/iguration bijection for general cases should be an
intriguing and meaningful future research topic.

This paper is organized as follows.
In Section~\ref{sec:background}, we collect necessary facts from the crystal bases and the tableau representations.
In Section~\ref{sec:RC_foudations}, we give the def\/inition of the rigged conf\/igurations and the rigged conf\/iguration
bijection and explain their basic properties (convexity relation of the vacancy numbers).
We also introduce the classical Kashiwara operators on the rigged conf\/igurations and show that they satisfy the axiom of
the classical crystals.
In Section~\ref{sec:statement} we give the statement of the main result and clarify what are the essential points to be
shown.
Main proofs are given in Sections~\ref{sec:main},~\ref{sec:main2} and~\ref{sec:main2_e}.
During the proof we have tried to avoid unnecessarily long arguments.
Nevertheless we have to deal with many subtle cases and some of them look very special.
Thus some readers might wonder whether such case does exist and is worthy of careful analysis.
In order to clarify that point, we include several examples for seemingly subtle cases.

\section{Background on crystals and tableaux}\label{sec:background}

\subsection{Crystal bases}\label{se:crystal}

Let us brief\/ly recall some def\/initions about the crystal bases theory of Kashiwara~\cite{Kashiwara:1991}.
For more detailed introduction, see, for example,~\cite{Kashiwara:book}.
Let $\mathfrak{g}$ be an af\/f\/ine Kac--Moody algebra, $\mathfrak{g}'$ be its derived subalgebra, $\mathfrak{g}_0$ be the
corresponding f\/inite dimensional simple Lie algebra obtained by removing the 0 node of the Dynkin diagram of
$\mathfrak{g}$.
Here we follow Kac' convention~\cite{Kac} for the labeling of the Dynkin nodes and denote by $I$ the vertex set of the
Dynkin diagram of $\mathfrak{g}$ and $I_0:=I
\setminus \{0\}$.
Let $U_q(\mathfrak{g})$, $U_q'(\mathfrak{g})$ and $U_q(\mathfrak{g}_0)$ be the quantized universal enveloping algebra
corresponding to $\mathfrak{g}$, $\mathfrak{g}'$ and $\mathfrak{g}_0$, respectively.
In this paper we will mainly focus on the case $\mathfrak{g}=D^{(1)}_n$ which has the following Dynkin diagram
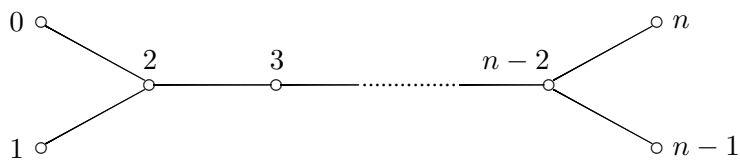
\begin{center}
\unitlength 12pt
\begin{picture}(23,4)
\put(0,-0.3){1}
\put(0,3.7){0}
\multiput(1,0)(0,4){2}{\circle{0.3}}
\qbezier(1.15,0.15)(1.15,0.15)(4.25,1.85)
\qbezier(1.15,3.85)(1.15,3.85)(4.25,2.15)
\multiput(4.4,2)(4,0){2}{\circle{0.3}}
\put(4.2,2.5){2}
\qbezier(4.55,2)(4.55,2)(8.25,2)
\put(8.2,2.5){3}
\qbezier(8.55,2)(8.55,2)(11,2)
\multiput(11.2,2)(0.2,0){15}{\circle*{0.1}}
\qbezier(14.2,2)(14.2,2)(16.85,2)
\put(17,2){\circle{0.3}}
\put(14.9,2.5){$n-2$}
\multiput(20.4,0)(0,4){2}{\circle{0.3}}
\qbezier(17.15,1.85)(17.15,1.85)(20.25,0.15)
\qbezier(17.15,2.15)(17.15,2.15)(20.25,3.85)
\put(20.9,-0.2){$n-1$}
\put(20.9,3.8){$n$}
\end{picture}
\end{center}

Let $\alpha_i$, $h_i$, $\Lambda_i$ ($i\in I$) be the simple roots, simple coroots and fundamental weights of
$\mathfrak{g}$, $\bar{\Lambda}_i$ ($i\in I_0$) be the fundamental weights of $\mathfrak{g}_0$.
In the present case $\mathfrak{g}_0=D_n$, the simple roots are
\begin{gather}\label{simple_root}
\alpha_i=\epsilon_i-\epsilon_{i+1},
\quad
1\leq i\leq n-1,
\qquad
\alpha_n=\epsilon_{n-1}+\epsilon_n,
\end{gather}
and the fundamental weights are
\begin{gather}
\bar{\Lambda}_i=\epsilon_1+\dots+\epsilon_i,
\qquad
1\leq i\leq n-2,
\label{fundamental_weight}
\\
\bar{\Lambda}_{n-1}=(\epsilon_1+\dots+\epsilon_{n-1}-\epsilon_n)/2,
\qquad
\bar{\Lambda}_{n}=(\epsilon_1+\dots+\epsilon_{n-1}+\epsilon_n)/2,
\nonumber
\end{gather}
where $\epsilon_i\in\mathbb{Z}^n$ is the $i$th standard unit vector.
Let $Q$, $Q^\vee$, $P$ be the root, coroot and weight lattices of $\mathfrak{g}$.
Let $\langle\cdot,\cdot\rangle:Q^\vee\otimes P\rightarrow\mathbb{Z}$ be the pairing such that $\langle
h_i,\Lambda_j\rangle=\Lambda_j(h_i)=\delta_{ij}$.
Note that we have $\langle h_i,\alpha_j\rangle=\alpha_j(h_i)=A_{ij}$ where $A_{ij}$ is the Cartan matrix of
$\mathfrak{g}=D^{(1)}_n$.

Now we give the axiomatic def\/inition of $U_q'(\mathfrak{g})$-crystals.
\begin{Definition}\label{def:crystal}
$U_q'(\mathfrak{g})$-{\it crystal} is a~nonempty set $B$ equipped with maps $\mathrm{wt}:B\rightarrow P$,
$\varepsilon_i,\varphi_i:B\rightarrow\mathbb{Z}\cup\{-\infty\}$ and the {\it Kashiwara operators}
$\widetilde{e}_i,\widetilde{f}_i:B\rightarrow B\cup\{0\}$ for all $i\in I$ such that
\begin{enumerate}\itemsep=0pt
\item[(1)]
$\langle h_i,\mathrm{wt}(b)\rangle =\varphi_i(b)-\varepsilon_i(b)$,
\item[(2)]
$\mathrm{wt}\bigl(\widetilde{e}_i(b)\bigr)=\mathrm{wt}(b)+\alpha_i$ if $\widetilde{e}_i(b)\in B$,
\item[(3)]
$\mathrm{wt}\bigl(\widetilde{f}_i(b)\bigr)=\mathrm{wt}(b)-\alpha_i$ if $\widetilde{f}_i(b)\in B$,
\item[(4)]
$\varepsilon_i\bigl(\widetilde{e}_i(b)\bigr)=\varepsilon_i(b)-1$,
$\varphi_i\bigl(\widetilde{e}_i(b)\bigr)=\varphi_i(b)+1$ if $\widetilde{e}_i(b)\in B$,
\item[(5)]
$\varepsilon_i\bigl(\widetilde{f}_i(b)\bigr)=\varepsilon_i(b)+1$,
$\varphi_i\bigl(\widetilde{f}_i(b)\bigr)=\varphi_i(b)-1$ if $\widetilde{f}_i(b)\in B$,
\item[(6)]
for $b,b'\in B$, $\widetilde{f}_i(b)=b'$ if and only if $\widetilde{e}_i(b')=b$,
\item[(7)]
if $\varphi_i(b)=-\infty$ for $b\in B$, then $\widetilde{e}_i(b)=\widetilde{f}_i(b)=0$.
\end{enumerate}
\end{Definition}
In this paper, we only consider the case when the maps $\varepsilon_i,\varphi_i$ are def\/ined by
\begin{gather*}
\varepsilon_i(b)= \max\big\{m\geq 0\,\big|\,\widetilde{e}^m_i(b)\neq0\big\},
\qquad
\varphi_i(b)= \max\big\{m\geq 0\,\big|\,\widetilde{f}^m_i(b)\neq 0\big\}
\end{gather*}
for $b\in B$.
In particular, we have $0\leq \varepsilon_i(b),\varphi_i(b)<\infty$ for all $i\in I$ and $b\in B$.
Note that the crystals with these conditions are called the regular crystals.
If we have $\widetilde{f}_i(b)=b'$ for $b,b'\in B$, we write an arrow with color (or label) $i$ from~$b$ to $b'$.
In this way, the crystal $B$ can be regarded as the colored oriented graph whose vertices are the elements of $B$.
We call such graph the {\it crystal graph}.

One of the nice properties of the crystal bases is that it behaves nicely with respect to the tensor product.
Let $B_2\otimes B_1$ be the tensor product of two crystals $B_1$ and $B_2$.
As the set, it coincides with the Cartesian product $B_2\times B_1$.
The action of the Kashiwara operators $\widetilde{e}_i$, $\widetilde{f}_i$ on an element $b_2\otimes b_1\in B_2\otimes B_1$
is given explicitly as follows:
\begin{gather*}
\widetilde{e}_i(b_2\otimes b_1)=
\begin{cases}
\widetilde{e}_i(b_2)\otimes b_1\quad\text{if}\ \ \varepsilon_i(b_2)>\varphi_i(b_1),
\\
b_2\otimes \widetilde{e}_i(b_1)\quad\text{if}\ \ \varepsilon_i(b_2)\leq\varphi_i(b_1),
\end{cases}
\\
\widetilde{f}_i(b_2\otimes b_1)=
\begin{cases}
\widetilde{f}_i(b_2)\otimes b_1\quad\text{if}\ \ \varepsilon_i(b_2)\geq\varphi_i(b_1),
\\
b_2\otimes \widetilde{f}_i(b_1)\quad\text{if}\ \ \varepsilon_i(b_2)<\varphi_i(b_1),
\end{cases}
\end{gather*}
where the result is declared to be $0$ if either of its tensor factors are $0$.
We also have $\mathrm{wt}(b_2\otimes b_1)=\mathrm{wt}(b_2)+\mathrm{wt}(b_1)$.

\begin{Remark}
We use the opposite of the Kashiwara's tensor product convention.
\end{Remark}

In order to deal with multiple tensor products, it is convenient to use the signature rule.
Let $B$ be the tensor product of crystals $B=B_L\otimes\cdots\otimes B_2\otimes B_1$.
For $i\in I$ and $b=b_L\otimes\dots\otimes b_2\otimes b_1\in B$, the {\it $i$-signature} of~$b$ is the following
sequence of the symbols $+$ and $-$
\begin{gather*}
\underbrace{+\cdots\cdots\cdots +}_{\varphi_i(b_L)\text{ times}} \underbrace{-\cdots\cdots\cdots -}_{\varepsilon_i(b_L)\text{times}}
\cdots\cdots \underbrace{+\cdots\cdots\cdots +}_{\varphi_i(b_1)\text{ times}} \underbrace{-\cdots\cdots\cdots -}_{\varepsilon_i(b_1)\text{ times}}.
\end{gather*}
The {\it reduced $i$-signature} of~$b$ is obtained by removing subsequence $-+$ of the $i$-signature of~$b$ repeatedly
until it becomes the following form
\begin{gather*}
\underbrace{+\cdots\cdots\cdots +}_{\varphi_i(b)\text{ times}} \underbrace{-\cdots\cdots\cdots -}_{\varepsilon_i(b)\text{times}}.
\end{gather*}
By using this information the action of $\widetilde{e}_i$ and $\widetilde{f}_i$ on~$b$ can be described as follows.
If $\varepsilon_i(b)=0$, then $\widetilde{e}_i(b)=0$.
Otherwise
\begin{gather*}
\widetilde{e}_i(b_L\otimes\dots\otimes b_1)= b_L\otimes\dots\otimes b_{j+1}\otimes \widetilde{e}_i(b_j)\otimes
b_{j-1}\otimes \dots\otimes b_1,
\end{gather*}
where the leftmost $-$ in the reduced $i$-signature of~$b$ comes from $b_j$.
Similarly, if $\varphi_i(b)=0$, then $\widetilde{f}_i(b)=0$.
Otherwise
\begin{gather*}
\widetilde{f}_i(b_L\otimes\dots\otimes b_1)= b_L\otimes\dots\otimes b_{j+1}\otimes \widetilde{f}_i(b_j)\otimes
b_{j-1}\otimes \dots\otimes b_1,
\end{gather*}
where the rightmost $+$ in the reduced $i$-signature of~$b$ comes from $b_j$.
Finally the actions of $\widetilde{e}_i$ and $\widetilde{f}_i$ on each tensor factor $b_j$ are described by tableau
representations of $b_j$ described in the sequel.

\subsection{Kirillov--Reshetikhin tableaux}

In order to perform explicit manipulations on the crystal bases, it is convenient to use an explicit realization of the
elements of crystals.
For the present purpose, it is more convenient to use a~new kind of tableau representation which we call the
Kirillov--Reshetikhin (KR) tableaux~\cite{OSS:2012} than the usual Kashiwara--Nakashima (KN) tableaux~\cite{KN:1994}.

\subsubsection{Reviews on tableau representations}\label{subsec:tableaux}

Let us brief\/ly recall tableau representations of crystals.
Let $B(\Lambda)$ be the crystal associated to the highest weight representation of highest weight $\Lambda$ of
$U_q(\mathfrak{g}_0)$.
Then in the present case $\mathfrak{g}=D^{(1)}_n$, the $U_q'(\mathfrak{g})$ crystal $B^{r,s}$ has the following
decomposition under the restriction to its f\/inite dimensional subalgebra $U_q(\mathfrak{g}_0)$.
If $r<n-1$, we have
\begin{gather}
\label{eq:classical_decomposition_1}
B^{r,s}
\simeq\bigoplus_\lambda B(\lambda)
\qquad
\text{as }
\mathfrak{g}_0\text{ crystals},
\end{gather}
and if $r=n-1,n$, we have
\begin{gather}
\label{eq:classical_decomposition_2}
B^{r,s}\simeq B(s\bar{\Lambda}_r)
\qquad
\text{as }\mathfrak{g}_0\text{ crystals}.
\end{gather}
Here in the f\/irst case $r<n-1$, $\lambda$ in the summation are determined as follows.
For this we identify $\lambda=
\sum\limits_ik_i\bar{\Lambda}
_i$ with the Young diagram which has $k_i$ columns of height $i$ for all $i$.
Then one of $\lambda$ in decomposition~\eqref{eq:classical_decomposition_1} is shape $(s^r)$ rectangular diagram and the
remaining $\lambda$ are obtained by removing the vertical dominoes $\Yboxdim5pt\yng(1,1)$ from $(s^r)$ rectangular
diagram in all possible ways.
Then if $b\in B^{r,s}$ belongs to some $B(\lambda)$ in the right hand side of the
decomposition~\eqref{eq:classical_decomposition_1} or~\eqref{eq:classical_decomposition_2}, then the KN tableau
representation of~$b$ has shape $\lambda$.

{\bf KN tableaux for $\boldsymbol{B^{r,s}}$ when $\boldsymbol{r\leq n-2}$.} For $J
\subset I$, we say~$b$ is {\it $J$-highest element}
if we have $\widetilde{e}_i(b)=0$ for all $i\in J$.
Suppose that $u_\lambda\in B^{r,s}$ ($r\leq n-2$) is the $I_0$-highest element of $B(\lambda)$ which appears in the
decomposition~\eqref{eq:classical_decomposition_1}.
Then the corresponding KN tableau is obtained by f\/illing the letter $j$ to all the boxes at $j$th row of the Young
diagram $\lambda$.
In order to obtain the other tableaux, we def\/ine the action of the Kashiwara operators on tableaux as follows.
The crystal graph of $B^{1,1}$ for $\mathfrak{g}_0=D_n$ is as follows
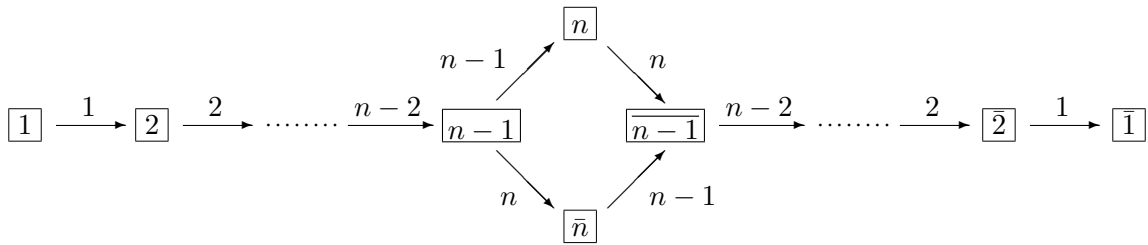
\begin{center}
\unitlength 12pt
\begin{picture}(36,7)(0,-0.2)
\put(0,3){$\Yvcentermath1\young(1)$}
\put(1.5,3.2){\vector(1,0){2.2}}
\put(2.3,3.5){1}
\put(4,3){$\Yvcentermath1\young(2)$}
\put(5.5,3.2){\vector(1,0){2.2}}
\put(6.3,3.5){2}
\multiput(8.2,3.2)(0.3,0){8}{\circle*{0.12}}
\put(10.7,3.2){\vector(1,0){2.7}}
\put(10.9,3.5){$n-2$}
\put(13.7,2.7){\frame{$\,n-1\,$\rule{0pt}{11pt}}}
\put(15.4,4.0){\vector(1,1){1.8}}
\put(13.6,5){$n-1$}
\put(15.4,2.4){\vector(1,-1){1.8}}
\put(15.5,0.7){$n$}
\put(17.5,6.2){$\Yvcentermath1\young(n)$}
\put(18.9,0.6){\vector(1,1){1.8}}
\put(18.9,5.7){\vector(1,-1){1.8}}
\put(20.2,5){$n$}
\put(20.2,0.7){$n-1$}
\put(19.5,2.7){\frame{$\,\overline{n-1}\,$\rule{0pt}{11pt}}}
\put(17.5,-0.2){$\Yvcentermath1\young(\mn)$}
\put(22.4,3.2){\vector(1,0){2.7}}
\put(22.6,3.5){$n-2$}
\multiput(25.6,3.2)(0.3,0){8}{\circle*{0.12}}
\put(28.1,3.2){\vector(1,0){2.2}}
\put(28.9,3.5){2}
\put(30.7,3){$\Yvcentermath1\young(\mtwo)$}
\put(32.2,3.2){\vector(1,0){2.2}}
\put(33.0,3.5){1}
\put(34.8,3){$\Yvcentermath1\young(\mone)$}
\end{picture}
\end{center}
Weight of each vertex is $\mathrm{wt}\Bigl(\Yvcentermath1\young(i)\Bigr)=\epsilon_i$ and
$\mathrm{wt}\Bigl(\Yvcentermath1\young(\minusi)\Bigr)=-\epsilon_i$, respectively.
For a~given tableau $t\in B^{r,s}$ ($r\leq n-2$),
\begin{center}
\unitlength 12pt
\begin{picture}(13,10)
\put(0,5){$t=$}
\put(2,3){
\put(7,6.7){\line(1,0){4}}
\put(7,5.7){\line(1,0){4}}
\put(7,4.7){\line(1,0){4}}
\put(7,2.7){\line(1,0){4}}
\put(7,1.7){\line(1,0){4}}
\put(7,0.7){\line(1,0){2.7}}
\put(7,-0.3){\line(1,0){2.7}}
\put(11,6.7){\line(0,-1){5}}
\put(9.7,6.7){\line(0,-1){7}}
\put(7.8,6.7){\line(0,-1){7}}
\put(10,6){$t_1$}
\put(10,5){$t_2$}
\multiput(10.3,4.4)(0,-0.2){8}{\circle*{0.1}}
\put(10,2){$t_r$}
\put(8,6){$t_{r+1}$}
\put(8,5){$t_{r+2}$}
\multiput(8.3,4.4)(0,-0.2){8}{\circle*{0.1}}
\put(8,2){$t_{2r}$}
\multiput(8.3,1.5)(0,-0.2){4}{\circle*{0.1}}
\put(8,0){$t_{r'}$}
\multiput(2.4,6.7)(0.2,0){23}{\circle*{0.1}}
\put(0,6.7){\line(1,0){2.3}}
\put(0,6.7){\line(0,-1){9}}
\put(0,-2.3){\line(1,0){1.3}}
\put(0,-1.3){\line(1,0){1.3}}
\put(1.3,-2.3){\line(0,1){9}}
\multiput(0.6,-1.0)(0,0.2){10}{\circle*{0.1}}
\put(0.2,-2){$t_N$}
\multiput(1.4,-2.3)(0.2,0.07){28}{\circle*{0.1}}
}
\end{picture}
\end{center}
we introduce the so-called Japanese reading word and regard $t$ as the element of $(B^{1,1})^{\otimes N}$ as follows
\begin{gather*}
\Yvcentermath1
t\longmapsto\young(\tN)\otimes\cdots\otimes
\young(\ttwo)\otimes\young(\tone)
\in\big(B^{1,1}\big)^{\otimes N}.
\end{gather*}
Via this identif\/ication we can introduce the Kashiwara operators on tableaux and obtain the KN tableau representation
for the general case.

{\bf KN tableaux for $\boldsymbol{B^{n,2s+1}}$.} Let us consider the basic case $B^{n,1}\simeq B(\bar{\Lambda}_n)$.
This is the so-called spin representation of dimension $2^{n-1}$ whose basis are $(s_1,s_2,\ldots,s_n)$ where $s_i=\pm$
and $s_1s_2\cdots s_n=1$.
The classical Kashiwara operators act on each basis as follows:
\begin{gather}\label{e_on_spin}
\widetilde{e}_i(s_1,\ldots,s_n)=
\begin{cases}
(s_1,\ldots,+,-,\ldots,s_n), & \text{if}\ \ i\neq n,\ (s_i,s_{i+1})=(-,+),
\\
(s_1,\ldots,s_{n-2},+,+), & \text{if}\ \ i=n,\ (s_{n-1},s_n)=(-,-),
\\
0, &  \text{otherwise},
\end{cases}
\\
\widetilde{f}_i(s_1,\ldots,s_n)=
\begin{cases}
(s_1,\ldots,-,+,\ldots,s_n), & \text{if}\ \ i\neq n,\ (s_i,s_{i+1})=(+,-),
\\
(s_1,\ldots,s_{n-2},-,-), & \text{if}\ \ i=n,\ (s_{n-1},s_n)=(+,+),
\\
0, & \text{otherwise}.
\end{cases}
\label{f_on_spin}
\end{gather}
The weight of each basis is
\begin{gather}
\label{wt_on_spin}
\mathrm{wt}\bigl((s_1,\ldots,s_n)\bigr)= \frac{1}{2}(s_1\epsilon_1+\dots+s_n\epsilon_n).
\end{gather}
Thus we can represent each element $\epsilon_i$ by half width tableau with entry $i$.
We arrange the entries within a~column according to the following order on the letters (smaller ones at higher places)
\begin{gather*}
1\prec 2\prec\dots\prec n-1\prec\begin{matrix}n \\ \overline{n}\end{matrix}
\prec \overline{n-1}\prec\dots\prec\overline{2}\prec\overline{1}.
\end{gather*}
Here we do not introduce the order between $n$ and $\overline{n}$.

For the general case, we start by the highest weight element $u_\lambda\in B^{n,2s+1}$ where
$\lambda=(2s+1)\bar{\Lambda}_n$.
Then represent $u_\lambda$ by the spin column $(+,\ldots,+)$ and $s$ columns of $1,2,\ldots,n$ such as
\begin{center}
\unitlength 12pt
\begin{picture}(13,6)(0.5,0)
\put(0.5,2.8){$u_{(2s+1)\bar{\Lambda}_{n}}=$}
\put(6.2,0){\line(0,1){6}}
\put(7,0){\line(0,1){6}}
\multiput(8.5,0)(1.5,0){4}{\line(0,1){6}}
\multiput(6.2,0)(0,1.5){5}{\line(1,0){6.8}}
\put(6.35,4.9){1}
\put(6.35,3.4){2}
\multiput(6.6,2.6)(0,-0.37){3}{\circle*{0.1}}
\put(6.35,0.5){$n$}
\put(1.15,0){
\put(6.35,4.9){1}
\put(6.35,3.4){2}
\multiput(6.6,2.6)(0,-0.37){3}{\circle*{0.1}}
\put(6.35,0.5){$n$}
}
\put(2.7,0){
\put(6.35,4.9){1}
\put(6.35,3.4){2}
\multiput(6.6,2.6)(0,-0.37){3}{\circle*{0.1}}
\put(6.35,0.5){$n$}
}
\put(5.6,0){
\put(6.35,4.9){1}
\put(6.35,3.4){2}
\multiput(6.6,2.6)(0,-0.37){3}{\circle*{0.1}}
\put(6.35,0.5){$n$}
}
\multiput(10.35,0.7)(0.37,0){3}{\circle*{0.1}}
\multiput(10.35,2.6)(0.37,-0.37){3}{\circle*{0.1}}
\multiput(10.35,3.7)(0.37,0){3}{\circle*{0.1}}
\multiput(10.35,5.2)(0.37,0){3}{\circle*{0.1}}
\end{picture}
\end{center}
In order to obtain the other tableaux, we embed $u_\lambda$ into $B^{r,1}\otimes (B^{1,1})^{\otimes rs}$ by using the
Japanese reading word for the right $s$ columns and apply the Kashiwara operators $\widetilde{f}_i$ in all possible
ways.

We can obtain the KN tableaux for the case $B^{n,2s}$ similarly if we ignore the leftmost spin column in $B^{n,2s+1}$
case.

{\bf KN tableaux for $\boldsymbol{B^{n-1,2s+1}}$.} In the case $B^{n-1,1}
\simeq B(\bar{\Lambda}
_{n-1})$, we can represent the basis as $(s_1,s_2,\ldots,s_n)$ where $s_i=\pm$ and $s_1s_2\cdots s_n=-1$.
We can use the same formulae~\eqref{e_on_spin} and~\eqref{f_on_spin} for the Kashiwara operators and~\eqref{wt_on_spin}
for the weight.
For the general case $B^{n-1,2s+1}$, we represent the highest weight element $u_{(2s+1)\bar{\Lambda}_{n-1}}\in
B^{n-1,2s+1}$~by
\begin{center}
\unitlength 12pt
\begin{picture}(13,6)(0.0,0)
\put(0.0,2.8){$u_{(2s+1)\bar{\Lambda}_{n-1}}=$}
\put(6.2,0){\line(0,1){6}}
\put(7,0){\line(0,1){6}}
\multiput(8.5,0)(1.5,0){4}{\line(0,1){6}}
\multiput(6.2,0)(0,1.5){5}{\line(1,0){6.8}}
\put(6.35,4.9){1}
\put(6.35,3.4){2}
\multiput(6.6,2.6)(0,-0.37){3}{\circle*{0.1}}
\put(6.35,0.4){$\mn$}
\put(1.15,0){
\put(6.35,4.9){1}
\put(6.35,3.4){2}
\multiput(6.6,2.6)(0,-0.37){3}{\circle*{0.1}}
\put(6.35,0.4){$\mn$}
}
\put(2.7,0){
\put(6.35,4.9){1}
\put(6.35,3.4){2}
\multiput(6.6,2.6)(0,-0.37){3}{\circle*{0.1}}
\put(6.35,0.4){$\mn$}
}
\put(5.6,0){
\put(6.35,4.9){1}
\put(6.35,3.4){2}
\multiput(6.6,2.6)(0,-0.37){3}{\circle*{0.1}}
\put(6.35,0.4){$\mn$}
}
\multiput(10.35,0.7)(0.37,0){3}{\circle*{0.1}}
\multiput(10.35,2.6)(0.37,-0.37){3}{\circle*{0.1}}
\multiput(10.35,3.7)(0.37,0){3}{\circle*{0.1}}
\multiput(10.35,5.2)(0.37,0){3}{\circle*{0.1}}
\end{picture}
\end{center}
Here we obtain the tableau for $u_{(2s+1)\bar{\Lambda}_{n-1}}$ by replacing all $n$ of the bottom row of
$u_{(2s+1)\bar{\Lambda}_{n}}$ by~$\mn$.
Then we can use the Kashiwara operators to obtain all the tableaux of $B^{n-1,2s+1}$.
We can obtain tableaux for the case $B^{n-1,2s}$ if we ignore the leftmost spin column.

\subsubsection{Filling map}

The purpose to introduce the KR tableaux is to give a~rectangular presentation of elements of crystals as opposed to the
KN tableaux whose shapes depend on the classical decompositions given in~\eqref{eq:classical_decomposition_1}
and~\eqref{eq:classical_decomposition_2}.
The def\/inition of the KR tableaux is given explicitly for the highest weight element of $B^{r,s}$ and extended to the
general case by using the Kashiwara operators $\widetilde{f}_i$ as in the KN tableaux case.

To begin with we give a~def\/inition of the KR tableaux for the case $B^{r,s}$ ($r\leq n-2$).
Let $u_\lambda$ be the $I_0$-highest element of $B^{r,s}$ with weight $\lambda=
k_r\bar{\Lambda}_r+k_{r-2}\bar{\Lambda}_{r-2}+\cdots$.
We denote by $\mathrm{f\/ill}(u_\lambda)$ the KR tableau representation of $u_\lambda$.
Now we describe the algorithm to obtain $\mathrm{f\/ill}(u_\lambda)$ which we call the {\it filling map}.
Let $k_c$ be the f\/irst odd integer in the sequence $k_{r-2},k_{r-4},\dots$.
If there is no such $k_c$, set $k_c=k_{-1}$, that is, $c=-1$.
Let $t$ be the KN tableau representation of $u_\lambda$.
We place $t$ at the top left corner of the $r$ by $s$ rectangle.
In order to f\/ill all the empty places, we start from the leftmost column of $t$ and move rightwards by the following
procedure.
\begin{enumerate}\itemsep=0pt
\item
We do nothing for the height $r$ columns.
If $c\geq 0$, remove one column from the height $c$ columns and move the columns that are shorter than $c$ to left by 1.
\item
For the columns whose heights are equal to or larger than $c$, put the following stuf\/f to the empty places as much as
possible.
If the height of the corresponding columns is $h$, put the transpose of
\begin{gather*}
\begin{array}{|c|c|c|c|}
\hline
\overline{r}&\overline{r-1}&\dots&\overline{h+1}\rule{0pt}{12pt}
\\
\hline
h+1&h+2&\dots&r
\\
\hline
\end{array}
\end{gather*}
If $c=-1$, we stop here since all empty places are f\/illed by this procedure.
\item
For the remaining columns except for the rightmost one, put the following stuf\/f recursively from left to right with
recursively redef\/ining an integer $x$.
As the initial condition, set $x=c+1$.

If the corresponding column has height $h$, put the transpose of
\begin{gather*}
\begin{array}{|c|c|c|c|c|c|c|c|}
\hline
y&\dots&r-1&r&\overline{r}&\dots&\overline{x+1}&\overline{x}\rule{0pt}{12pt}
\\
\hline
\end{array}
\end{gather*}
to the empty place.
Since the number of the empty places is $r-h$, we have $y=r-(x-h-2)$.
We redef\/ine $x=y$ and do the same procedure to the next column.
\item
By using the f\/inal output of $x$ from the last step, we put the transpose of the following to the rightmost column
\begin{gather*}
\begin{array}{|c|c|c|c|c|c|c|c|}
\hline
1&2&\dots&y&\overline{y}&\dots&\overline{x+1}&\overline{x}\rule{0pt}{12pt}
\\
\hline
\end{array}
\end{gather*}
Recall that the rightmost column is empty by Step 1 whenever $c\geq 0$.
Then in order to achieve the above pattern, we have $y=(r+x-1)/2$.
\end{enumerate}
\begin{Example}
Let us consider the $I_0$-highest element $u_\lambda$ of $B^{8,7}$ where
$\lambda=\bar{\Lambda}_8+2\bar{\Lambda}_6+\bar{\Lambda}_4+2\bar{\Lambda}_2$.
In this case, we have $(k_6,k_4,k_2,k_0)=(2,1,2,1)$, thus $c=4$.
We put the KN tableau representation of $u_\lambda$ into the 8 by 7 rectangle.
Then
\begin{gather*}
\begin{array}{|c|c|c|c|c|c|c|}
\hline
1&1&1&1&1&1&\rule{5pt}{0pt}\rule{0pt}{11pt}
\\
\hline
2&2&2&2&2&2&\rule{0pt}{11pt}
\\
\hline
3&3&3&3& & &\rule{0pt}{11pt}
\\
\hline
4&4&4&4& & &\rule{0pt}{11pt}
\\
\hline
5&5&5& & & &\rule{0pt}{11pt}
\\
\hline
6&6&6& & & &\rule{0pt}{11pt}
\\
\hline
7& & & & & &\rule{0pt}{11pt}
\\
\hline
8& & & & & &\rule{0pt}{11pt}
\\
\hline
\end{array}
\, \xrightarrow{\text{ Step 1 }}\,
\begin{array}{|c|c|c|c|c|c|c|}
\hline
1&1&1&1&1&\rule{5pt}{0pt}&\rule{5pt}{0pt}\rule{0pt}{11pt}
\\
\hline
2&2&2&2&2& &\rule{0pt}{11pt}
\\
\hline
3&3&3& & & &\rule{0pt}{11pt}
\\
\hline
4&4&4& & & &\rule{0pt}{11pt}
\\
\hline
5&5&5& & & &\rule{0pt}{11pt}
\\
\hline
6&6&6& & & &\rule{0pt}{11pt}
\\
\hline
7& & & & & &\rule{0pt}{11pt}
\\
\hline
8& & & & & &\rule{0pt}{11pt}
\\
\hline
\end{array}
\end{gather*}
Then we f\/ill the empty places as
\begin{gather*}
\begin{array}{|c|c|c|c|c|c|c|}
\hline
\rule{5pt}{0pt}& & & & &5&1\rule{0pt}{11pt}
\\
\hline
& & & & &6&2\rule{0pt}{11pt}
\\
\hline
& & &7&5&7&3\rule{0pt}{11pt}
\\
\hline
& & &8&6&8&4\rule{0pt}{11pt}
\\
\hline
& & &\bar{8}&7&\bar{8}&5\rule{0pt}{11pt}
\\
\hline
& & &\bar{7}&8&\bar{7}&6\rule{0pt}{11pt}
\\
\hline
&\bar{8}&7&\bar{6}&\bar{8}&\bar{6}&\bar{6}\rule{0pt}{11pt}
\\
\hline
&\bar{7}&8&\mfive&\bar{7}&\mfive&\mfive\rule{0pt}{11pt}
\\
\hline
\end{array}
\qquad \text{thus,} \quad \mathrm{f\/ill}(u_\lambda)=
\begin{array}{|c|c|c|c|c|c|c|}
\hline
1&1&1&1&1&5&1\rule{0pt}{11pt}
\\
\hline
2&2&2&2&2&6&2\rule{0pt}{11pt}
\\
\hline
3&3&3&7&5&7&3\rule{0pt}{11pt}
\\
\hline
4&4&4&8&6&8&4\rule{0pt}{11pt}
\\
\hline
5&5&5&\bar{8}&7&\bar{8}&5\rule{0pt}{11pt}
\\
\hline
6&6&6&\bar{7}&8&\bar{7}&6\rule{0pt}{11pt}
\\
\hline
7&\bar{8}&7&\bar{6}&\bar{8}&\bar{6}&\bar{6}\rule{0pt}{11pt}
\\
\hline
8&\bar{7}&8&\mfive&\bar{7}&\mfive&\mfive\rule{0pt}{11pt}
\\
\hline
\end{array}
\end{gather*}
\end{Example}

In order to def\/ine the KR tableaux for the spin cases, we start from the following expression for the highest weight
elements (collection of $s$ spin columns)
\begin{center}
\unitlength 12pt
\begin{picture}(20,6)
\put(-0.2,2.8){$u_{s\bar{\Lambda}_{n-1}}=$}
\multiput(4.2,0)(0.8,0){3}{\line(0,1){6}}
\multiput(4.2,0)(0,1.5){5}{\line(1,0){3.9}}
\put(4.35,4.9){1}
\put(4.35,3.4){2}
\multiput(4.6,2.6)(0,-0.37){3}{\circle*{0.1}}
\put(4.35,0.45){$\mn$}
\put(0.8,0){
\put(4.35,4.9){1}
\put(4.35,3.4){2}
\multiput(4.6,2.6)(0,-0.37){3}{\circle*{0.1}}
\put(4.35,0.45){$\mn$}
}
\put(3.1,0){
\put(4.35,4.9){1}
\put(4.35,3.4){2}
\multiput(4.6,2.6)(0,-0.37){3}{\circle*{0.1}}
\put(4.35,0.45){$\mn$}
}
\multiput(7.3,0)(0.8,0){2}{\line(0,1){6}}
\multiput(6.15,0.7)(0.37,0){3}{\circle*{0.1}}
\multiput(6.15,2.6)(0.37,-0.37){3}{\circle*{0.1}}
\multiput(6.15,3.7)(0.37,0){3}{\circle*{0.1}}
\multiput(6.15,5.2)(0.37,0){3}{\circle*{0.1}}
\put(8.4,2.8){,}
%%%%%%%%%%%%%%%%%%%%%%%%%%%%%%%%%%%
\put(10,0){
\put(0.5,2.8){$u_{s\bar{\Lambda}_{n}}=$}
\multiput(4.2,0)(0.8,0){3}{\line(0,1){6}}
\multiput(4.2,0)(0,1.5){5}{\line(1,0){3.9}}
\put(4.35,4.9){1}
\put(4.35,3.4){2}
\multiput(4.6,2.6)(0,-0.37){3}{\circle*{0.1}}
\put(4.35,0.5){$n$}
\put(0.8,0){
\put(4.35,4.9){1}
\put(4.35,3.4){2}
\multiput(4.6,2.6)(0,-0.37){3}{\circle*{0.1}}
\put(4.35,0.5){$n$}
}
\put(3.1,0){
\put(4.35,4.9){1}
\put(4.35,3.4){2}
\multiput(4.6,2.6)(0,-0.37){3}{\circle*{0.1}}
\put(4.35,0.5){$n$}
}
\multiput(7.3,0)(0.8,0){2}{\line(0,1){6}}
\multiput(6.15,0.7)(0.37,0){3}{\circle*{0.1}}
\multiput(6.15,2.6)(0.37,-0.37){3}{\circle*{0.1}}
\multiput(6.15,3.7)(0.37,0){3}{\circle*{0.1}}
\multiput(6.15,5.2)(0.37,0){3}{\circle*{0.1}}
}
\end{picture}
\end{center}
and apply the Kashiwara operators $\widetilde{f}_i$ ($i\in I_0$) for all possible ways to obtain the rest of the
elements of $B^{n-1,s}$ and $B^{n,s}$.

From Remark 5.2 and Remark 5.4 in~\cite{OSS:2012}, it follows that the new Kirillov--Reshetikhin (KR) tableaux carry
a~crystal structure obtained by their reading word which gives a~natural embedding of $B(\lambda)$ into
$B(\Lambda_1)^{\otimes rs}$.
Since each classical factor in $B^{r,s}$ occurs with multiplicity 1, we can use the Kirillov--Reshetikhin tableaux
representation in place of the usual Kashiwara--Nakashima tableaux representation.

For the later purposes, we introduce several operations on the KR tableaux representation of crystals.
\begin{Definition}
\label{def:operations_on_tableaux}
We use the KR tableaux representation for the element $b=b_L\otimes b_{L-1}\otimes\dots\otimes b_1\in B^{r_L,s_L}\otimes
B^{r_{L-1},s_{L-1}}\otimes\dots\otimes B^{r_1,s_1}$.
Then we def\/ine the following operations.
\begin{enumerate}\itemsep=0pt
\item[(1)] Suppose that $B^{r_L,s_L}=B^{1,1}$.
Then we def\/ine the operation called {\it left-hat} by
\begin{gather*}
\!\!\!\!\!\! \mathrm{lh}(b)=b_{L-1}\otimes\dots\otimes b_1\in B^{r_{L-1},s_{L-1}}\otimes\dots\otimes B^{r_1,s_1}.
\end{gather*}
\item[(2)] Suppose that $B^{r_L,s_L}=B^{r_L,1}$ and $r_L>1$.
Thus $b_L$ has the form $b_L=
\begin{array}{|c|}
\hline
t_1
\\
\hline
\vdots
\\
\hline
t_{r_L-1}
\\
\hline
t_{r_L}
\\
\hline
\end{array}
$.
Then we def\/ine the operation called {\it left-box} by
\begin{gather*}
\!\!\!\!\!\!\mathrm{lb}(b)=
\begin{array}{|c|}
\hline
t_{r_L}
\\
\hline
\end{array}
\otimes
\begin{array}{|c|}
\hline
t_1
\\
\hline
\vdots
\\
\hline
t_{r_L-1}
\\
\hline
\end{array}
\otimes b_{L-1}\otimes\cdots\otimes b_1\in B^{1,1}\otimes B^{r_L-1,1}\otimes B^{r_{L-1},s_{L-1}}\otimes\cdots\otimes B^{r_1,s_1}.
\end{gather*}
\item[(3)] Suppose that $s_L>1$.
Let $c_i$ be the $i$th column of $b_L$ such that $b_L=c_1c_2\cdots c_{s_L}$.
Then we def\/ine the operation called {\it left-split} by
\begin{gather*}
\!\!\!\!\!\!\mathrm{ls}(b)= c_1 \otimes c_2\cdots c_{s_L} \otimes b_{L-1}\otimes\dots\otimes b_1\in B^{r_L,1}\otimes
B^{r_L,s_L-1}\otimes B^{r_{L-1},s_{L-1}}\otimes\dots\otimes B^{r_1,s_1}.
\end{gather*}
\end{enumerate}
\end{Definition}

\section{Rigged conf\/igurations and the main bijection}\label{sec:RC_foudations}

\subsection{Definition of the rigged conf\/igurations}

Our strategy to def\/ine the rigged conf\/igurations of type $\mathfrak{g}=D^{(1)}_n$ has three steps.
First we def\/ine the highest weight rigged conf\/igurations.
Next we def\/ine the classical Kashiwara operators $\widetilde{e}_i$ and $\widetilde{f}_i$ for $i\in I_0$ on the rigged
conf\/igurations.
Finally, we def\/ine the general rigged conf\/igurations by all possible applications of the operators $\widetilde{f}_i$ on
highest weight rigged conf\/igurations.

Let us prepare several notation to describe the rigged con\-f\/i\-gu\-ra\-tions.
The rigged conf\/igura\-tions are combinatorial objects made up with partitions $\nu^{(a)}$ and integers $J^{(a)}$ for $a\in
I_0$.
$\nu^{(a)}=\big(\nu^{(a)}_1,\ldots,\nu^{(a)}_k\big)$ is a~partition called the {\it configuration} and
$J^{(a)}=\big(J^{(a)}_1,\ldots,J^{(a)}_k\big)$ is a~sequence of integers called the {\it rigging} which are associated with rows
of the corresponding conf\/iguration.
We should regard that the pairs $\big(\nu^{(a)},J^{(a)}\big)$ are located on the vertices $I_0$ of the Dynkin diagram of
$D^{(1)}_n$.
We denote by $m^{(a)}_l(\nu)$ the number of length $l$ rows of $\nu^{(a)}$.
We will abbreviate the pair $\big(\nu^{(a)},J^{(a)}\big)$ by $(\nu,J)^{(a)}$ and call the components of
$(\nu,J)^{(a)}=\big\{\big(\nu^{(a)}_1,J^{(a)}_1\big),\ldots,\big(\nu^{(a)}_k,J^{(a)}_k\big)\big\}$ the strings.

The rigged conf\/igurations also depend on the shape of the tensor product $B=B^{r_L,s_L}\otimes\dots\otimes B^{r_1,s_1}$.
Let $L^{(a)}_l$ be the number of the component $B^{a,l}$ within $B$.
Then we def\/ine the partition~$\mu^{(a)}$ such that the number of length $l$ rows of $\mu^{(a)}$ is equal to $L^{(a)}_l$.
Thus the rigged conf\/iguration~$(\nu,J)$ is described by the sequence of conf\/igurations and the riggings
\begin{gather*}
(\nu,J)= \big( (\nu,J)^{(1)},\ldots,(\nu,J)^{(n-1)},(\nu,J)^{(n)} \big)
\end{gather*}
together with additional data $\mu=\big(\mu^{(1)},\ldots,\mu^{(n-1)},\mu^{(n)}\big)$ under certain constraints to be described
below.
We denote by $Q^{(a)}_l(\nu)$ (resp.
$Q^{(a)}_l(\mu)$) the number of boxes in the f\/irst $l$ columns of $\nu^{(a)}$ (resp.
$\mu^{(a)}$);
\begin{gather*}
Q^{(a)}_l(\nu)=
\sum\limits_{j>0}
\min\{l,j\}m^{(a)}_j(\nu),
\qquad
Q^{(a)}_l(\mu)=
\sum\limits_{j>0}
\min\{l,j\}L^{(a)}_j.
\end{gather*}
Therefore we have $Q^{(a)}_\infty(\nu)=|\nu^{(a)}|$ and $Q^{(a)}_\infty(\mu)=|\mu^{(a)}|$ where $|\nu^{(a)}|$ and
$|\mu^{(a)}|$ are the total number of boxes in $\nu^{(a)}$ and $\mu^{(a)}$, respectively.
From the data $\mu$ and $\nu$ we def\/ine the {\it vacancy number} $P^{(a)}_l(\mu,\nu)$ (hereafter we will abbreviate as
$P^{(a)}_l(\nu)$) by the formula
\begin{gather}
P^{(a)}_l(\nu)=Q^{(a)}_l(\mu)+
\sum\limits_{b\in I_0}
A_{ab}Q^{(b)}_l(\nu)\nonumber\\
\hphantom{P^{(a)}_l(\nu)}{}
=Q^{(a)}_l(\mu)-2Q^{(a)}_l(\nu)+
\sum\limits_{b\in I_0 \atop b\sim a}
Q^{(b)}_l(\nu),\label{eq:def_vacancy}
\end{gather}
where $A_{ab}$ is the Cartan matrix of $\mathfrak{g}_0=D_n$ and $a
\sim b$ means that the vertices~$a$ and~$b$ are connected by a~single edge on the Dynkin diagram.

\begin{Definition}
The data $\mu$ and $(\nu,J)$ is the highest weight rigged conf\/iguration if all the strings $\big(\nu^{(a)}_i,J^{(a)}_i\big)$ and
the corresponding vacancy numbers satisfy the following condition
\begin{gather*}
0\leq J^{(a)}_i\leq P^{(a)}_{\nu^{(a)}_i}(\nu).
\end{gather*}
\end{Definition}

\begin{Remark}
If $\nu^{(n-1)}=\nu^{(n)}=\varnothing$, then we can identify the rigged conf\/iguration as type $A^{(1)}_{n-2}$.
In this way we can recover proofs for type $A^{(1)}_n$ from those for type $D^{(1)}_n$.
\end{Remark}

\begin{Example}
The following object is a~highest weight rigged conf\/iguration corresponding to the tensor product $B^{3,2}\otimes
B^{3,1}\otimes B^{2,2}\otimes B^{1,2}\otimes B^{1,1}$ of type $D^{(1)}_5$:
\begin{center}
\unitlength 12pt
\begin{picture}(25,4)
\put(0,2){
\put(-0.7,0.1){0}
\put(-0.7,1.1){0}
\put(0,0){\Yboxdim12pt\yng(2,2)}
\put(2.2,0.1){0}
\put(2.2,1.1){0}
}
\put(5,1){
\put(-0.7,0.1){1}
\put(-0.7,1.1){2}
\put(-0.7,2.1){1}
\put(0,0){\Yboxdim12pt\yng(3,2,1)}
\put(1.2,0.1){1}
\put(2.2,1.1){1}
\put(3.2,2.1){0}
}
\put(11,0){
\put(-0.7,0.1){0}
\put(-0.7,1.1){0}
\put(-0.7,2.1){1}
\put(-0.7,3.1){0}
\put(0,0){\Yboxdim12pt\yng(3,2,1,1)}
\put(1.2,0.1){0}
\put(1.2,1.1){0}
\put(2.2,2.1){0}
\put(3.2,3.1){0}
}
\put(17,2){
\multiput(-0.7,0.1)(0,1){2}{0}
\put(0,0){\Yboxdim12pt\yng(2,1)}
\put(1.2,0.1){0}
\put(2.2,1.1){0}
}
\put(22,3){
\put(-0.7,0.1){2}
\put(0,0){\Yboxdim12pt\yng(2)}
\put(2.2,0.1){1}
}
\end{picture}
%\{\{2,2\},\{3,2,1\},\{3,2,1,1\},\{2,1\},\{2\}\}\\
%\{\{0,0\},\{0,1,1\},\{0,0,0,0\},\{0,0\},\{1\}\}
\end{center}
Here we put the vacancy number (resp.\
rigging) on the left (resp.\
right) of the corresponding row of the conf\/iguration represented by a~Young diagram.
By the rigged conf\/iguration bijection~$\Phi$ (see Section~\ref{sec:def_RCbijection}) the above rigged conf\/iguration
$(\nu,J)$ corresponds to the following highest weight tensor product:
\begin{gather*}
\Phi(\nu,J)= \Yvcentermath1 \young(13,2\mfive,45)\otimes\young(1,4,5)\otimes\young(21,3\mone)
\otimes\young(12)\otimes\young(1)
\end{gather*}
\end{Example}

The weight $\lambda$ of the rigged conf\/iguration is def\/ined by the relation (sometimes called the
$(L,\lambda)$-conf\/iguration condition)
\begin{gather*}%\label{wt_RC_1}
\sum\limits_{a\in I_0,i>0}
im^{(a)}_i\alpha_a=
\sum\limits_{a\in I_0,i>0}
iL^{(a)}_i\bar{\Lambda}_a-\lambda.
\end{gather*}
If we expand $\lambda$ by the basis $\epsilon_i$, we can rewrite as
\begin{gather}
\label{wt_RC_2}
\lambda=
\sum\limits_i\lambda_i\epsilon_i= \sum\limits_{a\in I_0}
\big(\big|\mu^{(a)}\big|\bar{\Lambda}_a-\big|\nu^{(a)}\big|\alpha_a\big).
\end{gather}
Then we can use the expressions~\eqref{simple_root} and~\eqref{fundamental_weight} to obtain the explicit expressions
for the weight~$\lambda_i$.
We write the weight of the rigged conf\/iguration by $\mathrm{wt}(\nu,J)$.

Following~\cite{Sch:2006} we introduce the classical Kashiwara operators on the rigged conf\/igurations and use them to
def\/ine the general rigged conf\/igurations.
For the string $(l,x)$ of $(\nu,J)^{(i)}$, we call the quantity $P^{(i)}_l(\nu)-x$ the {\it corigging}.

\begin{Definition}
The rigged conf\/igurations are obtained by all possible applications of the Kashiwara operators $\widetilde{f}_i$ ($i\in
I_0$) on the highest weight rigged conf\/igurations.
Here the Kashiwara operators on the rigged conf\/igurations are def\/ined as follows.
Let $x_\ell$ be the smallest rigging of~$(\nu,J)^{(i)}$.
\begin{enumerate}\itemsep=0pt
\item[(1)] Let $\ell$ be the minimal length of the strings of $(\nu,J)^{(i)}$ with the rigging $x_\ell$.
If $x_\ell\geq 0$, def\/ine $\widetilde{e}_i(\nu,J)=0$.
Otherwise $\widetilde{e}_i(\nu,J)$ is obtained by replacing the string $(\ell,x_\ell)$ by $(\ell-1,x_\ell+1)$ while
changing all other riggings to keep coriggings f\/ixed.
\item[(2)] Let $\ell$ be the maximal length of the strings of $(\nu,J)^{(i)}$ with the rigging $x_\ell$.
$\widetilde{f}_i(\nu,J)$ is obtained by the following procedure.
If $x_\ell>0$, add a~string $(1,-1)$ to $(\nu,J)^{(i)}$.
Otherwise replace the string $(\ell,x_\ell)$ by $(\ell+1,x_\ell-1)$.
Change other riggings to keep coriggings f\/ixed.
If the new rigging is strictly larger than the corresponding new vacancy number, def\/ine $\widetilde{f}_i(\nu,J)=0$.
\end{enumerate}
\end{Definition}

By def\/inition every rigging is always less than or equal to the corresponding vacancy number.

\begin{Proposition}
\label{prop:check_crystal}
The above defined Kashiwara operators satisfy the following relations
\begin{gather*}
\mathrm{wt}\big(\widetilde{e}_i(\nu,J)\big)=\mathrm{wt}(\nu,J)+\alpha_i,
\qquad
\mathrm{wt}\big(\widetilde{f}_i(\nu,J)\big)=\mathrm{wt}(\nu,J)-\alpha_i.
\end{gather*}
Moreover $\widetilde{f}_i(\nu,J)\neq 0\Longrightarrow\widetilde{e}_i\widetilde{f}_i(\nu,J)=(\nu,J)$ and
$\widetilde{e}_i(\nu,J)\neq 0\Longrightarrow\widetilde{f}_i\widetilde{e}_i(\nu,J)=(\nu,J)$.
\end{Proposition}
\begin{proof}
The relations for wt follow from~\eqref{wt_RC_2}.
In order to show $\widetilde{e}_i\widetilde{f}_i(\nu,J)=(\nu,J)$, suppose that $\widetilde{f}_i$ creates the string
$(\ell+1,x_\ell-1)$.
Let $(\widetilde{\nu},\widetilde{J})=\widetilde{f}_i(\nu,J)$.
Suppose that $0<\ell$.
\begin{enumerate}\itemsep=0pt
\item[(a)] Let $(k,x_k)$ be a~string of $(\nu,J)^{(i)}$ which satisf\/ies $k\leq\ell$.
Since $\widetilde{f}_i$ acts on the string with smallest rigging, we have $x_k\geq x_\ell$.
Recall that we have $P^{(i)}_k(\widetilde{\nu})=P^{(i)}_k(\nu)$ since $\widetilde{f}_i$ adds a~box to the $(\ell+1)$-th
column.
Thus the string $(k,x_k)$ remains as it is after $\widetilde{f}_i$ and its rigging satisf\/ies $x_k>x_\ell-1$.
\item[(b)] Let $(k,x_k)$ be a~string of $(\nu,J)^{(i)}$ which satisf\/ies $\ell<k$.
Since $\widetilde{f}_i$ acts on the longest string with rigging $x_\ell$, we have $x_\ell<x_k$.
Recall that we have $P^{(i)}_k(\widetilde{\nu})=P^{(i)}_k(\nu)-2$ in this case.
Thus the string $(k,x_k)$ becomes $(k,x_k-2)$ after $\widetilde{f}_i$.
Then its rigging satisf\/ies $x_k-2\geq x_\ell-1$.
\end{enumerate}
Thus the string $(\ell+1,x_\ell-1)$ is the shortest string among the strings of the smallest rigging of
$(\widetilde{\nu},\widetilde{J})^{(i)}$, hence $\widetilde{e}_i$ will act on it.
Therefore we obtain $\widetilde{e}_i\widetilde{f}_i(\nu,J)=(\nu,J)$.
Let us consider the case $\ell=0$.
Then we can use the same arguments as before if we put $(\ell,x_\ell)=(0,0)$.

The opposite relation $\widetilde{f}_i\widetilde{e}_i(\nu,J)=(\nu,J)$ can be shown similarly.
\end{proof}

Finally, let us def\/ine the maps $\varepsilon_i,\varphi_i:(\nu,J)\rightarrow\mathbb{Z}$ by
\begin{gather*}
\varepsilon_i(\nu,J)= \max\big\{m\geq 0\,\big|\,\widetilde{e}^m_i(\nu,J)\neq 0\big\},
\qquad
\varphi_i(\nu,J)= \max\big\{m\geq 0\,\big|\,\widetilde{f}^m_i(\nu,J)\neq 0\big\}.
\end{gather*}

\subsection{Convexity relations of the vacancy numbers and their applications}\label{sec:crystal_axiom_RC}

In this subsection, we introduce a~fundamental property of the vacancy numbers called the convexity relation.
As an application, we provide a~direct proof of the fact that the crystal structure on the set of the rigged
conf\/igurations indeed satisf\/ies the axiom of the $U_q(\mathfrak{g}_0)$-crystals.
Finally, we prepare ref\/ined estimates for the vacancy numbers which are necessary in the later arguments.

One of the basic properties of the vacancy number is the following convexity relation.

\begin{Proposition}[convexity]
\label{lem:convexity1}
Suppose that $m^{(a)}_l(\nu)=0$ for some $l\geq 1$.
Then the vacancy numbers satisfy the following convexity relation
\begin{gather*}
2P^{(a)}_l(\nu)\geq P^{(a)}_{l-1}(\nu)+P^{(a)}_{l+1}(\nu).
\end{gather*}
\end{Proposition}
\begin{proof}
Recall that the functions $Q^{(a)}_k(\mu)$ and $Q^{(b)}_k(\nu)$ for $b\neq a$ are all upper convex functions of $k$.
On the other hand, since $m^{(a)}_l(\nu)=0$, the function $Q^{(a)}_k(\nu)$ is a~linear function between $l-1\leq k\leq
l+1$.
Thus the combination in~\eqref{eq:def_vacancy} gives the upper convex relation for the vacancy numbers.
\end{proof}

By repeated use of the above convexity relation, we obtain several useful criteria.
\begin{Corollary}
Suppose that $m^{(a)}_k(\nu)=0$ for all $l_1<k<l_2$.
\begin{enumerate}\itemsep=0pt
\item[$(1)$] $P^{(a)}_k(\nu)\geq\min\big\{P^{(a)}_{l_1}(\nu),P^{(a)}_{l_2}(\nu)\big\}$.
\item[$(2)$] For some $k$ satisfying $l_1<k<l_2$ suppose that we have $P^{(a)}_{l_1}(\nu)\geq P^{(a)}_k(\nu)
<P^{(a)}_{l_2}(\nu)$ or $P^{(a)}_{l_1}(\nu)>P^{(a)}_k(\nu)\leq P^{(a)}_{l_2}(\nu)$.
Then the corresponding rigged configuration is forbidden.
\item[$(3)$] For some $k$ satisfying $l_1<k<l_2$ suppose that we have $P^{(a)}_{l_1}(\nu)\geq P^{(a)}_k(\nu)\leq
P^{(a)}_{l_2}(\nu)$.
Then the only admissible situation is $P^{(a)}_{l_1}(\nu)=P^{(a)}_{l_1+1}(\nu)=\dots=P^{(a)}_{l_2}(\nu)$.
\end{enumerate}
\end{Corollary}
Since we will use these properties so many times in the rest of the paper, we will sometimes use these relations without
giving an explicit reference.

Let us give the f\/irst application of the convexity relations.

\begin{Theorem}
\label{prop:phi_RC}
Let $x_\ell$ be the smallest rigging of $(\nu,J)^{(i)}$ and define $s=\min\{0,x_\ell\}$.
If $\nu^{(i)}=\varnothing$, set $s=0$.
Then we have
\begin{enumerate}\itemsep=0pt
\item[$(1)$] $\varepsilon_i(\nu,J)=-s$, \item[$(2)$] $\varphi_i(\nu,J)=P^{(i)}_\infty(\nu)-s$.
\end{enumerate}
\end{Theorem}
\begin{proof}
(1) We proceed by induction on $-s$.
Let $(\widetilde{\nu},\widetilde{J}):=\widetilde{e}_i(\nu,J)$, $\widetilde{x}_{\widetilde{\ell}}$ be the smallest
rigging of $(\widetilde{\nu},\widetilde{J})^{(i)}$ and def\/ine $\widetilde{s}=\min\{0,\widetilde{x}_{\widetilde{\ell}}\}$.

Suppose that $-s=0$.
This implies that $x_\ell\geq 0$.
Then by the def\/inition of $\widetilde{e}_i$ we have $\widetilde{e}_i(\nu,J)=0$ as requested.
Suppose that $-s>0$.
Then we have $s=x_\ell<0$.
Suppose that $\widetilde{e}_i$ acts on the string $(\ell,x_\ell)$.
Let us analyze the behaviors of the riggings by $\widetilde{e}_i$.
\begin{enumerate}\itemsep=0pt
\item[(a)] The string $(\ell,x_\ell)$ becomes $(\ell-1,x_\ell+1)$.
\item[(b)] Let $(k,x_k)$ be an arbitrary string of $(\nu,J)^{(i)}$ satisfying $k<\ell$.
Since $\widetilde{e}_i$ acts on the shortest string with rigging $x_\ell$, we have $x_k>x_\ell$.
From $k<\ell$, we see that $\widetilde{e}_i$ will not change corigging of $(k,x_k)$, thus the string $(k,x_k)$ remains
as it is in $(\widetilde{\nu},\widetilde{J})$.
\item[(c)] Let $(k,x_k)$ be an arbitrary string of $(\nu,J)^{(i)}$ satisfying $\ell\leq k$.
Since $x_\ell$ is the minimal rigging of $(\nu,J)^{(i)}$, we have $x_\ell\leq x_k$.
Recall that we have $P^{(i)}_k(\widetilde{\nu})=P^{(i)}_k(\nu)+2$ since $\widetilde{e}_i$ removes a~box from $\ell$th
column of $\nu^{(i)}$.
Thus $\widetilde{e}_i$ makes the string $(k,x_k)$ into $(k,x_k+2)$, in particular, its rigging satisf\/ies
$x_\ell+1<x_k+2$.
\end{enumerate}
To summarize, we have $\widetilde{x}_{\widetilde{\ell}}=x_\ell+1$ and thus $\widetilde{s}=s+1$, that is,
$-\widetilde{s}=\varepsilon_i(\nu,J)-1=\varepsilon_i(\widetilde{\nu},\widetilde{J})$ as requested.

(2) This is proved in~\cite[Lemma 3.6]{Sch:2006}.
However, for the sake of the completeness, we include a~proof here.
We prove by induction on $P^{(i)}_\infty(\nu)-s$.
Suppose that $\widetilde{f}_i$ acting on $(\nu,J)$ creates a~length $\ell+1$ string.
By def\/inition of $s$, $\widetilde{f}_i$ creates the string $(\ell+1,s-1)$ of
$(\widetilde{\nu},\widetilde{J}):=\widetilde{f}_i(\nu,J)$.
Let $\widetilde{x}_{\widetilde{\ell}}$ be the smallest rigging of $(\widetilde{\nu},\widetilde{J})^{(i)}$ and def\/ine
$\widetilde{s}=\min\{0,\widetilde{x}_{\widetilde{\ell}}\}$.

{\bf Step 1.} Let us consider the case $P^{(i)}_\infty(\nu)-s=0$.
Suppose that $x_\ell\leq 0$.
Then we have $\ell>0$ and $s=x_\ell$, thus, $P^{(i)}_\infty(\nu)=x_\ell$.
Let $j$ be the maximal integer such that $\ell\leq j$ and $m^{(i)}_j(\nu)>0$.
Suppose if possible that $\ell<j$.
Let us consider the string $(j,x_j)$.
Since $\widetilde{f}_i$ acts on the longest string with rigging $x_\ell$, the assumption $\ell<j$ implies that
$x_\ell<x_j$.
Note that we have $x_j\leq P^{(i)}_j(\nu)$ since $(\nu,J)$ is the rigged conf\/iguration.
Thus $x_\ell<P^{(i)}_j(\nu)$.
Recall that the def\/inition~\eqref{eq:def_vacancy} implies that there exists a~large integer $L$ such that
$P^{(i)}_L(\nu)=P^{(i)}_{L+1}(\nu)=\dots=P^{(i)}_\infty(\nu)=x_\ell$.
Since $j$ is the maximal integer such that $m^{(i)}_j(\nu)>0$, we see that $P^{(i)}_k(\nu)$ satisfy the convexity
relation between $j\leq k\leq\infty$.
In particular, the relation $P^{(i)}_j(\nu)>P^{(i)}_L(\nu)=P^{(i)}_{L+1}(\nu)$ for $j<L<L+1$ is a~contradiction.
Therefore we have $\ell=j$ and thus there is no string of $(\nu,J)^{(i)}$ that is longer than $\ell$.
Then from the convexity relation of $P^{(i)}_k(\nu)$ between $\ell\leq k\leq\infty$, the only possibility that is
compatible with the requirement $P^{(i)}_\ell(\nu)\geq x_\ell$ is the relation
$P^{(i)}_{\ell}(\nu)=P^{(i)}_{\ell+1}(\nu)=\dots =x_\ell$.

Let us consider the string $(\ell+1,x_\ell-1)$ of $(\widetilde{\nu},\widetilde{J})^{(i)}$ created by $\widetilde{f}_i$.
Due to the extra box added by $\widetilde{f}_i$, we have
$P^{(i)}_{\ell+1}(\widetilde{\nu})=P^{(i)}_{\ell+1}(\nu)-2=x_\ell-2$.
Thus the string $(\ell+1,x_\ell-1)$ of $(\widetilde{\nu},\widetilde{J})^{(i)}$ has the rigging that is strictly larger
than the corresponding vacancy number.
Thus we have $\widetilde{f}_i(\nu,J)=0$ as requested.

Finally let us consider the case $x_\ell>0$.
Then we have $P^{(i)}_\infty(\nu)=0$ by $s=0$.
Let $j$ be the maximal integer such that $m^{(i)}_j(\nu)>0$ and consider the corresponding string $(j,x_j)$.
Since $x_\ell$ is the minimal rigging, we have $P^{(i)}_j(\nu)\geq x_j\geq x_\ell>0$.
Then for a~suf\/f\/iciently large integer $L$, we have $P^{(i)}_j(\nu)>P^{(i)}_L(\nu)=P^{(i)}_{L+1}(\nu)=\dots=0$.
This is a~contradiction since $P^{(i)}_k(\nu)$ must satisfy the convexity relation between $j\leq k\leq\infty$.
Therefore we see that there is no string in $(\nu,J)^{(i)}$.
This is a~contradiction since we assume that $x_\ell>0$.
Hence this case cannot happen.

{\bf Step 2.} Let us consider the case $P^{(i)}_\infty(\nu)-s>0$.
Suppose that $x_\ell\leq 0$.
Then we have $P^{(i)}_\infty(\nu)>s=x_\ell$.
Let $j$ be the minimal integer such that $\ell<j$ and $m^{(i)}_j(\nu)>0$.
If there is no such $j$, set $j=\infty$.
Then we have $x_\ell<P^{(i)}_j(\nu)$.
For, if $j<\infty$ we have $x_\ell<x_j\leq P^{(i)}_j(\nu)$ as in the previous step and if $j=\infty$ we have
$x_\ell<P^{(i)}_j(\nu)$ by the assumption.

Let us show that
\begin{gather}
\label{eq:phiRC_check}
x_\ell<P^{(i)}_{\ell+1}(\nu).
\end{gather}
If $j=\ell+1$ this relation is already conf\/irmed.
Thus suppose that $\ell+1<j$.
Suppose if possible that we have $x_\ell\geq P^{(i)}_{\ell+1}(\nu)$.
Then we have $P^{(i)}_{\ell}(\nu)\geq P^{(i)}_{\ell+1}(\nu)<P^{(i)}_{j}(\nu)$ by $x_\ell\leq P^{(i)}_{\ell}(\nu)$ and
$x_\ell<P^{(i)}_j(\nu)$.
This is a~contradiction since $P^{(i)}_{k}(\nu)$ must satisfy the convexity relation between $\ell\leq k\leq j$.
In conclusion, we have $x_\ell<P^{(i)}_{\ell+1}(\nu)$.

Since $\widetilde{f}_i$ adds a~box to the $(\ell+1)$-th column of $\nu^{(i)}$, we have
$P^{(i)}_{\ell+1}(\widetilde{\nu})=P^{(i)}_{\ell+1}(\nu)-2\geq x_\ell-1$ by~\eqref{eq:phiRC_check}.
Therefore the string $(\ell+1,x_\ell-1)$ of $(\widetilde{\nu},\widetilde{J})^{(i)}$ created by $\widetilde{f}_i$ has the
rigging which is smaller than or equal to the corresponding vacancy number.
Hence we have $\widetilde{f}_i(\nu,J)\neq 0$.

Let us determine $\widetilde{s}$.
\begin{enumerate}\itemsep=0pt
\item[(a)] $\widetilde{f}_i$ makes the string $(\ell+1,x_\ell-1)$.
Below we consider the remaining strings of $(\nu,J)^{(i)}$.
\item[(b)] Let $(k,x_k)$ be an arbitrary string of $(\nu,J)^{(i)}$ satisfying $\ell<k$.
Since $\widetilde{f}_i$ acts on the longest string with rigging $x_\ell$, $\ell<k$ implies that $x_\ell<x_k$.
Recall that we have $P^{(i)}_k(\widetilde{\nu})=P^{(i)}_k(\nu)-2$ since $\widetilde{f}_i$ adds a~box to the $(\ell+1)$th
column of $\nu^{(i)}$.
In order to keep the corigging, the string $(k,x_k)$ becomes $(k,x_k-2)$ after $\widetilde{f}_i$.
In particular, the new rigging satisf\/ies $x_\ell-1\leq x_k-2$.
\item[(c)] Let $(k,x_k)$ be an arbitrary string of $(\nu,J)^{(i)}$ satisfying $k\leq\ell$.
Then we have $x_k\geq x_\ell$.
Since~$\widetilde{f}_i$ does not change the corigging of the string, the string $(k,x_k)$ remains as it is after~$\widetilde{f}_i$.
\end{enumerate}
To summarize, we have $\widetilde{s}=x_\ell-1=s-1$.
Since we have $P^{(i)}_\infty(\widetilde{\nu})=P^{(i)}_\infty(\nu)-2$, we have
$P^{(i)}_\infty(\widetilde{\nu})-\widetilde{s}=P^{(i)}_\infty(\nu)-s-1
=\varphi_i(\nu,J)-1=\varphi_i(\widetilde{\nu},\widetilde{J})$ as requested.

Finally consider the case $x_\ell>0$.
In this case we have $s=0$.
Then $\widetilde{f}_i$ adds the string $(1,-1)$ to $(\nu,J)^{(i)}$.
If we show that $\widetilde{f}_i(\nu,J)\neq 0$, then we have $\widetilde{s}=-1=s-1$ as requested.
This can be done if we formally setting $(\ell,x_\ell)=(0,0)$ in the proof of~\eqref{eq:phiRC_check} to deduce that
$P^{(i)}_1(\widetilde{\nu})=P^{(i)}_1(\nu)-2\geq -1$.
\end{proof}

\begin{Corollary}
\label{cor:check_crystal}
Let $\lambda=\mathrm{wt}(\nu,J)$.
Then we have
\begin{gather*}
\langle h_i,\lambda\rangle=\varphi_i(\nu,J)-\varepsilon_i(\nu,J).
\end{gather*}
\end{Corollary}
\begin{proof}
From Theorem~\ref{prop:phi_RC}, we have
$\varphi_i(\nu,J)-\varepsilon_i(\nu,J)=\big(P^{(i)}_\infty(\nu)-s\big)-(-s)=P^{(i)}_\infty(\nu)$.
On the other hand, from~\eqref{wt_RC_2} we have (see also def\/initions in Section~\ref{se:crystal})
\begin{gather*}
\langle h_i,\lambda\rangle=
\sum\limits_{a\in I_0}
\left(\big|\mu^{(a)}\big|\langle h_i,\bar{\Lambda}_a\rangle -\big|\nu^{(a)}\big|\langle h_i,\alpha_a\rangle\right)
\\
\phantom{\langle h_i,\lambda\rangle}
=\big|\mu^{(i)}\big|-2\big|\nu^{(i)}\big|+
\sum\limits_{a\in I_0,\, a\sim i}
\big|\nu^{(a)}\big|
=P^{(i)}_\infty(\nu).
\end{gather*}
Hence we have $\langle h_i,\lambda\rangle=\varphi_i(\nu,J)-\varepsilon_i(\nu,J)$.
\end{proof}

\begin{Example}
Let us consider the following element of $(B^{1,1})^{\otimes 8}$
\begin{gather*}
b= \Yvcentermath1 \young(1)\otimes\young(2)\otimes\young(3)\otimes\young(1)\otimes
\young(1)\otimes\young(2)\otimes\young(1)\otimes\young(1)
\end{gather*}
We can see that $\varphi_1(b)=3$.
Then we can apply $\widetilde{f}_1$ as follows
\begin{center}
\unitlength 12pt
\begin{picture}(40,2)
%{{{1,1,1,1,1,1,1,1}},{{2,1},{1},{},{}},{{2,1},{0},{},{}}}
\put(0.3,1.1){5}
\put(0.3,2.1){3}
\put(1,1){\Yboxdim12pt\yng(2,1)}
\put(2.2,1.1){1}
\put(3.2,2.1){2}
\put(4.5,0){
\put(0.3,2.1){0}
\put(1,2){\Yboxdim12pt\yng(1)}
\put(2.2,2.1){0}
}
\put(7.6,1.1){$\xrightarrow{\widetilde{f}_1}$}
\put(9.2,0){
\put(0.3,0.1){3}
\put(0.3,1.1){3}
\put(0.3,2.1){1}
\put(1,0){\Yboxdim12pt\yng(2,1,1)}
\put(2.2,0.1){$-1$}
\put(2.2,1.1){$-1$}
\put(3.2,2.1){0}
}
\put(13.7,0){
\put(0.3,2.1){1}
\put(1,2){\Yboxdim12pt\yng(1)}
\put(2.2,2.1){1}
}
\put(16.8,1.1){$\xrightarrow{\widetilde{f}_1}$}
\put(19.0,0){
\put(0.3,0.1){3}
\put(-0.4,1.1){$-1$}
\put(-0.4,2.1){$-1$}
\put(1,0){\Yboxdim12pt\yng(2,2,1)}
\put(2.2,0.1){$-1$}
\put(3.2,1.1){$-2$}
\put(3.2,2.1){$-2$}
}
\put(24.0,0){
\put(0.3,2.1){1}
\put(1,2){\Yboxdim12pt\yng(1)}
\put(2.2,2.1){1}
}
\put(27.1,1.1){$\xrightarrow{\widetilde{f}_1}$}
\put(29.3,0){
\put(0.3,0.1){3}
\put(-0.4,1.1){$-1$}
\put(-0.4,2.1){$-3$}
\put(1,0){\Yboxdim12pt\yng(3,2,1)}
\put(2.2,0.1){$-1$}
\put(3.2,1.1){$-2$}
\put(4.2,2.1){$-3$}
}
\put(35.2,0){
\put(0.3,2.1){1}
\put(1,2){\Yboxdim12pt\yng(1)}
\put(2.2,2.1){1}
}
\end{picture}
\end{center}
The leftmost rigged conf\/iguration corresponds to the path~$b$.
For $(\nu,J)=\Phi^{-1}(b)$ we have $P^{(1)}_\infty(\nu)=3$ and $s$ in Theorem~\ref{prop:phi_RC} is $s=0$ since the
smallest rigging is 1.
Therefore we have $\varphi_1(\nu,J)=P^{(1)}_\infty(\nu)-s=3$.
\end{Example}

Combining Proposition~\ref{prop:check_crystal} and Corollary~\ref{cor:check_crystal}, we obtain the following result.
\begin{Theorem}%\label{th:axiom}
The maps $\varepsilon_i$, $\varphi_i$, $\mathrm{wt}$ and the Kashiwara operators $\widetilde{e}_i$, $\widetilde{f}_i$
defined for the rigged configurations satisfy the axioms for the $U_q(\mathfrak{g}_0)$-crystals for
$\mathfrak{g}=A^{(1)}_n$ or $D^{(1)}_n$ given in Definition~{\rm \ref{def:crystal}}.
\end{Theorem}

For the later purposes, let us prepare some ref\/ined estimates related with the convexity relations.
The f\/irst one is an intuitive estimate.

\begin{Lemma}
\label{lem:convexity2}
Let $m^{(i)}_l(\nu)=m$.
Then $P^{(i)}_{l-1}(\nu)$, $P^{(i)}_{l}(\nu)$, $P^{(i)}_{l+1}(\nu)-2m$ satisfy the convexity relation.
\end{Lemma}
\begin{proof}
Let $(\nu',J)$ be the f\/ictitious rigged conf\/iguration obtained by replacing all length $l$ strings of $(\nu,J)^{(i)}$ by
length $l+1$ strings.
From the similar argument of the proof of Proposition~\ref{lem:convexity1}, we can show that $P^{(i)}_{l-1}(\nu')$,
$P^{(i)}_{l}(\nu')$, $P^{(i)}_{l+1}(\nu')$ satisfy the convexity relation.
Since $P^{(i)}_{l-1}(\nu')=P^{(i)}_{l-1}(\nu)$, $P^{(i)}_{l}(\nu')=P^{(i)}_{l}(\nu)$ and
$P^{(i)}_{l+1}(\nu')=P^{(i)}_{l+1}(\nu)-2m$ we obtain the assertion.
\end{proof}

For some purposes, we need a~more ref\/ined estimate as follows.

\begin{Lemma}
\label{lem:convexity3}
Suppose that $l\geq 1$.
For $1\leq a< n-2$ we have
\begin{gather*}
-P^{(a)}_{l-1}(\nu)+2P^{(a)}_l(\nu)-P^{(a)}_{l+1}(\nu)\geq m^{(a-1)}_l(\nu)-2m^{(a)}_l(\nu)+m^{(a+1)}_l(\nu),
\end{gather*}
and for $n-2\leq a$ we have
\begin{gather*}
-P^{(n-2)}_{l-1}(\nu)+2P^{(n-2)}_l(\nu)-P^{(n-2)}_{l+1}(\nu)\geq
m^{(n-3)}_l(\nu)-2m^{(n-2)}_l(\nu)+m^{(n-1)}_l(\nu)+m^{(n)}_l(\nu),
\\
-P^{(n-1)}_{l-1}(\nu)+2P^{(n-1)}_l(\nu)-P^{(n-1)}_{l+1}(\nu)\geq m^{(n-2)}_l(\nu)-2m^{(n-1)}_l(\nu),
\\
-P^{(n)}_{l-1}(\nu)+2P^{(n)}_l(\nu)-P^{(n)}_{l+1}(\nu)\geq m^{(n-2)}_l(\nu)-2m^{(n)}_l(\nu).
\end{gather*}
\end{Lemma}
\begin{proof}
Since $Q^{(a)}_l(\nu)$ counts the number of boxes of the left $l$ columns of $\nu$, we have
\begin{gather*}
-Q^{(a)}_{l-1}(\nu)+2Q^{(a)}_l(\nu)-Q^{(a)}_{l+1}(\nu)
=\big\{Q^{(a)}_l(\nu)-Q^{(a)}_{l-1}(\nu)\big\}
-\big\{Q^{(a)}_{l+1}(\nu)-Q^{(a)}_{l}(\nu)\big\}
\\
\qquad
=\sum\limits_{k\geq l}
m^{(a)}_k(\nu)-
\sum\limits_{k\geq l+1}
m^{(a)}_k(\nu)
=m^{(a)}_l(\nu).
\end{gather*}
Similarly, we obtain
\begin{gather*}
-Q^{(a)}_{l-1}(\mu)+2Q^{(a)}_l(\mu)-Q^{(a)}_{l+1}(\mu)=L^{(a)}_l.
\end{gather*}
Then for $a<n-2$ we have
\begin{gather*}
-P^{(a)}_{l-1}(\nu)+2P^{(a)}_l(\nu)-P^{(a)}_{l+1}(\nu)= L^{(a)}_l+m^{(a-1)}_l(\nu)-2m^{(a)}_l(\nu)+m^{(a+1)}_l(\nu).
\end{gather*}
Since $L^{(a)}_l\geq 0$, we obtain the assertion.
Similarly we can show the assertions for $n-2\leq a$.
\end{proof}

A nice property about these convexity relations is that they do not involve explicitly the information on $L^{(a)}_l$ or
$\mu$.
Indeed, we will not refer to $L^{(a)}_l$ or $\mu$ explicitly during the arguments of Sections~\ref{sec:main},~\ref{sec:main2} and~\ref{sec:main2_e}.

\subsection{Rigged conf\/iguration bijection}\label{sec:def_RCbijection}

Let us def\/ine the bijection between the rigged conf\/igurations and tensor products of crystals expressed by the KR
tableaux which we call {\it paths}.
For the map from the rigged conf\/iguration to paths, the basic operation is called~$\delta$
\begin{gather*}
\delta: \ (\nu,J)\longmapsto \{(\nu',J'),k\},
\end{gather*}
where $(\nu,J)$ and $(\nu',J')$ are rigged conf\/igurations and
$k\in\{1,2,\ldots,n,\mn,\ldots,\mtwo,\mone\}$.
In the description, we call a~string $(l,x_l)$ of $(\nu,J)^{(a)}$ {\it singular} if we have $x_l=P^{(a)}_l(\nu)$, that
is, the rigging takes the maximal possible value.
To begin with we consider the generic case (non-spin case).
Necessary modif\/ications related with the spin cases will be given in Section~\ref{sec:spin}.

\begin{Definition}%\label{def:delta}
Let us consider a~rigged conf\/iguration $(\nu,J)$ corresponding to the tensor product of the form $B^{a,l}\otimes B'$
(equivalently, $\mu^{(a)}$ has length $l$ row).
For $1\leq a\leq n-2$, the map~$\delta^{(a)}_l$
\begin{gather*}
\delta^{(a)}_l: \ (\nu,J)\longmapsto \{(\nu',J'),k\}
\end{gather*}
is def\/ined by the following procedure.
Set $\ell^{(a-1)}=l$.
\begin{enumerate}\itemsep=0pt
\item[(1)] For $a\leq i\leq n-2$, assume that $\ell^{(i-1)}$ is already determined.
Then we search for the shortest singular string in $(\nu,J)^{(i)}$ that is longer than or equal to $\ell^{(i-1)}$.
\begin{enumerate}\itemsep=0pt
\item
If there exists such a~string, set $\ell^{(i)}$ to be the length of the selected string and continue the process
recursively.
If there is more than one such string, choose any of them.
\item
If there is no such string, set $\ell^{(i)}=\infty$, $k=i$ and stop.
\end{enumerate}
\item[(2)] Suppose that $\ell^{(n-2)}<\infty$.
Then we search for the shortest singular string in $(\nu,J)^{(n-1)}$ (resp.
$(\nu,J)^{(n)}$) that is longer than or equal to $\ell^{(n-2)}$ and def\/ine $\ell^{(n-1)}$ (resp.
$\ell^{(n)}$) similarly.
\begin{enumerate}\itemsep=0pt
\item
If $\ell^{(n-1)}=\infty$ and $\ell^{(n)}=\infty$, set $k=n-1$ and stop.
\item
If $\ell^{(n-1)}<\infty$ and $\ell^{(n)}=\infty$, set $k=n$ and stop.
\item
If $\ell^{(n-1)}=\infty$ and $\ell^{(n)}<\infty$, set $k=\mn$ and stop.
\item
If $\ell^{(n-1)}<\infty$ and $\ell^{(n)}<\infty$, set $\ell_{(n-1)}=\max\{\ell^{(n-1)},\ell^{(n)}\}$ and continue.
\end{enumerate}
\item[(3)] For $1\leq i\leq n-2$, assume that $\ell_{(i+1)}$ is already def\/ined.
Then we search for the shortest singular string in $(\nu,J)^{(i)}$ that is longer than or equal to $\ell_{(i+1)}$ and
has not yet been selected as $\ell^{(i)}$.
Def\/ine $\ell_{(i)}$ similarly.
If $\ell_{(i)}=\infty$, set $k=\overline{i+1}$ and stop.
Otherwise continue.
If $\ell_{(1)}<\infty$, set $k=\mone$ and stop.
\item[(4)] Once the process has stopped, remove the rightmost box of each selected row specif\/ied by~$\ell^{(i)}$ or~$\ell_{(i)}$.
The result gives the output $\nu'$.
\item[(5)] Def\/ine the new riggings $J'$ as follows.
For the rows that are not selected by $\ell^{(i)}$ or $\ell_{(i)}$, take the corresponding riggings from $J$.
In order to def\/ine the remaining riggings, replace one $B^{a,l}$ in $B$ by $B^{a-1,1}\otimes B^{a,l-1}$ (equivalently,
replace one of the length $l$ row of $\mu^{(a)}$ by a~length $(l-1)$ row and add a~length $1$ row to $\mu^{(a-1)}$).
Denote the result by $B'$.
Use $B'$ to compute all the vacancy numbers for $\nu'$.
Then the remaining riggings are def\/ined so that all the corresponding rows become singular with respect to the new
vacancy number.
\end{enumerate}
\end{Definition}

We remark that the resulting rigged conf\/iguration $(\nu',J')$ is associated with the tensor pro\-duct~$B'$.
We write $\delta_2\delta_1(\nu,J)$ etc.
for repeated applications of~$\delta$ on the rigged conf\/igurations.

\begin{Definition}
For a~given rigged conf\/iguration $(\nu,J)$ corresponding to the tensor product of the form $B=B^{r_1,s_1}\otimes
B^{r_2,s_2}\otimes \dots\otimes B^{r_L,s_L}$, def\/ine the map $\Phi_B$ (sometimes also just denoted~$\Phi$)
\begin{gather*}
\Phi_B: \ (\nu,J)\longmapsto b
\end{gather*}
as follows.
\begin{enumerate}\itemsep=0pt
\item[(1)] Suppose that $\delta^{(1)}_1\cdots\delta^{(r_1-1)}_1\delta^{(r_1)}_{s_1}(\nu,J)=(\nu',J')$ yields the sequence
of letters $k^{(r_1)}$, $k^{(r_1-1)}$, $\ldots$, $k^{(1)}$ ($k^{(a)}$ corresponds to $\delta^{(a)}$).
Put the transpose of the row
\begin{gather*}
\begin{array}{|c|c|c|c|}
\hline
k^{(1)}& k^{(2)}& \dots& k^{(r_1)}\rule{0pt}{12pt}
\\
\hline
\end{array}
\end{gather*}
as the leftmost column of the rectangle $(s_1^{r_1})$.
\item[(2)] Continue the previous step for
$\delta^{(1)}_1\cdots\delta^{(r_1-1)}_1\delta^{(r_1)}_{s_1-1}(\nu',J')=(\nu'',J'')$ and f\/ill the second column for the
rectangle $(s_1^{r_1})$ with the produced letters.
Repeat the process until all places of $(s_1^{r_1})$ are f\/illed.
\item[(3)] Repeat the previous two steps for the remaining rectangles $(s_2^{r_2}),(s_3^{r_3}),\ldots,(s_L^{r_L})$.
\end{enumerate}
\end{Definition}

\begin{Remark}
\label{rem:delta_inverse}
The inverse procedure $\Phi^{-1}$ is obtained by reversing all procedure in~$\Phi$ step by step.
As an example let us give a~sketch of the algorithm for the operation $\big(\delta^{(1)}_1\big)^{-1}$.
Suppose that we start from an unbarred letter $k$ in a~KR tableau.
Then we look for the largest singular string of $\nu^{(k)}$.
If there is no such string, add a~length one row to the bottom of $\nu^{(k)}$.
Otherwise add a~box to the longest singular row of $\nu^{(k)}$.
Suppose that we add the box to the $\ell^{(k)}$th column of $\nu^{(k)}$.
Then we look for the singular string of $\nu^{(k-1)}$ which is strictly shorter than $\ell^{(k)}$.
We continue this process up to $\nu^{(1)}$.
Finally change the riggings for the modif\/ied strings according to the new vacancy numbers.
Here remind that we have to add a~length one row to $\mu^{(1)}$ for the computation of the vacancy numbers.
\end{Remark}

The following fundamental result is conjectured in~\cite{OSS:2012}.

\begin{Conjecture}\label{conj:well_def}
Consider a~tensor product of crystals $B=B^{r_1,s_1}\otimes B^{r_2,s_2}\otimes \dots\otimes B^{r_L,s_L}$.
Then the map $\Phi_B$ gives a~well-defined bijection between the rigged configuration of type $B$ and the elements of
$B$ expressed via the KR tableaux representation.
\end{Conjecture}

Currently this conjecture is verif\/ied in the following cases:
\begin{itemize}\itemsep=0pt
\item
$B=B^{r_1,s_1}\otimes B^{r_2,s_2}\otimes \dots\otimes B^{r_L,s_L}$ type tensor products of $A^{(1)}_n$~\cite{KSS:2002},
\item
$B=B^{r_1,1}\otimes B^{r_2,1}\otimes \dots\otimes B^{r_L,1}$ type tensor products of $D^{(1)}_n$~\cite{S:2005},
\item
$B=B^{1,s_1}\otimes B^{1,s_2}\otimes \dots\otimes B^{1,s_L}$ type tensor products of $D^{(1)}_n$~\cite{SS},
\item
$B=B^{r,s}$ ($r<n-1$) case of $D^{(1)}_n$~\cite{OSS:2012}.
\end{itemize}
A proof for the general case $B=B^{r_1,s_1}\otimes B^{r_2,s_2}\otimes \dots\otimes B^{r_L,s_L}$ of type $D^{(1)}_n$ is
the subject of~\cite{OSSS:2012}.
However, since the paper~\cite{OSSS:2012} is yet not available for public, we will not assume this result in this paper.

If there is no afraid of confusion, we use the abbreviation~$\Phi$ instead of $\Phi_B$.

\begin{Example}
Let us consider the following rigged conf\/iguration corresponding to the tensor product $B^{2,2}\otimes B^{3,2}$ of type
$D^{(1)}_5$:
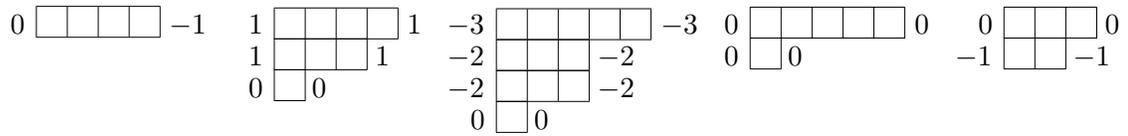
\begin{center}
\unitlength 12pt
\begin{picture}(35,4)(-1,0)
\put(0,3){
\put(-0.8,0.1){0}
\put(0,0){\Yboxdim12pt\yng(4)}
\put(4.2,0.1){$-1$}
}
\put(7.5,1){
\put(-0.8,0.1){0}
\put(-0.8,1.1){1}
\put(-0.8,2.1){1}
\put(0,0){\Yboxdim12pt\yng(4,3,1)}
\put(1.2,0.1){0}
\put(3.2,1.1){1}
\put(4.2,2.1){1}
}
\put(14.5,0){
\put(-0.8,0.1){0}
\put(-1.53,1.1){$-2$}
\put(-1.53,2.1){$-2$}
\put(-1.53,3.1){$-3$}
\put(0,0){\Yboxdim12pt\yng(5,3,3,1)}
\put(1.2,0.1){0}
\put(3.2,1.1){$-2$}
\put(3.2,2.1){$-2$}
\put(5.2,3.1){$-3$}
}
\put(22.5,2){
\multiput(-0.8,0.1)(0,1){2}{0}
\put(0,0){\Yboxdim12pt\yng(5,1)}
\put(1.2,0.1){0}
\put(5.2,1.1){0}
}
\put(30.5,2){
\put(-1.53,0.1){$-1$}
\put(-0.8,1.1){0}
\put(0,0){\Yboxdim12pt\yng(3,2)}
\put(2.2,0.1){$-1$}
\put(3.2,1.1){0}
}
\end{picture}
\end{center}
We start by removing $B^{2,2}$ part.
Then the f\/irst operation is $\delta^{(2)}_2$ which start by searching the shortest singular string of
$\nu^{(2)}=(4,3,1)$ whose length is equal to or larger than 2.
We work out all the procedure corresponding to $B^{2,2}$ in the following sequence of diagrams
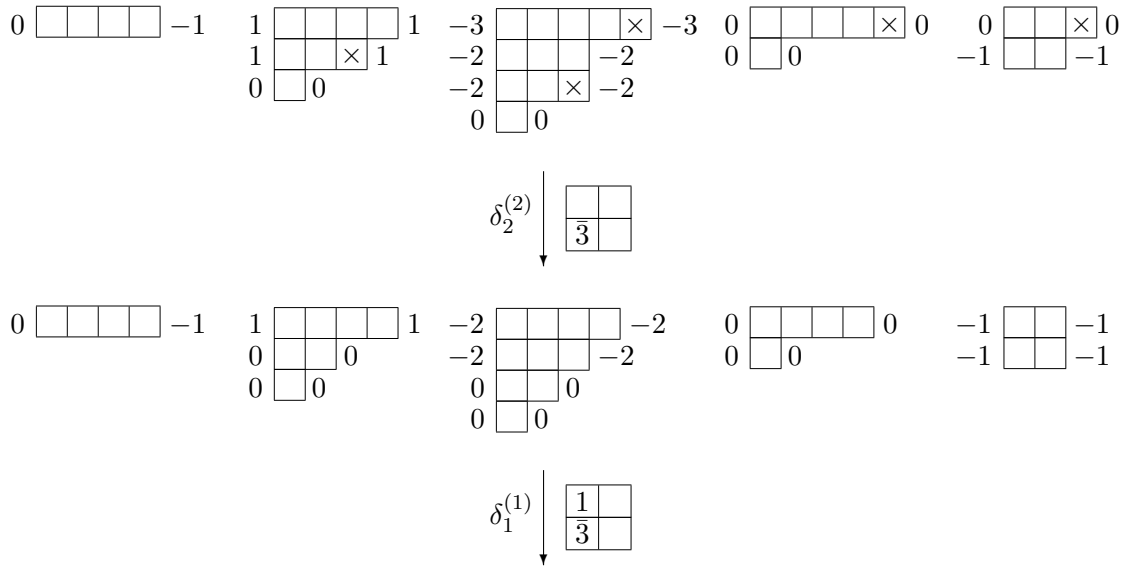
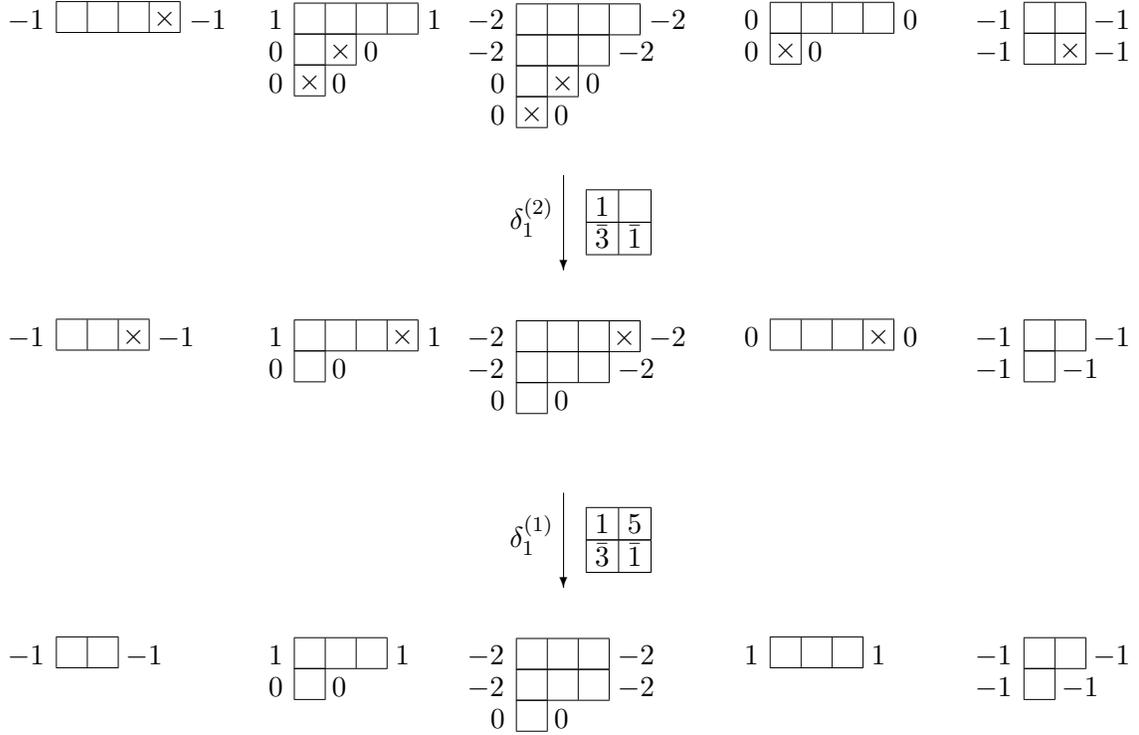
\begin{center}
\unitlength 12pt
\begin{picture}(35,4)(-1,0)
\put(0,3){
\put(-0.8,0.1){0}
\put(0,0){\Yboxdim12pt\yng(4)}
\put(4.2,0.1){$-1$}
}
\put(7.5,1){
\put(-0.8,0.1){0}
\put(-0.8,1.1){1}
\put(-0.8,2.1){1}
\put(0,0){$\Yboxdim12pt\young(\aki\aki\aki\aki,\aki\aki\times,\aki)$}
\put(1.2,0.1){0}
\put(3.2,1.1){1}
\put(4.2,2.1){1}
}
\put(14.5,0){
\put(-0.8,0.1){0}
\put(-1.53,1.1){$-2$}
\put(-1.53,2.1){$-2$}
\put(-1.53,3.1){$-3$}
\put(0,0){\Yboxdim12pt\young(\aki\aki\aki\aki\times,\aki\aki\aki,\aki\aki\times,\aki)}
\put(1.2,0.1){0}
\put(3.2,1.1){$-2$}
\put(3.2,2.1){$-2$}
\put(5.2,3.1){$-3$}
}
\put(22.5,2){
\multiput(-0.8,0.1)(0,1){2}{0}
\put(0,0){\Yboxdim12pt\young(\aki\aki\aki\aki\times,\aki)}
\put(1.2,0.1){0}
\put(5.2,1.1){0}
}
\put(30.5,2){
\put(-1.53,0.1){$-1$}
\put(-0.8,1.1){0}
\put(0,0){\Yboxdim12pt\young(\aki\aki\times,\aki\aki)}
\put(2.2,0.1){$-1$}
\put(3.2,1.1){0}
}
\end{picture}
%\{\{4\},\{4,3,1\},\{5,3,3,1\},\{5,1\},\{3,2\}\}\\
%\{\{-1\},\{1,1,0\},\{-3,-2,-2,0\},\{0,0\},\{0,-1\}\}
\end{center}
\begin{center}
\unitlength 12pt
\begin{picture}(5,3)
\put(0.3,1.3){$\delta^{(2)}_2$}
\put(2,3){\vector(0,-1){3}}
\put(2.7,0.5){\young(\aki\aki,\mthree\aki)}
\end{picture}
\end{center}
\begin{center}
\unitlength 12pt
\begin{picture}(35,4)(-1,0)
\put(0,3){
\put(-0.8,0.1){0}
\put(0,0){\Yboxdim12pt\yng(4)}
\put(4.2,0.1){$-1$}
}
\put(7.5,1){
\put(-0.8,0.1){0}
\put(-0.8,1.1){0}
\put(-0.8,2.1){1}
\put(0,0){$\Yboxdim12pt\young(\aki\aki\aki\aki,\aki\aki,\aki)$}
\put(1.2,0.1){0}
\put(2.2,1.1){0}
\put(4.2,2.1){1}
}
\put(14.5,0){
\put(-0.8,0.1){0}
\put(-0.8,1.1){0}
\put(-1.53,2.1){$-2$}
\put(-1.53,3.1){$-2$}
\put(0,0){\Yboxdim12pt\young(\aki\aki\aki\aki,\aki\aki\aki,\aki\aki,\aki)}
\put(1.2,0.1){0}
\put(2.2,1.1){0}
\put(3.2,2.1){$-2$}
\put(4.2,3.1){$-2$}
}
\put(22.5,2){
\multiput(-0.8,0.1)(0,1){2}{0}
\put(0,0){\Yboxdim12pt\young(\aki\aki\aki\aki,\aki)}
\put(1.2,0.1){0}
\put(4.2,1.1){0}
}
\put(30.5,2){
\put(-1.53,0.1){$-1$}
\put(-1.53,1.1){$-1$}
\put(0,0){\Yboxdim12pt\young(\aki\aki,\aki\aki)}
\put(2.2,0.1){$-1$}
\put(2.2,1.1){$-1$}
}
\end{picture}
%\{\{4\},\{4,2,1\},\{4,3,2,1\},\{4,1\},\{2,2\}\}\\
%\{\{-1\},\{1,0,0\},\{-2,-2,0,0\},\{0,0\},\{-1,-1\}\}
\end{center}
\begin{center}
\unitlength 12pt
\begin{picture}(5,3)
\put(0.3,1.3){$\delta^{(1)}_1$}
\put(2,3){\vector(0,-1){3}}
\put(2.7,0.5){\young(1\aki,\mthree\aki)}
\end{picture}
\end{center}
\begin{center}
\unitlength 12pt
\begin{picture}(35,4)(-1,0)
\put(0,3){
\put(-1.53,0.1){$-1$}
\put(0,0){\Yboxdim12pt\young(\aki\aki\aki\times)}
\put(4.2,0.1){$-1$}
}
\put(7.5,1){
\put(-0.8,0.1){0}
\put(-0.8,1.1){0}
\put(-0.8,2.1){1}
\put(0,0){$\Yboxdim12pt\young(\aki\aki\aki\aki,\aki\times,\times)$}
\put(1.2,0.1){0}
\put(2.2,1.1){0}
\put(4.2,2.1){1}
}
\put(14.5,0){
\put(-0.8,0.1){0}
\put(-0.8,1.1){0}
\put(-1.53,2.1){$-2$}
\put(-1.53,3.1){$-2$}
\put(0,0){\Yboxdim12pt\young(\aki\aki\aki\aki,\aki\aki\aki,\aki\times,\times)}
\put(1.2,0.1){0}
\put(2.2,1.1){0}
\put(3.2,2.1){$-2$}
\put(4.2,3.1){$-2$}
}
\put(22.5,2){
\multiput(-0.8,0.1)(0,1){2}{0}
\put(0,0){\Yboxdim12pt\young(\aki\aki\aki\aki,\times)}
\put(1.2,0.1){0}
\put(4.2,1.1){0}
}
\put(30.5,2){
\put(-1.53,0.1){$-1$}
\put(-1.53,1.1){$-1$}
\put(0,0){\Yboxdim12pt\young(\aki\aki,\aki\times)}
\put(2.2,0.1){$-1$}
\put(2.2,1.1){$-1$}
}
\end{picture}
%\{\{4\},\{4,2,1\},\{4,3,2,1\},\{4,1\},\{2,2\}\}\\
%\{\{-1\},\{1,0,0\},\{-2,-2,0,0\},\{0,0\},\{-1,-1\}\}
\end{center}
\begin{center}
\unitlength 12pt
\begin{picture}(5,3)
\put(0.3,1.3){$\delta^{(2)}_1$}
\put(2,3){\vector(0,-1){3}}
\put(2.7,0.5){\young(1\aki,\mthree\mone)}
\end{picture}
\end{center}
\begin{center}
\unitlength 12pt
\begin{picture}(35,4)(-1,0)
\put(0,3){
\put(-1.53,0.1){$-1$}
\put(0,0){\Yboxdim12pt\young(\aki\aki\times)}
\put(3.2,0.1){$-1$}
}
\put(7.5,1){
\put(-0.8,1.1){0}
\put(-0.8,2.1){1}
\put(0,1){$\Yboxdim12pt\young(\aki\aki\aki\times,\aki)$}
\put(1.2,1.1){0}
\put(4.2,2.1){1}
}
\put(14.5,1){
\put(-0.8,0.1){0}
\put(-1.53,1.1){$-2$}
\put(-1.53,2.1){$-2$}
\put(0,0){\Yboxdim12pt\young(\aki\aki\aki\times,\aki\aki\aki,\aki)}
\put(1.2,0.1){0}
\put(3.2,1.1){$-2$}
\put(4.2,2.1){$-2$}
}
\put(22.5,3){
\put(-0.8,0.1){0}
\put(0,0){\Yboxdim12pt\young(\aki\aki\aki\times)}
\put(4.2,0.1){0}
}
\put(30.5,2){
\put(-1.53,0.1){$-1$}
\put(-1.53,1.1){$-1$}
\put(0,0){\Yboxdim12pt\young(\aki\aki,\aki)}
\put(1.2,0.1){$-1$}
\put(2.2,1.1){$-1$}
}
\end{picture}
%\{\{3\},\{4,1\},\{4,3,1\},\{4\},\{2,1\}\}\\
%\{\{-1\},\{1,0\},\{-2,-2,0\},\{0\},\{-1,-1\}\}
\end{center}
\begin{center}
\unitlength 12pt
\begin{picture}(5,3)
\put(0.3,1.3){$\delta^{(1)}_1$}
\put(2,3){\vector(0,-1){3}}
\put(2.7,0.5){\young(15,\mthree\mone)}
\end{picture}
\end{center}
\begin{center}
\unitlength 12pt
\begin{picture}(35,4)(-1,0)
\put(0,3){
\put(-1.53,0.1){$-1$}
\put(0,0){\Yboxdim12pt\young(\aki\aki)}
\put(2.2,0.1){$-1$}
}
\put(7.5,1){
\put(-0.8,1.1){0}
\put(-0.8,2.1){1}
\put(0,1){$\Yboxdim12pt\young(\aki\aki\aki,\aki)$}
\put(1.2,1.1){0}
\put(3.2,2.1){1}
}
\put(14.5,1){
\put(-0.8,0.1){0}
\put(-1.53,1.1){$-2$}
\put(-1.53,2.1){$-2$}
\put(0,0){\Yboxdim12pt\young(\aki\aki\aki,\aki\aki\aki,\aki)}
\put(1.2,0.1){0}
\put(3.2,1.1){$-2$}
\put(3.2,2.1){$-2$}
}
\put(22.5,3){
\put(-0.8,0.1){1}
\put(0,0){\Yboxdim12pt\young(\aki\aki\aki)}
\put(3.2,0.1){1}
}
\put(30.5,2){
\put(-1.53,0.1){$-1$}
\put(-1.53,1.1){$-1$}
\put(0,0){\Yboxdim12pt\young(\aki\aki,\aki)}
\put(1.2,0.1){$-1$}
\put(2.2,1.1){$-1$}
}
\end{picture}
\end{center}
Here the boxes to be removed by $\delta^{(a)}_l$ indicated below the corresponding rigged conf\/igurations are marked by
``$\times$".
The output of each $\delta^{(a)}_l$ is given on the right of the corresponding arrows.
Note that after $\delta^{(2)}_2$ we recalculate all the vacancy numbers assuming that the resulting rigged conf\/iguration
corresponds to the tensor product of type $B^{1,1}\otimes B^{2,1}\otimes B^{3,2}$.
We remark that the f\/irst $\delta^{(1)}_1$ cannot remove a~box from $\nu^{(1)}$ since there is no singular string and,
$\delta^{(2)}_1$ starts by removing from $\nu^{(2)}$ and ends at $\nu^{(1)}$.
The f\/inal output for the given rigged conf\/iguration $(\nu,J)$ is as follows:
\begin{gather}
\label{eq:Phi_example_1}
\Phi_{B^{2,2}\otimes B^{3,2}}(\nu,J)=\Yvcentermath1\young(15,\mthree\mone)\otimes \young(1\mfive,4\mthree,5\mone)
\end{gather}
\end{Example}

\begin{Remark}
\label{rem:combR}
In the above example, we can reverse the order of the tensor product.
In this case, we obtain
\begin{gather}
\label{eq:Phi_example_2}
\Phi_{B^{3,2}\otimes B^{2,2}}(\nu,J)= \Yvcentermath1 \young(14,5\mtwo,\mthree\mone)\otimes
\young(2\mfive,5\mthree)
\end{gather}
Then the two tensor products in~\eqref{eq:Phi_example_1} and~\eqref{eq:Phi_example_2} are isomorphic under the
combinatorial $R$-matrix
\begin{gather*}
R:\Yvcentermath1\young(15,\mthree\mone)\otimes \young(1\mfive,4\mthree,5\mone) \longmapsto
\young(14,5\mtwo,\mthree\mone)\otimes \young(2\mfive,5\mthree)
\end{gather*}
Therefore in this case we have $R=\Phi_{B^{3,2}\otimes B^{2,2}}\circ\Phi_{B^{2,2}\otimes B^{3,2}}^{-1}$.
This relation is valid not only for two times tensor product but also for multiple applications of pairwise
combinatorial $R$-matrices for arbitrary tensor products.
\end{Remark}

\section{Rigged conf\/igurations and the Kashiwara operators}\label{sec:statement}

\subsection{Statement of the main result}

Let us state the main result of the present paper.

\begin{Theorem}
\label{th:main}
Suppose that the rigged configuration bijection~$\Phi$ is well-defined.
Then the bijection~$\Phi$ and the Kashiwara operators $\widetilde{e}_i$ and $\widetilde{f}_i$ $(i\in I_0)$ of types
$A^{(1)}_n$ and $D^{(1)}_n$ satisfy the following commutative diagrams:
\begin{gather*}
\begin{CD}
\mathrm{RC}(B) @>\Phi >> B
\\
@V\widetilde{f}_i VV @VV\widetilde{f}_iV
\\
\mathrm{RC}(B)\cup\{0\} @>>\Phi > B\cup\{0\}
\end{CD}
\qquad
\begin{CD}
\mathrm{RC}(B) @>\Phi >> B
\\
@V\widetilde{e}_i VV @VV\widetilde{e}_iV
\\
\mathrm{RC}(B)\cup\{0\} @>>\Phi > B\cup\{0\}
\end{CD}
\end{gather*}
\end{Theorem}
Here for the tensor product of crystals $B=B^{r_L,s_L}\otimes B^{r_{L-1},s_{L-1}}\otimes\dots\otimes B^{r_1,s_1}$, we
denote by $\mathrm{RC}(B)$ the rigged conf\/igurations which are mapped to elements of $B$ under $\Phi_B$, that is,
$\Phi_B:\mathrm{RC}(B)\rightarrow B$.

\begin{Remark}
Note that $\mathrm{RC}(B)$ does not fully depend on $B$ but depends only on the shape of~$B$ via the quantities
$L^{(a)}_l$ or $\mu$.
Thus it is possible to denote the set of rigged conf\/igurations as~$\mathrm{RC}(L(B))$.
\end{Remark}

Therefore we have the compatibility of the rigged conf\/iguration bijection and the classical Kashiwara operators.
We give a~few words on the compatibility for the af\/f\/ine Kashiwara ope\-rators $\widetilde{e}_0$ and $\widetilde{f}_0$ of
type $D^{(1)}_n$.
Recall that in~\cite{S:2008} the operators $\widetilde{e}_0$ and $\widetilde{f}_0$ are realized via the bijection
between the set of $J$-highest element of $B^{r,s}$ and the combinatorial objects called the plus-minus diagrams.
Here $J=\{2,3,\ldots,n\}$.
In~\cite[Section 4]{OSS:2012}, we introduce an analogue of the plus-minus diagrams on the set of the rigged
conf\/igurations and show that this indeed yields the af\/f\/ine isomorphism with the KR crystal $B^{r,s}$.

The main points of the arguments in~\cite{OSS:2012} are as follows.
In~\cite[Theorem 5.9]{OSS:2012} we prove the combinatorial bijection~$\Phi$ for the $I_0$-highest element of $B^{r,s}$
directly.
Then in~\cite[Proposi\-tion~4.3]{OSS:2012} we show that the same sequence of the classical Kashiwara operators
$\widetilde{f}_i$ gives the plus-minus diag\-ram and the corresponding element in the set of rigged conf\/iguration.
Therefore Theorem~\ref{th:main} asserts that the rigged conf\/iguration analogue of the plus-minus diag\-rams introduced
in~\cite[Section~4]{OSS:2012} are indeed the image of the plus-diagrams under the rigged conf\/iguration bijection~$\Phi$.
Hence, in view of the construction of~\cite{S:2008}, we have the compatibility of the af\/f\/ine Kashiwara operators
$\widetilde{e}_0$, $\widetilde{f}_0$ and the rigged conf\/iguration bijection for the case~$B^{r,s}$.
However, the compatibility of the af\/f\/ine Kashiwara operators and the combinatorial bijection~$\Phi$ for more than two
tensor factors is still an open problem.

The rest of the paper is devoted to the proof of Theorem~\ref{th:main}.

\subsection{Preliminary steps}

Our main strategy of the proof is to reduce the essential combinatorial aspects of Theorem~\ref{th:main} to the
fundamental operation $\delta=\delta^{(1)}_1$.
For this purpose we decompose the original bijection~$\Phi$ into several steps according to the operations lh, lb and ls
def\/ined on the KR tableau representation of the crystals.
In the following, we def\/ine the corresponding operations (summarized in the following table) on the set of the rigged
conf\/igurations.
\begin{center}
\begin{tabular}
{|c|c|}
\hline
$B$&$\mathrm{RC}(B)$
\\
\hline
\hline
ls&$\gamma$
\\
lb&$\beta$
\\
lh&$\delta$
\\
\hline
\end{tabular}
\end{center}

First, let us remind that the following symbols bear specif\/ic meanings throughout the paper.
\begin{itemize}\itemsep=0pt
\item
$\nu$: shapes of the conf\/igurations,
\item
$J$: riggings,
\item
$\mu$: shape of the tensor product.
\end{itemize}

\begin{Definition}\label{def:RC_operations}
Let $B=B^{r_L,s_L}\otimes B^{r_{L-1},s_{L-1}}\otimes\dots\otimes B^{r_1,s_1}$.
\begin{enumerate}\itemsep=0pt
\item[(1)] Corresponding to $B^{r_L,s_L}=B^{a,l}$, $\gamma$ replaces a~length $l$ row of $\mu^{(a)}$ by rows of lengths $l-1$
and $1$.
\item[(2)] Suppose that $B^{r_L,s_L}=B^{r,1}$ where $r>1$.
If $r<n-1$, then $\beta$ removes a~length one row from $\mu^{(r)}$, adds a~length one row to each of $\mu^{(1)}$ and
$\mu^{(r-1)}$ and adds a~length one singular string to each of $(\nu,J)^{(a)}$ for $a<r$.
\item[(3)] Suppose that $B^{r_L,s_L}=B^{1,1}$.
Then def\/ine $\delta=\delta^{(1)}_1$.
\end{enumerate}
\end{Definition}

\begin{Proposition}
\label{prop:RC_decomp}
We can decompose the map~$\Phi$ into the following three steps.
\begin{enumerate}\itemsep=0pt
\item[$(1)$] Let $B=B^{r,s}\otimes \bar{B}$ with $s>1$.
Then we reduce the problem to the case of $B^{r,1}$ by
\begin{gather*}
\begin{CD}
\mathrm{RC}(B) @>\Phi >> B
\\
@V\gamma VV @VV\mathrm{ls}V
\\
\mathrm{RC}(B') @>>\Phi > B'
\end{CD}
\end{gather*}
where $B'=B^{r,1}\otimes B^{r,s-1}\otimes \bar{B}$.
\item[$(2)$] Let $B=B^{r,1}\otimes \bar{B}$ with $r>1$.
Then we reduce the problem to the case of $B^{1,1}$ by
\begin{gather*}
\begin{CD}
\mathrm{RC}(B) @>\Phi >> B
\\
@V\beta VV @VV\mathrm{lb}V
\\
\mathrm{RC}(B') @>>\Phi > B'
\end{CD}
\end{gather*}
where $B'=B^{1,1}\otimes B^{r-1,1}\otimes \bar{B}$.
\item[$(3)$] Let $B=B^{1,1}\otimes \bar{B}$.
Then we have to deal with the fundamental operation~$\delta$ as
\begin{gather*}
\begin{CD}
\mathrm{RC}(B) @>\Phi >> B
\\
@V\delta VV @VV\mathrm{lh} V
\\
\mathrm{RC}(\bar{B}) @>>\Phi > \bar{B}
\end{CD}
\end{gather*}
\end{enumerate}
\end{Proposition}

\begin{proof}
(1) Suppose that~$\Phi$ maps $(\nu,J)\in\mathrm{RC}(B)$ to $b\otimes\bar{b}\in B$.
Recall that $\gamma$ replaces a~length~$s$ row of $\mu^{(r)}$ by rows of lengths $s-1$ and~$1$.
Thus after application of $\gamma$ there is no singular string of $\nu^{(r)}$ which are strictly shorter than~$s$.
Then the length~$1$ row of~$\mu^{(r)}$ corresponds to the same column which appear as the leftmost column of the tableau~$b$.

(2) If we apply $\delta^{(1)}_1$ for $\mathrm{RC}(B')$, it automatically selects the length 1 singular strings of
$\nu^{(a)}$ ($a<r$) created by $\beta$.
As for $\nu^{(r)}$, since $\beta$ removes length 1 row from $\mu^{(r)}$ and adds a~length 1 row to $\nu^{(r-1)}$, all
the coriggings for $\nu^{(r)}$ do not change after $\beta$.
Thus $\delta^{(1)}_1$ selects the same string which is selected by $\delta^{(r)}_1$ before $\beta$.

(3) This follows from the def\/inition.
\end{proof}

\subsection{Reduction steps}

We show that the operations $\gamma$ and $\beta$ commutes with the Kashiwara operators.
\begin{Proposition}
\label{prop:ls_commute}
Let $B=B^{r,s}\otimes \bar{B}$ where $s>1$.
Then we have the following commutative diagrams
\begin{gather*}
\begin{CD}
\mathrm{RC}(B) @>\gamma >> \mathrm{RC}(B')
\\
@V\widetilde{e}_i VV @VV\widetilde{e}_i V
\\
\mathrm{RC}(B) @>>\gamma > \mathrm{RC}(B')
\end{CD}
\qquad
\begin{CD}
\mathrm{RC}(B) @>\gamma >> \mathrm{RC}(B')
\\
@V\widetilde{f}_i VV @VV\widetilde{f}_i V
\\
\mathrm{RC}(B) @>>\gamma > \mathrm{RC}(B')
\end{CD}
\end{gather*}
where $B'=B^{r,1}\otimes B^{r,s-1}\otimes \bar{B}$ and $i\in I_0$.
\end{Proposition}
\begin{proof}
For the rigged conf\/iguration $\{\mu,(\nu,J)\}$, $\widetilde{e}_i$ and~$\widetilde{f}_i$ change~$(\nu,J)$ part only and~$\gamma$ changes~$\mu$ part only.
Thus they commute.

We remark that if we have $\widetilde{e}_i(\nu,J)=0$ or $\widetilde{f}_i(\nu,J)=0$, we have
$\widetilde{e}_i\circ\gamma(\nu,J)=0$ or $\widetilde{f}_i\circ\gamma(\nu,J)=0$ by Theorem~\ref{prop:phi_RC} since
$\gamma$ does not alter the vacancy numbers nor riggings.
\end{proof}

\begin{Proposition}
\label{prop:lb_commute}
Let $B=B^{r,1}\otimes \bar{B}$ where $r>1$.
Then we have the following commutative diagrams
\begin{gather*}
\begin{CD}
\mathrm{RC}(B) @>\beta >> \mathrm{RC}(B')
\\
@V\widetilde{e}_i VV @VV\widetilde{e}_i V
\\
\mathrm{RC}(B) @>>\beta > \mathrm{RC}(B')
\end{CD}
\qquad
\begin{CD}
\mathrm{RC}(B) @>\beta >> \mathrm{RC}(B')
\\
@V\widetilde{f}_i VV @VV\widetilde{f}_i V
\\
\mathrm{RC}(B) @>>\beta > \mathrm{RC}(B')
\end{CD}
\end{gather*}
where $B'=B^{1,1}\otimes B^{r-1,1}\otimes \bar{B}$ and $i\in I_0$.
\end{Proposition}

\begin{proof}
Let us prove the $\widetilde{f}_i$ part f\/irst.
During the proof we denote that $\widetilde{f}_i(\nu,J)=(\widetilde{\nu},\widetilde{J})$ and
$\beta(\nu,J)=(\dot{\nu},\dot{J})$.

Suppose that $i>r$.
Recall that the operation $\beta$ does nothing on $(\nu,J)^{(a)}$ for $a\geq r$.
Thus $\beta$ and $\widetilde{f}_i$ commute since they do not have interaction.

Suppose that $i=r$.
Since the operation $\beta$ adds a~length one string to $(\nu,J)^{(r-1)}$, the interaction between $\beta$ and
$\widetilde{f}_{r}$ can occur if $\widetilde{f}_{r}$ adds the string $(1,-1)$ to $(\nu,J)^{(r)}$.
In this case, since $\beta$ does nothing on $(\nu,J)^{(r)}$, $\widetilde{f}_r$ adds the same string $(1,-1)$ after
$\beta$.
Let us consider $(\nu,J)^{(r-1)}$.
Then we have to check the coincidence of the string added by $\beta$:
\begin{gather*}
\begin{array}{@{}lllll}
\varnothing&\overset{\widetilde{f}_r}{\longmapsto}&\varnothing
&\overset{\beta}{\longmapsto}&\big(1,P^{(r-1)}_1(\widetilde{\nu})\big),
\\
\varnothing&\overset{\beta}{\longmapsto}&\big(1,P^{(r-1)}_1(\nu)\big)
&\overset{\widetilde{f}_r}{\longmapsto}&\big(1,P^{(r-1)}_1(\widetilde{\nu})\big).
\end{array}
\end{gather*}
Here we use the fact that $\widetilde{f}_r$ changes the riggings so as to keep the coriggings.

Finally, suppose that $i<r$.
Let~$x_k$ be the smallest rigging of length $k$ strings of $(\nu,J)^{(i)}$.
Let~$x_\ell$ be the smallest rigging of $(\nu,J)^{(i)}$ where $\ell$ is the length of the largest string with rigging~$x_\ell$.

Suppose that $\nu^{(i)}\neq\varnothing$ and $x_\ell\leq 0$.
Recall that $\beta$ adds the singular string $(1,P^{(i)}_1(\dot{\nu}))$ to~$(\nu,J)^{(i)}$.
If we show
\begin{gather*}
P^{(i)}_1(\dot{\nu})\geq x_\ell,
\end{gather*}
then we see that $\widetilde{f}_i$ chooses the same string $(\ell,x_\ell)$ even after $\beta$.
Let $j$ be the length of the shortest string of $(\nu,J)^{(i)}$.
Suppose that $j=1$.
Then we have $P^{(i)}_1(\dot{\nu})=P^{(i)}_1(\nu)$ by def\/inition of~$\beta$, $P^{(i)}_1(\nu)\geq x_1$ by the def\/inition
of the rigged conf\/igurations and $x_1\geq x_\ell$ by the minimality of~$x_\ell$.
Thus we have $P^{(i)}_1(\dot{\nu})\geq x_\ell$ in this case.
Next, suppose that $j>1$.
By the convexity relation of $P^{(i)}_k(\nu)$ between $0\leq k\leq j$ we have
\begin{gather*}
P^{(i)}_1(\dot{\nu})=P^{(i)}_1(\nu)\geq\min\big\{P^{(i)}_0(\nu),P^{(i)}_j(\nu)\big\}=\min\big\{0,P^{(i)}_j(\nu)\big\}.
\end{gather*}
If we have $P^{(i)}_j(\nu)\geq 0$, we have $P^{(i)}_1(\dot{\nu})\geq 0\geq x_\ell$.
On the other hand, suppose that $0>P^{(i)}_j(\nu)$.
Then we have $P^{(i)}_1(\dot{\nu})\geq P^{(i)}_j(\nu)$ by the above inequality, $P^{(i)}_j(\nu)\geq x_j$ by the
def\/inition of the rigged conf\/igurations and $x_j\geq x_\ell$ by the minimality of $x_\ell$.
Therefore we have $P^{(i)}_1(\dot{\nu})\geq x_\ell$.

Let us consider the case $\nu^{(i)}=\varnothing$ or $x_\ell>0$.
In both cases, if we have
\begin{gather*}
P^{(i)}_1(\dot{\nu})>0,
\end{gather*}
then we have the following behaviors of the two strings of $(\nu,J)^{(i)}$:
\begin{gather*}
\begin{array}{@{}lllll}
\begin{array}{l}
\varnothing
\\
\varnothing
\end{array}
&\overset{\beta}{\longmapsto}&
\begin{array}{l}
\big(1,P^{(i)}_1(\dot{\nu})\big)
\\
\varnothing
\end{array}
&\overset{\widetilde{f}_i}{\longmapsto}&
\begin{array}{l}
\big(1,P^{(i)}_1(\widetilde{\dot{\nu}})\big)
\\
(1,-1)
\end{array}
\end{array}
\end{gather*}
and, on the other hand,
\begin{gather*}
\begin{array}{@{}lllll}
\begin{array}{l}
\varnothing
\\
\varnothing
\end{array}
&\overset{\widetilde{f}_i}{\longmapsto}&
\begin{array}{l}
\varnothing
\\
(1,-1)
\end{array}
&\overset{\beta}{\longmapsto}&
\begin{array}{l}
(1,P^{(i)}_1(\dot{\widetilde{\nu}}))
\\
(1,-1)
\end{array}
.
\end{array}
\end{gather*}
Here we denote by $\widetilde{f}_i\circ\beta(\nu,J)=(\widetilde{\dot{\nu}},\widetilde{\dot{J}})$ and
$\beta\circ\widetilde{f}_i(\nu,J)=(\dot{\widetilde{\nu}},\dot{\widetilde{J}})$.
Since $\beta$ does not change vacancy numbers, we have
$P^{(i)}_1(\widetilde{\dot{\nu}})=P^{(i)}_1(\dot{\widetilde{\nu}})$.
Thus $\widetilde{f}_i$ and $\beta$ commute in this case.

Now let us consider the case $P^{(i)}_1(\nu)=P^{(i)}_1(\dot{\nu})\leq 0$.
\begin{itemize}\itemsep=0pt
\item
Suppose that we have $\nu^{(i)}=\varnothing$.
Then we have $P^{(i)}_\infty(\nu)\geq 0$ since we have $Q^{(i)}_\infty(\nu)=0$.
Now we invoke the convexity relation of $P^{(i)}_k(\nu)$ between $0\leq k\leq\infty$ with $0=P^{(i)}_0(\nu)\geq
P^{(i)}_1(\nu)\leq P^{(i)}_\infty(\nu)$ to deduce that $P^{(i)}_0(\nu)=P^{(i)}_1(\nu)=\dots =P^{(i)}_\infty(\nu)=0$.
Since $\nu^{(i)}=\varnothing$, this relation also implies that $\nu^{(i-1)}=\nu^{(i+1)}=\mu^{(i)}=\varnothing$.\\
Consider $\widetilde{f}_i$ acting on $(\nu,J)$.
Suppose that $i<r-1$.
Then it adds the string $(1,-1)$ to $(\nu,J)^{(i)}$ whereas the corresponding vacancy number is
$P^{(i)}_1(\widetilde{\nu})=P^{(i)}_1(\nu)-2=-2$.
Thus we have $\widetilde{f}_i(\nu,J)=0$.
On the other hand, $\beta$ acting on $(\nu,J)$ adds the length $1$ singular string and gives
$(\dot{\nu},\dot{J})^{(i)}=\{(1,P^{(i)}_1(\dot{\nu}))\}$.
The next $\widetilde{f}_i$ selects this string since we have $P^{(i)}_1(\dot{\nu})=P^{(i)}_1(\nu)=0$.
Then we have $(\widetilde{\dot{\nu}},\widetilde{\dot{J}})^{(i)}=\{(2,-1)\}$.
In particular we have $Q^{(i)}_2(\widetilde{\dot{\nu}})=2$.
Since $\nu^{(i-1)}=\nu^{(i+1)}=\varnothing$, both $\widetilde{\dot{\nu}}^{(i-1)}$ and $\widetilde{\dot{\nu}}^{(i+1)}$
consist of length $1$ string.
Thus $Q^{(i-1)}_2(\widetilde{\dot{\nu}})=Q^{(i+1)}_2(\widetilde{\dot{\nu}})=1$.
Since we also know $\mu^{(i)}=\varnothing$, we have $Q^{(i)}_2(\widetilde{\dot{\mu}})=0$.
Hence $P^{(i)}_2(\widetilde{\dot{\nu}})=Q^{(i)}_2(\widetilde{\dot{\mu}})
-2Q^{(i)}_2(\widetilde{\dot{\nu}})+Q^{(i-1)}_2(\widetilde{\dot{\nu}})+Q^{(i+1)}_2(\widetilde{\dot{\nu}})=-2$.
Therefore we have $\widetilde{f}_i\circ\beta(\nu,J)=0$.
Hence $\widetilde{f}_i$ and $\beta$ commute in this case.
If $i=r-1$, we can use the parallel argument if we replace $\nu^{(i+1)}$ by $\mu^{(i)}$.

\item
Since we have completed the proof for the case $\nu^{(i)}=\varnothing$, we can assume that $\nu^{(i)}\neq\varnothing$.
Let $j$ be the length of the smallest string of $(\nu,J)^{(i)}$.
Suppose that $j=1$.
Since we are considering the case $\nu^{(i)}=\varnothing$ or $x_\ell>0$, the condition $\nu^{(i)}\neq\varnothing$ implies
that $x_\ell>0$.
On the other hand we have $0\geq P^{(i)}_1(\nu)\geq x_1$ by the assumption and the def\/inition of the rigged
conf\/igurations.
This is in contradiction to the relation $x_1\geq x_\ell$ which follows from the minimality of $x_\ell$.
Hence this case cannot happen.

\item
Suppose that $j>1$.
Again we have $x_\ell>0$.
By the convexity relation of $P^{(i)}_k(\nu)$ between $0\leq k\leq j$ and the def\/inition of the rigged conf\/igurations we
have $0=P^{(i)}_0(\nu)\geq P^{(i)}_1(\nu)\geq P^{(i)}_j(\nu)\geq x_j$.
This is in contradiction to the relation $x_j\geq x_\ell$.
Hence this case cannot occur.
\end{itemize}
Thus we have checked the commutativity of $\widetilde{f}_i$ and $\beta$ in all possible cases.

Finally let us consider the proof for $\widetilde{e}_i$.
If $\widetilde{e}_i(\nu,J)\neq 0$, then the commutativity of~$\widetilde{e}_i$ and~$\beta$ follows from the result of
$\widetilde{f}_i$ since $\widetilde{e}_i$ and $\widetilde{f}_i$ are mutually inverse in this case.
Suppose that $\widetilde{e}_i(\nu,J)=0$.
Then all the riggings of $(\nu,J)^{(i)}$ are non-negative.
Recall that $\beta$ acting on $(\nu,J)$ adds the string $(1,P^{(i)}_1(\dot{\nu}))$.
Then if we show $P^{(i)}_1(\dot{\nu})=P^{(i)}_1(\nu)\geq 0$ we have $\widetilde{e}_i\circ\beta(\nu,J)=0$.
Let $j$ be the length of the shortest string of $(\nu,J)^{(i)}$.
If $\nu^{(i)}=\varnothing$, set $j=\infty$.
In both cases we have $P^{(i)}_j(\nu)\geq 0$ since we have $P^{(i)}_j(\nu)\geq x_j\geq 0$ if $j<\infty$ and
$Q^{(i)}_j(\nu)=0$ if $j=\infty$.
Now we invoke the convexity relation of $P^{(i)}_k(\nu)$ between $0\leq k\leq j$.
Then we have $P^{(i)}_1(\nu)\geq\min\big\{P^{(i)}_0(\nu),P^{(i)}_j(\nu)\big\} =\min\big\{0,P^{(i)}_j(\nu)\big\}=0$ as desired.
\end{proof}

\subsection{Statements to be proved}\label{subsec:notation}

Due to the arguments in the previous section, the essential points of the proof of Theorem~\ref{th:main} are reduced to
the following two assertions on the case $B=B^{1,1}\otimes \bar{B}$.
Up to Section~\ref{sec:main2} from here we will use the following notation when we consider the case $B=B^{1,1}\otimes
\bar{B}$.
Let $(\nu,J)\in\mathrm{RC}(B)$.
Denote the rigged conf\/igurations obtained by the actions of $\widetilde{f}_i$ or~$\delta$ as in the following diagram:
\begin{gather*}
\begin{CD}
(\nu,J) @>\delta >> (\bar{\nu},\bar{J})
\\
@V\widetilde{f}_iVV @VV\widetilde{f}_iV
\\
(\widetilde{\nu},\widetilde{J}) @>>\delta >
(\bar{\widetilde{\nu}},\bar{\widetilde{J}}),(\widetilde{\bar{\nu}},\widetilde{\bar{J}})
\end{CD}
\end{gather*}
Here we denote as $\delta:(\widetilde{\nu},\widetilde{J})\mapsto (\bar{\widetilde{\nu}},\bar{\widetilde{J}})$ and
$\widetilde{f}_i:(\bar{\nu},\bar{J})\mapsto (\widetilde{\bar{\nu}},\widetilde{\bar{J}})$.
We denote the smallest rigging associated with length $k$ rows of $\nu^{(i)}$ by $x_k$ and express such string as
$(k,x_k)$.
For later purposes, let us prepare several notation for the specif\/ic strings as follows:
\begin{itemize}\itemsep=0pt
\item
$(\ell,x_\ell)$ is the string of $(\nu,J)$ that is selected by $\widetilde{f}_i$,
\item
$(\bar{\ell},x_{\bar{\ell}})$ is the string of $(\bar{\nu},\bar{J})$ that is selected by $\widetilde{f}_i$,
\item
$\ell^{(a)}$ and $\ell_{(a)}$ are the lengths of the singular strings of $(\nu,J)^{(a)}$ which are selected by~$\delta$
under the condition $\ell^{(a)}\leq\ell_{(a)}$,
\item
$\widetilde{\ell}^{(a)}$ and $\widetilde{\ell}_{(a)}$ are the lengths of the singular strings of
$(\widetilde{\nu},\widetilde{J})^{(a)}$ which are selected by~$\delta$ under the condition
$\widetilde{\ell}^{(a)}\leq\widetilde{\ell}_{(a)}$.
\end{itemize}

Then our main claims to be proved are the following two assertions which are proved in Sections~\ref{sec:main},~\ref{sec:main2} and~\ref{sec:main2_e}.

\begin{Proposition}
\label{th:core1}
Let $B=B^{1,1}\otimes \bar{B}$ and $i\in I_0$.
Then we have the following commutative diagrams
\begin{gather*}
\begin{CD}
\mathrm{RC}(B) @>\delta >> \mathrm{RC}(\bar{B})
\\
@V\widetilde{e}_i VV @VV\widetilde{e}_i V
\\
\mathrm{RC}(B) @>>\delta > \mathrm{RC}(\bar{B})
\end{CD}
\qquad
\begin{CD}
\mathrm{RC}(B) @>\delta >> \mathrm{RC}(\bar{B})
\\
@V\widetilde{f}_i VV @VV\widetilde{f}_i V
\\
\mathrm{RC}(B) @>>\delta > \mathrm{RC}(\bar{B})
\end{CD}
\end{gather*}
if all $\widetilde{e}_i$ and $\widetilde{f}_i$ in the above diagrams are defined.
\end{Proposition}

\begin{Proposition}
\label{th:core2}
Let $B=B^{1,1}\otimes \bar{B}$, $i\in I_0$, $(\nu,J)\in \mathrm{RC}(B)$ and $b=b'\otimes\bar{b}\in B^{1,1}\otimes
\bar{B}$.
Suppose that $\Phi:(\nu,J)\longmapsto b$.
Then we have the following properties:
\begin{enumerate}\itemsep=0pt
\item[$(1)$] $\widetilde{f}_i(\nu,J)\neq 0$ and $\widetilde{f}_i(\bar{\nu},\bar{J})=0$ $\Longrightarrow$
$\widetilde{f}_i(b)\neq 0$, $\widetilde{f}_i(\bar{b})=0$ and
$\Phi(\widetilde{f}_i(\nu,J))=\widetilde{f}_i(\Phi(\nu,J))$,
\item[$(2)$] $\widetilde{f}_i(b)\neq 0$ and $\widetilde{f}_i(\bar{b})=0$ $\Longrightarrow$ $\widetilde{f}_i(\nu,J)\neq 0$,
$\widetilde{f}_i(\bar{\nu},\bar{J})=0$ and $\Phi^{-1}(\widetilde{f}_i(b))=\widetilde{f}_i(\Phi^{-1}(b))$,
\item[$(3)$] $\widetilde{f}_i(\nu,J)= 0$ and $\widetilde{f}_i(\bar{\nu},\bar{J})\neq 0$ $\Longrightarrow$
$\widetilde{f}_i(b)=0$, $\widetilde{f}_i(\bar{b})\neq 0$ and
$\Phi(\widetilde{f}_i(\nu,J))=\widetilde{f}_i(\Phi(\nu,J))$,
\item[$(4)$] $\widetilde{f}_i(b)=0$ and $\widetilde{f}_i(\bar{b})\neq 0$ $\Longrightarrow$ $\widetilde{f}_i(\nu,J)=0$,
$\widetilde{f}_i(\bar{\nu},\bar{J})\neq 0$ and $\Phi^{-1}(\widetilde{f}_i(b))=\widetilde{f}_i(\Phi^{-1}(b))$,
\item[$(5)$] $\widetilde{e}_i(\nu,J)\neq 0$ and $\widetilde{e}_i(\bar{\nu},\bar{J})=0$ $\Longrightarrow$
$\widetilde{e}_i(b)\neq 0$, $\widetilde{e}_i(\bar{b})=0$ and
$\Phi(\widetilde{e}_i(\nu,J))=\widetilde{e}_i(\Phi(\nu,J))$,
\item[$(6)$] $\widetilde{e}_i(b)\neq 0$ and $\widetilde{e}_i(\bar{b})=0$ $\Longrightarrow$ $\widetilde{e}_i(\nu,J)\neq 0$,
$\widetilde{e}_i(\bar{\nu},\bar{J})=0$ and $\Phi^{-1}(\widetilde{e}_i(b))=\widetilde{e}_i(\Phi^{-1}(b))$,
\item[$(7)$] the situation $\widetilde{e}_i(\nu,J)=0$ and $\widetilde{e}_i(\bar{\nu},\bar{J})\neq 0$ cannot happen,
\item[$(8)$] the situation $\widetilde{e}_i(b)=0$ and $\widetilde{e}_i(\bar{b})\neq 0$ cannot happen.
\end{enumerate}
Here from $(1)$ to $(4)$ we assume the commutativity of $\widetilde{f}_i$ with~$\Phi$ for $\bar{B}$ and from $(5)$ to $(8)$ we
assume the commutativity of $\widetilde{e}_i$ and $\widetilde{f}_i$ with~$\Phi$ for $\bar{B}$.
\end{Proposition}

Let us show how these properties lead to the proof of Theorem~\ref{th:main}.
In the proof, we will use the diagram of the following kind as in~\cite{KSS:2002}:
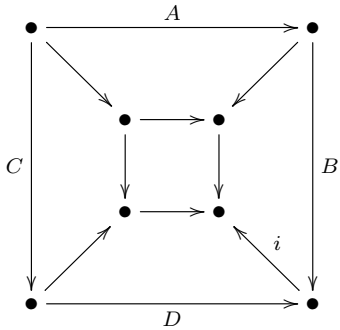
\begin{gather*}
\xymatrix{ {\bullet} \ar[rrr]^{A} \ar[ddd]_{C} \ar[dr] & & & {\bullet} \ar[ddd]^{B} \ar[dl]
\\
& {\bullet} \ar[r] \ar[d] & {\bullet} \ar[d] &
\\
& {\bullet} \ar[r] & {\bullet} &
\\
{\bullet} \ar[rrr]_{D} \ar[ur] & & & {\bullet} \ar[ul]_{i} }
\end{gather*}
We regard this as a~cube with front face given by the large square.
Suppose that the square diagrams given by the faces of the cube except for the front face commute and $i$ is the
injective map.
Then the front face also commutes since we have
\begin{gather*}
i\circ B\circ A=i\circ D\circ C
\end{gather*}
by a~diagram chasing.

\begin{proof}
[Proof of Theorem~\ref{th:main}] First we prove the $\widetilde{f}_i$ case.
Then we can prove the $\widetilde{e}_i$ case by a~parallel argument.

We follow the decomposition of~$\Phi$ as described in Proposition~\ref{prop:RC_decomp}.
The simplest case $B^{1,1}$ can be shown directly.
We give a~list of the explicit rigged conf\/igurations in this case:
\begin{center}
\begin{tabular}
{|l|l|l|}
\hline
$b\in B^{1,1}$&$(\nu,J)^{(a)}$ of $\Phi^{-1}(b)$&values for~$a$
\\
\hline
\hline
$\Yvcentermath1\young(i)$\rule{0pt}{13pt} & $\{(1,0)\}$ & $a<i-1$,
\\
&$\{(1,-1)\}$&$a=i-1$
\\
&$\varnothing$ & $i\leq a$.
\\
\hline
$\Yvcentermath1\young(\mn)$\rule{0pt}{13pt} & $\{(1,0)\}$ & $a\leq n-2$,
\\
&$\varnothing$&$a=n-1$
\\
&$\{(1,-1)\}$&$a=n$
\\
\hline
$\Yvcentermath1\young(\minusi)$\rule{0pt}{13pt} & $\{(1,0)\}$ & $a\leq i-2$,
\\
($i<n$)&$\{(1,1)\}$&$a=i-1$
\\
&$\{(1,-1),(1,-1)\}$&$a=i$
\\
&$\{(1,0),(1,0)\}$ & $i+1\leq a\leq n-2$,
\\
&$\{(1,0)\}$ & $a= n-1,n$
\\
\hline
\end{tabular}
\end{center}
By using the crystal graph of $B^{1,1}$ presented in Section~\ref{subsec:tableaux} we can verify the commutativity of
$\widetilde{f}_i$ and~$\Phi$.\footnote{In $B^{1,1}$ case~$\Phi$ coincides with~$\delta$.
In view of the crystal graph of $B^{1,1}$ and the def\/inition of $\widetilde{f}_i$ on the rigged conf\/igurations we see
a~connection between the crystal graph and~$\delta$ in the case $B^{1,1}$.
However this is only the tip of an iceberg.
Indeed, in $E^{(1)}_6$ case~\cite{OSano},~$\delta$ is described by the crystal graph of $B^{1,1}$ as in the present case
and the same algorithm can be applied for general tensor products of the form $(B^{1,1})^{L}$.}

Let us consider the case $B=B^{1,1}\otimes\bar{B}$.
This case corresponds to an addition of a~new tensor factor.
For this we use the following diagram:
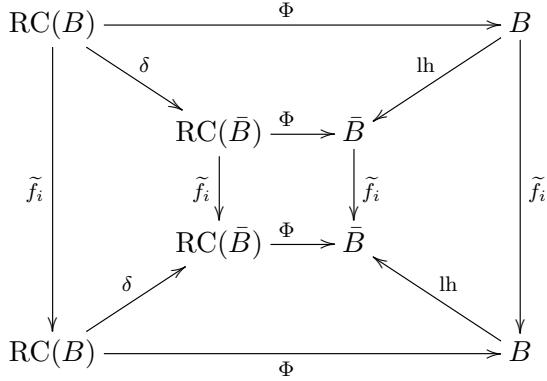
\begin{gather*}
\xymatrix{ {\mathrm{RC}(B)} \ar[rrr]^{\Phi} \ar[ddd]_{\widetilde{f}_i} \ar[dr] ^{\delta}& & & {B}
\ar[ddd]^{\widetilde{f}_i} \ar[dl]_{\mathrm{lh}}
\\
& {\mathrm{RC}(\bar{B})} \ar[r]^{\Phi} \ar[d]_{\widetilde{f}_i} & {\bar{B}} \ar[d] ^{\widetilde{f}_i}&{\hspace{20mm}}
\\
& {\mathrm{RC}(\bar{B})} \ar[r]^{\Phi} & {\bar{B}} &
\\
{\mathrm{RC}(B)} \ar[rrr]_{\Phi} \ar[ur] ^{\delta}& & &{B} \ar[ul]_{\mathrm{lh}} }
\end{gather*}
Suppose that all $\widetilde{f}_i$ which appear in the above diagram are def\/ined.
The top face and the bottom face commute by Proposition~\ref{prop:RC_decomp}.
The left face commutes by Proposition~\ref{th:core1}.
The right face commutes by def\/inition of $\mathrm{lh}$ on tensor products of crystals.
The back face commutes by the assumption.
Finally the map $\mathrm{lh}$ is injective if we consider the weight of the element of $B^{1,1}$ subtracted by the
operation.
Hence the front face commutes.
This proves the commutativity of~$\Phi$ and $\widetilde{f}_i$ when both $\widetilde{f}_i$ before or after~$\delta$ (or
$\mathrm{lh}$) are def\/ined.
On the other hand we know from Proposition~\ref{th:core2} that~$\Phi$ and $\widetilde{f}_i$ also commute even if one of
$\widetilde{f}_i$ before or after~$\delta$ (or $\mathrm{lh}$) is undef\/ined.
By exclusion we have the commutativity for the case when both $\widetilde{f}_i$ before and after~$\delta$ (or
$\mathrm{lh}$) are undef\/ined.

Let us consider the the case $B=B^{r,1}\otimes\bar{B}$ with $r>1$.
Let $B'=B^{1,1}\otimes B^{r-1,1}\otimes\bar{B}$.
Consider the following diagram:
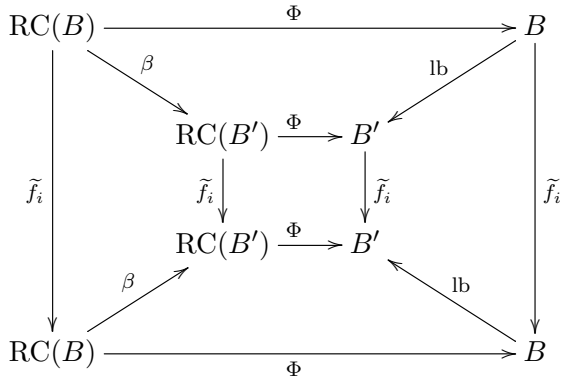
\begin{gather*}
\xymatrix{ {\mathrm{RC}(B)} \ar[rrr]^{\Phi} \ar[ddd]_{\widetilde{f}_i} \ar[dr] ^{\beta}& & & {B}
\ar[ddd]^{\widetilde{f}_i} \ar[dl]_{\mathrm{lb}}
\\
& {\mathrm{RC}(B')} \ar[r]^{\Phi} \ar[d]_{\widetilde{f}_i} & {B'} \ar[d] ^{\widetilde{f}_i}&{\hspace{20mm}}
\\
& {\mathrm{RC}(B')} \ar[r]^{\Phi} & {B'} &
\\
{\mathrm{RC}(B)} \ar[rrr]_{\Phi} \ar[ur] ^{\beta}& & &{B} \ar[ul]_{\mathrm{lb}} }
\end{gather*}

\looseness=1
The top face and the bottom face commute by Proposition~\ref{prop:RC_decomp}.
The left face commutes by Proposition~\ref{prop:lb_commute}.
The right face commutes by def\/inition of $\mathrm{lb}$ on tensor products of crystals.
The back face commutes by the assumption.
Since $\mathrm{lb}$ is injective, the front face commutes.

\pagebreak
Finally let us consider the the case $B=B^{r,s}\otimes\bar{B}$ with $s>1$.
Let $B'=B^{r,1}\otimes B^{r,s-1}\otimes\bar{B}$.
Consider the following diagram:
\begin{gather*}
\xymatrix{ {\mathrm{RC}(B)} \ar[rrr]^{\Phi} \ar[ddd]_{\widetilde{f}_i} \ar[dr] ^{\gamma}& & & {B}
\ar[ddd]^{\widetilde{f}_i} \ar[dl]_{\mathrm{ls}}
\\
& {\mathrm{RC}(B')} \ar[r]^{\Phi} \ar[d]_{\widetilde{f}_i} & {B'} \ar[d] ^{\widetilde{f}_i}&{\hspace{20mm}}
\\
& {\mathrm{RC}(B')} \ar[r]^{\Phi} & {B'} &
\\
{\mathrm{RC}(B)} \ar[rrr]_{\Phi} \ar[ur] ^{\gamma}& & &{B} \ar[ul]_{\mathrm{ls}} }
\end{gather*}
The top face and the bottom face commute by Proposition~\ref{prop:RC_decomp}.
The left face commutes by Proposition~\ref{prop:ls_commute}.
The right face commutes by the def\/inition of $\mathrm{ls}$ on tensor products of crystals.
The back face commutes by the assumption.
Since $\mathrm{ls}$ is injective, the front face commutes.
This completes the proof of Theorem~\ref{th:main}.
\end{proof}

\subsection{Spin cases}\label{sec:spin}

Let us consider the arguments which are needed to treat the spin cases.
Here we will concentrate on the case $B^{n,s}$.
The other case $B^{n-1,s}$ can be obtained by replacing the Dynkin nodes $n$ and $n-1$ in the following description.
To begin with, following~\cite{S:2005}, we introduce the three sets $\hat{B}^{n-1,1}$, $\hat{B}^{n,1}$ and
$\hat{\bar{B}}^{n,1}$ which are generated by
\begin{gather*}
\hat{u}_{n-1}=
\begin{array}{|c|}
\hline
1
\\
\hline
2
\\
\hline
\vdots
\\
\hline
n-1
\\
\hline
\end{array}
,
\qquad
\hat{u}_{n}=
\begin{array}{|c|}
\hline
1
\\
\hline
2
\\
\hline
\vdots
\\
\hline
n-1
\\
\hline
n
\\
\hline
\end{array}
,
\qquad
\hat{\bar{u}}_{n}=
\begin{array}{|c|}
\hline
1
\\
\hline
2
\\
\hline
\vdots
\\
\hline
n-1
\\
\hline
\mn
\\
\hline
\end{array}
\end{gather*}
Here we regard them as the usual single columns and apply the Kashiwara operators $\widetilde{f}_i$ $(i\in I_0)$ as in
the case for non-spin KN tableaux to obtain the whole elements of $\hat{B}^{n-1,1}$, $\hat{B}^{n,1}$ and~$\hat{\bar{B}}^{n,1}$.
On these sets, we def\/ine the operation left-box analogously to Def\/inition~\ref{def:operations_on_tableaux}.
To be more explicit, we have $\mathrm{lb}:\hat{B}^{n,1}\longrightarrow B^{1,1}\otimes\hat{B}^{n-1,1}$ and
$\mathrm{lb}:\hat{B}^{n-1,1}\longrightarrow B^{1,1}\otimes B^{n-2,1}$.

Now we describe the algorithm for the rigged conf\/iguration bijection for the spin cases $B^{n-1,l}$ and $B^{n,l}$
following~\cite{S:2005}.
To begin with, we introduce the embedding of the rigged conf\/igurations
\begin{gather*}
\mathrm{emb}:\{\mu,(\nu,J)\}\longmapsto\{\mu',(\nu',J')\},
\end{gather*}
where the quantities with primes means that the everything is doubled as $\mu_i'^{(a)}=2\mu^{(a)}_i$,
$\nu'^{(a)}_i=2\nu^{(a)}_i$ and $J'^{(a)}_i=2J^{(a)}_i$.
Similarly, if the rigged conf\/iguration $\{\mu,(\nu,J)\}$ is composed of even integers, we can def\/ine $\mathrm{emb}^{-1}$
by dividing all its components by 2.
Then the operation to obtain the leftmost column of an element of $B^{n,l}$ is
\begin{gather*}
\mathrm{emb}^{-1}\circ\delta^{(1)}_1\cdots\delta^{(n-2)}_1
\hat{\delta}^{(n-1)}_1\hat{\delta}^{(n)}_{2}\circ\mathrm{emb}\circ\gamma(\nu,J).
\end{gather*}
Here $\gamma$ is the same as Def\/inition~\ref{def:RC_operations}; $\gamma$ replaces one of length $l$ rows of $\mu^{(n)}$
by two rows of length $l-1$ and $1$.
The new operations $\hat{\delta}^{(n)}_{2}$ and $\hat{\delta}^{(n-1)}_1$ are def\/ined as follows.
Let us denote by $(\nu',J')=\mathrm{emb}\circ\gamma(\nu,J)$.
\begin{itemize}\itemsep=0pt
\item
$\hat{\delta}^{(n)}_{2}:(\nu',J')\longmapsto\{(\nu'',J''),k\}$ is def\/ined as follows.
We look for the shortest singular string of $(\nu',J')^{(n)}$ which is longer than or equal to $2$.
If there is no such string, set $\ell^{(n)}=\infty$, $k=n$ and stop.
Otherwise we set $\ell_{(n-1)}$ to be the length of the selected string of $(\nu',J')^{(n)}$ and do the same procedure
as in the usual $\delta^{(a)}_l$.\\
The def\/inition of the new rigged conf\/iguration $(\nu'',J'')$ is the same as in the usual $\delta^{(a)}_l$ except for the
following point.
We replace one of the length $2$ rows of $\mu'^{(n)}$ by length $1$ row and add a~length one row to~$\mu'^{(n-1)}$.

\item
$\hat{\delta}^{(n-1)}_{1}:(\nu',J')\longmapsto\{(\nu'',J''),k\}$ is def\/ined as follows.
Set $\ell^{(n-2)}=1$.
Then, as in the usual $\delta^{(a)}_l$ we begin to determine $\ell^{(n-1)}$ and $\ell^{(n)}$ and proceed similarly to
the usual case.
The new rigged conf\/iguration $(\nu'',J'')$ is def\/ined similarly as in the usual case except for the following point.
We remove a~length one row from each of $\mu'^{(n-1)}$ and $\mu'^{(n)}$ and add a~length one row to $\mu'^{(n-2)}$.
\end{itemize}

We regard the output of the sequence of~$\delta$'s as the element of $\hat{B}^{n,1}$.
Corresponding to the f\/inal $\mathrm{emb}^{-1}$, we halve the width of the column with keeping all its contents.
We regard the f\/inal output as the element of $B^{n,1}$.

\begin{Remark}
In the above algorithm, some readers might feel that the change in $\mu$ is peculiar.
To get a~better understanding of this point, let us take a~look at the following computations
\begin{gather*}
2\bar{\Lambda}_{n}=\epsilon_1+\dots +\epsilon_{n-2}+\epsilon_{n-1}+\epsilon_n,
\\
\bar{\Lambda}_{n}+\bar{\Lambda}_{n-1}=\epsilon_1+\dots +\epsilon_{n-2}+\epsilon_{n-1},
\\
\bar{\Lambda}_{n-2}=\epsilon_1+\dots +\epsilon_{n-2}.
\end{gather*}
Let us compare with the special case $\nu^{(a)}=\varnothing$ for all $a\in I_0$ with weight $2\bar{\Lambda}_{n}$
(see~\eqref{wt_RC_2}).
In this case we should obtain a~column f\/illed by $1,2,\ldots,n$ from top to bottom.
Note that each letter $i$ in this column corresponds to $\epsilon_i$.
\end{Remark}

For the proof of Theorem~\ref{th:main}, we have to modify $\beta$ in Def\/inition~\ref{def:RC_operations} as follows.
$\hat{\beta}^{(n)}$ replaces a~length two row of $\mu^{(n)}$ by a~length one row, adds a~length one row to each of
$\mu^{(1)}$ and $\mu^{(n-1)}$ and adds a~length one singular string to each of $(\nu,J)^{(a)}$ for $a\leq n-1$.
$\hat{\beta}^{(n-1)}$ removes a~length one row from each of $\mu^{(n-1)}$ and $\mu^{(n)}$, adds a~length one row to each
of $\mu^{(1)}$ and $\mu^{(n-2)}$ and adds a~length one singular string to each of $(\nu,J)^{(a)}$ for $a\leq n-2$.
Note that we have the following property.
\begin{Proposition}
On the rigged configurations, we have
\begin{gather*}
\mathrm{emb}\circ\widetilde{e}_i=\widetilde{e}_i^2\circ\mathrm{emb},
\qquad
\mathrm{emb}\circ\widetilde{f}_i=\widetilde{f}_i^2\circ\mathrm{emb}.
\end{gather*}
\end{Proposition}
\begin{proof}
We prove the $\widetilde{e}_i$ part.
Let us consider the case $\widetilde{e}_i(\nu,J)\neq 0$.
Suppose that $\widetilde{e}_i$ acts on the string $(\ell,x_\ell)$ of $(\nu,J)^{(i)}$.
Then $\widetilde{e}_i$ makes the string $(2\ell,2x_\ell)$ of $\mathrm{emb}(\nu,J)^{(i)}$ into $(2\ell-1,2x_\ell+1)$.
Since $x_\ell<0$, we have $2x_\ell+1<0$ and thus $\widetilde{e}_i^2\circ\mathrm{emb}(\nu,J)\neq 0$.
We show that the second $\widetilde{e}_i$ acts on the string $(2\ell-1,2x_\ell+1)$.
\begin{enumerate}\itemsep=0pt
\item[(a)] Suppose that there is a~string $(k,x_k)$ with $k<\ell$ in $(\nu,J)^{(i)}$.
By the minimality of $\ell$ we have $x_k>x_\ell$.
Therefore we have $2x_k>2x_\ell+1$.
\item[(b)] Suppose that there is a~string $(k,x_k)$ with $\ell\leq k$ in $(\nu,J)^{(i)}$.
Then in the $\widetilde{e}_i\circ\mathrm{emb}(\nu,J)^{(i)}$, the string becomes $(2k,x_k+2)$.
By the minimality of $x_\ell$, we have $x_k+2>x_\ell+1$.
\end{enumerate}
Therefore the second $\widetilde{e}_i$ acts on the string $(2\ell-1,2x_\ell+1)$ to make $(2\ell-2,2x_\ell+2)$.
Thus $\mathrm{emb}\circ\widetilde{e}_i=\widetilde{e}_i^2\circ\mathrm{emb}$ in this case.
On the other hand, if $\widetilde{e}_i(\nu,J)=0$, we have $\widetilde{e}_i\circ\mathrm{emb}(\nu,J)=0$ by
Theorem~\ref{prop:phi_RC}.

Similarly we can show $\mathrm{emb}\circ\widetilde{f}_i=\widetilde{f}_i^2\circ\mathrm{emb}$.
\end{proof}

Then we can use the same arguments of the proof of Proposition~\ref{prop:lb_commute} to show that~$\hat{\beta}^{(n-1)}$
and~$\hat{\beta}^{(n)}$ commute with~$\widetilde{e}_i$ and~$\widetilde{f}_i$.
Thus by using the same arguments we can show Theorem~\ref{th:main} in this case.

\section{Proof of Proposition~\ref{th:core1}}\label{sec:main}

\subsection{Outline}

We concentrate on the proof for $\widetilde{f}_i$ part since the other part $\widetilde{e}_i$ follows from the former case.
We divide the proof into the following f\/ive cases which exhaust all the possibilities:
\begin{alignat*}{3}
& \text{Case A:}\quad && \ell+1<\ell^{(i-1)},&
\\
& \text{Case B:}\quad && \ell^{(i-1)}\leq\ell+1<\ell^{(i)}=\infty,&
\\
& \text{Case C:}\quad && \ell^{(i-1)}\leq \ell+1\leq \ell^{(i)}<\infty,&
\\
& \text{Case D:}\quad && \ell^{(i)}=\ell,&
\\
& \text{Case E:}\quad && \ell^{(i)}<\ell.&
\end{alignat*}

The actual proof for each case is based on the following behaviors of the vacancy numbers induced by~$\delta$.

\begin{Lemma}\label{lem:vacancy_delta}
For $1\leq i\leq n$ we have
\begin{gather*}
P^{(i)}_{\ell^{(i)}-1}(\bar{\nu})=
\left\{\begin{array}{@{}lll}
P^{(i)}_{\ell^{(i)}-1}(\nu)& (\text{if}\ \ \ell^{(i-1)}=\ell^{(i)}),&\mathrm{(I)}
\vspace{1mm}\\
P^{(i)}_{\ell^{(i)}-1}(\nu)-1 & (\text{if}\ \ \ell^{(i-1)}<\ell^{(i)}).&\mathrm{(II)}
\end{array}\right.
\end{gather*}
Here if $i=n$ we understand $\ell^{(i-1)}$ as $\ell^{(n-2)}$.
For $1\leq i\leq n-3$ we have
\begin{gather*}
P^{(i)}_{\ell_{(i)}-1}(\bar{\nu})=
\left\{
\begin{array}{@{}lll}
P^{(i)}_{\ell_{(i)}-1}(\nu)& (\text{if}\ \  \ell^{(i-1)}=\ell^{(i)}=\dots =\ell_{(i)}),&\mathrm{(III)}
\vspace{1mm}\\
P^{(i)}_{\ell_{(i)}-1}(\nu)-1& (\text{if}\ \ \ell^{(i-1)}<\ell^{(i)}=\dots =\ell_{(i)}),&\mathrm{(IV)}
\vspace{1mm}\\
P^{(i)}_{\ell_{(i)}-1}(\nu)+1& (\text{if} \ \ \ell^{(i)}<\ell^{(i+1)}=\dots =\ell_{(i)}),&\mathrm{(V)}
\vspace{1mm}\\
P^{(i)}_{\ell_{(i)}-1}(\nu)& (\text{if}\ \ \ell^{(i+1)}<\ell_{(i+1)}=\ell_{(i)}),&\mathrm{(VI)}
\vspace{1mm}\\
P^{(i)}_{\ell_{(i)}-1}(\nu)-1& (\text{if} \ \ \ell_{(i+1)}<\ell_{(i)}),&\mathrm{(VII)}
\end{array}\right.
\end{gather*}
and for $i=n-2$ we have
\begin{gather*}
P^{(n-2)}_{\ell_{(n-2)}-1}(\bar{\nu}) =
\left\{
\begin{array}{@{}lll}
P^{(n-2)}_{\ell_{(n-2)}-1}(\nu)& (\text{if}\ \  \ell^{(n-3)}=\ell^{(n-2)}=\ell^{(n-1)}=\ell^{(n)}=\ell_{(n-2)}),
&\mathrm{(VIII)}
\vspace{1mm}\\
P^{(n-2)}_{\ell_{(n-2)}-1}(\nu)-1& (\text{if}\ \  \ell^{(n-3)}<\ell^{(n-2)}=\ell^{(n-1)}=\ell^{(n)}=\ell_{(n-2)}),
&\mathrm{(IX)}
\vspace{1mm}\\
P^{(n-2)}_{\ell_{(n-2)}-1}(\nu)+1& (\text{if}\ \ \ell^{(n-2)}<\ell^{(n-1)}=\ell^{(n)}=\ell_{(n-2)}), &\mathrm{(X)}
\vspace{1mm}\\
P^{(n-2)}_{\ell_{(n-2)}-1}(\nu)& (\text{if}\ \ \min\{\ell^{(n-1)},\ell^{(n)}\}<\ell_{(n-1)}=\ell_{(n-2)}), &\mathrm{(XI)}
\vspace{1mm}\\
P^{(n-2)}_{\ell_{(n-2)}-1}(\nu)-1 & (\text{if}\ \ \ell_{(n-1)}<\ell_{(n-2)}).&\mathrm{(XII)}
\end{array}\right.
\end{gather*}
\end{Lemma}

During the proof, the cases (VIII) to (XII) can be treated as the special case of the cases (III) to (VII).
Indeed, we can assume that $\ell^{(n-1)}\leq\ell^{(n)}$ without loss of the generality.
Then we have $\ell_{(n-1)}=\ell^{(n)}$ and regard $\ell^{(n-1)}$ and $\ell_{(n-1)}$ as lengths of strings of the
concatenation $\left((\nu^{(n-1)},\nu^{(n)}),(J^{(n-1)},J^{(n)})\right)$ of $(\nu,J)^{(n-1)}$ and $(\nu,J)^{(n)}$.

For the reader's convenience, we give a~summary of all major subcases here.
\begin{enumerate}\itemsep=0pt
\item[\fbox{A}] $\ell+1<\ell^{(i-1)}$.
\item[\fbox{B}] $\ell^{(i-1)}\leq\ell+1<\ell^{(i)}=\infty$,
\begin{enumerate}\itemsep=0pt
\item[(1)] $\widetilde{f}_i$ creates a~non-singular string, \item[(2)] $\widetilde{f}_i$ creates a~singular string.
\end{enumerate}
\item[\fbox{C}] $\ell^{(i-1)}\leq \ell+1\leq \ell^{(i)}<\infty$,
\begin{enumerate}\itemsep=0pt
\item[(1)] $m^{(i)}_{\ell+1}(\nu)>0$ and $\widetilde{f}_i$ creates a~singular string, or
$\ell+1=\ell^{(i)}$ and $\widetilde{f}_i$ creates a~non-singular string,
\item[(2)] $m^{(i)}_{\ell+1}(\nu)=0$ and $\widetilde{f}_i$ creates a~singular string, \item[(3)] $\ell+1<\ell^{(i)}$ and
$\widetilde{f}_i$ creates a~non-singular string.
\end{enumerate}
\item[\fbox{D}] $\ell^{(i)}=\ell$,
\begin{enumerate}\itemsep=0pt
\item[(1)] $\ell=\ell^{(i)}<\ell_{(i)}$ and $m^{(i)}_\ell(\nu)>1$, \item[(2)] $\ell=\ell^{(i)}=\ell_{(i)}$ and
$m^{(i)}_\ell(\nu)>2$, \item[(3)] $\ell=\ell^{(i)}<\ell_{(i)}$ and $m^{(i)}_\ell(\nu)=1$, \item[(4)]
$\ell=\ell^{(i)}=\ell_{(i)}$ and $m^{(i)}_\ell(\nu)=2$.
\end{enumerate}
\item[\fbox{E}] $\ell^{(i)}<\ell$,
\begin{enumerate}\itemsep=0pt
\item[(1)] $\ell^{(i)}<\ell$ and $\ell+1<\ell_{(i+1)}$, \item[(2)] $\ell^{(i)}<\ell$ and
$\ell_{(i+1)}\leq\ell+1<\ell_{(i)}=\infty$, \item[(3)] $\ell^{(i)}<\ell$ and
$\ell_{(i+1)}\leq\ell+1\leq\ell_{(i)}<\infty$, \item[(4)] $\ell^{(i)}<\ell$ and $\ell_{(i)}=\ell$, \item[(5)]
$\ell^{(i)}<\ell$ and $\ell_{(i)}<\ell$.
\end{enumerate}
\end{enumerate}

Finally, let us prepare the following lemma which gives an useful criterion.

\begin{Lemma}
\label{lem:f_singular}
Suppose that $\widetilde{f}_i$ acts on a~string $(\ell,x_\ell)$ of $(\nu,J)^{(i)}$ and creates the new string
$(\ell+1,x_\ell-1)$.
\begin{enumerate}\itemsep=0pt
\item[$(1)$] The string $(\ell+1,x_\ell-1)$ is singular if and only if $P^{(i)}_{\ell+1}(\nu)=x_\ell+1$.
\item[$(2)$] Suppose that $\widetilde{f}_i$ creates a~singular string.
If there is a~string $(l,x_l)$ of $(\nu,J)^{(i)}$ that satisfies $\ell<l$, then we have $P^{(i)}_{\ell+1}(\nu)\leq x_l$.
\end{enumerate}
\end{Lemma}
\begin{proof}
(1) Suppose that the string $(\ell+1,x_\ell-1)$ is singular.
Then it satisf\/ies $P^{(i)}_{\ell+1}(\widetilde{\nu})=x_\ell-1$.
Since $\widetilde{f}_i$ adds the box to the $(\ell+1)$-th column of $\nu^{(i)}$, we have
$P^{(i)}_{\ell+1}(\widetilde{\nu})=P^{(i)}_{\ell+1}(\nu)-2$.
Thus we have $P^{(i)}_{\ell+1}(\nu)= P^{(i)}_{\ell+1}(\widetilde{\nu})+2=x_\ell+1$.
We can reverse the arguments to obtain the if part.

(2) If $x_\ell\geq x_l$, then $\widetilde{f}_i$ will act on $(l,x_l)$ instead of $(\ell,x_\ell)$ since $\ell<l$.
Thus we have $x_\ell<x_l$.
Then we have $P^{(i)}_{\ell+1}(\nu)=x_\ell+1\leq x_l$.
\end{proof}

\subsection{Proof for Case A}

The proof of this case depends on the following fundamental properties.

\begin{Proposition}
\label{lem:case_A}
If $\ell<\ell^{(i-1)}-1$, we have
\begin{gather*}
P^{(i)}_{\ell^{(i)}-1}(\bar{\nu})>x_\ell,%\label{eq:caseA}
\qquad
P^{(i)}_{\ell_{(i)}-1}(\bar{\nu})>x_\ell.%\label{eq:caseA2}
\end{gather*}
\end{Proposition}
\begin{proof}
The proof proceeds by case by case analysis according to the classif\/ication (I) to (VII) in
Lemma~\ref{lem:vacancy_delta}.
We remark that the condition $\ell<\ell^{(i-1)}-1$ automatically implies that $\ell<\ell_{(i+1)}-1$.

{\bf Case (I).} We have to show $P^{(i)}_{\ell^{(i)}-1}(\nu)>x_\ell$ by Lemma~\ref{lem:vacancy_delta}.
Suppose if possible that $P^{(i)}_{\ell^{(i)}-1}(\nu)\leq x_\ell$.
We have $P^{(i)}_{\ell^{(i)}}(\nu)\geq x_{\ell^{(i)}}>x_\ell$ since the string $(\ell^{(i)},x_{\ell^{(i)}})$ is not
selected by~$\widetilde{f}_i$ whereas $\ell<\ell^{(i)}$.
Let $(j,x_j)$ be the longest string of $\nu^{(i)}$ that is strictly shorter than $\ell^{(i)}$.
Since $m^{(i)}_\ell (\nu)>0$ we have $\ell\leq j$.
If $j=\ell^{(i)}-1$, we have $\ell<j$ (since $\ell<\ell^{(i-1)}-1\leq \ell^{(i)}-1$) and $x_j\leq P^{(i)}_{j}(\nu)\leq
x_\ell$, which contradict the def\/inition of $(\ell,x_\ell)$.
If $j<\ell^{(i)}-1$, combining $P^{(i)}_{\ell^{(i)}}(\nu)>x_\ell$ with the assumption $P^{(i)}_{\ell^{(i)}-1}(\nu)\leq
x_\ell$, we have $P^{(i)}_{\ell^{(i)}-1}(\nu)<P^{(i)}_{\ell^{(i)}}(\nu)$.
Then by the convexity of $P^{(i)}_k(\nu)$ between $j\leq k\leq\ell^{(i)}$, we have $x_j\leq
P^{(i)}_j(\nu)<P^{(i)}_{\ell^{(i)}-1}(\nu)\leq x_\ell$, which contradicts the minimality of~$x_\ell$.
Thus we have proved $P^{(i)}_{\ell^{(i)}-1}(\nu)>x_\ell$.

{\bf Case (II).} We have to show $P^{(i)}_{\ell^{(i)}-1}(\nu)>x_\ell +1$ by Lemma~\ref{lem:vacancy_delta}.
Suppose if possible that we have $P^{(i)}_{\ell^{(i)}-1}(\nu)\leq x_\ell +1$.
Recall that we have $x_\ell<x_{\ell^{(i)}}\leq P^{(i)}_{\ell^{(i)}}(\nu)$ by the minimality of $x_\ell$ and $\ell
<\ell^{(i)}$.
Thus we have
\begin{gather}
\label{eq:caseA(II)_1}
P^{(i)}_{\ell^{(i)}-1}(\nu)\leq P^{(i)}_{\ell^{(i)}}(\nu).
\end{gather}
Let $j$ be the length of the longest string of $\nu^{(i)}$ under the condition $j<\ell^{(i)}$.
Since $m^{(i)}_\ell (\nu)>0$ we have $\ell\leq j$.

Let us show that $j<\ell^{(i-1)}$.
Suppose if possible that $\ell^{(i-1)}\leq j$.
In particular, we have $\ell<j$ by $\ell<\ell^{(i-1)}-1$.
Then we shall show that $P^{(i)}_{j}(\nu)\leq x_\ell +1$.
\begin{enumerate}\itemsep=0pt
\item[(i)] The case $j=\ell^{(i)}-1$.
Then $P^{(i)}_{j}(\nu)\leq x_\ell +1$ follows from the assumption.
\item[(ii)] The case $j<\ell^{(i)}-1$.
We can apply the convexity relation of $P^{(i)}_k(\nu)$ between $j\leq k\leq\ell^{(i)}$ with~\eqref{eq:caseA(II)_1} to
deduce that $P^{(i)}_{j}(\nu)\leq P^{(i)}_{\ell^{(i)}-1}(\nu)\leq x_\ell+1$.
\end{enumerate}
On the other hand, since $\ell<j$, we have the condition $x_\ell<x_j\leq P^{(i)}_{j}(\nu)$.
Combining this with the above result $P^{(i)}_{j}(\nu)\leq x_\ell +1$ we deduce that $P^{(i)}_{j}(\nu)=x_\ell +1$ and
thus $x_j=P^{(i)}_{j}(\nu)$.
But then the string $(j,x_j)$ is the singular string whose length satisf\/ies $\ell^{(i-1)}\leq j<\ell^{(i)}$, which is in
contradiction to the def\/inition of $\ell^{(i)}$.
Hence we have proved $j<\ell^{(i-1)}$.

Then we can apply the convexity relation of $P^{(i)}_k(\nu)$ between $\ell^{(i-1)}\leq k\leq\ell^{(i)}$
with~\eqref{eq:caseA(II)_1} to deduce that $P^{(i)}_{\ell^{(i-1)}}(\nu)\leq P^{(i)}_{\ell^{(i)}-1}(\nu)$.
Since $m^{(i-1)}_{\ell^{(i-1)}}(\nu)>0$ the convexity relation of $P^{(i)}_k(\nu)$ between $\ell^{(i-1)}-1\leq
k\leq\ell^{(i-1)}+1$ is strict.
Thus we have
\begin{gather*}
P^{(i)}_j(\nu)<P^{(i)}_{\ell^{(i-1)}}(\nu) \leq P^{(i)}_{\ell^{(i)}-1}(\nu)\leq x_\ell+1,
\end{gather*}
in particular, $x_j\leq P^{(i)}_j(\nu)\leq x_\ell$.
If $\ell<j$, this contradicts the requirement $x_\ell <x_j$.
Hence we are left with the case $j=\ell$.
Since we are assuming $\ell<\ell^{(i-1)}-1$ we can ref\/ine the above relation as
\begin{gather*}
P^{(i)}_\ell(\nu)<P^{(i)}_{\ell^{(i-1)}-1}(\nu)<P^{(i)}_{\ell^{(i-1)}}(\nu) \leq P^{(i)}_{\ell^{(i)}-1}(\nu)\leq
x_\ell+1,
\end{gather*}
that is, $P^{(i)}_{\ell}(\nu)<x_\ell$.
But this contradicts the requirement that the string $(\ell,x_\ell)$ must satisfy $x_\ell\leq P^{(i)}_\ell(\nu)$.
Thus we have proved $P^{(i)}_{\ell^{(i)}-1}(\nu)>x_\ell+1$.

{\bf Case (III).} We have to show $P^{(i)}_{\ell_{(i)}-1}(\nu)>x_\ell$ by Lemma~\ref{lem:vacancy_delta}.
Since $\ell_{(i)}=\ell^{(i)}$ we can use the same arguments of~(I) to obtain the result.

{\bf Case (IV).} We have to show $P^{(i)}_{\ell_{(i)}-1}(\nu)>x_\ell+1$ by Lemma~\ref{lem:vacancy_delta}.
Since $\ell_{(i)}=\ell^{(i)}$ we can use the same arguments of (II) to obtain the result.

{\bf Case (V).} We have to show $P^{(i)}_{\ell_{(i)}-1}(\nu)\geq x_\ell$ by Lemma~\ref{lem:vacancy_delta}.
Suppose if possible that $P^{(i)}_{\ell_{(i)}-1}(\nu)<x_\ell$.
Let $j$ be the length of the longest string of $\nu^{(i)}$ under the condition $j<\ell_{(i)}$.
Since $m^{(i)}_{\ell^{(i)}}(\nu)>0$ we have $\ell^{(i)}\leq j$.
If $j=\ell_{(i)}-1$ we have $x_j\leq P^{(i)}_j(\nu)<x_\ell$ which is in contradiction to the minimality of $x_\ell$.
Thus assume that $j<\ell_{(i)}-1$.
Then by the convexity relation of $P^{(i)}_k(\nu)$ between $j\leq k\leq\ell_{(i)}$ we have
$P^{(i)}_{\ell_{(i)}-1}(\nu)\geq\min\big\{P^{(i)}_j(\nu),P^{(i)}_{\ell_{(i)}}(\nu)\big\}$.
However this implies that $x_\ell >x_j$ or $x_\ell >x_{\ell_{(i)}}$ both of which are in contradiction to the minimality
of~$x_\ell$.
Thus we conclude that $P^{(i)}_{\ell_{(i)}-1}(\nu)\geq x_\ell$.

{\bf Case (VI).} We have to show $P^{(i)}_{\ell_{(i)}-1}(\nu)>x_\ell$ by Lemma~\ref{lem:vacancy_delta}.
Suppose if possible that $P^{(i)}_{\ell_{(i)}-1}(\nu)\leq x_\ell$.
Let $j$ be the length of the longest string of $\nu^{(i)}$ under the condition $j<\ell_{(i)}$.
Since $m^{(i)}_{\ell^{(i)}}(\nu)>0$ we have $\ell^{(i)}\leq j$.
In particular, we have $\ell<j$ by $\ell<\ell^{(i-1)}-1$.
If $j=\ell_{(i)}-1$ we have $x_j\leq x_\ell$ which contradicts the relation $x_j>x_\ell$ required by $\ell<j$.
So suppose that $j<\ell_{(i)}-1$.
Again we have $P^{(i)}_{\ell_{(i)}-1}(\nu)\geq\min\big\{P^{(i)}_j(\nu),P^{(i)}_{\ell_{(i)}}(\nu)\big\}$.
Then we have $x_\ell\geq x_j$ or $x_\ell\geq x_{\ell_{(i)}}$.
This is a~contradiction since $\ell<j$ and $\ell<\ell^{(i)}$.

{\bf Case (VII).} We have to show $P^{(i)}_{\ell_{(i)}-1}(\nu)>x_\ell+1$ by Lemma~\ref{lem:vacancy_delta}.
We can prove this by the same arguments of the case (II) if we replace $\ell^{(i)}$ and $\ell^{(i-1)}$ there by
$\ell_{(i)}$ and $\ell_{(i+1)}$ respectively.
\end{proof}

\begin{Proposition}
If $\ell<\ell^{(i-1)}-1$, we have the following identities:
\begin{gather*}
 (1)\quad \bar{\widetilde{\nu}}=\widetilde{\bar{\nu}}, \qquad (2) \quad \bar{\widetilde{J}}=\widetilde{\bar{J}}.
\end{gather*}
\end{Proposition}

\begin{proof}
(1) {\bf Step 1.} Let us consider the case $\ell^{(i-1)}=\infty$.
Then~$\delta$ will not choose any string from $(\nu,J)^{(i-1)}$ so that $\widetilde{f}_i$ chooses the same string before
and after~$\delta$.
Next, recall that $\widetilde{f}_i$ does not change coriggings of untouched strings.
Then we have $\widetilde{\ell}^{(a)}=\ell^{(a)}$ for all $a\leq i-1$.
In particular, we have $\widetilde{\ell}^{(i-1)}=\infty$.
Thus we have $\widetilde{\ell}^{(a)}=\ell^{(a)}=\infty$ for all $a\geq i-1$.

{\bf Step 2.} Based on the above observation, let us assume that $\ell^{(i-1)}<\infty$ in the sequel.
Let us show that~$\delta$ chooses the same strings before and after $\widetilde{f}_i$.
Again we have $\widetilde{\ell}^{(a)}=\ell^{(a)}$ for all $a\leq i-1$ since $\widetilde{f}_i$ does not change coriggings
of untouched strings.
Recall that we have $\ell +1<\ell^{(i-1)}=\widetilde{\ell}^{(i-1)}$ by the assumption.
Thus~$\delta$ cannot choose the string $(\ell+1,x_\ell-1)$ created by $\widetilde{f}_i$.
Thus we have $\widetilde{\ell}^{(i)}=\ell^{(i)}$ which implies that $\widetilde{\ell}^{(a)}=\ell^{(a)}$ for all $a\geq
i$ and $\widetilde{\ell}_{(a)}=\ell_{(a)}$ for all~$a$.

{\bf Step 3.} Let us show that $\widetilde{f}_i$ chooses the same string before and after~$\delta$.
Recall that~$\delta$ creates the strings $\big(\ell^{(i)}-1,P^{(i)}_{\ell^{(i)}-1}(\bar{\nu})\big)$ and
$\big(\ell_{(i)}-1,P^{(i)}_{\ell_{(i)}-1}(\bar{\nu})\big)$ of $(\bar{\nu},\bar{J})^{(i)}$.
Note that the string $(\ell,x_\ell)$ remains as it is in $(\bar{\nu},\bar{J})^{(i)}$ since~$\delta$ acting on $(\nu,J)$
does not touch the string by the assumption $\ell<\ell^{(i-1)}$.
By Proposition~\ref{lem:case_A} we have $P^{(i)}_{\ell^{(i)}-1}(\bar{\nu})>x_\ell$ and
$P^{(i)}_{\ell_{(i)}-1}(\bar{\nu})>x_\ell$.
Therefore $\widetilde{f}_i$ acts on the same string $(\ell,x_\ell)$ before and after the application of~$\delta$.

To summarize, we have $\bar{\widetilde{\nu}}=\widetilde{\bar{\nu}}$.

(2)~It is enough to consider the strings $\big(\ell^{(a)},P^{(a)}_{\ell^{(a)}}(\nu)\big)$ of $(\nu,J)^{(a)}$ for $a=i-1,i,i+1$,
since~$\widetilde{f}_i$ adds a~box to the same place before and after~$\delta$ and since~$\delta$ does not change
riggings for untouched strings.
For example, the string for $a=i-1$ case behaves as follows:
\begin{gather*}
\begin{array}{@{}lllll}
\big(\ell^{(i-1)},P^{(i-1)}_{\ell^{(i-1)}}(\nu)\big) &\overset{\widetilde{f}_i}{\longmapsto}&
\big(\ell^{(i-1)},P^{(i-1)}_{\ell^{(i-1)}}(\widetilde{\nu})\big) &\overset{\delta}{\longmapsto}&
\big(\ell^{(i-1)}-1,P^{(i-1)}_{\ell^{(i-1)}-1}(\bar{\widetilde{\nu}})\big),
\\
\big(\ell^{(i-1)},P^{(i-1)}_{\ell^{(i-1)}}(\nu)\big) &\overset{\delta}{\longmapsto}&
\big(\ell^{(i-1)}-1,P^{(i)}_{\ell^{(i-1)}-1}(\bar{\nu})\big) &\overset{\widetilde{f}_i}{\longmapsto}&
\big(\ell^{(i-1)}-1,P^{(i)}_{\ell^{(i-1)}-1}(\widetilde{\bar{\nu}})\big).
\end{array}
\end{gather*}
In all cases, the riggings are determined so that $\widetilde{f}_i$ and~$\delta$ create the singular strings.
Since we have $\bar{\widetilde{\nu}}=\widetilde{\bar{\nu}}$, we obtain
$P^{(i-1)}_{\ell^{(i-1)}-1}(\bar{\widetilde{\nu}})= P^{(i)}_{\ell^{(i-1)}-1}(\widetilde{\bar{\nu}})$.
The other cases are similar.
Thus we conclude that $\bar{\widetilde{J}}=\widetilde{\bar{J}}$.
\end{proof}

\subsection{Proof for Case B}

Let us show the commutativity of $\widetilde{f}_i$ and~$\delta$ for the Case B, that is,
$\ell^{(i-1)}\leq\ell+1<\ell^{(i)}=\infty$.
We divide the proof into the two cases.

\subsubsection[Case 1: $\widetilde{f}_i$ creates a~non-singular string]{Case 1: $\boldsymbol{\widetilde{f}_i}$ creates a~non-singular string}

\begin{Proposition}
Suppose that $\ell^{(i-1)}\leq\ell+1<\ell^{(i)}=\infty$.
If $\widetilde{f}_i$ acting on $(\nu,J)$ creates a~non-singular string, then we have the following identities:
\begin{gather*}
 (1)\quad \bar{\widetilde{\nu}}=\widetilde{\bar{\nu}}, \qquad (2) \quad \bar{\widetilde{J}}=\widetilde{\bar{J}}.
\end{gather*}
\end{Proposition}

\begin{proof}
(1) Since $\widetilde{f}_i$ does not change coriggings,~$\delta$ chooses the same strings in $(\nu,J)^{(a)}$ and
$(\widetilde{\nu},\widetilde{J})^{(a)}$ for all $a<i$.
Also, since $\widetilde{f}_i$ creates a~non-singular string,~$\delta$ does not choose a~string from
$(\widetilde{\nu},\widetilde{J})^{(i)}$, which coincides with the action of~$\delta$ on $(\nu,J)^{(i)}$.
Next, since~$\delta$ does not touch $(\nu,J)^{(i)}$, we have $\bar{J}^{(i)}=J^{(i)}$.
Thus $\widetilde{f}_i$ chooses the same string in both $(\nu,J)^{(i)}$ and $(\bar{\nu},\bar{J})^{(i)}$.
To summarize,~$\delta$ and $\widetilde{f}_i$ choose the same strings in $\delta\circ\widetilde{f}_i$ and
$\widetilde{f}_i\circ\delta$ which implies $\bar{\widetilde{\nu}}=\widetilde{\bar{\nu}}$.

(2) We only need to consider $J^{(i-1)}$ and $J^{(i)}$ in this case.
If the string is not touched by~$\delta$ and $\widetilde{f}_i$, the riggings after $\delta\circ\widetilde{f}_i$ and
$\widetilde{f}_i\circ\delta$ coincide since~$\delta$ and $\widetilde{f}_i$ choose the same strings in both cases
respectively.
The string $\big(\ell^{(i-1)},P^{(i-1)}_{\ell^{(i-1)}}(\nu)\big)$ becomes
$\big(\ell^{(i-1)}-1,P^{(i-1)}_{\ell^{(i-1)}-1}(\bar{\widetilde{\nu}})\big)$ and
$\big(\ell^{(i-1)}-1,P^{(i-1)}_{\ell^{(i-1)}-1}(\widetilde{\bar{\nu}})\big)$ after application of $\delta\circ\widetilde{f}_i$
and $\widetilde{f}_i\circ\delta$ respectively.
Since $\bar{\widetilde{\nu}}=\widetilde{\bar{\nu}}$ we see the coincidence of the riggings.
As for the string $(\ell,x_\ell)$, we see that both $\delta\circ\widetilde{f}_i$ and $\widetilde{f}_i\circ\delta$ give
the string $(\ell+1,x_\ell-1)$.
Hence we have $\bar{\widetilde{J}}=\widetilde{\bar{J}}$.
\end{proof}

\subsubsection[Case 2: $\widetilde{f}_i$ creates a~singular string]{Case 2: $\boldsymbol{\widetilde{f}_i}$ creates a~singular string}

This situation does not satisfy the condition of Proposition~\ref{th:core1} due to the following fact.
\begin{Proposition}
\label{f_undefined}
Suppose that $\ell^{(i-1)}\leq\ell+1<\ell^{(i)}=\infty$.
If $\widetilde{f}_i$ acting on $(\nu,J)$ creates a~singular string, then we have
\begin{gather*}
\widetilde{f}_i\circ\delta (\nu,J)=0.
\end{gather*}
\end{Proposition}
\begin{proof}
Since $\ell^{(i)}=\infty$, we have $(\nu,J)^{(i)}=(\bar{\nu},\bar{J})^{(i)}$ and thus $\bar{\ell}=\ell$.
Then the string $(\ell,x_\ell)$ in $(\nu,J)^{(i)}$ will become $(\ell+1,x_\ell-1)$ in $\widetilde{f}_i\circ\delta
(\nu,J)$.
Let us compute $P^{(i)}_{\ell+1}(\widetilde{\bar{\nu}})$.
From Lemma~\ref{lem:f_singular}(1), we have $P^{(i)}_{\ell+1}(\nu)=x_\ell+1$.
By the assumption $\ell^{(i-1)}\!\leq\!\ell+1<\ell^{(i)}$ we have
\mbox{$P^{(i)}_{\ell+1}(\bar{\nu})=P^{(i)}_{\ell+1}(\nu)\!-\!1=x_\ell$}.
Therefore we have $P^{(i)}_{\ell+1}(\widetilde{\bar{\nu}})= P^{(i)}_{\ell+1}(\bar{\nu})-2=x_\ell-2$.
But then the string $(\ell+1,x_\ell-1)$ in $\widetilde{f}_i\circ\delta (\nu,J)$ satisf\/ies
$P^{(i)}_{\ell+1}(\widetilde{\bar{\nu}})<x_\ell-1$ which implies that $\widetilde{f}_i\circ\delta (\nu,J)$ is not
a~valid rigged conf\/iguration.
Thus we have $\widetilde{f}_i\circ\delta (\nu,J)=0$.
\end{proof}
\begin{Example}
Consider the following rigged conf\/iguration $(\nu,J)$ of type $(B^{1,1})^{\otimes 4}\otimes B^{1,3}\otimes
B^{2,1}\otimes B^{2,2}\otimes B^{3,1}$ of $D^{(1)}_5$
\begin{center}
\unitlength 12pt
\begin{picture}(32,6)
\put(0,1){
\multiput(-0.8,0.1)(0,1){4}{1}
\put(-0.8,4.1){0}
\put(0,0){\Yboxdim12pt\yng(4,2,1,1,1)}
\put(1.2,0.1){0}
\put(1.2,1.1){0}
\put(1.2,2.1){0}
\put(2.2,3.1){1}
\put(4.2,4.1){$-1$}
}
\put(7.5,0){
\multiput(-0.8,0.1)(0,1){3}{1}
\multiput(-0.8,3.1)(0,1){3}{0}
\put(0,0){\Yboxdim12pt\yng(4,2,2,1,1,1)}
\put(1.2,0.1){0}
\put(1.2,1.1){1}
\put(1.2,2.1){1}
\put(2.2,3.1){$-1$}
\put(2.2,4.0){$-1$}
\put(4.2,5.0){$-1$}
}
\put(14.5,0){
\multiput(-0.8,0.1)(0,1){4}{0}
\put(-0.8,4.1){1}
\put(-0.8,5.1){2}
\put(0,0){\Yboxdim12pt\yng(4,2,1,1,1,1)}
\put(1.2,0.1){$-1$}
\put(1.2,1.1){$-1$}
\put(1.2,2.1){0}
\put(1.2,3.1){0}
\put(2.2,4.1){0}
\put(4.2,5.1){0}
}
\put(21.5,3){
\multiput(-0.8,0.1)(0,1){2}{0}
\put(-1.53,2.1){$-1$}
\put(0,0){\Yboxdim12pt\yng(3,1,1)}
\put(1.2,0.1){0}
\put(1.2,1.1){0}
\put(3.2,2.1){$-1$}
}
\put(27.5,4){
\put(-0.8,0.1){2}
\put(-0.8,1.1){0}
\put(0,0){\Yboxdim12pt\yng(4,1)}
\put(1.2,0.1){2}
\put(4.2,1.1){0}
}
\end{picture}
\end{center}
The corresponding tensor product $\Phi(\nu,J)$ is
\begin{gather*}
\Yboxdim14pt \Yvcentermath1 \young(2)\otimes \young(\mone)\otimes \young(1)\otimes \young(3)\otimes \young(112)\otimes
\young(4,\mfour)\otimes \young(13,\mfive\mone)\otimes \young(4,5,\mfour)
\end{gather*}
$\widetilde{f}_2$ acts on the string $(4,-1)$ of $(\nu,J)^{(2)}$ and makes it into $(5,-2)$ of
$(\widetilde{\nu},\widetilde{J})^{(2)}$.
Since we have $P^{(2)}_5(\widetilde{\nu})=-2$, the latter string is singular.
In particular, we see that $\widetilde{f}_2(\nu,J)\neq 0$.
Now $(\bar{\nu},\bar{J})$ is
\begin{center}
\unitlength 12pt
\begin{picture}(33,6)
\put(0,1){
\multiput(-0.8,0.1)(0,1){4}{0}
\put(-0.8,4.1){1}
\put(0,0){\Yboxdim12pt\yng(4,1,1,1,1)}
\put(1.2,0.1){0}
\put(1.2,1.1){0}
\put(1.2,2.1){0}
\put(1.2,3.1){0}
\put(4.2,4.1){$-1$}
}
\put(7.5,0){
\multiput(-0.8,0.1)(0,1){3}{1}
\multiput(-1.53,3.1)(0,1){3}{$-1$}
\put(0,0){\Yboxdim12pt\yng(4,2,2,1,1,1)}
\put(1.2,0.1){0}
\put(1.2,1.1){1}
\put(1.2,2.1){1}
\put(2.2,3.1){$-1$}
\put(2.2,4.0){$-1$}
\put(4.2,5.0){$-1$}
}
\put(14.5,0){
\multiput(-0.8,0.1)(0,1){4}{0}
\put(-0.8,4.1){1}
\put(-0.8,5.1){2}
\put(0,0){\Yboxdim12pt\yng(4,2,1,1,1,1)}
\put(1.2,0.1){$-1$}
\put(1.2,1.1){$-1$}
\put(1.2,2.1){0}
\put(1.2,3.1){0}
\put(2.2,4.1){0}
\put(4.2,5.1){0}
}
\put(21.5,3){
\multiput(-0.8,0.1)(0,1){2}{0}
\put(-1.53,2.1){$-1$}
\put(0,0){\Yboxdim12pt\yng(3,1,1)}
\put(1.2,0.1){0}
\put(1.2,1.1){0}
\put(3.2,2.1){$-1$}
}
\put(27.5,4){
\put(-0.8,0.1){2}
\put(-0.8,1.1){0}
\put(0,0){\Yboxdim12pt\yng(4,1)}
\put(1.2,0.1){2}
\put(4.2,1.1){0}
}
\end{picture}
\end{center}
Thus we have $\ell^{(1)}=2<\ell+1=5<\ell^{(2)}=\infty$ (see Proposition~\ref{f_undefined}).
Indeed, we have $P^{(2)}_5(\widetilde{\bar{\nu}})=-3$ and thus $\widetilde{f}_2(\bar{\nu},\bar{J})=0$.
\end{Example}

\subsection{Proof for Case C: preliminary steps}%\label{sec:caseC_easy}

\subsubsection{Classif\/ication}

Recall that the def\/ining condition of the Case C is $\ell^{(i-1)}\leq \ell+1\leq \ell^{(i)}<\infty$.
We divide the proof into the three cases according to the following extra conditions:
\begin{enumerate}\itemsep=0pt
\item[(1)] $m^{(i)}_{\ell+1}(\nu)>0$ and $\widetilde{f}_i$ creates a~singular string, or\\
$\ell+1=\ell^{(i)}$ and $\widetilde{f}_i$ creates a~non-singular string,
\item[(2)] $m^{(i)}_{\ell+1}(\nu)=0$ and $\widetilde{f}_i$ creates a~singular string,
 \item[(3)] $\ell+1<\ell^{(i)}$ and
$\widetilde{f}_i$ creates a~non-singular string.
\end{enumerate}
We give several remarks on the classif\/ication.
In case (2), the condition $m^{(i)}_{\ell+1}(\nu)=0$ automa\-ti\-cal\-ly implies that $\ell+1<\ell^{(i)}$ by the assumption
$\ell+1\leq \ell^{(i)}$.
On the other hand, in the f\/irst case of case (1), we will show that $\ell+1=\ell^{(i)}$.
We remark that we may have $m^{(i)}_{\ell+1}(\nu)>0$ in case (3).

\subsubsection{A common property for Case C}

In this case, we can show the following property.
Note that the corresponding analysis for the case $\ell^{(i)}=\ell_{(i)}$ is included in the analysis of
$P^{(i)}_{\ell^{(i)}-1}(\bar{\nu})$ which will be given later in this section.

\begin{Proposition}
\label{prop:C_prep}
Suppose that $\ell+1\leq \ell^{(i)}<\ell_{(i)}<\infty$.
Then we have
\begin{gather*}
P^{(i)}_{\ell_{(i)}-1}(\bar{\nu})>x_\ell.
\end{gather*}
\end{Proposition}
\begin{proof}
We follow the classif\/ication in Lemma~\ref{lem:vacancy_delta}.
Since $\ell^{(i)}<\ell_{(i)}$, we have to consider Cases~(V), (VI) and (VII).
Let $j$ be the largest integer satisfying $j<\ell_{(i)}$ and $m^{(i)}_j(\nu)>0$.
From $\ell^{(i)}<\ell_{(i)}$, we see that $\ell^{(i)}\leq j$ and in particular we have $\ell<j$.
In fact we only need the latter relation $\ell<j$ in the following proof\footnote{We will use this property in the
proof of Proposition~\ref{prop:D_1}.}.

{\bf Case (V).} According to Lemma~\ref{lem:vacancy_delta} we have to show $P^{(i)}_{\ell_{(i)}-1}(\nu)\geq
x_\ell$ under the condition $\ell^{(i)}<\ell^{(i+1)}=\dots =\ell_{(i)}$.
Suppose if possible that $P^{(i)}_{\ell_{(i)}-1}(\nu)<x_\ell$.
If $j=\ell_{(i)}-1$, then we have $x_{\ell_{(i)}-1}\leq P^{(i)}_{\ell_{(i)}-1}(\nu)<x_\ell$ which contradicts the
minimality of $x_\ell$.
So suppose that $j<\ell_{(i)}-1$.
By the assumption we see that $\ell<\ell_{(i)}$.
Then we have $P^{(i)}_{\ell_{(i)}}(\nu)\geq x_{\ell_{(i)}}>x_\ell$.
Thus we have $P^{(i)}_{\ell_{(i)}-1}(\nu)<P^{(i)}_{\ell_{(i)}}(\nu)$.
By the convexity of $P^{(i)}_k(\nu)$ between $j\leq k\leq\ell_{(i)}$, we have
$P^{(i)}_j(\nu)<P^{(i)}_{\ell_{(i)}-1}(\nu)<x_\ell$.
Since $x_j\leq P^{(i)}_j(\nu)$, this contradicts the minimality of $x_\ell$.
Thus we conclude that $P^{(i)}_{\ell_{(i)}-1}(\nu)\geq x_\ell$.

{\bf Case (VI).} According to Lemma~\ref{lem:vacancy_delta} we have to show $P^{(i)}_{\ell_{(i)}-1}(\nu)>x_\ell$
under the condition $\ell^{(i+1)}<\ell_{(i+1)}=\ell_{(i)}$.
Suppose if possible that $P^{(i)}_{\ell_{(i)}-1}(\nu)\leq x_\ell$.
If $j=\ell_{(i)}-1$, then we have $x_{j}\leq P^{(i)}_{j}(\nu)\leq x_\ell$.
However this relation contradicts the requirement $x_{j}>x_\ell$ which follows from $\ell<j$.
Next, suppose that $j<\ell_{(i)}-1$.
Then, as in the Case (V), we have $P^{(i)}_{\ell_{(i)}}(\nu)>x_\ell$ by $\ell<\ell_{(i)}$.
Thus we have $P^{(i)}_{\ell_{(i)}-1}(\nu)<P^{(i)}_{\ell_{(i)}}(\nu)$.
By the similar arguments in the previous case, we can deduce the contradiction.
Therefore we obtain $P^{(i)}_{\ell_{(i)}-1}(\nu)>x_\ell$ in this case.

{\bf Case (VII).} According to Lemma~\ref{lem:vacancy_delta} we have to show
$P^{(i)}_{\ell_{(i)}-1}(\nu)>x_\ell+1$ under the condition $\ell_{(i+1)}<\ell_{(i)}$.
Suppose if possible that $P^{(i)}_{\ell_{(i)}-1}(\nu)\leq x_\ell+1$.
If $j=\ell_{(i)}-1$, then we have $P^{(i)}_j(\nu)\geq x_j>x_\ell$ by $j>\ell$.
Then the only possibility that is compatible with the assumption is $P^{(i)}_j(\nu)=x_j=x_\ell +1$.
In particular, the string $(j,x_j)$ is singular whose length satisf\/ies $\ell_{(i+1)}\leq j<\ell_{(i)}$.
However this is in contradiction to the def\/inition of~$\ell_{(i)}$.

If $j<\ell_{(i)}-1$, again we have $P^{(i)}_{\ell_{(i)}}(\nu)>x_\ell$ as in the Case~(V).
Similarly, by $\ell<j$ we have $P^{(i)}_{j}(\nu)>x_\ell$.
Then by the convexity of $P^{(i)}_k(\nu)$ between $j\leq k\leq \ell_{(i)}$, the only possibility that is compatible with
$P^{(i)}_{\ell_{(i)}-1}(\nu)\leq x_\ell+1$ is $P^{(i)}_{j}(\nu)=\dots =P^{(i)}_{\ell_{(i)}}(\nu)=x_\ell+1$.
Since $P^{(i)}_{j}(\nu)\geq x_j>x_\ell$ it has to be $x_j=x_\ell+1$, in particular, the string $(j,x_j)$ is singular.
If $\ell_{(i+1)}\leq j$, this contradicts the def\/inition of $\ell_{(i)}$.
On the other hand, if $j<\ell_{(i+1)}$, $m^{(i+1)}_{\ell_{(i+1)}}(\nu)>0$ means that the convexity of $P^{(i)}_k(\nu)$
between $\ell_{(i+1)}-1\leq k\leq \ell_{(i+1)}+1$ is strict.
This is in contradiction to the relation $P^{(i)}_{j}(\nu)=\dots =P^{(i)}_{\ell_{(i)}}(\nu)=x_\ell+1$.
\end{proof}

\subsection{Proof for Case C (1)}

\begin{Proposition}
\label{prop:C(1)}
Suppose that $\ell^{(i-1)}\leq \ell+1\leq \ell^{(i)}<\infty$, $m^{(i)}_{\ell+1}(\nu)>0$ and $\widetilde{f}_i$ creates
a~singular string.
Then we have the following identities:
\begin{gather*}
 (1)\quad \bar{\widetilde{\nu}}=\widetilde{\bar{\nu}} , \qquad (2) \quad  \bar{\widetilde{J}}=\widetilde{\bar{J}}.
\end{gather*}
\end{Proposition}

\begin{proof}
(1) To begin with let us show that $\ell+1=\ell^{(i)}$.
By the assumption $m^{(i)}_{\ell+1}(\nu)>0$ there exists the string $(\ell+1,x_{\ell+1})$ of $(\nu,J)^{(i)}$.
We show that the string $(\ell+1,x_{\ell+1})$ is singular.
For this we derive the following three relations:
\begin{itemize}\itemsep=0pt
\item
$P^{(i)}_{\ell+1}(\nu)=x_\ell+1$ by Lemma~\ref{lem:f_singular},
\item
$x_\ell+1\leq x_{\ell+1}$, that is, $x_\ell<x_{\ell+1}$, since the string $(\ell+1,x_{\ell+1})$ of $(\nu,J)$ is not
selected by~$\widetilde{f}_i$ although it is longer than~$\ell$,
\item
$x_{\ell+1}\leq P^{(i)}_{\ell+1}(\nu)$ by def\/inition of the rigged conf\/igurations.
\end{itemize}
Combining these three inequalities, we obtain $P^{(i)}_{\ell+1}(\nu)=x_{\ell+1}$ and thus the string
$(\ell+1,x_{\ell+1})$ is singular.
Suppose if possible that $\ell+1<\ell^{(i)}$.
Then the string $(\ell+1,x_{\ell+1})$ is singular whose length satisf\/ies $\ell^{(i-1)}\leq\ell+1$, contradicting the
def\/inition of $\ell^{(i)}$.
Thus $\ell+1=\ell^{(i)}$.
Then we can show the commutativity of $\widetilde{f}_i$ and~$\delta$ as follows.

Let us consider the case $\delta\circ\widetilde{f}_i$.
It is enough to analyze the strings $(\ell,x_\ell)$ and $\big(\ell+1,P^{(i)}_{\ell+1}(\nu)\big)$ in $(\nu,J)^{(i)}$:
\begin{gather*}
\begin{array}{@{}l}
(\ell,x_\ell)
\\
\big(\ell+1,P^{(i)}_{\ell+1}(\nu)\big)
\end{array}
\overset{\widetilde{f}_i}{\longmapsto}
\begin{array}
{l} (\ell+1,x_{\ell}-1)
\\
\big(\ell+1,P^{(i)}_{\ell+1}(\widetilde{\nu})\big)
\end{array}
\overset{\delta}{\longmapsto}
\begin{array}
{l} (\ell+1,x_{\ell}-1)
\\
\big(\ell,P^{(i)}_{\ell}(\bar{\widetilde{\nu}})\big)
\end{array}.
\end{gather*}
Recall that $\widetilde{f}_i$ does not change coriggings other than the touched string.
In particular, we have $\widetilde{\ell}^{(i-1)}=\ell^{(i-1)}$.
Thus the next~$\delta$ can act on one of the length $\ell+1$ singular strings.

Let us consider the case $\widetilde{f}_i\circ\delta$.
Recall that if $\ell^{(i)}<\ell_{(i)}$ we have $P^{(i)}_{\ell_{(i)}-1}(\bar{\nu})>x_\ell$ by
Proposition~\ref{prop:C_prep}.
Therefore $\widetilde{f}_i$ will not act on the length $\ell_{(i)}-1$ string of $(\bar{\nu},\bar{J})^{(i)}$ that has
been created by~$\delta$.
On the other hand, if $\ell^{(i)}=\ell_{(i)}$ we do not need to take care of this issue.
Thus it is enough to analyze the behaviors of the lengths $\ell$ and $\ell+1$ strings:
\begin{gather*}
\begin{array}
{l} (\ell,x_\ell)
\\
\big(\ell+1,P^{(i)}_{\ell+1}(\nu)\big)
\end{array}
\overset{\delta}{\longmapsto}
\begin{array}
{l} (\ell,x_\ell)
\\
\big(\ell,P^{(i)}_{\ell}(\bar{\nu})\big)
\end{array}
\overset{\widetilde{f}_i}{\longmapsto}
\begin{array}
{l} (\ell+1,x_\ell-1)
\\
\big(\ell,P^{(i)}_{\ell}(\widetilde{\bar{\nu}})\big)
\end{array}.
\end{gather*}
Recall that $\ell+1=\ell^{(i)}$ so that~$\delta$ acts on the string $\big(\ell+1,P^{(i)}_{\ell+1}(\nu)\big)$.
Also, since $(\bar{\nu},\bar{J})$ is the valid rigged conf\/iguration, we have the requirement
$P^{(i)}_{\ell}(\bar{\nu})\geq x_\ell$.
Thus the string $(\ell,x_\ell)$ has the largest length among the strings of the smallest rigging within
$(\bar{\nu},\bar{J})^{(i)}$.

To summarize, we are left with lengths $\ell$ and $\ell+1$ strings in both $\delta\circ\widetilde{f}_i$ and
$\widetilde{f}_i\circ\delta$.
Thus we have $\bar{\widetilde{\nu}}=\widetilde{\bar{\nu}}$.

(2) We have $P^{(i)}_{\ell}(\bar{\widetilde{\nu}})=P^{(i)}_{\ell}(\widetilde{\bar{\nu}})$ by
$\bar{\widetilde{\nu}}=\widetilde{\bar{\nu}}$.
Thus $\bar{\widetilde{J}}=\widetilde{\bar{J}}$.
\end{proof}

\begin{Proposition}
Suppose that we have $\ell+1=\ell^{(i)}$ and $\widetilde{f}_i$ acting on $(\nu,J)$ creates the non-singular string.
Then we have the following identities:
\begin{gather*}
 (1)\quad \bar{\widetilde{\nu}}=\widetilde{\bar{\nu}}, \qquad (2) \quad \bar{\widetilde{J}}=\widetilde{\bar{J}}.
\end{gather*}
\end{Proposition}
\begin{proof}
Since $\ell+1=\ell^{(i)}$ there is a~singular string $\big(\ell+1,P^{(i)}_{\ell+1}(\nu)\big)$ in $(\nu,J)^{(i)}$.
Then this string and $(\ell,x_\ell)$ behave just as in the proof of the previous proposition which conf\/irms the claim.
\end{proof}

\subsection{Proof for Case C (2)}%\label{sec:caseC(2)}

\subsubsection{Outline}

As the second step, we treat the case when $m^{(i)}_{\ell+1}(\nu)=0$ and $\widetilde{f}_i$ creates a~singular string.
In this case, we will show the following assertion by case by case analysis.
\begin{Proposition}
\label{prop:caseC1}
Suppose that $\ell^{(i-1)}\leq \ell+1<\ell^{(i)}<\infty$ and $m^{(i)}_{\ell+1}(\nu)=0$.
If $\widetilde{f}_i$ acting on~$(\nu,J)$ creates a~singular string, then the following two conditions are satisfied:
\begin{gather}
P^{(i)}_{\ell^{(i)}-1}(\bar{\nu})= x_\ell,\label{eq:caseC1}
\\
m_k^{(i+1)}(\nu)=0\qquad \text{for}\ \ \ell<k<\ell^{(i)}.\label{eq:caseC2}
\end{gather}
\end{Proposition}

Note that we have $\ell+1\neq\ell^{(i)}$ by $m^{(i)}_{\ell+1}(\nu)=0$.
To begin with, we explain how these properties provide the proof of $\bar{\widetilde{\nu}}=\widetilde{\bar{\nu}}$.

\begin{Proposition}
\label{prop:caseC(2)1}
Suppose that $\ell^{(i-1)}\leq \ell+1<\ell^{(i)}<\infty$, $m^{(i)}_{\ell+1}(\nu)=0$ and $\widetilde{f}_i$ acting on
$(\nu,J)$ creates a~singular string.
Then we have $\bar{\widetilde{\nu}}=\widetilde{\bar{\nu}}$.
\end{Proposition}

\begin{proof}
Let us analyze $\widetilde{f}_i\circ\delta$ and $\delta\circ\widetilde{f}_i$ respectively.

{\bf Step 1.} Let us consider the case $\widetilde{f}_i\circ\delta$.
Let us consider the strings $\big(\ell^{(i)}-1,P^{(i)}_{\ell^{(i)}-1}(\bar{\nu})\big)$ and
$\big(\ell_{(i)}-1,P^{(i)}_{\ell_{(i)}-1}(\bar{\nu})\big)$ of $(\bar{\nu},\bar{J})^{(i)}$ created by~$\delta$.
By~\eqref{eq:caseC1} we see that the rigging of the string $\big(\ell^{(i)}-1,P^{(i)}_{\ell^{(i)}-1}(\bar{\nu})\big)$ of
$(\bar{\nu},\bar{J})$ is $x_\ell$.
On the other hand, if $\ell^{(i)}<\ell_{(i)}$ we have $P^{(i)}_{\ell_{(i)}-1}(\bar{\nu})>x_\ell$ by
Proposition~\ref{prop:C_prep}.
Since we are assuming that $\ell<\ell^{(i)}-1$, the string $\big(\ell^{(i)}-1,P^{(i)}_{\ell^{(i)}-1}(\bar{\nu})\big)$ has the
largest length among the strings of smallest riggings within $(\bar{\nu},\bar{J})^{(i)}$.
Thus $\widetilde{f}_i$ adds a~box to the length $\ell^{(i)}-1$ string that has been shortened by~$\delta$.

{\bf Step 2.} Let us consider the case $\delta\circ\widetilde{f}_i$.
Recall that $\widetilde{f}_i$ does not change coriggings of untouched strings.
By the assumption $\widetilde{f}_i$ creates the singular string $(\ell+1,x_\ell-1)$ whose length satisf\/ies
$\ell^{(i-1)}\leq \ell+1<\ell^{(i)}$.
Hence~$\delta$ choose this string and obtain $\widetilde{\ell}^{(i)}=\ell+1$.
By def\/inition of~$\delta$ we have $\widetilde{\ell}^{(i)}=\ell+1\leq\widetilde{\ell}^{(i+1)}$ and by~\eqref{eq:caseC2}
we conclude that $\ell^{(i)}\leq\widetilde{\ell}^{(i+1)}$.
Thus~$\delta$ choose a~string of length $\ell^{(i+1)}$ from $(\widetilde{\nu},\widetilde{J})^{(i+1)}$.
Therefore we obtain $\widetilde{\ell}^{(a)}=\ell^{(a)}$ for all $a\neq i$ and $\widetilde{\ell}_{(a)}=\ell_{(a)}$ for
all~$a$.

{\bf Step 3.} To summarize, we have $\bar{\widetilde{\nu}}^{(a)}=\widetilde{\bar{\nu}}^{(a)}=\bar{\nu}^{(a)}$ for
all $a\neq i$ and $\bar{\widetilde{\nu}}^{(i)}=\widetilde{\bar{\nu}}^{(i)}$ is obtained by removing a~box from length
$\ell_{(i)}$ string of $\nu^{(i)}$.
\end{proof}

The remaining property $\bar{\widetilde{J}}=\widetilde{\bar{J}}$ will be checked from case by case analysis and it will
be done along with the proof of Proposition~\ref{prop:caseC1}.

\subsubsection{Classif\/ication}

We divide the proofs of Proposition~\ref{prop:caseC1} as well as the identity $\bar{\widetilde{J}}=\widetilde{\bar{J}}$
into the two cases according to the following classif\/ication.
\begin{Lemma}
\label{lem:C_1}
Suppose that $\ell^{(i-1)}\leq\ell+1<\ell^{(i)}<\infty$.
Then there are only two possibilities as follows:
\begin{gather}
(\ell,x_\ell)\text{ is singular and }\ell=\ell^{(i-1)}-1,
\label{caseC_1}
\\
(\ell,x_\ell)\text{ is non-singular.}
\label{caseC_2}
\end{gather}
\end{Lemma}
\begin{proof}
Suppose that the string $(\ell,x_\ell)$ is singular.
If its length satisf\/ies $\ell^{(i-1)}\leq\ell<\ell^{(i)}$, then it is in contradiction to the def\/inition of
$\ell^{(i)}$.
Therefore we must have $\ell=\ell^{(i-1)}-1$.
\end{proof}

\subsubsection{Proof for the case~(\ref{caseC_1})}

\label{sec:caseC(2)3}
Under the present assumptions, let us show the following property.
\begin{Proposition}
\label{lem:C_2}
Suppose that we have $\ell^{(i-1)}=\ell+1<\ell^{(i)}<\infty$, $m^{(i)}_{\ell+1}(\nu)=0$, the string $(\ell,x_\ell)$ is
singular and $\widetilde{f}_i$ creates a~singular string.
Then we have
\begin{gather}
P^{(i)}_\ell(\nu)=x_\ell,
\qquad
P^{(i)}_{\ell+1}(\nu)=P^{(i)}_{\ell+2}(\nu) =\dots =P^{(i)}_{\ell^{(i)}}(\nu)=x_\ell+1.
\label{eq:C_4}
\end{gather}
\end{Proposition}

\begin{proof}
Since $\widetilde{f}_i$ creates a~singular string, we have
\begin{gather}\label{eq:C_2}
P^{(i)}_{\ell+1}(\nu)\overset{\text{by Lemma~\ref{lem:f_singular}(1)}}{=}x_\ell+1
\overset{\text{since the string $(\ell,x_\ell)$ is singular}}{=}
P^{(i)}_\ell(\nu)+1  .
\end{gather}
Due to the existence of the length $\ell^{(i-1)}=\ell+1$ string at $\nu^{(i-1)}$ and the assumption
$m^{(i)}_{\ell+1}(\nu)=0$, the vacancy number $P^{(i)}_k(\nu)$ is strictly convex between the interval $\ell\leq k\leq
\ell+2$.
Then the relation~\eqref{eq:C_2} implies $P^{(i)}_{\ell+1}(\nu)\geq P^{(i)}_{\ell+2}(\nu)$.
Let us show that this inequality is in fact equality.
For this let $j$ be the length of the smallest string of $\nu^{(i)}$ that is longer than $\ell$.
Since $m^{(i)}_{\ell+1}(\nu)=0$ we have $\ell^{(i-1)}=\ell+1<j\leq\ell^{(i)}$.
By Lemma~\ref{lem:f_singular}(2) we have
\begin{gather}
P^{(i)}_{\ell+1}(\nu)\leq x_j\leq P^{(i)}_j(\nu).
\label{eq:C_3}
\end{gather}
Then by the convexity of $P^{(i)}_k(\nu)$ between $\ell\leq k\leq j$, the relations $P^{(i)}_{\ell+1}(\nu)\geq
P^{(i)}_{\ell+2}(\nu)$ and $P^{(i)}_{\ell+1}(\nu)\leq P^{(i)}_j(\nu)$ imply the relation
$P^{(i)}_{\ell+1}(\nu)=P^{(i)}_{\ell+2}(\nu)=\dots =P^{(i)}_{j}(\nu)$.
Combining this equality with~\eqref{eq:C_3} we obtain $x_j=P^{(i)}_j(\nu)$, thus the string $(j,x_j)$ is singular.
By the relation $\ell^{(i-1)}<j\leq\ell^{(i)}$ we deduce that $j=\ell^{(i)}$.
\end{proof}

\begin{Corollary}
\label{cor:C_1}
Under the same assumptions of Proposition~{\rm \ref{lem:C_2}}, we have
\begin{enumerate}\itemsep=0pt
\item[$(1)$] $m^{(i-1)}_{\ell+1}(\nu)=1$, \item[$(2)$] $m^{(i-1)}_k(\nu)=0$ for $\ell+1<k<\ell^{(i)}$, \item[$(3)$]
$m^{(i)}_k(\nu)=0$ for $\ell+1\leq k<\ell^{(i)}$, \item[$(4)$] $m^{(i+1)}_k(\nu)=0$ for $\ell+1\leq k<\ell^{(i)}$.
\end{enumerate}
\end{Corollary}

\begin{proof}
The assertions for $m^{(i)}_k(\nu)$ given in (3) are proved at the end of the proof of Proposition~\ref{lem:C_2}.

Let us show the assertions for $m^{(i+1)}_k(\nu)$.
If $m^{(i+1)}_k(\nu)>0$ for some $\ell+1<k<\ell^{(i)}$, the relation~\eqref{eq:C_4} must be strictly convex around such
$k$, which is a~contradiction.
Therefore $m^{(i+1)}_k(\nu)=0$ for all $\ell+1<k<\ell^{(i)}$.
By using~\eqref{eq:C_4} and the assumption $m^{(i)}_{\ell+1}(\nu)=0$, the estimate of Lemma~\ref{lem:convexity3} with
$l=\ell+1$ gives
\begin{gather*}
-P^{(i)}_\ell(\nu)+2P^{(i)}_{\ell+1}(\nu)-P^{(i)}_{\ell+2}(\nu)=1
\\
\qquad{}
\geq m^{(i-1)}_{\ell+1}(\nu)-2m^{(i)}_{\ell+1}(\nu)+m^{(i+1)}_{\ell+1}(\nu)
=m^{(i-1)}_{\ell+1}(\nu)+m^{(i+1)}_{\ell+1}(\nu).
\end{gather*}
Recall that we have $m^{(i-1)}_{\ell+1}(\nu)>0$ by def\/inition of $\ell^{(i-1)}(=\ell+1)$.
Then the only possibility that is compatible with the above inequality is $m^{(i-1)}_{\ell+1}(\nu)=1$ and
$m^{(i+1)}_{\ell+1}(\nu)=0$.
Hence we have conf\/irmed the relations (4) and (1).

Finally, since $m^{(i)}_k(\nu)=0$ and $m^{(i+1)}_k(\nu)=0$ for $\ell+1<k<\ell^{(i)}$, the relation~\eqref{eq:C_4}
implies that $m^{(i-1)}_k(\nu)=0$ for $\ell+1<k<\ell^{(i)}$, since otherwise the relation would be strictly convex.
Thus we have conf\/irmed the relation (2).
\end{proof}

\begin{proof}
[Proof of Proposition~\ref{prop:caseC1} for the case~(\ref{caseC_1})] The assertions $m^{(i+1)}_k(\nu)=0$ for $\ell<
k<\ell^{(i)}$ are already proved in Corollary~\ref{cor:C_1}.

From~\eqref{eq:C_4} we have $P^{(i)}_{\ell^{(i)}-1}(\nu)=x_\ell+1$ since we have $\ell<\ell^{(i)}-1$ by the assumption
$(\ell^{(i-1)}=)$ $\ell+1<\ell^{(i)}$.
Since $\ell^{(i-1)}\leq\ell^{(i)}-1$ we have $P^{(i)}_{\ell^{(i)}-1}(\bar{\nu})= P^{(i)}_{\ell^{(i)}-1}(\nu)-1=x_\ell$.
This completes the proof of Proposition~\ref{prop:caseC1} for the case~\eqref{caseC_1}.
\end{proof}

Thus we have conf\/irmed $\bar{\widetilde{\nu}}=\widetilde{\bar{\nu}}$ in this case.
In order to prove $\bar{\widetilde{J}}=\widetilde{\bar{J}}$, we need additional relations obtained in
Corollary~\ref{cor:C_1}.

\begin{Proposition}
\label{prop:C_3}
Suppose that we have $\ell^{(i-1)}=\ell+1<\ell^{(i)}<\infty$, $m^{(i)}_{\ell+1}(\nu)=0$, the string $(\ell,x_\ell)$ is
singular and $\widetilde{f}_i$ creates a~singular string.
Then we have $\bar{\widetilde{J}}=\widetilde{\bar{J}}$.
\end{Proposition}
\begin{proof}
From Corollary~\ref{cor:C_1} and Proposition~\ref{prop:caseC(2)1} we see that there are no strings of length from
$\ell+1$ to $\ell^{(i)}-1$ that are not touched by~$\delta$ nor $\widetilde{f}_i$ during both
$\delta\circ\widetilde{f}_i$ and $\widetilde{f}_i\circ\delta$.
Thus it is enough to analyze the strings $(\ell,x_\ell)$ and $(\ell^{(i)},x_{\ell^{(i)}})$ of $(\nu,J)^{(i)}$.

Let us analyze the behavior of the string $(\ell,x_\ell)$.
From the proof of Proposition~\ref{prop:caseC(2)1} we have
\begin{gather*}
\begin{array}{@{}lllll}
(\ell,x_\ell)&\overset{\widetilde{f}_i}{\longmapsto}&(\ell+1,x_\ell-1)
&\overset{\delta}{\longmapsto}&\big(\ell,P^{(i)}_\ell(\bar{\widetilde{\nu}})\big),
\\
(\ell,x_\ell)&\overset{\delta}{\longmapsto}&(\ell,x_\ell) &\overset{\widetilde{f}_i}{\longmapsto}&(\ell,x_\ell).
\end{array}
\end{gather*}
By $\ell<\ell^{(i-1)}$ we have $P^{(i)}_\ell(\bar{\widetilde{\nu}})=P^{(i)}_\ell(\nu)$.
Recall that we have $P^{(i)}_\ell(\nu)=x_\ell$ by~\eqref{eq:C_4}.
Hence we obtain $P^{(i)}_\ell(\bar{\widetilde{\nu}})=x_\ell$.

Next let us analyze the behavior of the string $(\ell^{(i)},x_{\ell^{(i)}})$.
From the proof of Proposition~\ref{prop:caseC(2)1} we have
\begin{gather*}
\begin{array}{@{}lllll}
(\ell^{(i)},x_{\ell^{(i)}})&\overset{\widetilde{f}_i}{\longmapsto}& (\ell^{(i)},x_{\ell^{(i)}}-2)
&\overset{\delta}{\longmapsto}&(\ell^{(i)},x_{\ell^{(i)}}-2),
\\
(\ell^{(i)},x_{\ell^{(i)}})&\overset{\delta}{\longmapsto}&\big(\ell^{(i)}-1,P^{(i)}_{\ell^{(i)}-1}(\bar{\nu})\big)
&\overset{\widetilde{f}_i}{\longmapsto}&\big(\ell^{(i)},P^{(i)}_{\ell^{(i)}-1}(\bar{\nu})-1\big).
\end{array}
\end{gather*}
Recall that we have $\ell+1<\ell^{(i)}$.
Thus the f\/irst $\widetilde{f}_i$ of $\delta\circ\widetilde{f}_i$ decreases the rigging of the string
$(\ell^{(i)},x_{\ell^{(i)}})$ by 2.

Let us check the coincidence of the f\/inal rigging of this case.
Recall that we have $x_{\ell^{(i)}}=P^{(i)}_{\ell^{(i)}}(\nu)$ since the string $(\ell^{(i)},x_{\ell^{(i)}})$ is
singular by def\/inition of $\ell^{(i)}$.
We also have $P^{(i)}_{\ell^{(i)}}(\nu)=x_\ell+1$ by~\eqref{eq:C_4}.
Therefore we have $x_{\ell^{(i)}}-2=x_\ell-1$.

On the other hand, since $\ell^{(i-1)}\leq\ell^{(i)}-1$ we have
$P^{(i)}_{\ell^{(i)}-1}(\bar{\nu})=P^{(i)}_{\ell^{(i)}-1}(\nu)-1$.
From~\eqref{eq:C_4} we have $P^{(i)}_{\ell^{(i)}-1}(\nu)=x_\ell+1$.
Therefore we have $P^{(i)}_{\ell^{(i)}-1}(\bar{\nu})-1=x_\ell-1$.
Hence we conclude that $x_{\ell^{(i)}}-2=P^{(i)}_{\ell^{(i)}-1}(\bar{\nu})-1$.
\end{proof}

\begin{Example}
Consider the following rigged conf\/iguration $(\nu,J)$ of type $(B^{1,1})^{\otimes 3}\otimes B^{1,3}\otimes
B^{2,1}\otimes B^{2,2}\otimes B^{3,1}$ of $D^{(1)}_5$
\begin{center}
\unitlength 12pt
\begin{picture}(36,6)
\put(0,0){
\multiput(-0.8,2.1)(0,1){3}{1}
\put(-0.8,5.1){0}
\put(0,2){\Yboxdim12pt\yng(5,3,1,1)}
\put(1.2,2.1){0}
\put(1.2,3.1){0}
\put(3.2,4.1){1}
\put(5.2,5.1){0}
}
\put(7.5,1){
\put(-0.8,0.1){1}
\multiput(-1.53,1.1)(0,1){2}{$-1$}
\multiput(-0.8,3.1)(0,1){2}{0}
\put(0,0){\Yboxdim12pt\yng(5,4,2,2,1)}
\put(1.2,0.1){1}
\put(2.2,1.1){$-1$}
\put(2.2,2.1){$-1$}
\put(4.2,3.1){0}
\put(5.2,4.1){0}
}
\put(15.2,1){
\multiput(-0.8,0.1)(0,1){2}{0}
\multiput(-0.8,2.1)(0,1){3}{2}
\put(0,0){\Yboxdim12pt\yng(5,4,4,1,1)}
\put(1.2,0.1){$-1$}
\put(1.2,1.1){0}
\put(4.2,2.1){2}
\put(4.2,3.1){2}
\put(5.2,4.1){2}
}
\put(23.0,4){
\multiput(-1.53,0.1)(0,1){2}{$-2$}
\put(0,0){\Yboxdim12pt\yng(4,4)}
\put(4.2,0.1){$-2$}
\put(4.2,1.1){$-2$}
}
\put(30.5,4){
\put(-1.53,0.1){$-2$}
\put(-1.53,1.1){$-3$}
\put(0,0){\Yboxdim12pt\yng(5,4)}
\put(4.2,0.1){$-2$}
\put(5.2,1.1){$-3$}
}
\end{picture}
%\{\{5,3,1,1\},\{5,4,2,2,1\},\{5,4,4,1,1\},\{4,4\},\{5,4\}\},
%\{\{0,1,0,0\},\{0,0,-1,-1,1\},\{2,2,2,0,-1\},\{-2,-2\},\{-3,-2\}\}
\end{center}
The corresponding tensor product $\Phi(\nu,J)$ is
\begin{gather*}
\Yboxdim14pt \Yvcentermath1 \young(\mone)\otimes \young(\mone)\otimes \young(\mtwo)\otimes
\young(\mfive\mthree\mone)\otimes \young(1,\mfour)\otimes \young(13,2\mfour)\otimes \young(1,4,\mfour)
\end{gather*}
We see that $\widetilde{f}_2$ acts on the string $(\ell,x_\ell)=(2,-1)$ of $(\nu,J)^{(2)}$.
Then $(\bar{\nu},\bar{J})$ is
\begin{center}
\unitlength 12pt
\begin{picture}(36,6)
\put(0,0){
\multiput(-0.8,2.1)(0,1){2}{0}
\multiput(-0.8,4.1)(0,1){2}{1}
\put(0,2){\Yboxdim12pt\yng(4,2,1,1)}
\put(1.2,2.1){0}
\put(1.2,3.1){0}
\put(2.2,4.1){1}
\put(4.2,5.1){1}
}
\put(7.5,1){
\put(-0.8,0.1){1}
\multiput(-1.53,1.1)(0,1){4}{$-1$}
\put(0,0){\Yboxdim12pt\yng(4,3,2,2,1)}
\put(1.2,0.1){1}
\put(2.2,1.1){$-1$}
\put(2.2,2.1){$-1$}
\put(3.2,3.1){$-1$}
\put(4.2,4.1){$-1$}
}
\put(15.2,1){
\multiput(-0.8,0.1)(0,1){2}{0}
\multiput(-0.8,2.1)(0,1){3}{2}
\put(0,0){\Yboxdim12pt\yng(5,3,3,1,1)}
\put(1.2,0.1){$-1$}
\put(1.2,1.1){0}
\put(3.2,2.1){2}
\put(3.2,3.1){2}
\put(5.2,4.1){2}
}
\put(23,4){
\put(-1.53,0.1){$-1$}
\put(-1.53,1.1){$-2$}
\put(0,0){\Yboxdim12pt\yng(4,3)}
\put(3.2,0.1){$-1$}
\put(4.2,1.1){$-2$}
}
\put(30.5,4){
\put(-1.53,0.1){$-1$}
\put(-1.53,1.1){$-3$}
\put(0,0){\Yboxdim12pt\yng(5,3)}
\put(3.2,0.1){$-1$}
\put(5.2,1.1){$-3$}
}
\end{picture}
%\{\{4,2,1,1\},\{4,3,2,2,1\},\{5,3,3,1,1\},\{4,3\},\{5,3\}\},\\
%\{\{1,1,0,0\},\{-1,-1,-1,-1,1\},\{2,2,2,0,-1\},\{-2,-1\},\{-3,-1\}\}
\end{center}
and $(\widetilde{\nu},\widetilde{J})$ is
\begin{center}
\unitlength 12pt
\begin{picture}(36,6)
\put(0,0){
\multiput(-0.8,2.1)(0,1){2}{1}
\put(-0.8,4.1){2}
\put(-0.8,5.1){1}
\put(0,2){\Yboxdim12pt\yng(5,3,1,1)}
\put(1.2,2.1){0}
\put(1.2,3.1){0}
\put(3.2,4.1){2}
\put(5.2,5.1){1}
}
\put(7.5,1){
\put(-0.8,0.1){1}
\put(-1.53,1.1){$-1$}
\multiput(-1.53,2.1)(0,1){3}{$-2$}
\put(0,0){\Yboxdim12pt\yng(5,4,3,2,1)}
\put(1.2,0.1){1}
\put(2.2,1.1){$-1$}
\put(3.2,2.1){$-2$}
\put(4.2,3.1){$-2$}
\put(5.2,4.1){$-2$}
}
\put(15.2,1){
\multiput(-0.8,0.1)(0,1){2}{0}
\multiput(-0.8,2.1)(0,1){3}{3}
\put(0,0){\Yboxdim12pt\yng(5,4,4,1,1)}
\put(1.2,0.1){$-1$}
\put(1.2,1.1){0}
\put(4.2,2.1){3}
\put(4.2,3.1){3}
\put(5.2,4.1){3}
}
\put(23,4){
\multiput(-1.53,0.1)(0,1){2}{$-2$}
\put(0,0){\Yboxdim12pt\yng(4,4)}
\put(4.2,0.1){$-2$}
\put(4.2,1.1){$-2$}
}
\put(30.5,4){
\put(-1.53,0.1){$-2$}
\put(-1.53,1.1){$-3$}
\put(0,0){\Yboxdim12pt\yng(5,4)}
\put(4.2,0.1){$-2$}
\put(5.2,1.1){$-3$}
}
\end{picture}
%\{\{5,3,1,1\},\{5,4,3,2,1\},\{5,4,4,1,1\},\{4,4\},\{5,4\}\},
%\{\{1,2,0,0\},\{-2,-2,-2,-1,1\},\{3,3,3,0,-1\},\{-2,-2\},\{-3,-2\}\}
\end{center}
Note that we have $\ell^{(1)}=3=\ell+1<\ell^{(2)}=4$ and $m^{(2)}_3(\nu)=0$.
Moreover, we see that the string $(\ell,x_\ell)$ is singular since we have $P^{(2)}_2=-1$, and $\widetilde{f}_2$ makes
it into the singular string $(3,-2)$ of $(\widetilde{\nu},\widetilde{J})^{(2)}$ since we have
$P^{(2)}_3(\widetilde{\nu})=-2$.
Hence this is the example for the present case.
Then we have $P^{(2)}_2(\nu)=-1$ and $P^{(2)}_3(\nu)=P^{(2)}_4(\nu)=0$ where $\ell^{(2)}=4$ (see
Proposition~\ref{lem:C_2}).
We have $m^{(1)}_3(\nu)=1$ and $m^{(2)}_3(\nu)=m^{(3)}_3(\nu)=0$ (see Corollary~\ref{cor:C_1}).
Finally $(\bar{\widetilde{\nu}},\bar{\widetilde{J}})=(\widetilde{\bar{\nu}},\widetilde{\bar{J}})$ is
\begin{center}
\unitlength 12pt
\begin{picture}(36,6)
\put(0,0){
\multiput(-0.8,2.1)(0,1){2}{0}
\multiput(-0.8,4.1)(0,1){2}{1}
\put(0,2){\Yboxdim12pt\yng(4,2,1,1)}
\put(1.2,2.1){0}
\put(1.2,3.1){0}
\put(2.2,4.1){1}
\put(4.2,5.1){1}
}
\put(7.5,1){
\put(-0.8,0.1){1}
\multiput(-1.53,1.1)(0,1){3}{$-1$}
\put(-1.53,4.1){$-2$}
\put(0,0){\Yboxdim12pt\yng(5,3,2,2,1)}
\put(1.2,0.1){1}
\put(2.2,1.1){$-1$}
\put(2.2,2.1){$-1$}
\put(3.2,3.1){$-1$}
\put(5.2,4.1){$-2$}
}
\put(15.2,1){
\multiput(-0.8,0.1)(0,1){2}{0}
\multiput(-0.8,2.1)(0,1){2}{2}
\put(-0.8,4.1){3}
\put(0,0){\Yboxdim12pt\yng(5,3,3,1,1)}
\put(1.2,0.1){$-1$}
\put(1.2,1.1){0}
\put(3.2,2.1){2}
\put(3.2,3.1){2}
\put(5.2,4.1){3}
}
\put(23,4){
\put(-1.53,0.1){$-1$}
\put(-1.53,1.1){$-2$}
\put(0,0){\Yboxdim12pt\yng(4,3)}
\put(3.2,0.1){$-1$}
\put(4.2,1.1){$-2$}
}
\put(30.5,4){
\put(-1.53,0.1){$-1$}
\put(-1.53,1.1){$-3$}
\put(0,0){\Yboxdim12pt\yng(5,3)}
\put(3.2,0.1){$-1$}
\put(5.2,1.1){$-3$}
}
\end{picture}
%\{\{4,2,1,1\},\{5,3,2,2,1\},\{5,3,3,1,1\},\{4,3\},\{5,3\}\},
%\{\{1,1,0,0\},\{-2,-1,-1,-1,1\},\{3,2,2,0,-1\},\{-2,-1\},\{-3,-1\}\}
\end{center}
\end{Example}

\subsubsection{Proof for the case~(\ref{caseC_2})}\label{sec:caseC(2)4}

Under the present assumptions, let us show the following property.
\begin{Proposition}
\label{lem:C_3}
Suppose that we have $\ell^{(i-1)}\leq\ell+1<\ell^{(i)}<\infty$, $m^{(i)}_{\ell+1}(\nu)=0$, the string $(\ell,x_\ell)$
is non-singular and $\widetilde{f}_i$ creates a~singular string.
Then we have
\begin{gather}
P^{(i)}_\ell(\nu)=P^{(i)}_{\ell+1}(\nu) =\dots =P^{(i)}_{\ell^{(i)}}(\nu)=x_\ell+1.
\label{eq:C_5}
\end{gather}
\end{Proposition}

\begin{proof}
Let $j$ be the length of the shortest string of $\nu^{(i)}$ that satisf\/ies $\ell<j$.
By the assumption $m^{(i)}_{\ell +1}(\nu)=0$ we have $\ell+1<j\leq\ell^{(i)}$.
From the convexity of $P^{(i)}_k(\nu)$ between $\ell\leq k\leq j$ we have
\begin{gather}
\label{eq:C_1}
P^{(i)}_{\ell+1}(\nu)\geq \min\big\{P^{(i)}_\ell(\nu),P^{(i)}_j(\nu)\big\}.
\end{gather}
In order to evaluate the minimum in~\eqref{eq:C_1}, let us suppose if possible that $P^{(i)}_\ell(\nu)>P^{(i)}_j(\nu)$.
Then the inequality~\eqref{eq:C_1} reads $P^{(i)}_{\ell+1}(\nu)\geq P^{(i)}_j(\nu)$.
Recall from Lemma~\ref{lem:f_singular}(2) that we have $P^{(i)}_{\ell+1}(\nu)\leq x_j$ and by def\/inition of the rigged
conf\/iguration we have $x_j\leq P^{(i)}_j(\nu)$, that is, $P^{(i)}_{\ell+1}(\nu)\leq P^{(i)}_j(\nu)$.
Combining the two inequalities we obtain $P^{(i)}_{\ell+1}(\nu)=P^{(i)}_j(\nu)$.
By the assumption, we conclude that $P^{(i)}_\ell(\nu)>P^{(i)}_{\ell+1}(\nu)=P^{(i)}_j(\nu)$ which is in contradiction
to the convexity relation of $P^{(i)}_k(\nu)$ between $\ell\leq k\leq j$.

Thus we assume that $P^{(i)}_\ell(\nu)\leq P^{(i)}_j(\nu)$.
Then~\eqref{eq:C_1} gives $P^{(i)}_{\ell}(\nu)\leq P^{(i)}_{\ell+1}(\nu)$.
We also have
\begin{gather*}
P^{(i)}_{\ell+1}(\nu)\overset{\text{by Lemma~\ref{lem:f_singular}(1)}}{=}x_\ell+1
\overset{\text{since the string $(\ell,x_\ell)$ is non-singular}}{\leq} P^{(i)}_\ell(\nu).
\end{gather*}
Thus we have $P^{(i)}_{\ell}(\nu)=P^{(i)}_{\ell+1}(\nu)$.
Combining this with the inequality $P^{(i)}_\ell(\nu)\leq P^{(i)}_j(\nu)$ and using the convexity of $P^{(i)}_k(\nu)$
between $\ell\leq k\leq j$, we obtain
\begin{gather}
P^{(i)}_{\ell}(\nu)=P^{(i)}_{\ell+1}(\nu)=\cdots =P^{(i)}_{j}(\nu)=x_\ell+1.
\label{eq:C_5.5}
\end{gather}
Suppose if possible that $j<\ell^{(i)}$.
Since $\ell<j$ we have $P^{(i)}_{\ell+1}(\nu)\leq x_j$ by Lemma~\ref{lem:f_singular}(2).
Then we have $P^{(i)}_{\ell+1}(\nu)\leq x_j\leq P^{(i)}_{j}(\nu)$.
By~\eqref{eq:C_5.5} we obtain $x_j=P^{(i)}_{j}(\nu)$, that is, the string $(j,x_j)$ is singular.
Since $\ell^{(i-1)}\leq j<\ell^{(i)}$ this contradicts the def\/inition of $\ell^{(i)}$.
Hence we conclude that $j=\ell^{(i)}$.
\end{proof}

\begin{proof}[Proof of Proposition~\ref{prop:caseC1} for the case~(\ref{caseC_2})] Since $\ell^{(i-1)}<\ell^{(i)}$ we have
\begin{gather*}
P^{(i)}_{\ell^{(i)}-1}(\bar{\nu})=P^{(i)}_{\ell^{(i)}-1}(\nu)-1=x_\ell,
\end{gather*}
where we have used~\eqref{eq:C_5}.
From the result at the end of the proof of the previous lemma, we see that $m^{(i)}_k(\nu)=0$ for $\ell<k<\ell^{(i)}$.
Then the relation~\eqref{eq:C_5} implies that $m^{(i+1)}_k(\nu)=0$ for $\ell<k<\ell^{(i)}$ since the existence of such
a~string would imply that the relation~\eqref{eq:C_5} have to be strictly convex.
Hence we complete the proof of Proposition~\ref{prop:caseC1}.
\end{proof}

By using the relation $\bar{\widetilde{\nu}}=\widetilde{\bar{\nu}}$ we show $\bar{\widetilde{J}}=\widetilde{\bar{J}}$.
For this purpose we note the following facts which are the consequences of the proof of Proposition~\ref{lem:C_3}.

\begin{Corollary}
\label{cor:C_2}
Under the same assumptions with the previous proposition, we have the following relations.
\begin{enumerate}\itemsep=0pt
\item[$(1)$] $m^{(a)}_k(\nu)=0$ for all $\ell<k<\ell^{(i)}$ and $a=i-1,i,i+1$, \item[$(2)$] $\ell^{(i-1)}\leq\ell$.
\end{enumerate}
\end{Corollary}
\begin{proof}
(1) From $m^{(i)}_k(\nu)=m^{(i+1)}_k(\nu)=0$ for $\ell<k<\ell^{(i)}$, we see that if $m^{(i-1)}_k(\nu)>0$ for some
$\ell<k<\ell^{(i)}$, then $P^{(i)}_k(\nu)$ becomes strictly convex around such the $k$, which is in contradiction to the
relation~\eqref{eq:C_5}.
Thus we obtain the assertion.

(2) If $\ell<\ell^{(i-1)}$, we have $m^{(i-1)}_{\ell^{(i-1)}}(\nu)>0$ which is in contradiction to the previous
assertion since $\ell^{(i-1)}<\ell^{(i)}$.
Thus we have $\ell^{(i-1)}\leq\ell$.
\end{proof}

\begin{Proposition}%\label{prop:C_4}
Suppose that we have $\ell^{(i-1)}\leq\ell+1<\ell^{(i)}<\infty$, $m^{(i)}_{\ell+1}(\nu)=0$, the string $(\ell,x_\ell)$
is non-singular and $\widetilde{f}_i$ creates a~singular string.
Then we have $\bar{\widetilde{J}}=\widetilde{\bar{J}}$.
\end{Proposition}
\begin{proof}
By Proposition~\ref{prop:caseC(2)1} and Corollary~\ref{cor:C_2}(1) we see that it is enough to check the coincidence of
riggings for the strings $(\ell,x_\ell)$ and $(\ell^{(i)},x_{\ell^{(i)}})$ that are touched by~$\delta$ and
$\widetilde{f}_i$.

Let us analyze the change of the string $(\ell,x_\ell)$.
Since $\ell+1<\ell^{(i)}$ we have
\begin{gather*}
\begin{array}{@{}lllll}
(\ell,x_\ell)&\overset{\widetilde{f}_i}{\longmapsto}&(\ell+1,x_\ell-1)
&\overset{\delta}{\longmapsto}&\big(\ell,P^{(i)}_\ell(\bar{\widetilde{\nu}})\big),
\\
(\ell,x_\ell)&\overset{\delta}{\longmapsto}&(\ell,x_\ell) &\overset{\widetilde{f}_i}{\longmapsto}&(\ell,x_\ell).
\end{array}
\end{gather*}
Since $\ell^{(i-1)}\leq\ell$ by Corollary~\ref{cor:C_2}(2) we have
$P^{(i)}_\ell(\bar{\widetilde{\nu}})=P^{(i)}_\ell(\nu)-1$.
From~\eqref{eq:C_5}, we have $P^{(i)}_\ell(\nu)=x_\ell+1$.
Hence we obtain $P^{(i)}_\ell(\bar{\widetilde{\nu}})=x_\ell$.

The analysis of the string $(\ell^{(i)},x_{\ell^{(i)}})$ is the same with the corresponding part of the proof of
Proposition~\ref{prop:C_3} if we replace~\eqref{eq:C_4} there by~\eqref{eq:C_5}.
\end{proof}

\subsection{Proof for Case C (3)}\label{sec:caseC(3)}

As the f\/inal step of the proof for Case C, let us treat the case where $\ell+1<\ell^{(i)}$ and $\widetilde{f}_i$ creates
a~non-singular string.
\begin{Proposition}
\label{prop:caseC2}
Assume that $\ell^{(i-1)}\leq \ell+1<\ell^{(i)}<\infty$.
If $\widetilde{f}_i$ acting on $(\nu,J)$ creates a~non-singular string, then the following relation is satisfied:
\begin{gather}
P^{(i)}_{\ell^{(i)}-1}(\bar{\nu})>x_\ell.
\label{eq:caseC3}
\end{gather}
\end{Proposition}

\begin{proof}
In this case we have
\begin{gather}
P^{(i)}_{\ell^{(i)}-1}(\bar{\nu})=P^{(i)}_{\ell^{(i)}-1}(\nu)-1
\label{eq:C_6}
\end{gather}
since $\ell^{(i-1)}\leq\ell^{(i)}-1$.

Throughout the proof of this proposition, let $j$ be the largest integer such that $j<\ell^{(i)}$ and
$m^{(i)}_j(\nu)>0$.
Since $m^{(i)}_\ell(\nu)>0$, we have $\ell\leq j\leq\ell^{(i)}-1$.
We divide the proof into four cases.

{\bf Case $\boldsymbol{j=\ell^{(i)}-1}$.} Let us denote the corresponding string as $(\ell^{(i)}-1,x_{\ell^{(i)}-1})$.
Since $\ell^{(i-1)}<\ell^{(i)}$ we have $\ell^{(i-1)}\leq\ell^{(i)}-1$ which implies that the string
$(\ell^{(i)}-1,x_{\ell^{(i)}-1})$ is non-singular by def\/inition of $\ell^{(i)}$.
Thus we have $x_{\ell^{(i)}-1}<P^{(i)}_{\ell^{(i)}-1}(\nu)$.
On the other hand, from the assumption $\ell+1<\ell^{(i)}$, we have $x_\ell<x_{\ell^{(i)}-1}$ since $\widetilde{f}_i$
chooses the string $(\ell,x_\ell)$ which is strictly shorter than the string $(\ell^{(i)}-1,x_{\ell^{(i)}-1})$.
Combining both inequalities we obtain $x_\ell<P^{(i)}_{\ell^{(i)}-1}(\nu)-1$.
Hence we obtain $x_\ell<P^{(i)}_{\ell^{(i)}-1}(\bar{\nu})$.

{\bf Case $\boldsymbol{\ell<j<\ell^{(i)}-1}$.} In this case the string $(j,x_j)$ have to be non-singular since $j$ satisf\/ies
$\ell^{(i-1)}\leq\ell+1\leq j<\ell^{(i)}$.
Then we have
\begin{gather*}
P^{(i)}_{j}(\nu)\overset{\text{since the string $(j,x_j)$ is non-singular}}{>}x_j
\overset{\text{since $\widetilde{f}_i$ acts on the string $(\ell,x_\ell)$ although we have $j>\ell$}}{>}x_\ell .
\end{gather*}
In particular we have $P^{(i)}_j(\nu)>x_\ell+1$.
Similarly let us consider the string $\big(\ell^{(i)},P^{(i)}_{\ell^{(i)}}(\nu)\big)$ whose existence is assured by the
def\/inition of $\ell^{(i)}$.
Then its rigging must satisfy $P^{(i)}_{\ell^{(i)}}(\nu)>x_\ell$ since its length satisf\/ies $\ell^{(i)}>\ell$.
From the convexity relation of $P^{(i)}_k(\nu)$ between $j\leq k\leq\ell^{(i)}$ we have
\begin{gather*}
P^{(i)}_{\ell^{(i)}-1}(\nu)\geq \min\big\{P^{(i)}_j(\nu),P^{(i)}_{\ell^{(i)}}(\nu)\big\}>x_\ell.
\end{gather*}
Suppose if possible that $P^{(i)}_{\ell^{(i)}-1}(\nu)=x_\ell+1$.
Then we would have $P^{(i)}_j(\nu)>P^{(i)}_{\ell^{(i)}-1}(\nu) \leq P^{(i)}_{\ell^{(i)}}(\nu)$ which is forbidden by the
convexity of the vacancy numbers.
Thus we have $P^{(i)}_{\ell^{(i)}-1}(\nu)>x_\ell+1$.
By~\eqref{eq:C_6} we have $P^{(i)}_{\ell^{(i)}-1}(\bar{\nu})>x_\ell$.

{\bf Case $\boldsymbol{j=\ell}$ and $\boldsymbol{(\ell,x_\ell)}$ is non-singular.} In this case we have $P^{(i)}_\ell(\nu)>x_\ell$ since the
string $(\ell,x_\ell)$ is non-singular.
Let us consider the string $\big(\ell^{(i)},P^{(i)}_{\ell^{(i)}}(\nu)\big)$.
Then its rigging must satisfy $P^{(i)}_{\ell^{(i)}}(\nu)>x_\ell$ since its length satisfy $\ell^{(i)}>\ell$.
Now we invoke the convexity of $P^{(i)}_k(\nu)$ between $\ell=j\leq k\leq\ell^{(i)}$ to conclude that
$P^{(i)}_{\ell^{(i)}-1}(\nu)>x_\ell$.

Suppose if possible that $P^{(i)}_{\ell^{(i)}-1}(\nu)=x_\ell+1$.
Again by the convexity relation the only possibility is $P^{(i)}_\ell(\nu)=\dots=P^{(i)}_{\ell^{(i)}}=x_\ell+1$.
In particular, we have $P^{(i)}_{\ell+1}(\nu)=x_\ell+1$ and by Lemma~\ref{lem:f_singular} (1) this implies that
$\widetilde{f}_i$ makes a~singular string.
This is in contradiction to the assumption that $\widetilde{f}_i$ makes a~non-singular string.
Thus we have $P^{(i)}_{\ell^{(i)}-1}(\nu)>x_\ell+1$.
By~\eqref{eq:C_6} we have $P^{(i)}_{\ell^{(i)}-1}(\bar{\nu})>x_\ell$.

{\bf Case $\boldsymbol{j=\ell}$ and $\boldsymbol{(\ell,x_\ell)}$ is singular.} In this case we have $P^{(i)}_\ell(\nu)=x_\ell$ since the
string $(\ell,x_\ell)$ is singular.
Also the rigging for the string $\big(\ell^{(i)},P^{(i)}_{\ell^{(i)}}(\nu)\big)$ must satisfy $P^{(i)}_{\ell^{(i)}}(\nu)>x_\ell$
since $\ell^{(i)}>\ell$.
Then from the convexity relation of $P^{(i)}_k(\nu)$ between $\ell\leq k\leq \ell^{(i)}$ we have
\begin{gather*}
P^{(i)}_{\ell^{(i)}-1}(\nu)\geq \min\big\{P^{(i)}_\ell(\nu),P^{(i)}_{\ell^{(i)}}(\nu)\big\}=x_\ell.
\end{gather*}

If $P^{(i)}_{\ell^{(i)}-1}(\nu)=x_\ell$ we have
$P^{(i)}_{\ell}(\nu)=P^{(i)}_{\ell^{(i)}-1}(\nu)<P^{(i)}_{\ell^{(i)}}(\nu)$ which is in contradiction to the convexity
relation since we have $\ell<\ell^{(i)}-1$ by the assumption $\ell+1<\ell^{(i)}$.
Suppose if possible that $P^{(i)}_{\ell^{(i)}-1}(\nu)=x_\ell+1$.
Since $x_\ell<P^{(i)}_{\ell^{(i)}}(\nu)$ we have $P^{(i)}_{\ell^{(i)}-1}(\nu)\leq P^{(i)}_{\ell^{(i)}}(\nu)$.
Then the only possibility that is compatible with the convexity relation is
$x_\ell=P^{(i)}_\ell(\nu)<P^{(i)}_{\ell+1}(\nu)=\dots=P^{(i)}_{\ell^{(i)}}(\nu)=x_\ell+1$.
In particular, we have $P^{(i)}_{\ell+1}(\nu)=x_\ell+1$.
However, by Lemma~\ref{lem:f_singular}(1), this implies that $\widetilde{f}_i$ creates a~singular string, which is the
contradiction.
Thus we conclude that $P^{(i)}_{\ell^{(i)}-1}(\nu)>x_\ell+1$.
By~\eqref{eq:C_6} we have $P^{(i)}_{\ell^{(i)}-1}(\bar{\nu})>x_\ell$, which completes the whole proof of
Proposition~\ref{prop:caseC2}.
\end{proof}

\begin{Proposition}%\label{prop:C_5}
Assume that $\ell^{(i-1)}\leq \ell+1<\ell^{(i)}<\infty$ and $\widetilde{f}_i$ acting on $(\nu,J)$ creates a~non-singular
string.
Then we have the following identities:
\begin{gather*}
 (1)\quad \bar{\widetilde{\nu}}=\widetilde{\bar{\nu}} , \qquad (2)\quad \bar{\widetilde{J}}=\widetilde{\bar{J}}.
\end{gather*}
\end{Proposition}
\begin{proof}
(1) We show that $\widetilde{f}_i$ acts on the same string $(\ell,x_\ell)$ in both $(\nu,J)$ and $(\bar{\nu},\bar{J})$.
For~$\delta$ creates the strings $\big(\ell^{(i)}-1,P^{(i)}_{\ell^{(i)}-1}(\bar{\nu})\big)$ and
$\big(\ell_{(i)}-1,P^{(i)}_{\ell_{(i)}-1}(\bar{\nu})\big)$.
Then by~\eqref{eq:caseC3} we have $P^{(i)}_{\ell^{(i)}-1}(\bar{\nu})>x_\ell$ in this case.
Recall also that we have $P^{(i)}_{\ell_{(i)}-1}(\bar{\nu})>x_\ell$ if $\ell^{(i)}<\ell_{(i)}$ by
Proposition~\ref{prop:C_prep}.
Therefore the string $(\ell,x_\ell)$ remains as the string with the smallest rigging of the largest length in
$(\bar{\nu},\bar{J})$.
Thus $\widetilde{f}_i$ acts on the string $(\ell,x_\ell)$ of $(\bar{\nu},\bar{J})$.

Also, since $\widetilde{f}_i$ creates the non-singular string and does not change other coriggings,~$\delta$ chooses the
same strings for both $(\nu,J)$ and $(\widetilde{\nu},\widetilde{J})$.
Thus we have $\bar{\widetilde{\nu}}=\widetilde{\bar{\nu}}$.

(2) Since $\widetilde{f}_i$ acts on the string $(\ell,x_\ell)$ in both cases, we have to analyze the strings
$(\ell,x_\ell)$, $\big(\ell^{(i)},P^{(i)}_{\ell^{(i)}}(\nu)\big)$ and $\big(\ell_{(i)},P^{(i)}_{\ell_{(i)}}(\nu)\big)$.
The string $(\ell,x_\ell)$ behaves as follows:
\begin{gather*}
\begin{array}{@{}lllll}
(\ell,x_\ell)&\overset{\widetilde{f}_i}{\longmapsto}& (\ell+1,x_\ell-1)&\overset{\delta}{\longmapsto}&
(\ell+1,x_\ell-1),
\\
(\ell,x_\ell)&\overset{\delta}{\longmapsto}& (\ell,x_\ell)&\overset{\widetilde{f}_i}{\longmapsto}& (\ell+1,x_\ell-1).
\end{array}
\end{gather*}
Hence the assertion follows.
On the other hand, the string $\big(\ell^{(i)},P^{(i)}_{\ell^{(i)}}(\nu)\big)$ behaves as follows:
\begin{gather*}
\begin{array}{@{}lllll}
\big(\ell^{(i)},P^{(i)}_{\ell^{(i)}}(\nu)\big)&\overset{\widetilde{f}_i}{\longmapsto}&
\big(\ell^{(i)},P^{(i)}_{\ell^{(i)}}(\widetilde{\nu})\big)&\overset{\delta}{\longmapsto}&
\big(\ell^{(i)}-1,P^{(i)}_{\ell^{(i)}-1}(\bar{\widetilde{\nu}})\big),
\\
\big(\ell^{(i)},P^{(i)}_{\ell^{(i)}}(\nu)\big)&\overset{\delta}{\longmapsto}&
\big(\ell^{(i)}-1,P^{(i)}_{\ell^{(i)}-1}(\bar{\nu})\big)&\overset{\widetilde{f}_i}{\longmapsto}&
\big(\ell^{(i)}-1,P^{(i)}_{\ell^{(i)}-1}(\widetilde{\bar{\nu}})\big).
\end{array}
\end{gather*}
By $\bar{\widetilde{\nu}}=\widetilde{\bar{\nu}}$ we obtain the desired fact.
The analysis of $\big(\ell_{(i)},P^{(i)}_{\ell_{(i)}}(\nu)\big)$ is similar.
Hence we have $\bar{\widetilde{J}}=\widetilde{\bar{J}}$.
\end{proof}

\subsection{Proof for Case~D: preliminary steps}

\subsubsection{Classif\/ication}

Recall that the def\/ining condition for Case D is $\ell^{(i)}=\ell$.
For the proof, it is convenient to further divide the case into the following four cases:
\begin{enumerate}\itemsep=0pt
\item[(1)] $\ell=\ell^{(i)}<\ell_{(i)}$ and $m^{(i)}_\ell(\nu)>1$, \item[(2)] $\ell=\ell^{(i)}=\ell_{(i)}$ and
$m^{(i)}_\ell(\nu)>2$, \item[(3)] $\ell=\ell^{(i)}<\ell_{(i)}$ and $m^{(i)}_\ell(\nu)=1$, \item[(4)]
$\ell=\ell^{(i)}=\ell_{(i)}$ and $m^{(i)}_\ell(\nu)=2$.
\end{enumerate}
Note that $\ell^{(i)}=\ell_{(i)}$ automatically requires $m^{(i)}_{\ell^{(i)}}(\nu)\geq 2$.

\subsubsection{A common property for Case D}

For the proof of Case D, we will need the following result.
\begin{Proposition}
\label{prop:caseD1}
Assume that $\ell^{(i)}=\ell$.
Then the following relation holds:
\begin{gather*}%\label{eq:caseD1}
P^{(i)}_{\ell-1}(\bar{\nu})\geq x_\ell.
\end{gather*}
\end{Proposition}

\begin{proof}
Throughout the proof of this proposition, let $j$ be the largest integer such that $j<\ell$ and $m^{(i)}_j(\nu)>0$.
If there is no such $j$, set $j=0$.
We divide the proof into the following four cases:
\begin{enumerate}\itemsep=0pt
\item[(a)] $j=\ell-1$, \item[(b)] $j<\ell-1$, $\ell^{(i-1)}=\ell$, \item[(c)] $j<\ell-1$, $\ell^{(i-1)}<\ell$ and
$\ell^{(i-1)}\leq j$, \item[(d)] $j<\ell-1$, $\ell^{(i-1)}<\ell$ and $j<\ell^{(i-1)}$.
\end{enumerate}

{\bf Case (a).} Denote the corresponding string of $(\nu,J)^{(i)}$ as $(\ell-1,x_{\ell-1})$.
Then we have $x_\ell\leq x_{\ell-1}$ by def\/inition of $\ell$.
Suppose that $\ell^{(i-1)}<\ell$.
Since we have $\ell^{(i-1)}\leq\ell-1<\ell^{(i)}$ the string $(\ell-1,x_{\ell-1})$ cannot be singular.
Thus we have $P^{(i)}_{\ell-1}(\nu)>x_\ell$ by $P^{(i)}_{\ell-1}(\nu)>x_{\ell-1}\geq x_\ell$.
Thus we obtain $P^{(i)}_{\ell-1}(\bar{\nu})\geq x_\ell$ since $P^{(i)}_{\ell-1}(\bar{\nu})=P^{(i)}_{\ell-1}(\nu)-1$ by
$\ell^{(i-1)}\leq\ell-1$.

Next suppose that $\ell^{(i-1)}=\ell$.
In this case the string $(\ell-1,x_{\ell-1})$ can be singular.
Therefore we have $P^{(i)}_{\ell-1}(\nu)\geq x_\ell$ by $P^{(i)}_{\ell-1}(\nu)\geq x_{\ell-1}\geq x_\ell$.
From $\ell-1<\ell^{(i-1)}$ we have $P^{(i)}_{\ell-1}(\bar{\nu})=P^{(i)}_{\ell-1}(\nu)$.
Thus we obtain $P^{(i)}_{\ell-1}(\bar{\nu})\geq x_\ell$.

{\bf Case (b).} By $j<\ell$ we have $P^{(i)}_j(\nu)\geq x_j\geq x_\ell$.
Note that this relation is valid even if $j=0$ since we have $P^{(i)}_0(\nu)=0\geq x_\ell$ by $\ell=\ell^{(i)}>0$.
Recall that we also have $P^{(i)}_\ell(\nu)\geq x_\ell$.
Thus by the convexity of $P^{(i)}_k(\nu)$ between $j\leq k\leq\ell$ we have $P^{(i)}_{\ell-1}(\nu)\geq x_\ell$.
Since $\ell-1<\ell^{(i-1)}$ we conclude that $P^{(i)}_{\ell-1}(\bar{\nu})\geq x_\ell$.

{\bf Case (c).} We can write the def\/ining condition of case (c) as $\ell^{(i-1)}\leq j<\ell-1$.
Since $\ell^{(i-1)}\leq j<\ell^{(i)}$ the string $(j,x_j)$ is non-singular.
Thus we have $P^{(i)}_j(\nu)>x_\ell$ by $P^{(i)}_j(\nu)>x_j\geq x_\ell$.
We also have $P^{(i)}_\ell(\nu)\geq x_\ell$.
Suppose if possible that $P^{(i)}_{\ell-1}(\nu)=x_\ell$.
Then we have $P^{(i)}_j(\nu)>P^{(i)}_{\ell-1}(\nu)\leq P^{(i)}_\ell(\nu)$ which is in contradiction to the convexity
relation of $P^{(i)}_k(\nu)$ between $j\leq k\leq\ell$.
Thus we have $P^{(i)}_{\ell-1}(\nu)>x_\ell$.
By $\ell^{(i-1)}<\ell-1$ and $\ell=\ell^{(i)}$ we have $P^{(i)}_{\ell-1}(\bar{\nu})\geq x_\ell$.

{\bf Case (d).} We can write the def\/ining condition of case (d) as $j<\ell^{(i-1)}<\ell$.
In this case the string $(j,x_j)$ can be singular.
Thus we have $P^{(i)}_j(\nu)\geq x_\ell$ by $P^{(i)}_j(\nu)\geq x_j\geq x_\ell$.
We also have $P^{(i)}_\ell(\nu)\geq x_\ell$.
Suppose if possible that $P^{(i)}_{\ell-1}(\nu)=x_\ell$.
Then by the convexity relation of $P^{(i)}_k(\nu)$ between $j\leq k\leq\ell$, the only possibility is
$P^{(i)}_j(\nu)=\dots=P^{(i)}_\ell(\nu)=x_\ell$.
However this is in contradiction to the fact that $P^{(i)}_k(\nu)$ is a~strictly convex function between
$\ell^{(i-1)}-1\leq k\leq\ell^{(i-1)}+1$ due to the existence of the length $\ell^{(i-1)}$ row at $\nu^{(i-1)}$.
Thus we have $P^{(i)}_{\ell-1}(\nu)>x_\ell$.
By $\ell^{(i-1)}\leq\ell-1$ we have $P^{(i)}_{\ell-1}(\bar{\nu})\geq x_\ell$.
\end{proof}

\subsection{Proof for Case D (1)}

In this section, let us consider the f\/irst case $\ell=\ell^{(i)}<\ell_{(i)}$ and $m^{(i)}_\ell(\nu)>1$.
Throughout this subsection, let $j$ be the largest integer satisfying $j<\ell_{(i)}$ and $m^{(i)}_j(\nu)>0$.
From $\ell=\ell^{(i)}<\ell_{(i)}$, we see that $\ell=\ell^{(i)}\leq j$.
The situation depends on whether $\ell<j$ or $\ell=j$.

\subsubsection[The case $\ell<j$]{The case $\boldsymbol{\ell<j}$}

\begin{Proposition}
\label{prop:D_1}
Assume that $\ell=\ell^{(i)}<\ell_{(i)}$, $m^{(i)}_\ell(\nu)>1$ and $\ell<j$.
Then we have the following identities:
\begin{gather*}
 (1)\quad \bar{\widetilde{\nu}}=\widetilde{\bar{\nu}} , \qquad (2)\quad \bar{\widetilde{J}}=\widetilde{\bar{J}} .
\end{gather*}
\end{Proposition}
\begin{proof}
(1) Let us analyze the action of $\widetilde{f}_i$ before and after the application of~$\delta$.
In this case we can choose two distinct length $\ell$ strings $(\ell,x_\ell)$ and $\big(\ell,P^{(i)}_\ell(\nu)\big)$ of
$(\nu,J)^{(i)}$ where $\widetilde{f}_i$ will act on the former one and~$\delta$ will act on the latter one.
Here note that if $x_\ell=P^{(i)}_\ell(\nu)$ then all the riggings for the length $\ell$ strings are the same since the
minimal value of the corresponding riggings is $x_\ell$ and the maximal one is $P^{(i)}_\ell(\nu)$.

By Proposition~\ref{prop:caseD1} we have $P^{(i)}_{\ell^{(i)}-1}(\bar{\nu})\geq x_\ell$.
Thus the string $\big(\ell^{(i)}-1,P^{(i)}_{\ell^{(i)}-1}(\bar{\nu})\big)$ is shorter than $(\ell,x_\ell)$ of
$(\bar{\nu},\bar{J})^{(i)}$ and its rigging is larger than or equal to $x_\ell$.
Also, by $\ell<j$, we can use the same arguments of Proposition~\ref{prop:C_prep} to show
$P^{(i)}_{\ell_{(i)}-1}(\bar{\nu})>x_\ell$.
Thus the string $\big(\ell_{(i)}-1,P^{(i)}_{\ell_{(i)}-1}(\bar{\nu})\big)$ has the rigging that is strictly larger than
$x_\ell$.
Therefore $\widetilde{f}_i$ will act on the same string $(\ell,x_\ell)$ in both $(\nu,J)$ and $(\bar{\nu},\bar{J})$.

Let us analyze~$\delta$.
Since $\widetilde{f}_i$ will not change the corigging of the string $\big(\ell^{(i)},P^{(i)}_{\ell^{(i)}}(\nu)\big)$ of
$(\nu,J)^{(i)}$, it remains as the shortest possible string in $(\widetilde{\nu},\widetilde{J})^{(i)}$ starting from
$\widetilde{\ell}^{(i-1)}=\ell^{(i-1)}$.
So~$\delta$ chooses this string after the application of $\widetilde{f}_i$ too.
Let us show that the string $\big(\ell_{(i)},P^{(i)}_{\ell_{(i)}}(\nu)\big)$ is chosen by~$\delta$ in both before and after the
application of $\widetilde{f}_i$.
Recall that we have two length $\ell$ strings $(\ell,x_\ell)$ and $\big(\ell,P^{(i)}_\ell(\nu)\big)$ of $(\nu,J)^{(i)}$ where
$\widetilde{f}_i$ acts on the former one and~$\delta$ acts on the latter one.

(i)~Consider the case $P^{(i)}_\ell(\nu)=x_\ell$, that is, the string $(\ell,x_\ell)$ is singular.
Then the assumptions $\ell^{(i)}<\ell_{(i)}$ and $m^{(i)}_{\ell^{(i)}}(\nu)\geq 2$ implies that
$\ell^{(i)}<\ell_{(i+1)}$, as otherwise we would have $\ell^{(i)}=\ell_{(i)}$.
Then a~subtlety could occur only if $\ell_{(i)}>\ell_{(i+1)}=\ell+1$ and $\widetilde{f}_i$ creates a~singular string of
length $\ell+1$.
So suppose if possible that such a~situation happens.
In this case, we have $P^{(i)}_{\ell+1}(\nu)=x_\ell+1$ by Lemma~\ref{lem:f_singular}.

Suppose if possible that $m^{(i)}_{\ell+1}(\nu)>0$.
Then the riggings for such string have to be $x_\ell+1$ by the minimality of $x_\ell$.
Thus the length $\ell+1$ strings are singular.
However this is in contradiction to the present assumption $\ell_{(i)}>\ell+1$.
Thus we have $m^{(i)}_{\ell+1}(\nu)=0$.

Next we show $m^{(i+1)}_{\ell+1}(\nu)=1$.
Let $j'$ be the smallest integer such that $\ell+1\leq j'$ and $m^{(i)}_{j'}(\nu)>0$.
Since $m^{(i)}_{\ell+1}(\nu)=0$ and $\ell_{(i)}>\ell+1$ we have $\ell+1<j'\leq\ell_{(i)}$.
Also we have $j'\leq j$ by the assumption $\ell<j$.
Recall that we have $P^{(i)}_\ell(\nu)=x_\ell$ and $P^{(i)}_{\ell+1}(\nu)=x_\ell+1$.
By the minimality of $x_\ell$ we have $P^{(i)}_{\ell+2}(\nu)\geq x_\ell+1$ as otherwise we would have
$P^{(i)}_{j'}(\nu)\leq x_\ell$.
Let us write $P^{(i)}_{\ell+2}(\nu)=x_\ell+1+\varepsilon$ with $\varepsilon\geq 0$.
Then the inequality of Lemma~\ref{lem:convexity3} with $l=\ell+1$ reads
\begin{gather*}
-P^{(i)}_\ell(\nu)+2P^{(i)}_{\ell+1}(\nu)-P^{(i)}_{\ell+2}(\nu)=1-\varepsilon
\\
\qquad
{}\geq m^{(i-1)}_{\ell+1}(\nu)-2m^{(i)}_{\ell+1}(\nu)+m^{(i+1)}_{\ell+1}(\nu)
=m^{(i-1)}_{\ell+1}(\nu)+m^{(i+1)}_{\ell+1}(\nu).
\end{gather*}
Since we have $m^{(i+1)}_{\ell+1}(\nu)>0$ by $\ell_{(i+1)}=\ell+1$, the only possibility is $\varepsilon=0$,
$m^{(i-1)}_{\ell+1}(\nu)=0$ and $m^{(i+1)}_{\ell+1}(\nu)=1$.

Since $\ell<j'$ we have $P^{(i)}_{j'}(\nu)\geq x_{j'}>x_\ell$.
Then from the convexity relation of $P^{(i)}_k(\nu)$ between $\ell\leq k\leq j'$, the only possibility that is
compatible with $P^{(i)}_{\ell+1}(\nu)=P^{(i)}_{\ell+2}(\nu)=x_\ell+1$ is $P^{(i)}_{\ell+1}(\nu)=\dots
=P^{(i)}_{j'}(\nu)=x_\ell+1$.
Since $x_{j'}>x_\ell$, we see that $x_{j'}=P^{(i)}_{j'}(\nu)$, in particular, the corresponding string is singular.
Since $j'\geq\ell+1=\ell_{(i+1)}$, we conclude that $j'=\ell_{(i)}$.
However this relation contradicts the requirement $j'\leq j<\ell_{(i)}$.

 (ii)  Consider the case $P^{(i)}_\ell(\nu)>x_\ell$.
Recall that we have $\ell^{(i)}\leq\ell_{(i+1)}$.
Again a~subtlety could occur only if $\ell_{(i)}>\ell_{(i+1)}=\ell+1$ and $\widetilde{f}_i$ creates a~singular string of
length $\ell+1$.
So suppose if possible that such a~situation happens.
By Lemma~\ref{lem:f_singular}, we have $P^{(i)}_{\ell+1}(\nu)=x_\ell+1$.
As in the previous case, we have $m^{(i)}_{\ell+1}(\nu)=0$.
Let $j'$ be the smallest integer such that $\ell+1\leq j'$ and $m^{(i)}_{j'}(\nu)>0$.
Then, as in the previous case, we have $\ell+1<j'\leq\ell_{(i)}$ and $j'\leq j$.
Then from the convexity relation of $P^{(i)}_k(\nu)$ between $\ell\leq k\leq j'$ the only possibility that is compatible
with the relations $P^{(i)}_\ell(\nu)>x_\ell$, $P^{(i)}_{\ell+1}(\nu)=x_\ell+1$ and $P^{(i)}_{j'}(\nu)>x_\ell$ is
$P^{(i)}_\ell(\nu)=\dots =P^{(i)}_{j'}(\nu)=x_\ell+1$.
On the other hand, by $\ell<j'$ we have $P^{(i)}_{j'}(\nu)\geq x_{j'}>x_\ell$.
Therefore we have $P^{(i)}_{j'}(\nu)=x_{j'}=x_\ell+1$ and, in particular, the corresponding string is singular.
Since $\ell+1=\ell_{(i+1)}\leq j'$ we conclude that $j'=\ell_{(i)}$.
This is a~contradiction since we have $j'\leq j<\ell_{(i)}$.

Therefore~$\delta$ will act on the same strings of lengths $\ell^{(i)}$ and $\ell_{(i)}$ in both $(\nu,J)^{(i)}$ and
$(\widetilde{\nu},\widetilde{J})^{(i)}$.
Thus we conclude that $\bar{\widetilde{\nu}}=\widetilde{\bar{\nu}}$ in this case.

(2) It is enough to consider the three strings $(\ell,x_\ell)$, $\big(\ell^{(i)},P^{(i)}_{\ell^{(i)}}(\nu)\big)$ and
$\big(\ell_{(i)},P^{(i)}_{\ell_{(i)}}(\nu)\big)$.
As for $(\ell,x_\ell)$, we have
\begin{gather*}
\begin{array}{@{}lllll}
(\ell,x_\ell)&\overset{\widetilde{f}_i}{\longmapsto}& (\ell+1,x_\ell-1)&\overset{\delta}{\longmapsto}&(\ell+1,x_\ell-1),
\\
(\ell,x_\ell)&\overset{\delta}{\longmapsto}& (\ell,x_\ell)&\overset{\widetilde{f}_i}{\longmapsto}&(\ell+1,x_\ell-1),
\end{array}
\end{gather*}
and as for $\big(\ell^{(i)},P^{(i)}_{\ell^{(i)}}(\nu)\big)$ we have
\begin{gather*}
\begin{array}{@{}lllll}
\big(\ell^{(i)},P^{(i)}_{\ell^{(i)}}(\nu)\big)&\overset{\widetilde{f}_i}{\longmapsto}&
\big(\ell^{(i)},P^{(i)}_{\ell^{(i)}}(\widetilde{\nu})\big)&\overset{\delta}{\longmapsto}&
\big(\ell^{(i)}-1,P^{(i)}_{\ell^{(i)}-1}(\bar{\widetilde{\nu}})\big),
\\
\big(\ell^{(i)},P^{(i)}_{\ell^{(i)}}(\nu)\big)&\overset{\delta}{\longmapsto}&
\big(\ell^{(i)}-1,P^{(i)}_{\ell^{(i)}-1}(\bar{\nu})\big)&\overset{\widetilde{f}_i}{\longmapsto}&
\big(\ell^{(i)}-1,P^{(i)}_{\ell^{(i)}-1}(\widetilde{\bar{\nu}})\big).
\end{array}
\end{gather*}
The situation for $\big(\ell_{(i)},P^{(i)}_{\ell_{(i)}}(\nu)\big)$ is similar.
Thus we have $\bar{\widetilde{J}}=\widetilde{\bar{J}}$.
\end{proof}

\subsubsection[The case $\ell=j$]{The case $\boldsymbol{\ell=j}$}

We follow the classif\/ication of Lemma~\ref{lem:vacancy_delta}.

\begin{Proposition}
Assume that $\ell=\ell^{(i)}<\ell^{(i+1)}=\dots=\ell_{(i)}$, $m^{(i)}_\ell(\nu)>1$ and $\ell=j$.
Then we have the following identities:
\begin{gather*}
 (1)\quad \bar{\widetilde{\nu}}=\widetilde{\bar{\nu}}, \qquad (2) \quad \bar{\widetilde{J}}=\widetilde{\bar{J}}.
\end{gather*}
\end{Proposition}
\begin{proof}
(1) In this case, we can use the same arguments of the proof of Case~(V) of Proposition~\ref{prop:C_prep} to show
$P^{(i)}_{\ell_{(i)}-1}(\bar{\nu})>x_\ell$.
We also have $P^{(i)}_{\ell^{(i)}-1}(\bar{\nu})\geq x_\ell$ by Proposition~\ref{prop:caseD1}.
Thus $\widetilde{f}_i$ acts on the same string before and after the application of~$\delta$.
We notice that~$\delta$ acts on the same string before and after the application of $\widetilde{f}_i$ even if
$\widetilde{f}_i$ creates a~singular string since we are assuming that $\ell_{(i+1)}=\ell_{(i)}$.
Here recall that $\widetilde{\ell}_{(i)}$ is determined as the length of the minimal possible string compared with the
length $\widetilde{\ell}_{(i+1)}=\ell_{(i+1)}$ and the length $\ell_{(i)}$ singular string remains singular after the
application of $\widetilde{f}_i$.

Proof of (2) is the same as Proposition~\ref{prop:D_1}.
\end{proof}

\begin{Proposition}
Assume that $\ell=\ell^{(i)}\leq\ell^{(i+1)}<\ell_{(i+1)}=\ell_{(i)}$, $m^{(i)}_\ell(\nu)>1$ and $\ell=j$.
Then we have the following identities:
\begin{gather*}
 (1) \quad \bar{\widetilde{\nu}}=\widetilde{\bar{\nu}}, \qquad (2) \quad \bar{\widetilde{J}}=\widetilde{\bar{J}}.
\end{gather*}
\end{Proposition}
\begin{proof}
(1) We divide the proof into two cases according to whether $\ell=\ell_{(i)}-1$ or $\ell<\ell_{(i)}-1$.

{\bf Case 1.} Let us consider the case $\ell=\ell_{(i)}-1$.
Then we have at least three strings $(\ell,x_\ell)$, $\big(\ell,P^{(i)}_\ell(\nu)\big)$
and $\big(\ell+1,P^{(i)}_{\ell+1}(\nu)\big)$ in
$(\nu,J)^{(i)}$.
Let us analyze the action of $\widetilde{f}_i$ before and after the application of~$\delta$.
In this case, we show $P^{(i)}_{\ell_{(i)}-1}(\bar{\nu})\geq x_\ell$.
According to Lemma~\ref{lem:vacancy_delta} (VI) we have to show $P^{(i)}_{\ell_{(i)}-1}(\nu)\geq x_\ell$, which is the
consequence of $\ell=\ell_{(i)}-1$.
Recall also that $P^{(i)}_{\ell-1}(\bar{\nu})\geq x_\ell$ by Proposition~\ref{prop:caseD1}.

If $P^{(i)}_{\ell}(\bar{\nu})>x_\ell$ there is no problem since the string $(\ell,x_\ell)$ of
$(\bar{\nu},\bar{J})^{(i)}$ becomes non-singular.
So suppose that $P^{(i)}_{\ell}(\bar{\nu})=x_\ell$.
Then there are three strings $\big(\ell-1,P^{(i)}_{\ell-1}(\bar{\nu})\big)$, $(\ell,x_\ell)$, and $(\ell,x_\ell)$ in
$(\bar{\nu},\bar{J})^{(i)}$.
Thus we can choose the strings such that $\widetilde{f}_i$ acts on the same string before and after the application of
$\delta$.
Also, since we are assuming that $\ell_{(i+1)}=\ell_{(i)}$,~$\delta$ acts on the same string before and after the
application of $\widetilde{f}_i$ even if $\widetilde{f}_i$ creates a~singular string.
Thus we have $\bar{\widetilde{\nu}}=\widetilde{\bar{\nu}}$ in this case.

{\bf Case 2.} Let us consider the case $\ell<\ell_{(i)}-1$.
In this case, we can use the same arguments of Case (VI) of the proof of Proposition~\ref{prop:C_prep} to show
$P^{(i)}_{\ell_{(i)}-1}(\bar{\nu})>x_\ell$.
Thus we see that $\widetilde{f}_i$ acts on the same string before and after the application of~$\delta$.
Also, by the assumption $\ell_{(i+1)}=\ell_{(i)}$, we see that~$\delta$ acts on the same string before and after the
application of $\widetilde{f}_i$.
Thus we have $\bar{\widetilde{\nu}}=\widetilde{\bar{\nu}}$ in this case.

Proof of (2) is the same as Proposition~\ref{prop:D_1}.
\end{proof}

\begin{Proposition}
\label{prop:D(1)}
Assume that $\ell=\ell^{(i)}\leq\ell_{(i+1)}<\ell_{(i)}$, $m^{(i)}_\ell(\nu)>1$ and $\ell=j$.
Then we have the following identities:
\begin{gather*}
 (1)\quad \bar{\widetilde{\nu}}=\widetilde{\bar{\nu}}, \qquad
 (2) \quad \bar{\widetilde{J}}=\widetilde{\bar{J}}.
\end{gather*}
\end{Proposition}
\begin{proof} {\bf Step 1.} Suppose that we have $P^{(i)}_{\ell_{(i)}-1}(\bar{\nu})>x_\ell$.
Since we have $P^{(i)}_{\ell-1}(\bar{\nu})\geq x_\ell$ by Proposition~\ref{prop:caseD1}, we see that $\widetilde{f}_i$
acts on the same string before and after~$\delta$.
Let us consider the behavior of~$\delta$ before and after $\widetilde{f}_i$.
Then a~subtlety could occur only if $\widetilde{f}_i$ creates a~singular string which satisf\/ies
$\ell_{(i)}>\ell+1\geq\ell_{(i+1)}$.

Thus suppose if possible that $P^{(i)}_{\ell_{(i)}-1}(\bar{\nu})>x_\ell$, $\widetilde{f}_i$ creates a~singular string of
length $\ell+1$ and $\ell_{(i)}>\ell+1\geq\ell_{(i+1)}$.
From Lemma~\ref{lem:vacancy_delta} (VII) we have $P^{(i)}_{\ell_{(i)}-1}(\nu)>x_\ell+1$ and from
Lemma~\ref{lem:f_singular} we have $P^{(i)}_{\ell+1}(\nu)=x_\ell+1$.
If $\ell_{(i)}-1=\ell+1$ this is a~contradiction.
Thus suppose that $\ell_{(i)}-1>\ell+1$.
Recall that we have $P^{(i)}_\ell(\nu)\geq x_\ell$ by def\/inition of the rigged conf\/igurations.
Suppose if possible that $P^{(i)}_\ell(\nu)>x_\ell$.
Then we have $P^{(i)}_\ell(\nu)\geq P^{(i)}_{\ell+1}(\nu)<P^{(i)}_{\ell_{(i)}-1}(\nu)$ which violates the convexity
relation of $P^{(i)}_k(\nu)$ between $\ell\leq k\leq \ell_{(i)}$.
Thus we have $P^{(i)}_\ell(\nu)=x_\ell$.
Suppose if possible that $\ell_{(i+1)}=\ell$.
Since $P^{(i)}_\ell(\nu)=x_\ell$ we see that all length $\ell$ strings of $(\nu,J)^{(i)}$ are singular.
This is a~contradiction since we have $m^{(i)}_\ell(\nu)>1$ and $\ell_{(i)}>\ell_{(i+1)}$.
Thus we have $\ell_{(i+1)}=\ell+1$.
From $P^{(i)}_{\ell+1}(\nu)=x_\ell+1$ and $P^{(i)}_{\ell_{(i)}-1}(\nu)>x_\ell+1$ we have
\begin{gather*}
P^{(i)}_{\ell+2}(\nu)\geq\min\big\{P^{(i)}_{\ell+1}(\nu),P^{(i)}_{\ell_{(i)}-1}(\nu)\big\}=x_\ell+1.
\end{gather*}
Let us write $P^{(i)}_{\ell+2}(\nu)=x_\ell+1+\varepsilon$ where $\varepsilon\geq 0$.
From Lemma~\ref{lem:convexity3} with $l=\ell+1$, we have
\begin{gather*}
-P^{(i)}_{\ell}(\nu)+2P^{(i)}_{\ell+1}(\nu)-P^{(i)}_{\ell+2}(\nu)=1-\varepsilon
\\
\qquad
{}\geq m^{(i-1)}_{\ell+1}(\nu)-2m^{(i)}_{\ell+1}(\nu)+m^{(i+1)}_{\ell+1}(\nu)
=
m^{(i-1)}_{\ell+1}(\nu)+m^{(i+1)}_{\ell+1}(\nu).
\end{gather*}
Since we have $m^{(i+1)}_{\ell+1}(\nu)>0$ by $\ell_{(i+1)}=\ell+1$, the only possibility is $\varepsilon=0$,
$m^{(i-1)}_{\ell+1}(\nu)=0$ and $m^{(i+1)}_{\ell+1}(\nu)=1$.
In particular, we have $P^{(i)}_{\ell+2}(\nu)=x_\ell+1$.
However this is a~contradiction since we have $P^{(i)}_{\ell+1}(\nu)=P^{(i)}_{\ell+2}(\nu)<P^{(i)}_{\ell_{(i)}-1}(\nu)$
which violates the convexity relation of $P^{(i)}_k(\nu)$ between $\ell\leq k\leq \ell_{(i)}$.
Hence this case cannot happen.
Therefore we have $\bar{\widetilde{\nu}}=\widetilde{\bar{\nu}}$ and $\bar{\widetilde{J}}=\widetilde{\bar{J}}$ in this
case.

{\bf Step 2.} Suppose that we have $\ell=\ell_{(i)}-1$.
By def\/inition of the rigged conf\/iguration we have $P^{(i)}_\ell(\nu)\geq x_\ell$.
Then by Lemma~\ref{lem:vacancy_delta} (VII) we have $P^{(i)}_\ell(\bar{\nu})\geq x_\ell-1$.
Suppose if possible that $P^{(i)}_\ell(\bar{\nu})=x_\ell-1$.
Recall that we have the three strings $(\ell,x_\ell)$, $\big(\ell,P^{(i)}_\ell(\nu)\big)$ and $\big(\ell+1,P^{(i)}_{\ell+1}(\nu)\big)$
of $(\nu,J)^{(i)}$.
Then we have the remaining string $(\ell,x_\ell)$ in $(\bar{\nu},\bar{J})^{(i)}$.
This is a~contradiction since we have $P^{(i)}_\ell(\bar{\nu})<x_\ell$.

Thus suppose that $P^{(i)}_\ell(\bar{\nu})\geq x_\ell$.
As in Case 1 of the proof of the previous proposition we can show that $\widetilde{f}_i$ chooses the same string before
and after~$\delta$.
Let us consider the behavior of~$\delta$ before and after $\widetilde{f}_i$.
In this case, there is no problem even if $\widetilde{f}_i$ creates a~singular string since we have $\ell+1=\ell_{(i)}$.
Hence we have $\bar{\widetilde{\nu}}=\widetilde{\bar{\nu}}$ and $\bar{\widetilde{J}}=\widetilde{\bar{J}}$.

{\bf Step 3.} We assume that $\ell<\ell_{(i)}-1$ and $P^{(i)}_{\ell_{(i)}-1}(\bar{\nu})\leq x_\ell$ in the
following proof.
According to Lemma~\ref{lem:vacancy_delta} (VII) we have $P^{(i)}_{\ell_{(i)}-1}(\nu)\leq x_\ell+1$.
Recall that we have $P^{(i)}_{\ell_{(i)}}(\nu)>x_\ell$ by $\ell_{(i)}>\ell$ and $P^{(i)}_{\ell}(\nu)\geq x_\ell$ by
def\/inition of the rigged conf\/igurations.
Since we are assuming that $\ell=\ell^{(i)}=j<\ell_{(i)}-1$, the convexity relation of the vacancy numbers allows the
following three possibilities:
\begin{enumerate}\itemsep=0pt
\item[(1)] $\ell=\ell_{(i)}-2$, $P^{(i)}_{\ell}(\nu)=x_\ell$, $P^{(i)}_{\ell_{(i)}-1}(\nu)=x_\ell+1$ and
$P^{(i)}_{\ell_{(i)}}(\nu)=x_\ell+2$, \item[(2)] $P^{(i)}_{\ell}(\nu)=x_\ell$ and $P^{(i)}_{\ell+1}(\nu)=\dots
=P^{(i)}_{\ell_{(i)}}(\nu)=x_\ell+1$, \item[(3)] $P^{(i)}_{\ell}(\nu)=P^{(i)}_{\ell+1}(\nu)=\dots
=P^{(i)}_{\ell_{(i)}}(\nu)=x_\ell+1$.
\end{enumerate}
In the following we consider them case by case.

{\bf Case (1).} According to Lemma~\ref{lem:convexity3} with $l=\ell_{(i)}-1$, we have
\begin{gather*}
-P^{(i)}_{\ell}(\nu)+2P^{(i)}_{\ell_{(i)}-1}(\nu)-P^{(i)}_{\ell_{(i)}}(\nu)=0
\\
\qquad
{}\geq m^{(i-1)}_{\ell_{(i)}-1}(\nu)-2m^{(i)}_{\ell_{(i)}-1}(\nu)+m^{(i+1)}_{\ell_{(i)}-1}(\nu)
=
m^{(i-1)}_{\ell_{(i)}-1}(\nu)+m^{(i+1)}_{\ell_{(i)}-1}(\nu)
\end{gather*}
where we have used $m^{(i)}_{\ell_{(i)}-1}(\nu)=0$ by $j=\ell$.
Thus we have $m^{(i-1)}_{\ell_{(i)}-1}(\nu)=m^{(i+1)}_{\ell_{(i)}-1}(\nu)=0$.
Since we are assuming that $\ell_{(i+1)}<\ell_{(i)}$, we conclude that $\ell^{(i+1)}=\ell_{(i+1)}=\ell$.

Now recall that we are assuming that $m^{(i)}_\ell(\nu)\geq 2$.
Then the minimality of $x_\ell$ and the assumption $P^{(i)}_{\ell}(\nu)=x_\ell$ imply that all the length $\ell$ strings
of $(\nu,J)^{(i)}$ are $\big(\ell,P^{(i)}_{\ell}(\nu)\big)$, in particular, they are singular.
However this is in contradiction to the assumption $\ell_{(i+1)}<\ell_{(i)}$ since we know that $\ell_{(i+1)}=\ell$
which forces that $\ell_{(i)}=\ell$.
Thus this case cannot happen.

{\bf Case (2).} The assumption $P^{(i)}_{\ell+1}(\nu)=\dots =P^{(i)}_{\ell_{(i)}}(\nu)$ implies that
$m^{(i-1)}_k(\nu)=m^{(i+1)}_k(\nu)=0$ for all $\ell+2\leq k\leq\ell_{(i)}-1$ as otherwise the above relation for the
vacancy numbers would be strictly convex.
Next, we apply Lemma~\ref{lem:convexity3} with $l=\ell +1$.
Then we have
\begin{gather*}
-P^{(i)}_{\ell}(\nu)+2P^{(i)}_{\ell+1}(\nu)-P^{(i)}_{\ell+2}(\nu)=1
\\
\qquad
{}\geq m^{(i-1)}_{\ell+1}(\nu)-2m^{(i)}_{\ell+1}(\nu)+m^{(i+1)}_{\ell+1}(\nu)
=m^{(i-1)}_{\ell+1}(\nu)+m^{(i+1)}_{\ell+1}(\nu),
\end{gather*}
where we have used $m^{(i)}_{\ell+1}(\nu)=0$ by $j=\ell$.

Suppose if possible that $m^{(i-1)}_{\ell+1}(\nu)=1$.
Then we have $m^{(i+1)}_{\ell+1}(\nu)=0$.
Since we are assuming that $\ell_{(i+1)}<\ell_{(i)}$, $m^{(i+1)}_k(\nu)=0$ for all $\ell+1\leq k\leq\ell_{(i)}-1$
implies that $\ell_{(i+1)}=\ell$.
However, as in the previous case, the assumption $P^{(i)}_{\ell}(\nu)=x_\ell$ implies that all the length $\ell$ strings
of $(\nu,J)^{(i)}$ are singular.
Thus we must have $\ell_{(i)}=\ell$ by $\ell_{(i+1)}=\ell$ and $m^{(i)}_\ell(\nu)>1$.
This is a~contradiction.
Hence this case cannot happen.

Therefore we assume that $m^{(i-1)}_{\ell+1}(\nu)=0$ and $m^{(i+1)}_{\ell+1}(\nu)=1$ in the sequel.
Then we can show that $\ell^{(i+1)}=\ell$ and $\ell_{(i+1)}=\ell+1$.
Since $m^{(i+1)}_k(\nu)=0$ for all $\ell+2\leq k\leq\ell_{(i)}-1$, we have $\ell_{(i+1)}=\ell$ or $\ell+1$ by
$\ell\leq\ell_{(i+1)}$.
Suppose if possible that $\ell_{(i+1)}=\ell$.
Since we are assuming that $P^{(i)}_\ell(\nu)=x_\ell$, the length $\ell=\ell^{(i)}$ strings of $(\nu,J)^{(i)}$ are
always equal to $(\ell,x_\ell)$ by the minimality of $x_\ell$.
Then by the assumption $m^{(i)}_\ell(\nu)\geq 2$, we must have $\ell_{(i)}=\ell$ by $\ell_{(i+1)}=\ell$.
This is in contradiction to the other assumption $\ell^{(i)}<\ell_{(i)}$.
Thus we have $\ell_{(i+1)}=\ell+1$.
Then from $m^{(i+1)}_{\ell+1}(\nu)=1$, we conclude that $\ell^{(i+1)}=\ell$ since there is the requirement
$\ell\leq\ell^{(i+1)}\leq\ell_{(i+1)}$.

Let us analyze $\widetilde{f}_i\circ\delta$ f\/irst.
Since $\widetilde{f}_i$ acts on the $i$-th place, we have $\widetilde{\bar{\nu}}^{(a)}=\bar{\nu}^{(a)}$ for all $a\neq
i$.
Next we consider $\widetilde{\bar{\nu}}^{(i)}$.
Recall that we have $P^{(i)}_{\ell^{(i)}-1}(\bar{\nu})\geq x_\ell$ by Proposition~\ref{prop:caseD1}.
Also, since we have $P^{(i)}_{\ell_{(i)}-1}(\nu)=x_\ell+1$, we have $P^{(i)}_{\ell_{(i)}-1}(\bar{\nu})=x_\ell$ by
Lemma~\ref{lem:vacancy_delta} (VII).
Therefore the string $\big(\ell_{(i)}-1,P^{(i)}_{\ell_{(i)}-1}(\bar{\nu})\big)$ of $(\bar{\nu},\bar{J})^{(i)}$ has the largest
length among the strings with smallest rigging, hence $\widetilde{f}_i$ acts on it.
To summarize, $\widetilde{\bar{\nu}}^{(i)}$ is obtained from $\nu^{(i)}$ by removing one box from a~length~$\ell^{(i)}$
row.

Next let us analyze $\delta\circ\widetilde{f}_i$.
Since $\widetilde{f}_i$ does not change the coriggings of untouched strings, and since we are assuming that
$m^{(i)}_\ell(\nu)\geq 2$, we have $\widetilde{\ell}^{(a)}=\ell^{(a)}$ for all~$a$ and
$\widetilde{\ell}_{(a)}=\ell_{(a)}$ for all $a\geq i+1$.
Let us consider $\widetilde{\ell}_{(i)}$.
Recall that we have $\widetilde{\ell}_{(i+1)}=\ell_{(i+1)}=\ell+1$ in this case.
Since we have $P^{(i)}_{\ell+1}(\nu)=x_\ell+1$, $\widetilde{f}_i$ creates a~singular string of length $\ell+1$ by
Lemma~\ref{lem:f_singular}.
Thus the length $\ell+1$ string created by $\widetilde{f}_i$ is the shortest possible string starting from
$\widetilde{\ell}_{(i+1)}$.
Therefore we have $\widetilde{\ell}_{(i)}=\ell+1$.
Since we have $m^{(i-1)}_k(\nu)=0$ for all $\ell+1\leq k\leq\ell_{(i)}-1$ and since $\widetilde{f}_i$ does not change
the coriggings of untouched strings, we have $\widetilde{\ell}_{(i-1)}=\ell_{(i-1)}$ and thus
$\widetilde{\ell}_{(a)}=\ell_{(a)}$ for all $a\leq i-1$.
To summarize, $\bar{\widetilde{\nu}}^{(i)}$ is obtained from $\nu^{(i)}$ by removing one box from a~length $\ell^{(i)}$
row and $\bar{\widetilde{\nu}}^{(a)}=\bar{\nu}^{(a)}$ for all $a\neq i$.
Hence we have $\bar{\widetilde{\nu}}=\widetilde{\bar{\nu}}$.

Finally, let us show $\bar{\widetilde{J}}=\widetilde{\bar{J}}$.
For this it is enough to check the three strings $\big(\ell^{(i)},P^{(i)}_{\ell^{(i)}}(\nu)\big)$, $(\ell,x_\ell)$ and
$\big(\ell_{(i)},P^{(i)}_{\ell_{(i)}}(\nu)\big)$ of $(\nu,J)^{(i)}$.
Here note that we have in fact $\big(\ell^{(i)},P^{(i)}_{\ell^{(i)}}(\nu)\big)=(\ell,x_\ell)$ under the present assumption.
To begin with, the string $\big(\ell^{(i)},P^{(i)}_{\ell^{(i)}}(\nu)\big)$ behaves as follows:
\begin{gather*}
\begin{array}{@{}lllll}
\big(\ell^{(i)},P^{(i)}_{\ell^{(i)}}(\nu)\big)&\overset{\widetilde{f}_i}{\longmapsto}&
\big(\ell^{(i)},P^{(i)}_{\ell^{(i)}}(\widetilde{\nu})\big)&\overset{\delta}{\longmapsto}&
\big(\ell^{(i)}-1,P^{(i)}_{\ell^{(i)}-1}(\bar{\widetilde{\nu}})\big),
\\
\big(\ell^{(i)},P^{(i)}_{\ell^{(i)}}(\nu)\big)&\overset{\delta}{\longmapsto}&
\big(\ell^{(i)}-1,P^{(i)}_{\ell^{(i)}-1}(\bar{\nu})\big)&\overset{\widetilde{f}_i}{\longmapsto}&
\big(\ell^{(i)}-1,P^{(i)}_{\ell^{(i)}-1}(\widetilde{\bar{\nu}})\big).
\end{array}
\end{gather*}
By $\bar{\widetilde{\nu}}=\widetilde{\bar{\nu}}$ we see the coincidence of the riggings.
Next, let us consider the string $(\ell,x_\ell)$:
\begin{gather*}
\begin{array}{@{}lllll}
(\ell,x_\ell)&\overset{\widetilde{f}_i}{\longmapsto}& (\ell+1,x_\ell-1)&\overset{\delta}{\longmapsto}&
\big(\ell,P^{(i)}_{\ell}(\bar{\widetilde{\nu}})\big),
\\
(\ell,x_\ell)&\overset{\delta}{\longmapsto}& (\ell,x_\ell)&\overset{\widetilde{f}_i}{\longmapsto}&(\ell,x_\ell).
\end{array}
\end{gather*}
Recall that we have $\widetilde{\ell}^{(i-1)}\leq\widetilde{\ell}^{(i)}=\widetilde{\ell}^{(i+1)}=\ell$ and
$\ell+1=\widetilde{\ell}_{(i+1)}=\widetilde{\ell}_{(i)}<\widetilde{\ell}_{(i-1)}$.
Therefore we have $P^{(i)}_{\ell}(\bar{\widetilde{\nu}})=P^{(i)}_\ell(\nu)=x_\ell$.
Let us consider the remaining string $\big(\ell_{(i)},P^{(i)}_{\ell_{(i)}}(\nu)\big)$.
Then we have
\begin{gather*}
\begin{array}{@{}lllll}
\big(\ell_{(i)},P^{(i)}_{\ell_{(i)}}(\nu)\big)&\overset{\widetilde{f}_i}{\longmapsto}&
\big(\ell_{(i)},P^{(i)}_{\ell_{(i)}}(\widetilde{\nu})\big)&\overset{\delta}{\longmapsto}&
\big(\ell_{(i)},P^{(i)}_{\ell_{(i)}}(\widetilde{\nu})\big),
\\
\big(\ell_{(i)},P^{(i)}_{\ell_{(i)}}(\nu)\big)&\overset{\delta}{\longmapsto}&
\big(\ell_{(i)}-1,P^{(i)}_{\ell_{(i)}-1}(\bar{\nu})\big)&\overset{\widetilde{f}_i}{\longmapsto}&
\big(\ell_{(i)},P^{(i)}_{\ell_{(i)}-1}(\bar{\nu})-1\big).
\end{array}
\end{gather*}
Since we have $\ell+1\leq\ell_{(i)}$ by the assumption, we have
$P^{(i)}_{\ell_{(i)}}(\widetilde{\nu})=P^{(i)}_{\ell_{(i)}}(\nu)-2=x_\ell-1$.
Also, since we are assuming that $\ell_{(i+1)}<\ell_{(i)}$, we have
$P^{(i)}_{\ell_{(i)}-1}(\bar{\nu})=P^{(i)}_{\ell_{(i)}-1}(\nu)-1$ by Lemma~\ref{lem:vacancy_delta} (VII).
Thus we have $P^{(i)}_{\ell_{(i)}-1}(\bar{\nu})-1=P^{(i)}_{\ell_{(i)}}(\nu)-2=x_\ell-1$.
Therefore we have checked $\bar{\widetilde{J}}=\widetilde{\bar{J}}$ in this case.

{\bf Case (3).} The assumption $P^{(i)}_\ell(\nu)=\dots =P^{(i)}_{\ell_{(i)}}(\nu)$ implies that
$m^{(i-1)}_k(\nu)=m^{(i+1)}_k(\nu)=0$ for all $\ell+1\leq k\leq\ell_{(i)}-1$.
Since we are assuming that $\ell_{(i+1)}<\ell_{(i)}$, the only possibility is $\ell_{(i+1)}=\ell$.

Let us analyze $\widetilde{f}_i\circ\delta$ f\/irst.
From Lemma~\ref{lem:vacancy_delta} (VII) we have $P^{(i)}_{\ell_{(i)}-1}(\bar{\nu})=x_\ell$ by
$P^{(i)}_{\ell_{(i)}-1}(\nu)=x_\ell+1$.
Recall that we have $P^{(i)}_{\ell^{(i)}-1}(\bar{\nu})\geq x_\ell$ by Proposition~\ref{prop:caseD1}.
From the assumption $\ell\leq\ell_{(i)}-1$, we see that the string $\big(\ell_{(i)}-1,P^{(i)}_{\ell_{(i)}-1}(\bar{\nu})\big)$
has the smallest rigging of the largest length.
Thus the next $\widetilde{f}_i$ acts on this string.
To summarize, we have $\widetilde{\bar{\nu}}^{(a)}=\bar{\nu}^{(a)}$ for all $a\neq i$ and $\widetilde{\bar{\nu}}^{(i)}$
is obtained by removing one box from the length $\ell^{(i)}$ row.

Next we analyze $\delta\circ\widetilde{f}_i$.
Since $\widetilde{f}_i$ does not change coriggings of untouched strings, and since we are assuming that
$m^{(i)}_\ell(\nu)>1$, we have $\widetilde{\ell}^{(a)}=\ell^{(a)}$ for all~$a$ and $\widetilde{\ell}_{(a)}=\ell_{(a)}$
for all $a\geq i+1$.
In particular, we have $\widetilde{\ell}_{(i+1)}=\ell_{(i+1)}=\ell$ by the assumption.
In order to determine $\widetilde{\ell}_{(i)}$, recall that the relation $P^{(i)}_{\ell+1}(\nu)=x_\ell+1$ implies that
$\widetilde{f}_i$ creates the singular string of length $\ell+1$ by Lemma~\ref{lem:f_singular}.
Thus we have $\widetilde{\ell}_{(i)}=\ell+1$.
Since we have $m^{(i-1)}_k(\nu)=0$ for all $\ell+1\leq k\leq\ell_{(i)}-1$, we have
$\widetilde{\ell}_{(i-1)}=\ell_{(i-1)}$ by $\ell_{(i-1)}\geq\ell_{(i)}$, thus $\widetilde{\ell}_{(a)}=\ell_{(a)}$ for
all $a\leq i-1$.
To summarize, we have $\bar{\widetilde{\nu}}^{(a)}=\bar{\nu}^{(a)}$ for all $a\neq i$ and $\bar{\widetilde{\nu}}^{(i)}$
is obtained by removing one box from the length $\ell^{(i)}$ row.
Hence we have $\widetilde{\bar{\nu}}=\bar{\widetilde{\nu}}$.

Let us show $\widetilde{\bar{J}}=\bar{\widetilde{J}}$ in this case.
Again it is enough to check the three strings $\big(\ell^{(i)},P^{(i)}_{\ell^{(i)}}(\nu)\big)$, $(\ell,x_\ell)$ and
$\big(\ell_{(i)},P^{(i)}_{\ell_{(i)}}(\nu)\big)$ of $(\nu,J)^{(i)}$.
The analysis for the string $\big(\ell^{(i)},P^{(i)}_{\ell^{(i)}}(\nu)\big)$ is the same with the previous Case (2).
Let us analyze the string $(\ell,x_\ell)$.
Then it behaves as
\begin{gather*}
\begin{array}{@{}lllll}
(\ell,x_\ell)&\overset{\widetilde{f}_i}{\longmapsto}& (\ell+1,x_\ell-1)&\overset{\delta}{\longmapsto}&
\big(\ell,P^{(i)}_{\ell}(\bar{\widetilde{\nu}})\big),
\\
(\ell,x_\ell)&\overset{\delta}{\longmapsto}& (\ell,x_\ell)&\overset{\widetilde{f}_i}{\longmapsto}& (\ell,x_\ell).
\end{array}
\end{gather*}
In this case, we have $\widetilde{\ell}^{(i-1)}\leq\widetilde{\ell}^{(i)}=\widetilde{\ell}^{(i+1)}=\ell$ and
$\ell=\widetilde{\ell}_{(i+1)}<\widetilde{\ell}_{(i)}\leq\widetilde{\ell}_{(i-1)}$.
Therefore we have $P^{(i)}_{\ell}(\bar{\widetilde{\nu}})=P^{(i)}_\ell(\nu)-1=x_\ell$.
The analysis of $\big(\ell_{(i)},P^{(i)}_{\ell_{(i)}}(\nu)\big)$ is the same with the previous Case~(2).
Thus we have $\widetilde{\bar{J}}=\bar{\widetilde{J}}$ in this case.
This completes the proof of proposition.
\end{proof}

\begin{Example}
Here we give an example for Case~(2) of Step~3 of the proof of Proposition~\ref{prop:D(1)}.
Consider the following rigged conf\/iguration $(\nu,J)$ of type $(B^{1,1})^{\otimes 3}\otimes B^{1,3}\otimes
B^{2,1}\otimes B^{2,2}\otimes B^{3,1}$ of $D^{(1)}_5$
\begin{center}
\unitlength 10pt
{\small
\begin{picture}(40,6)
\put(-0.8,2.1){2}
\multiput(-0.8,3.1)(0,1){2}{1}
\put(-1.6,5.1){$-1$}
\put(0,2){\Yboxdim10pt\yng(5,2,2,1)}
\put(1.2,2.1){$-1$}
\put(2.2,3.1){0}
\put(2.2,4.1){1}
\put(5.2,5.1){$-1$}
\put(8.5,0){
\put(-0.8,0.1){0}
\put(-0.8,1.1){0}
\multiput(-1.6,2.1)(0,1){3}{$-1$}
\put(-0.8,5.1){0}
\put(0,0){\Yboxdim10pt\yng(5,2,2,2,1,1)}
\put(1.2,0.1){0}
\put(1.2,1.1){0}
\put(2.2,2.1){$-1$}
\put(2.2,3.1){$-1$}
\put(2.2,4.0){$-1$}
\put(5.1,5.0){0}
}
\put(17,0){
\multiput(-0.8,0.1)(0,1){3}{0}
\put(-0.8,3.1){2}
\put(-0.8,4.1){1}
\put(-0.8,5.1){3}
\put(0,0){\Yboxdim10pt\yng(5,3,2,1,1,1)}
\put(1.2,0.1){0}
\put(1.2,1.1){0}
\put(1.2,2.1){0}
\put(2.2,3.1){2}
\put(3.2,4.1){1}
\put(5.2,5.1){3}
}
\put(26.0,3){
\put(-0.8,0.1){0}
\put(-1.6,1.1){$-1$}
\put(-1.6,2.1){$-3$}
\put(0,0){\Yboxdim10pt\yng(5,2,1)}
\put(1.2,0.1){0}
\put(2.2,1.1){$-1$}
\put(5.2,2.1){$-3$}
}
\put(35.0,4){
\put(-0.8,0.1){1}
\put(-1.6,1.1){$-1$}
\put(0,0){\Yboxdim10pt\yng(5,2)}
\put(2.2,0.1){1}
\put(5.2,1.1){$-2$}
}
\end{picture}}
%\{\{5,2,2,1\},\{5,2,2,2,1,1\},\{5,3,2,1,1,1\},\{5,2,1\},\{5,2\}\}\\
%\{\{-1,1,0,-1\},\{0,-1,-1,-1,0,0\},\{3,1,2,0,0,0\},\{-3,-1,0\},\{-2,1\}\}
\end{center}
The corresponding tensor product $\Phi^{-1}(\nu,J)$ is
\begin{gather*}
\Yboxdim14pt \Yvcentermath1 \young(\mone)\otimes \young(\mone)\otimes \young(5)\otimes
\young(1\mfour\mfour)\otimes \young(\mfive,\mthree)\otimes \young(15,3\mfour)\otimes \young(2,3,4)
\end{gather*}
Then $(\bar{\nu},\bar{J})$ is
\begin{center}
\unitlength 10pt
{\small
\begin{picture}(40,6)
\put(-0.8,2.1){1}
\multiput(-0.8,3.1)(0,1){2}{1}
\put(-0.8,5.1){0}
\put(0,2){\Yboxdim10pt\yng(4,2,1,1)}
\put(1.2,2.1){$-1$}
\put(1.2,3.1){1}
\put(2.2,4.1){0}
\put(4.2,5.1){0}
\put(8.5,0){
\put(-0.8,0.1){0}
\put(-0.8,1.1){0}
\put(-0.8,2.1){0}
\multiput(-1.6,3.1)(0,1){3}{$-1$}
\put(0,0){\Yboxdim10pt\yng(4,2,2,1,1,1)}
\put(1.2,0.1){0}
\put(1.2,1.1){0}
\put(1.2,2.1){0}
\put(2.2,3.1){$-1$}
\put(2.2,4.0){$-1$}
\put(4.1,5.0){$-1$}
}
\put(17,0){
\multiput(-0.8,0.1)(0,1){3}{0}
\put(-0.8,3.1){0}
\put(-0.8,4.1){1}
\put(-0.8,5.1){3}
\put(0,0){\Yboxdim10pt\yng(5,2,1,1,1,1)}
\put(1.2,0.1){0}
\put(1.2,1.1){0}
\put(1.2,2.1){0}
\put(1.2,3.1){0}
\put(2.2,4.1){1}
\put(5.2,5.1){3}
}
\put(26.0,3){
\put(-0.8,0.1){0}
\put(-0.8,1.1){0}
\put(-1.6,2.1){$-3$}
\put(0,0){\Yboxdim10pt\yng(5,1,1)}
\put(1.2,0.1){0}
\put(1.2,1.1){0}
\put(5.2,2.1){$-3$}
}
\put(35.0,4){
\put(-0.8,0.1){2}
\put(-1.6,1.1){$-1$}
\put(0,0){\Yboxdim10pt\yng(5,1)}
\put(1.2,0.1){2}
\put(5.2,1.1){$-2$}
}
\end{picture}}
%\{\{4,2,1,1\},\{4,2,2,1,1,1\},\{5,2,1,1,1,1\},\{5,1,1\},\{5,1\}\}\\
%\{\{0,0,1,-1\},\{-1,-1,-1,0,0,0\},\{3,1,0,0,0,0\},\{-3,0,0\},\{-2,2\}\}
\end{center}
and $(\widetilde{\nu},\widetilde{J})$ is
\begin{center}
\unitlength 10pt
{\small
\begin{picture}(40,6)
\put(-0.8,2.1){2}
\multiput(-0.8,3.1)(0,1){2}{1}
\put(-0.8,5.1){0}
\put(0,2){\Yboxdim10pt\yng(5,2,2,1)}
\put(1.2,2.1){$-1$}
\put(2.2,3.1){0}
\put(2.2,4.1){1}
\put(5.2,5.1){0}
\put(8.5,0){
\put(-0.8,0.1){0}
\put(-0.8,1.1){0}
\multiput(-1.6,2.1)(0,1){2}{$-1$}
\put(-1.6,4.1){$-2$}
\put(-1.6,5.1){$-2$}
\put(0,0){\Yboxdim10pt\yng(5,3,2,2,1,1)}
\put(1.2,0.1){0}
\put(1.2,1.1){0}
\put(2.2,2.1){$-1$}
\put(2.2,3.1){$-1$}
\put(3.2,4.0){$-2$}
\put(5.1,5.0){$-2$}
}
\put(17,0){
\multiput(-0.8,0.1)(0,1){3}{0}
\put(-0.8,3.1){2}
\put(-0.8,4.1){2}
\put(-0.8,5.1){4}
\put(0,0){\Yboxdim10pt\yng(5,3,2,1,1,1)}
\put(1.2,0.1){0}
\put(1.2,1.1){0}
\put(1.2,2.1){0}
\put(2.2,3.1){2}
\put(3.2,4.1){2}
\put(5.2,5.1){4}
}
\put(26.0,3){
\put(-0.8,0.1){0}
\put(-1.6,1.1){$-1$}
\put(-1.6,2.1){$-3$}
\put(0,0){\Yboxdim10pt\yng(5,2,1)}
\put(1.2,0.1){0}
\put(2.2,1.1){$-1$}
\put(5.2,2.1){$-3$}
}
\put(35.0,4){
\put(-0.8,0.1){1}
\put(-1.6,1.1){$-1$}
\put(0,0){\Yboxdim10pt\yng(5,2)}
\put(2.2,0.1){1}
\put(5.2,1.1){$-2$}
}
\end{picture}}
\end{center}
In this example, we have $\ell=\ell^{(2)}=2<\ell_{(3)}=3<\ell_{(2)}=5$, $m^{(2)}_2(\nu)=3$ and $m^{(2)}_k(\nu)=0$ for
$\ell^{(2)}<k<\ell_{(2)}$.
Moreover, we have $\ell<\ell_{(2)}-1$ and $P^{(2)}_{\ell_{(2)}-1}(\bar{\nu})=x_\ell$ where $x_\ell=-1$.
This example satisf\/ies the condition for Case (2)
\begin{gather*}
P^{(2)}_2(\nu)=x_\ell=-1,
\qquad
P^{(2)}_3(\nu)=P^{(2)}_4(\nu)=P^{(2)}_5(\nu)=x_\ell+1=0.
\end{gather*}
As discussed in the proof, we have $m^{(1)}_4(\nu)=m^{(3)}_4(\nu)=0$ where $\ell+2=\ell_{(2)}-1=4$, $m^{(1)}_3(\nu)=0$
and $m^{(3)}_3(\nu)=1$.
Finally $(\bar{\widetilde{\nu}},\bar{\widetilde{J}})=(\widetilde{\bar{\nu}},\widetilde{\bar{J}})$ is
\begin{center}
\unitlength 10pt
{\small
\begin{picture}(40,6)
\put(-0.8,2.1){1}
\multiput(-0.8,3.1)(0,1){2}{1}
\put(-0.8,5.1){0}
\put(0,2){\Yboxdim10pt\yng(4,2,1,1)}
\put(1.2,2.1){$-1$}
\put(1.2,3.1){1}
\put(2.2,4.1){0}
\put(4.2,5.1){0}
\put(8.5,0){
\put(-0.8,0.1){0}
\put(-0.8,1.1){0}
\put(-0.8,2.1){0}
\multiput(-1.6,3.1)(0,1){2}{$-1$}
\put(-1.6,5.1){$-2$}
\put(0,0){\Yboxdim10pt\yng(5,2,2,1,1,1)}
\put(1.2,0.1){0}
\put(1.2,1.1){0}
\put(1.2,2.1){0}
\put(2.2,3.1){$-1$}
\put(2.2,4.0){$-1$}
\put(5.1,5.0){$-2$}
}
\put(17,0){
\multiput(-0.8,0.1)(0,1){3}{0}
\put(-0.8,3.1){0}
\put(-0.8,4.1){1}
\put(-0.8,5.1){4}
\put(0,0){\Yboxdim10pt\yng(5,2,1,1,1,1)}
\put(1.2,0.1){0}
\put(1.2,1.1){0}
\put(1.2,2.1){0}
\put(1.2,3.1){0}
\put(2.2,4.1){1}
\put(5.2,5.1){4}
}
\put(26.0,3){
\put(-0.8,0.1){0}
\put(-0.8,1.1){0}
\put(-1.6,2.1){$-3$}
\put(0,0){\Yboxdim10pt\yng(5,1,1)}
\put(1.2,0.1){0}
\put(1.2,1.1){0}
\put(5.2,2.1){$-3$}
}
\put(35.0,4){
\put(-0.8,0.1){2}
\put(-1.6,1.1){$-1$}
\put(0,0){\Yboxdim10pt\yng(5,1)}
\put(1.2,0.1){2}
\put(5.2,1.1){$-2$}
}
\end{picture}}
%\{\{4,2,1,1\},\{5,2,2,1,1,1\},\{5,2,1,1,1,1\},\{5,1,1\},\{5,1\}\}\\
%\{\{0,0,1,-1\},\{-2,-1,-1,0,0,0\},\{4,1,0,0,0,0\},\{-3,0,0\},\{-2,2\}\}
\end{center}
\end{Example}

\subsection{Proof for Case D (2)}

\begin{Proposition}
Suppose that $\ell=\ell^{(i)}=\ell_{(i)}$ and $m^{(i)}_\ell(\nu)>2$.
Then we have the following identities:
\begin{gather*}
 (1)\quad \bar{\widetilde{\nu}}=\widetilde{\bar{\nu}},\qquad
 (2)\quad \bar{\widetilde{J}}=\widetilde{\bar{J}}.
\end{gather*}
\end{Proposition}
\begin{proof}
By the assumption $m^{(i)}_\ell(\nu)\geq 3$ we can choose the three distinct strings
\begin{gather*}
\big(\ell^{(i)},P^{(i)}_{\ell^{(i)}}(\nu)\big)=\big(\ell_{(i)},P^{(i)}_{\ell_{(i)}}(\nu)\big)\text{ and }(\ell,x_\ell)
\end{gather*}
of the $i$-th conf\/iguration.
Here we may have $P^{(i)}_{\ell^{(i)}}(\nu)=x_\ell$.

Let us analyze $\widetilde{f}_i\circ\delta$.
After the application of~$\delta$ the above three strings become
\begin{gather*}
\big(\ell^{(i)}-1,P^{(i)}_{\ell^{(i)}-1}(\bar{\nu})\big)=\big(\ell_{(i)}-1,P^{(i)}_{\ell_{(i)}-1}(\bar{\nu})\big) \text{ and }(\ell,x_\ell).
\end{gather*}
By Proposition~\ref{prop:caseD1} we have $P^{(i)}_{\ell^{(i)}-1}(\bar{\nu})=P^{(i)}_{\ell_{(i)}-1}(\bar{\nu})\geq
x_\ell$.
Thus the string $(\ell,x_\ell)$ remains as the string with smallest rigging of the largest length.
Hence $\widetilde{f}_i$ will act on it.
To summarize $\widetilde{\bar{\nu}}^{(a)}=\bar{\nu}^{(a)}$ for all $a\neq i$ and $\widetilde{\bar{\nu}}^{(i)}$ is
obtained by removing a~box from each of the two length $\ell^{(i)}=\ell_{(i)}$ rows and adding a~box to the length
$\ell$ row.

Let us analyze $\delta\circ\widetilde{f}_i$.
Since $\widetilde{f}_i$ does not change coriggings of untouched strings, we have $\widetilde{\ell}^{(a)}=\ell^{(a)}$ for
all~$a$ and $\widetilde{\ell}_{(a)}=\ell_{(a)}$ for all $a\geq i+1$.
Since the string $\big(\ell_{(i)},P^{(i)}_{\ell_{(i)}}(\nu)\big)$ remains as the singular string which is longer than
$\widetilde{\ell}_{(i+1)}$, we have $\widetilde{\ell}_{(i)}=\ell_{(i)}$ and hence $\widetilde{\ell}_{(a)}=\ell_{(a)}$
for all $a\leq i$.
Thus we obtain $\bar{\widetilde{\nu}}=\widetilde{\bar{\nu}}$.
The remaining statement $\bar{\widetilde{J}}=\widetilde{\bar{J}}$ follows from
$\bar{\widetilde{\nu}}=\widetilde{\bar{\nu}}$.
\end{proof}

\subsection{Proof for Case D (3)}

\subsubsection{Outline}

Let us consider the third case.
We prove the following properties in the latter part of this subsection by case by case analysis.
\begin{Proposition}
\label{prop:caseD2}
Assume that $\ell^{(i)}=\ell$ and $m^{(i)}_\ell(\nu)=1$.
Then the following relations hold:
\begin{gather}
P^{(i)}_{\ell+1}(\nu)=x_\ell+1,
\label{eq:caseD2}
\\
P^{(i)}_{\ell-1}(\bar{\nu})=x_\ell,
\label{eq:caseD3}
\\
\ell+1\leq\ell^{(i+1)}.
\label{eq:caseD4}
\end{gather}
\end{Proposition}

Next let us show the following property:
\begin{Proposition}
\label{prop:caseD2_2}
Assume that $\ell=\ell^{(i)}<\ell_{(i)}$ and $m^{(i)}_\ell(\nu)=1$.
Then the following relation holds:
\begin{gather*}
P^{(i)}_{\ell_{(i)}-1}(\bar{\nu})>x_\ell.
\end{gather*}
\end{Proposition}
\begin{proof}
We follow the classif\/ication of Lemma~\ref{lem:vacancy_delta}.
Since we are assuming that $\ell^{(i)}<\ell_{(i)}$, we have to deal with Cases~(V), (VI) and (VII).
Throughout the proof of this proposition, let $j$ be the largest integer such that $j<\ell_{(i)}$ and
$m^{(i)}_j(\nu)>0$.
Since we are assuming that $\ell^{(i)}<\ell_{(i)}$, we have $\ell^{(i)}\leq j$.

{\bf Case (V).} In this case, we have to show that $P^{(i)}_{\ell_{(i)}-1}(\nu)\geq x_\ell$ under the assumption
$\ell^{(i)}<\ell^{(i+1)}=\dots=\ell_{(i)}$.
We can use the same argument of the corresponding part of Proposition~\ref{prop:C_prep} to show the assertion.

{\bf Case (VI).} In this case, we have to show that $P^{(i)}_{\ell_{(i)}-1}(\nu)>x_\ell$ under the assumption
$\ell^{(i+1)}<\ell_{(i+1)}=\ell_{(i)}$.
Suppose if possible that $P^{(i)}_{\ell_{(i)}-1}(\nu)\leq x_\ell$.

Let us consider the case $j=\ell_{(i)}-1$.
Note that by~\eqref{eq:caseD4}, which is shown independently to the present proposition, we have $\ell<\ell^{(i+1)}$.
Combining with the current assumption $\ell^{(i+1)}<\ell_{(i)}$ we see that $\ell=\ell^{(i)}<\ell_{(i)}-1$.
Thus the string $(\ell_{(i)}-1,x_{\ell_{(i)}-1})$ of $(\nu,J)^{(i)}$ is dif\/ferent from the string $(\ell,x_\ell)$.
Then we have $x_{\ell_{(i)}-1}\leq P^{(i)}_{\ell_{(i)}-1}(\nu)\leq x_\ell$.
However this is in contradiction to the relation $x_{\ell_{(i)}-1}>x_\ell$ which follows from the minimality of $x_\ell$
under the condition $\ell_{(i)}-1>\ell$.

Next let us consider the case $j<\ell_{(i)}-1$.
Since $\ell^{(i)}<\ell_{(i)}$ we have $P^{(i)}_{\ell_{(i)}}(\nu)\geq x_{\ell_{(i)}}>x_\ell$ by the minimality of
$x_\ell$.
Combining with the assumption $P^{(i)}_{\ell_{(i)}-1}(\nu)\leq x_\ell$ we have
$P^{(i)}_{\ell_{(i)}-1}(\nu)<P^{(i)}_{\ell_{(i)}}(\nu)$.
Then from the convexity relation of $P^{(i)}_k(\nu)$ between $j\leq k\leq\ell_{(i)}$, we have $x_j\leq
P^{(i)}_j(\nu)<P^{(i)}_{\ell_{(i)}-1}(\nu)\leq x_\ell$.
This is in contradiction to the relation $x_j\geq x_\ell$ that follows from the minimality of $x_\ell$.

{\bf Case (VII).} In this case, we have to show that $P^{(i)}_{\ell_{(i)}-1}(\nu)>x_\ell+1$ under the assumption
\mbox{$\ell_{(i+1)}<\ell_{(i)}$}.
Suppose if possible that $P^{(i)}_{\ell_{(i)}-1}(\nu)\leq x_\ell+1$.
Again recall that we have \mbox{$\ell<\ell^{(i+1)}$} from~\eqref{eq:caseD4} which is shown independently to the present
proposition.
Since we have \mbox{$\ell^{(i+1)}\leq\ell_{(i+1)}$} by def\/inition, we have the relations
\mbox{$\ell=\ell^{(i)}<\ell_{(i+1)}<\ell_{(i)}$}.

We divide the proof of Case (VII) into three steps according to the position of $j$.

{\bf Step 1.} First, let us consider the case $j=\ell_{(i)}-1$.
Then there is the corresponding string $(\ell_{(i)}-1,x_{\ell_{(i)}-1})$.
Since $\ell<\ell_{(i)}-1$, we have $x_\ell<x_{\ell_{(i)}-1}\leq P^{(i)}_{\ell_{(i)}-1}(\nu)$ by the minimality of~$x_\ell$.
Combining with the assumption $P^{(i)}_{\ell_{(i)}-1}(\nu)\leq x_\ell+1$ we conclude that
$x_{\ell_{(i)}-1}=P^{(i)}_{\ell_{(i)}-1}(\nu)$.
Then the string $(\ell_{(i)}-1,x_{\ell_{(i)}-1})$ is singular and its length satisf\/ies
$\ell_{(i+1)}\leq\ell_{(i)}-1<\ell_{(i)}$.
This contradicts the def\/inition of $\ell_{(i)}$.

{\bf Step 2.} Next, let us consider the case $\ell=\ell^{(i)}<j<\ell_{(i)}-1$.
Since we are assuming that $\ell^{(i)}<\ell_{(i)}$ we have $P^{(i)}_{\ell_{(i)}}(\nu)\geq x_{\ell_{(i)}}>x_\ell$ by the
minimality of $x_\ell$.
Similarly, from $\ell<j$, we have $P^{(i)}_j(\nu)>x_\ell$.
Then, by the convexity of $P^{(i)}_k(\nu)$ between $j\leq k\leq\ell_{(i)}$, the only possibility that is compatible with
the assumption $P^{(i)}_{\ell_{(i)}-1}(\nu)\leq x_\ell+1$ is the situation
$P^{(i)}_j(\nu)=\dots=P^{(i)}_{\ell_{(i)}}(\nu)=x_\ell+1$.
\begin{enumerate}\itemsep=0pt
\item
Consider the case $\ell_{(i+1)}\leq j$.
Then, combining with the relation $P^{(i)}_j(\nu)\geq x_j>x_\ell$, we obtain $P^{(i)}_j(\nu)=x_j$.
In particular, the string $(j,x_j)$ is singular and its length satisf\/ies that $\ell_{(i+1)}\leq j<\ell_{(i)}$.
This contradicts the def\/inition of $\ell_{(i)}$.
\item
Consider the case $j<\ell_{(i+1)}$.
Since we have $m^{(i+1)}_{\ell_{(i+1)}}(\nu)>0$, the relation for $P^{(i)}_k(\nu)$ between $\ell_{(i+1)}-1\leq
k\leq\ell_{(i+1)}+1$ must be strictly convex.
This contradicts the above relation.
\end{enumerate}

{\bf Step 3.} Finally, let us consider the case $\ell=\ell^{(i)}=j$.
Since we are assuming that $m^{(i)}_\ell(\nu)=1$, there is the only one singular string $(\ell,x_\ell)$.
In particular, we have $P^{(i)}_\ell(\nu)=x_\ell$.
Recall that we have $P^{(i)}_{\ell+1}(\nu)=x_\ell+1$ by~\eqref{eq:caseD2} which is shown independently to the present
proposition.
Then the convexity relation of $P^{(i)}_k(\nu)$ between $j\leq k\leq\ell_{(i)}$ under the condition
$P^{(i)}_{\ell_{(i)}}(\nu)>x_\ell$ admits the following two possibilities:
\begin{enumerate}\itemsep=0pt
\item
The case $\ell^{(i)}+2=\ell_{(i)}$, $P^{(i)}_\ell(\nu)=x_\ell$, $P^{(i)}_{\ell+1}(\nu)=x_\ell+1$ and
$P^{(i)}_{\ell_{(i)}}(\nu)=x_\ell+2$.
Since the relation for the vacancy numbers is linear, this situation can happen only if $m^{(i+1)}_{\ell+1}(\nu)=0$.
Recall that the relation $\ell=\ell^{(i)}<\ell_{(i+1)}<\ell_{(i)}$ in this case implies that $\ell_{(i+1)}=\ell+1$, in
particular, $m^{(i+1)}_{\ell+1}(\nu)>0$.
This is a~contradiction.
\item
The case $P^{(i)}_\ell(\nu)=x_\ell$ and $P^{(i)}_{\ell+1}(\nu)=\dots=P^{(i)}_{\ell_{(i)}}(\nu)=x_\ell+1$.
Since the relation for $P^{(i)}_k(\nu)$ between $\ell+1\leq k\leq \ell_{(i)}$ is linear, we have $m^{(i+1)}_k(\nu)=0$
for all $\ell+1<k<\ell_{(i)}$.
Next, applying Lemma~\ref{lem:convexity3} with $l=\ell+1$ we have
\begin{gather*}
-P^{(i)}_{\ell}(\nu)+2P^{(i)}_{\ell+1}(\nu)-P^{(i)}_{\ell+2}(\nu)=1
\\
\qquad
{}\geq m^{(i-1)}_{\ell+1}(\nu)-2m^{(i)}_{\ell+1}(\nu)+m^{(i+1)}_{\ell+1}(\nu)
=m^{(i-1)}_{\ell+1}(\nu)+m^{(i+1)}_{\ell+1}(\nu),
\end{gather*}
where we have used the fact that $m^{(i)}_{\ell+1}(\nu)=0$ which follows from the assumption $j=\ell$ and
$\ell+2\leq\ell_{(i)}$.
Since we know that $m^{(i+1)}_k(\nu)=0$ for all $\ell+1<k<\ell_{(i)}$, the relation
$\ell=\ell^{(i)}<\ell_{(i+1)}<\ell_{(i)}$ with $m^{(i+1)}_{\ell_{(i+1)}}(\nu)>0$ forces that $\ell_{(i+1)}=\ell+1$ and
$m^{(i+1)}_{\ell+1}(\nu)=1$.
Now recall the relation $\ell<\ell^{(i+1)}$ of~\eqref{eq:caseD4} which is shown independently to the present
proposition.
By the assumption $\ell_{(i+1)}<\ell_{(i)}$ and the relation $\ell^{(i+1)}\leq\ell_{(i+1)}$ that follows from the
def\/inition, we also have $\ell=\ell^{(i)}<\ell^{(i+1)}<\ell_{(i)}$ with $m^{(i+1)}_{\ell^{(i+1)}}(\nu)>0$.
By def\/inition of the operation~$\delta$, the two strings $(\ell^{(i+1)},P^{(i)}_{\ell^{(i+1)}}(\nu))$ and
$(\ell_{(i+1)},P^{(i)}_{\ell_{(i+1)}}(\nu))$ of $(\nu,J)^{(i+1)}$ must be distinct.
However, as we have seen in the previous paragraph, there is only one string $(k,x_k)$ of $(\nu,J)^{(i+1)}$ under the
restriction $\ell<k<\ell_{(i)}$.
This is a~contradiction.
\end{enumerate}
In conclusion, we obtain $P^{(i)}_{\ell_{(i)}-1}(\nu)>x_\ell+1$ in this Case (VII).
\end{proof}

\begin{proof}
[Proof of ``Propositions~\ref{prop:caseD2} \&~\ref{prop:caseD2_2}
$\boldsymbol{\Longrightarrow\bar{\widetilde{\nu}}=\widetilde{\bar{\nu}}}$''] $\mathstrut$
{\bf Step 1.} Let us consider the case $\widetilde{f}_i\circ\delta$.
The f\/irst~$\delta$ creates the strings $\big(\ell-1,P^{(i)}_{\ell-1}(\bar{\nu})\big)$ and
$\big(\ell_{(i)}-1,P^{(i)}_{\ell_{(i)}-1}(\bar{\nu})\big)$ of $(\bar{\nu},\bar{J})^{(i)}$.
Recall that the rigging of the latter string satisf\/ies $P^{(i)}_{\ell_{(i)}-1}(\bar{\nu})>x_\ell$ by
Proposition~\ref{prop:caseD2_2}.
Then, since the operation~$\delta$ does not change the riggings of unchanged strings, the relation~\eqref{eq:caseD3} and
$m^{(i)}_\ell(\nu)=1$ imply that the string $(\ell-1,x_\ell)$ is the longest strings among the string with rigging
$x_\ell$, hence $\widetilde{f}_i$ will act on this string.

  {\bf Step 2.} Let us consider the case $\delta\circ\widetilde{f}_i$.
Recall that by def\/inition of $\widetilde{f}_i$ we have $\widetilde{\ell}^{(a)}=\ell^{(a)}$ for all $a<i$.
From Lemma~\ref{lem:f_singular}(1) the relation~\eqref{eq:caseD2} implies that the string $(\ell+1,x_\ell-1)$ of
$(\widetilde{\nu},\widetilde{J})^{(i)}$ created by $\widetilde{f}_i$ is singular.
Since $m^{(i)}_\ell(\nu)=1$ we know that there is no singular string of $(\widetilde{\nu},\widetilde{J})^{(i)}$ that is
shorter than or equal to $\ell$ and longer than $\widetilde{\ell}^{(i-1)}$.
Thus we have $\widetilde{\ell}^{(i)}=\ell+1$.
Then, since~$\widetilde{f}_i$ does not change coriggings of $(\widetilde{\nu},\widetilde{J})^{(i+1)}$, the
relation~\eqref{eq:caseD4} implies that $\widetilde{\ell}^{(a)}=\ell^{(a)}$ for all $i<a$ as well as
$\widetilde{\ell}_{(a)}=\ell_{(a)}$ for all~$a$.

 {\bf Step 3.} To summarize, we have $\bar{\widetilde{\nu}}^{(a)}= \widetilde{\bar{\nu}}^{(a)}=\bar{\nu}^{(a)}$
for all $a\neq i$ and $\bar{\widetilde{\nu}}^{(i)}=\widetilde{\bar{\nu}}^{(i)}$ are obtained by removing a~box from
a~length $\ell_{(i)}$ singular string of $\nu^{(i)}$.
Hence we have $\bar{\widetilde{\nu}}=\widetilde{\bar{\nu}}$.
\end{proof}

Proofs of $\bar{\widetilde{J}}=\widetilde{\bar{J}}$ will be given by case by case analysis assuming the relation
$\bar{\widetilde{\nu}}=\widetilde{\bar{\nu}}$ along with the proof of Proposition~\ref{prop:caseD2}.

\subsubsection[Proof of Proposition~\ref{prop:caseD2} for the case $\ell^{(i-1)}<\ell$]{Proof
of Proposition~\ref{prop:caseD2} for the case $\boldsymbol{\ell^{(i-1)}<\ell}$}

We divide the proof into two cases depending on whether $\ell^{(i-1)}<\ell$ or
$\ell^{(i-1)}=\ell$.
In this subsection, we consider the case $\ell^{(i-1)}<\ell$.

To begin with, let us show the following property.
\begin{Proposition}
\label{prop:caseD(3)1}
Suppose that $\ell^{(i-1)}<\ell=\ell^{(i)}$ and $m^{(i)}_\ell(\nu)=1$.
Then we have
\begin{gather}
\label{eq:D_1}
P^{(i)}_{\ell-1}(\nu)=x_\ell+1,
\qquad
P^{(i)}_{\ell}(\nu)=x_\ell,
\qquad
P^{(i)}_{\ell+1}(\nu)=x_\ell+1.
\end{gather}
Moreover we have
\begin{gather}
\label{eq:D_1_1}
m^{(i-1)}_\ell(\nu)=m^{(i+1)}_\ell(\nu)=0.
\end{gather}
\end{Proposition}

For the proof, we prepare the following lemma.

\begin{Lemma}
\label{lem:caseD(3)1}
Suppose that $\ell^{(i-1)}<\ell=\ell^{(i)}$ and $m^{(i)}_\ell(\nu)=1$.
Then we have the following relations:
\begin{gather*}
 (1)\quad P^{(i)}_{\ell-1}(\nu)>P^{(i)}_\ell(\nu), \qquad (2) \quad P^{(i)}_{\ell}(\nu)<P^{(i)}_{\ell+1}(\nu).
\end{gather*}
\end{Lemma}
\begin{proof}
(1) Note that from the assumptions $m^{(i)}_\ell(\nu)=1$ and $\ell=\ell^{(i)}$, we see that the string $(\ell,x_\ell)$
is singular which implies that $x_\ell=P^{(i)}_\ell(\nu)$.
Let $j$ be the largest integer such that $j\leq\ell-1$ and $m^{(i)}_j(\nu)>0$.
If there is no such~$j$, set $j=0$ and $x_j=0$.

We divide the proof into two cases depending on whether $j=\ell-1$ or $j<\ell-1$.

 {\bf Step 1.} Let us consider the case $j=\ell-1$.
Then the string $(\ell-1,x_{\ell-1})$ is non-singular since its length satisf\/ies $\ell^{(i-1)}\leq\ell-1<\ell^{(i)}$.
Thus we have $P^{(i)}_{\ell-1}(\nu)>x_{\ell-1}$.
Recall that we have $x_{\ell-1}\geq x_\ell$ by the minimality of $x_\ell$.
Then $x_\ell=P^{(i)}_\ell(\nu)$ implies that $P^{(i)}_{\ell-1}(\nu)>P^{(i)}_\ell(\nu)$ in this case.

 {\bf Step 2.} Let us consider the case $j<\ell-1$.
Then we have $x_j\geq x_\ell$ by def\/inition of~$x_\ell$ and~$\ell>0$.
Thus we have $P^{(i)}_j(\nu)\geq x_\ell$ by $P^{(i)}_j(\nu)\geq x_j$.
Since $x_\ell=P^{(i)}_\ell(\nu)$ we obtain $P^{(i)}_j(\nu)\geq P^{(i)}_\ell(\nu)$.
Then by the convexity of $P^{(i)}_k(\nu)$ between $j\leq k\leq\ell$ we obtain
$P^{(i)}_{\ell-1}(\nu)\geq\min\{P^{(i)}_j(\nu),P^{(i)}_\ell(\nu)\} =P^{(i)}_\ell(\nu)$.
Suppose if possible that $P^{(i)}_{\ell-1}(\nu)=P^{(i)}_\ell(\nu)$.
Then by the convexity of the vacancy numbers the only possibility that is compatible with $P^{(i)}_j(\nu)\geq
P^{(i)}_\ell(\nu)$ is the case $P^{(i)}_{j}(\nu)=\dots=P^{(i)}_\ell(\nu)$.
\begin{enumerate}\itemsep=0pt
\item[(i)] Suppose that $\ell^{(i-1)}\leq j$.
Recall that the string $(j,x_j)$ satisf\/ies $x_j\leq P^{(i)}_j(\nu)=P^{(i)}_\ell(\nu)=x_\ell$.
On the other hand, by the minimality of $x_\ell$, we have $x_j\geq x_\ell$.
Thus $x_j=P^{(i)}_j(\nu)$.
Then the length of the singular string $(j,x_j)$ satisf\/ies $\ell^{(i-1)}\leq j<\ell^{(i)}$ which is in contradiction to
the def\/inition of $\ell^{(i)}$.
Thus we conclude that $P^{(i)}_{\ell-1}(\nu)>P^{(i)}_\ell(\nu)$ if $\ell^{(i-1)}\leq j$.
\item[(ii)] Suppose that $j<\ell^{(i-1)}$.
Then the relation of $P^{(i)}_k(\nu)$ between $\ell^{(i-1)}-1\leq k\leq\ell^{(i-1)}+1$ must be strictly convex due to
the existence of the length $\ell^{(i-1)}$ string at $\nu^{(i-1)}$.
This is in contradiction to the relation $P^{(i)}_j(\nu)=\dots=P^{(i)}_\ell(\nu)$.
\end{enumerate}
Thus we conclude that $P^{(i)}_{\ell-1}(\nu)>P^{(i)}_\ell(\nu)$.

(2) Since $\widetilde{f}_i$ does not vanish, the rigging for the new string created by $\widetilde{f}_i$ does not exceed
the corresponding vacancy number.
Therefore, from the proof of Lemma~\ref{lem:f_singular}(1), we have $P^{(i)}_{\ell+1}(\nu)\geq x_\ell+1$.
Since the string $(\ell,x_\ell)$ is singular we have $P^{(i)}_\ell(\nu)=x_\ell$.
Combining both relations we obtain $P^{(i)}_{\ell}(\nu)<P^{(i)}_{\ell+1}(\nu)$.
\end{proof}

\begin{proof}
[Proof of Proposition~\ref{prop:caseD(3)1}] We apply Lemma~\ref{lem:convexity2} with $m^{(i)}_\ell(\nu)=1$.
Then $P^{(i)}_{\ell-1}(\nu)$, $P^{(i)}_{\ell}(\nu)$ and $P^{(i)}_{\ell+1}(\nu)-2$ have to satisfy the convex relation.
The only possibility that is compatible with the relation
$P^{(i)}_{\ell-1}(\nu)>P^{(i)}_{\ell}(\nu)<P^{(i)}_{\ell+1}(\nu)$ is
$P^{(i)}_{\ell-1}(\nu)=P^{(i)}_{\ell}(\nu)+1=P^{(i)}_{\ell+1}(\nu)$.

From Lemma~\ref{lem:convexity3} we have
\begin{gather*}
-P^{(i)}_{\ell-1}(\nu)+2P^{(i)}_{\ell}(\nu)-P^{(i)}_{\ell+1}(\nu)=-2
\\
\qquad
{}\geq m^{(i-1)}_\ell(\nu)-2m^{(i)}_\ell(\nu)+m^{(i+1)}_\ell(\nu)
=m^{(i-1)}_\ell(\nu)-2+m^{(i+1)}_\ell(\nu),
\end{gather*}
that is, $0\geq m^{(i-1)}_\ell(\nu)+m^{(i+1)}_\ell(\nu)$.
Since $m^{(i-1)}_\ell(\nu)\geq 0$ and $m^{(i+1)}_\ell(\nu)\geq 0$ we conclude that
$m^{(i-1)}_\ell(\nu)=m^{(i+1)}_\ell(\nu)=0$.
\end{proof}

\begin{proof}
[Proof of Proposition~\ref{prop:caseD2} for the case $\boldsymbol{\ell^{(i-1)}<\ell}$] By~\eqref{eq:D_1} we have
$P^{(i)}_{\ell+1}(\nu)=x_\ell+1$, that is,~\eqref{eq:caseD2}.
By the assumption $\ell^{(i-1)}<\ell$ we see $\ell^{(i-1)}\leq\ell-1$.
Thus we obtain $P^{(i)}_{\ell-1}(\bar{\nu})= P^{(i)}_{\ell-1}(\nu)-1=x_\ell$, that is,~\eqref{eq:caseD3}.
Finally, from~\eqref{eq:D_1_1} we have $m^{(i+1)}_\ell(\nu)=0$.
Therefore we have $\ell<\ell^{(i+1)}$, that is,~\eqref{eq:caseD4}.
\end{proof}

In view of the proof of $\bar{\widetilde{\nu}}=\widetilde{\bar{\nu}}$, we can show
$\bar{\widetilde{J}}=\widetilde{\bar{J}}$ for this case as follows.

\begin{Proposition}%\label{prop:D_2}
Suppose that $\ell^{(i-1)}<\ell=\ell^{(i)}<\ell_{(i)}$ and $m^{(i)}_\ell(\nu)=1$.
Then we have $\bar{\widetilde{J}}=\widetilde{\bar{J}}$.
\end{Proposition}
\begin{proof}
From~\eqref{eq:D_1_1} it is enough to analyze the string $(\ell,x_\ell)$ of $(\nu,J)^{(i)}$ where
$x_\ell=P^{(i)}_\ell(\nu)$.
Then we have
\begin{gather*}
\begin{array}{@{}lllll}
(\ell,x_\ell)&\overset{\widetilde{f}_i}{\longmapsto}& (\ell+1,x_\ell-1)&\overset{\delta}{\longmapsto}&
\big(\ell,P^{(i)}_\ell(\bar{\widetilde{\nu}})\big),
\\
(\ell,x_\ell)&\overset{\delta}{\longmapsto}&
\big(\ell-1,P^{(i)}_{\ell-1}(\bar{\nu})\big)&\overset{\widetilde{f}_i}{\longmapsto}&\big(\ell,P^{(i)}_{\ell-1}(\bar{\nu})-1\big).
\end{array}
\end{gather*}
Recall that we have $\ell^{(i-1)}<\ell=\ell^{(i)}<\ell^{(i+1)}$.
Then we have $P^{(i)}_\ell(\bar{\widetilde{\nu}})=P^{(i)}_\ell(\nu)-1=x_\ell-1$ since we have $P^{(i)}_\ell(\nu)=x_\ell$
by~\eqref{eq:D_1}.
Similarly we have $P^{(i)}_{\ell-1}(\bar{\nu})-1=P^{(i)}_{\ell-1}(\nu)-2=x_\ell-1$ since we have
$P^{(i)}_{\ell-1}(\nu)=x_\ell+1$ by~\eqref{eq:D_1}.
\end{proof}

\subsubsection[Proof of Proposition~\ref{prop:caseD2} for the case $\ell^{(i-1)}=\ell$]{Proof
of Proposition~\ref{prop:caseD2} for the case $\boldsymbol{\ell^{(i-1)}=\ell}$}

We begin by showing the following property.
\begin{Proposition}
\label{prop:caseD(3)2}
Suppose that $\ell^{(i-1)}=\ell=\ell^{(i)}$ and $m^{(i)}_\ell(\nu)=1$.
Then we have
\begin{gather}
\label{eq:D_2}
P^{(i)}_{\ell-1}(\nu)=P^{(i)}_{\ell}(\nu)=x_\ell,
\qquad
P^{(i)}_{\ell+1}(\nu)=x_\ell+1.
\end{gather}
Moreover, we have
\begin{gather}
\label{eq:D_2_1}
m^{(i-1)}_\ell(\nu)=1,
\qquad
m^{(i+1)}_\ell(\nu)=0.
\end{gather}
\end{Proposition}

For the proof, we prepare the following lemma.

\begin{Lemma}
\label{lem:caseD(3)2}
Suppose that $\ell^{(i-1)}=\ell=\ell^{(i)}$ and $m^{(i)}_\ell(\nu)=1$.
Then we have the following relations:
\begin{gather*}
 (1) \quad P^{(i)}_{\ell-1}(\nu)\geq P^{(i)}_\ell(\nu), \qquad (2) \quad P^{(i)}_{\ell}(\nu)<P^{(i)}_{\ell+1}(\nu).
\end{gather*}
\end{Lemma}
\begin{proof}
(1) Let $j$ be the largest integer such that $j\leq\ell-1$ and $m^{(i)}_j(\nu)>0$.
If there is no such~$j$, set $j=0$.
We divide the proof into two cases depending on whether $j=\ell-1$ or $j<\ell-1$.

  {\bf Step 1.} Let us consider the case $j=\ell-1$.
Since the string $(\ell,x_\ell)$ is singular, we have $P^{(i)}_\ell(\nu)=x_\ell$, and by the minimality of $x_\ell$ we
have $x_\ell\leq x_{\ell-1}$.
Since $x_{\ell-1}\leq P^{(i)}_{\ell-1}$, we conclude that $P^{(i)}_{\ell-1}(\nu)\geq P^{(i)}_\ell(\nu)$.

 {\bf Step 2.} Let us consider the case $j<\ell-1$.
Recall that the assumptions $\ell=\ell^{(i)}$ and $m^{(i)}_\ell(\nu)=1$ imply that the string $(\ell,x_\ell)$ is
singular and thus $P^{(i)}_\ell(\nu)=x_\ell$.
If $j>0$, then by the minimality of $\ell$ with $j<\ell$ we have $P^{(i)}_j(\nu)\geq x_j\geq x_\ell =P^{(i)}_\ell(\nu)$.
If $j=0$ we have $P^{(i)}_0(\nu)=0$.
In this case, we also have $P^{(i)}_j(\nu)\geq P^{(i)}_\ell(\nu)$ since we have $P^{(i)}_\ell(\nu)=x_\ell\leq 0$ which
is the consequence of $\ell>0$.
In both cases, by the convexity of $P^{(i)}_k(\nu)$ between $j\leq k\leq\ell$ we have $P^{(i)}_{\ell-1}(\nu)\geq
\min\big\{P^{(i)}_j(\nu),P^{(i)}_\ell(\nu)\big\}=P^{(i)}_\ell(\nu)$.

Proof of (2) is the same with the previous case.
\end{proof}

\begin{proof}[Proof of Proposition~\ref{prop:caseD(3)2}] Since $\ell^{(i-1)}=\ell$ we have $m^{(i-1)}_\ell(\nu)>0$.
Then Lemma~\ref{lem:convexity2} with $m^{(i)}_\ell(\nu)=1$ claims that $P^{(i)}_{\ell-1}(\nu)$, $P^{(i)}_{\ell}(\nu)$
and $P^{(i)}_{\ell+1}(\nu)-2$ are strictly convex.
The only possibility that is compatible with the relation $P^{(i)}_{\ell-1}(\nu)\geq
P^{(i)}_{\ell}(\nu)<P^{(i)}_{\ell+1}(\nu)$ is $P^{(i)}_{\ell-1}(\nu)+1=P^{(i)}_{\ell}(\nu)+1=P^{(i)}_{\ell+1}(\nu)$.
Then from Lemma~\ref{lem:convexity3} we have
\begin{gather*}
-P^{(i)}_{\ell-1}(\nu)+2P^{(i)}_{\ell}(\nu)-P^{(i)}_{\ell+1}(\nu)=-1
\\
\qquad
{}\geq m^{(i-1)}_\ell(\nu)-2m^{(i)}_\ell(\nu)+m^{(i+1)}_\ell(\nu)
=m^{(i-1)}_\ell(\nu)-2+m^{(i+1)}_\ell(\nu),
\end{gather*}
that is, $1\geq m^{(i-1)}_\ell(\nu)+m^{(i+1)}_\ell(\nu)$.
Since we know that $m^{(i-1)}_\ell(\nu)>0$ and $m^{(i+1)}_\ell(\nu)\geq 0$, the only possibility is
$m^{(i-1)}_\ell(\nu)=1$ and $m^{(i+1)}_\ell(\nu)=0$.
\end{proof}

\begin{proof}[Proof of Proposition~\ref{prop:caseD2} for the case $\boldsymbol{\ell^{(i-1)}=\ell}$]\sloppy By~\eqref{eq:D_2} we have
$P^{(i)}_{\ell+1}(\nu)=x_\ell+1$, that is,~\eqref{eq:caseD2}.
By the assumption we have $\ell-1<\ell^{(i-1)}$.
Then from~\eqref{eq:D_2} we have \mbox{$P^{(i)}_{\ell-1}(\bar{\nu})=P^{(i)}_{\ell-1}(\nu)=x_\ell$}, which
implies~\eqref{eq:caseD3}.
Finally $m^{(i+1)}_\ell(\nu)=0$ implies that $\ell<\ell^{(i+1)}$, that is,~\eqref{eq:caseD4}.
\end{proof}
This completes the whole proof of Proposition~\ref{prop:caseD2}.

Again, using the proof of $\bar{\widetilde{\nu}}=\widetilde{\bar{\nu}}$, we can show
$\bar{\widetilde{J}}=\widetilde{\bar{J}}$ for this case as follows.

\begin{Proposition}%\label{prop:D_3}
Suppose that $\ell^{(i-1)}=\ell=\ell^{(i)}<\ell_{(i)}$ and $m^{(i)}_\ell(\nu)=1$.
Then we have $\bar{\widetilde{J}}=\widetilde{\bar{J}}$.
\end{Proposition}
\begin{proof}
By~\eqref{eq:D_2_1} it is enough to analyze the string $(\ell,x_\ell)$ of $(\nu,J)^{(i)}$.
Then we have
\begin{gather*}
\begin{array}{@{}lllll}
(\ell,x_\ell)&\overset{\widetilde{f}_i}{\longmapsto}& (\ell+1,x_\ell-1)&\overset{\delta}{\longmapsto}&
\big(\ell,P^{(i)}_\ell(\bar{\widetilde{\nu}})\big),
\\
(\ell,x_\ell)&\overset{\delta}{\longmapsto}&
\big(\ell-1,P^{(i)}_{\ell-1}(\bar{\nu})\big)&\overset{\widetilde{f}_i}{\longmapsto}& \big(\ell,P^{(i)}_{\ell-1}(\bar{\nu})-1\big).
\end{array}
\end{gather*}
From $\ell^{(i-1)}=\ell<\ell^{(i+1)}$ by Proposition~\ref{prop:caseD2}, $P^{(i)}_\ell(\nu)=x_\ell$ by~\eqref{eq:D_2} and
the fact that the length of the string $(\ell,x_\ell)$ is preserved under both $\delta\circ\widetilde{f}_i$ and
$\widetilde{f}_i\circ\delta$, we have $P^{(i)}_\ell(\bar{\widetilde{\nu}})=P^{(i)}_\ell(\nu)-1=x_\ell-1$.
On the other hand, we have $P^{(i)}_{\ell-1}(\bar{\nu})-1=P^{(i)}_{\ell-1}(\nu)-1=x_\ell-1$ since we have
$\ell-1<\ell^{(i-1)}$ and $P^{(i)}_{\ell-1}(\nu)=x_\ell$ by~\eqref{eq:D_2}.
\end{proof}

\begin{Example}
Here we give an example for Proposition~\ref{prop:caseD(3)2}.
Consider the following rigged conf\/iguration $(\nu,J)$ of type $(B^{1,1})^{\otimes 3}\otimes B^{1,3}\otimes
B^{2,1}\otimes B^{2,2}\otimes B^{3,1}$ of $D^{(1)}_5$
\begin{center}
\unitlength 10pt
{\small
\begin{picture}(43,7)
\multiput(-0.8,2.1)(0,1){4}{1}
\put(-1.7,6.1){$-2$}
\put(0,2){\Yboxdim10pt\yng(6,3,1,1,1)}
\put(1.2,2.1){0}
\put(1.2,3.1){0}
\put(1.2,4.1){0}
\put(3.2,5.1){1}
\put(6.1,6.1){$-2$}
\put(9.5,0){
\put(-0.8,0.1){0}
\put(-0.8,1.1){0}
\put(-0.8,2.1){0}
\multiput(-1.6,3.1)(0,1){3}{$-1$}
\put(-0.8,6.1){1}
\put(0,0){\Yboxdim10pt\yng(6,3,2,2,1,1,1)}
\put(1.2,0.1){0}
\put(1.2,1.1){0}
\put(1.2,2.1){0}
\put(2.2,3.1){$-1$}
\put(2.2,4.1){$-1$}
\put(3.2,5.0){$-1$}
\put(6.1,6.0){1}
}
\put(18.5,0){
\multiput(-0.8,0.1)(0,1){3}{0}
\put(-0.8,3.1){1}
\put(-0.8,4.1){1}
\put(-0.8,5.1){0}
\put(-0.8,6.1){0}
\put(0,0){\Yboxdim10pt\yng(6,5,2,2,1,1,1)}
\put(1.2,0.1){0}
\put(1.2,1.1){0}
\put(1.2,2.1){0}
\put(2.2,3.0){0}
\put(2.2,4.0){0}
\put(5.2,5.0){0}
\put(6.2,6.0){0}
}
\put(28.0,4){
\put(-1.6,0.1){$-1$}
\put(-1.6,1.1){$-1$}
\put(-1.6,2.1){$-1$}
\put(0,0){\Yboxdim10pt\yng(5,2,2)}
\put(2.2,0.1){$-1$}
\put(2.2,1.1){$-1$}
\put(5.2,2.1){$-1$}
}
\put(37,4){
\put(-0.8,0.1){1}
\put(-1.6,1.1){$-1$}
\put(-1.6,2.1){$-2$}
\put(0,0){\Yboxdim10pt\yng(6,3,1)}
\put(1.2,0.1){0}
\put(3.2,1.1){$-1$}
\put(6.1,2.1){$-2$}
}
\end{picture}}
%\{\{6,3,1,1,1\},\{6,3,2,2,1,1,1\},\{6,5,2,2,1,1,1\},\{5,2,2\},\{6,3,1\}\}\\
%\{\{-2,1,0,0,0\},\{1,-1,-1,-1,0,0,0\},\{0,0,0,0,0,0,0\},\{-1,-1,-1\},\{-2,-1,0\}\}
\end{center}
The corresponding tensor product $\Phi^{-1}(\nu,J)$ is
\begin{gather*}
\Yboxdim14pt \Yvcentermath1 \young(\mone)\otimes \young(\mone)\otimes \young(\mthree)\otimes
\young(\mfive\mfour\mone)\otimes \young(2,\mone)\otimes \young(12,\mfive\mtwo)\otimes \young(1,5,\mtwo)
\end{gather*}
Then $(\bar{\nu},\bar{J})$ is
\begin{center}
\unitlength 10pt
{\small
\begin{picture}(43,7)
\multiput(-0.8,2.1)(0,1){3}{0}
\put(-0.8,5.1){1}
\put(-1.7,6.1){$-1$}
\put(0,2){\Yboxdim10pt\yng(5,2,1,1,1)}
\put(1.2,2.1){0}
\put(1.2,3.1){0}
\put(1.2,4.1){0}
\put(2.2,5.1){1}
\put(5.1,6.1){$-1$}
\put(9.5,0){
\put(-0.8,0.1){0}
\put(-0.8,1.1){0}
\put(-0.8,2.1){0}
\multiput(-1.6,3.1)(0,1){3}{$-1$}
\put(-0.8,6.1){1}
\put(0,0){\Yboxdim10pt\yng(5,2,2,2,1,1,1)}
\put(1.2,0.1){0}
\put(1.2,1.1){0}
\put(1.2,2.1){0}
\put(2.2,3.1){$-1$}
\put(2.2,4.1){$-1$}
\put(2.2,5.0){$-1$}
\put(5.1,6.0){1}
}
\put(18.5,0){
\multiput(-0.8,0.1)(0,1){3}{0}
\put(-0.8,3.1){1}
\put(-0.8,4.1){1}
\put(-0.8,5.1){0}
\put(-0.8,6.1){0}
\put(0,0){\Yboxdim10pt\yng(5,4,2,2,1,1,1)}
\put(1.2,0.1){0}
\put(1.2,1.1){0}
\put(1.2,2.1){0}
\put(2.2,3.0){0}
\put(2.2,4.0){0}
\put(4.2,5.0){0}
\put(5.2,6.0){0}
}
\put(28.0,4){
\put(-1.6,0.1){$-1$}
\put(-1.6,1.1){$-1$}
\put(-1.6,2.1){$-1$}
\put(0,0){\Yboxdim10pt\yng(4,2,2)}
\put(2.2,0.1){$-1$}
\put(2.2,1.1){$-1$}
\put(4.2,2.1){$-1$}
}
\put(37,4){
\put(-0.8,0.1){1}
\put(-1.6,1.1){$-1$}
\put(-1.6,2.1){$-2$}
\put(0,0){\Yboxdim10pt\yng(5,3,1)}
\put(1.2,0.1){0}
\put(3.2,1.1){$-1$}
\put(5.1,2.1){$-2$}
}
\end{picture}}
%\{\{5,2,1,1,1\},\{5,2,2,2,1,1,1\},\{5,4,2,2,1,1,1\},\{4,2,2\},\{5,3,1\}\}\\
%\{\{-1,1,0,0,0\},\{1,-1,-1,-1,0,0,0\},\{0,0,0,0,0,0,0\},\{-1,-1,-1\},\{-2,-1,0\}\}
\end{center}
and $(\widetilde{\nu},\widetilde{J})$ is
\begin{center}
\unitlength 10pt
{\small
\begin{picture}(43,7)
\multiput(-0.8,2.1)(0,1){4}{1}
\put(-1.7,6.1){$-1$}
\put(0,2){\Yboxdim10pt\yng(6,3,1,1,1)}
\put(1.2,2.1){0}
\put(1.2,3.1){0}
\put(1.2,4.1){0}
\put(3.2,5.1){1}
\put(6.1,6.1){$-1$}
\put(9.5,0){
\put(-0.8,0.1){0}
\put(-0.8,1.1){0}
\put(-0.8,2.1){0}
\multiput(-1.6,3.1)(0,1){2}{$-1$}
\put(-1.6,5.1){$-2$}
\put(-1.6,6.1){$-1$}
\put(0,0){\Yboxdim10pt\yng(6,4,2,2,1,1,1)}
\put(1.2,0.1){0}
\put(1.2,1.1){0}
\put(1.2,2.1){0}
\put(2.2,3.1){$-1$}
\put(2.2,4.1){$-1$}
\put(4.2,5.0){$-2$}
\put(6.1,6.0){$-1$}
}
\put(18.5,0){
\multiput(-0.8,0.1)(0,1){3}{0}
\put(-0.8,3.1){1}
\put(-0.8,4.1){1}
\put(-0.8,5.1){1}
\put(-0.8,6.1){1}
\put(0,0){\Yboxdim10pt\yng(6,5,2,2,1,1,1)}
\put(1.2,0.1){0}
\put(1.2,1.1){0}
\put(1.2,2.1){0}
\put(2.2,3.0){0}
\put(2.2,4.0){0}
\put(5.2,5.0){1}
\put(6.2,6.0){1}
}
\put(28.0,4){
\put(-1.6,0.1){$-1$}
\put(-1.6,1.1){$-1$}
\put(-1.6,2.1){$-1$}
\put(0,0){\Yboxdim10pt\yng(5,2,2)}
\put(2.2,0.1){$-1$}
\put(2.2,1.1){$-1$}
\put(5.2,2.1){$-1$}
}
\put(37,4){
\put(-0.8,0.1){1}
\put(-1.6,1.1){$-1$}
\put(-1.6,2.1){$-2$}
\put(0,0){\Yboxdim10pt\yng(6,3,1)}
\put(1.2,0.1){0}
\put(3.2,1.1){$-1$}
\put(6.1,2.1){$-2$}
}
\end{picture}}
%\{\{6,3,1,1,1\},\{6,4,2,2,1,1,1\},\{6,5,2,2,1,1,1\},\{5,2,2\},\{6,3,1\}\}\\
%\{\{-1,1,0,0,0\},\{-1,-2,-1,-1,0,0,0\},\{1,1,0,0,0,0,0\},\{-1,-1,-1\},\{-2,-1,0\}\}
\end{center}
In this example, we have $\ell^{(1)}=\ell=\ell^{(2)}=3$, $m^{(2)}_3(\nu)=1$ and $x_\ell=-1$.
As proved in Proposition~\ref{prop:caseD(3)2}, we have $P^{(2)}_2(\nu)=P^{(2)}_3(\nu)=-1=x_\ell$,
$P^{(2)}_4(\nu)=0=x_\ell+1$, $m^{(1)}_3(\nu)=1$ and $m^{(3)}_3(\nu)=0$.
Finally $(\bar{\widetilde{\nu}},\bar{\widetilde{J}})=(\widetilde{\bar{\nu}},\widetilde{\bar{J}})$ is
\begin{center}
\unitlength 10pt
{\small
\begin{picture}(43,7)
\multiput(-0.8,2.1)(0,1){3}{0}
\put(-0.8,5.1){1}
\put(-0.8,6.1){0}
\put(0,2){\Yboxdim10pt\yng(5,2,1,1,1)}
\put(1.2,2.1){0}
\put(1.2,3.1){0}
\put(1.2,4.1){0}
\put(2.2,5.1){1}
\put(5.1,6.1){0}
\put(9.5,0){
\put(-0.8,0.1){0}
\put(-0.8,1.1){0}
\put(-0.8,2.1){0}
\multiput(-1.6,3.1)(0,1){2}{$-1$}
\put(-1.6,5.1){$-2$}
\put(-1.6,6.1){$-1$}
\put(0,0){\Yboxdim10pt\yng(5,3,2,2,1,1,1)}
\put(1.2,0.1){0}
\put(1.2,1.1){0}
\put(1.2,2.1){0}
\put(2.2,3.1){$-1$}
\put(2.2,4.1){$-1$}
\put(3.2,5.0){$-2$}
\put(5.1,6.0){$-1$}
}
\put(18.5,0){
\multiput(-0.8,0.1)(0,1){3}{0}
\put(-0.8,3.1){1}
\put(-0.8,4.1){1}
\put(-0.8,5.1){1}
\put(-0.8,6.1){1}
\put(0,0){\Yboxdim10pt\yng(5,4,2,2,1,1,1)}
\put(1.2,0.1){0}
\put(1.2,1.1){0}
\put(1.2,2.1){0}
\put(2.2,3.0){0}
\put(2.2,4.0){0}
\put(4.2,5.0){1}
\put(5.2,6.0){1}
}
\put(28.0,4){
\put(-1.6,0.1){$-1$}
\put(-1.6,1.1){$-1$}
\put(-1.6,2.1){$-1$}
\put(0,0){\Yboxdim10pt\yng(4,2,2)}
\put(2.2,0.1){$-1$}
\put(2.2,1.1){$-1$}
\put(4.2,2.1){$-1$}
}
\put(37,4){
\put(-0.8,0.1){1}
\put(-1.6,1.1){$-1$}
\put(-1.6,2.1){$-2$}
\put(0,0){\Yboxdim10pt\yng(5,3,1)}
\put(1.2,0.1){0}
\put(3.2,1.1){$-1$}
\put(5.1,2.1){$-2$}
}
\end{picture}}
%\{\{5,2,1,1,1\},\{5,3,2,2,1,1,1\},\{5,4,2,2,1,1,1\},\{4,2,2\},\{5,3,1\}\}\\
%\{\{0,1,0,0,0\},\{-1,-2,-1,-1,0,0,0\},\{1,1,0,0,0,0,0\},\{-1,-1,-1\},\{-2,-1,0\}\}
\end{center}
\end{Example}

\subsection{Proof for Case D (4)}

\subsubsection{Outline}

We will prove the following proposition in the latter part of this subsection.
\begin{Proposition}
\label{prop:caseD(4)1}
Suppose that $\ell=\ell^{(i)}=\ell_{(i)}$ and $m^{(i)}_\ell(\nu)=2$.
Then the following relations hold:
\begin{gather}
P^{(i)}_{\ell+1}(\nu)=x_\ell+1,
\label{eq:caseD(4)1}
\\
P^{(i)}_{\ell-1}(\bar{\nu})=x_\ell,
\label{eq:caseD(4)2}
\\
\ell+1\leq\ell_{(i-1)}.
\label{eq:caseD(4)3}
\end{gather}
\end{Proposition}

\begin{proof}
[Proof of ``Proposition~\ref{prop:caseD(4)1} $\boldsymbol{\Longrightarrow\bar{\widetilde{\nu}}=\widetilde{\bar{\nu}}}$'']
 {\bf Step 1.} Let us consider the case $\widetilde{f}_i\circ\delta$.
The f\/irst~$\delta$ creates the two strings $\big(\ell-1,P^{(i)}_{\ell-1}(\bar{\nu})\big)$ of $(\bar{\nu},\bar{J})^{(i)}$.
Since the operation~$\delta$ does not change the riggings of unchanged strings, the relation~\eqref{eq:caseD(4)2} and
$m^{(i)}_\ell(\nu)=2$ imply that the strings $(\ell-1,x_\ell)$ are the longest strings among the strings with rigging
$x_\ell$, hence $\widetilde{f}_i$ will act on one of these strings.

 {\bf Step 2.} Let us consider the case $\delta\circ\widetilde{f}_i$.
Recall that by def\/inition of $\widetilde{f}_i$ we have $\widetilde{\ell}^{(a)}=\ell^{(a)}$ for all $a<i$.
From Lemma~\ref{lem:f_singular}(1) the relation~\eqref{eq:caseD(4)1} implies that the string $(\ell+1,x_\ell-1)$ of~$(\widetilde{\nu},\widetilde{J})^{(i)}$ created by~$\widetilde{f}_i$ is singular.
Then from the assumption $m^{(i)}_\ell(\nu)=2$ we see that there is one singular string of length $\ell$ and at least
one singular string of length $\ell+1$ within $(\widetilde{\nu},\widetilde{J})^{(i)}$.
Thus we have $\widetilde{\ell}^{(a)}=\ell^{(a)}$ for all $i\leq a$, $\widetilde{\ell}_{(a)}=\ell_{(a)}$ for all $i<a$
and $\widetilde{\ell}_{(i)}=\ell+1$.
Then from~\eqref{eq:caseD(4)3} we have $\widetilde{\ell}_{(i-1)}=\ell_{(i-1)}$ and hence
$\widetilde{\ell}_{(a)}=\ell_{(a)}$ for all $a<i$.

 {\bf Step 3.} To summarize, we have $\bar{\widetilde{\nu}}^{(a)}= \widetilde{\bar{\nu}}^{(a)}=\bar{\nu}^{(a)}$
for all $a\neq i$ and $\bar{\widetilde{\nu}}^{(i)}=\widetilde{\bar{\nu}}^{(i)}$ are obtained by removing a~box from
a~length $\ell$ row of $\nu^{(i)}$.
Hence we have $\bar{\widetilde{\nu}}=\widetilde{\bar{\nu}}$.
\end{proof}

Proofs of $\bar{\widetilde{J}}=\widetilde{\bar{J}}$ will be given by case by case analysis along with the proof of
Proposition~\ref{prop:caseD(4)1} assuming the relation $\bar{\widetilde{\nu}}=\widetilde{\bar{\nu}}$ we have just
proved.

\subsubsection[Proof of Proposition~\ref{prop:caseD(4)1} for the case $\ell^{(i-1)}<\ell$]{Proof
of Proposition~\ref{prop:caseD(4)1} for the case $\boldsymbol{\ell^{(i-1)}<\ell}$}

We divide the proof of Proposition~\ref{prop:caseD(4)1} into two cases according to
whether $\ell^{(i-1)}<\ell$ or $\ell^{(i-1)}=\ell$.
In this subsection, we consider the case $\ell^{(i-1)}<\ell$.
Throughout the proof of Proposition~\ref{prop:caseD(4)1}, let $j$ be the largest integer such that $j\leq\ell-1$ and
$m^{(i)}_j(\nu)>0$.
If there is no such $j$, set $j=0$.

\begin{Proposition}
\label{prop:caseD(4)2}
Suppose that $\ell=\ell^{(i)}=\ell_{(i)}$, $m^{(i)}_\ell(\nu)=2$ and $\ell^{(i-1)}<\ell$.
Then we have
\begin{gather}
\label{eq:caseD(4)4}
P^{(i)}_{\ell-1}(\nu)=x_\ell+1,
\qquad
P^{(i)}_{\ell}(\nu)=x_\ell,
\qquad
P^{(i)}_{\ell+1}(\nu)=x_\ell+1.
\end{gather}
Moreover, we have
\begin{gather}
\label{eq:caseD(4)5}
m^{(i-1)}_\ell(\nu)=0,
\qquad
m^{(i+1)}_\ell(\nu)=2.
\end{gather}
\end{Proposition}

For the proof of this proposition, we prepare the following lemma.

\begin{Lemma}
Suppose that $\ell=\ell^{(i)}=\ell_{(i)}$, $m^{(i)}_\ell(\nu)=2$ and $\ell^{(i-1)}<\ell$.
Then we have
\begin{gather*}
 (1) \quad P^{(i)}_{\ell-1}(\nu)>P^{(i)}_\ell(\nu),\qquad
 (2) \quad P^{(i)}_\ell(\nu)<P^{(i)}_{\ell+1}(\nu).
\end{gather*}
\end{Lemma}
\begin{proof}
Recall that by def\/inition of~$\delta$, the strings $\big(\ell^{(i)},P^{(i)}_{\ell^{(i)}}(\nu)\big)$ and
$\big(\ell_{(i)},P^{(i)}_{\ell_{(i)}}(\nu)\big)$ of $(\nu,J)^{(i)}$ must be dif\/ferent.
Then the assumptions $\ell=\ell^{(i)}=\ell_{(i)}$ and $m^{(i)}_\ell(\nu)=2$ imply that the string $(\ell,x_\ell)$ must
be singular.
In particular, we have $P^{(i)}_\ell(\nu)=x_\ell$.
Then we can use the same arguments of Lemma~\ref{lem:caseD(3)1} to prove the assertions.
\end{proof}

\begin{proof}
[Proof of Proposition~\ref{prop:caseD(4)2}] According to the previous lemma, we can write $P^{(i)}_{\ell-1}(\nu)$ and
$P^{(i)}_{\ell+1}(\nu)$ as follows:
\begin{gather*}
P^{(i)}_{\ell-1}(\nu)=x_\ell+\varepsilon_1,
\qquad
P^{(i)}_{\ell+1}(\nu)=x_\ell+\varepsilon_2,
\end{gather*}
where $\varepsilon_1,\varepsilon_2\geq 1$.
Then Lemma~\ref{lem:convexity3} with $l=\ell$ reads
\begin{gather*}
-P^{(i)}_{\ell-1}(\nu)+2P^{(i)}_{\ell}(\nu)-P^{(i)}_{\ell+1}(\nu)=-\varepsilon_1-\varepsilon_2
\\
\qquad
{}\geq m^{(i-1)}_\ell(\nu)-2m^{(i)}_\ell(\nu)+m^{(i+1)}_\ell(\nu)
=m^{(i-1)}_\ell(\nu)-4+m^{(i+1)}_\ell(\nu),
\end{gather*}
that is,
\begin{gather*}
4-\varepsilon_1-\varepsilon_2\geq m^{(i-1)}_\ell(\nu)+m^{(i+1)}_\ell(\nu).
\end{gather*}
Note that the assumption $\ell^{(i)}=\ell_{(i)}$ implies that $\ell^{(i+1)}=\ell_{(i+1)}=\ell$ by def\/inition of
$\delta$.
In particular, we have $m^{(i+1)}_\ell(\nu)\geq 2$.
Then the only possibility that is compatible with the above inequality is $\varepsilon_1=\varepsilon_2=1$,
$m^{(i-1)}_\ell(\nu)=0$ and $m^{(i+1)}_\ell(\nu)=2$.
\end{proof}

\begin{proof}
[Proof of Proposition~\ref{prop:caseD(4)1} for the case $\boldsymbol{\ell^{(i-1)}<\ell}$] The relation~\eqref{eq:caseD(4)1} is proved
in~\eqref{eq:caseD(4)4}.
As for~\eqref{eq:caseD(4)2}, since $\ell^{(i-1)}<\ell=\ell^{(i)}$ we have
$P^{(i)}_{\ell-1}(\bar{\nu})=P^{(i)}_{\ell-1}(\nu)-1=x_\ell$ by~\eqref{eq:caseD(4)4}.
Finally, by~\eqref{eq:caseD(4)5} we have $m^{(i-1)}_\ell(\nu)=0$.
Recall that we have $\ell_{(i-1)}\geq\ell_{(i)}$ by def\/inition of~$\delta$.
Then, since we are assuming that $\ell_{(i)}=\ell$, we have $\ell_{(i-1)}>\ell$, that is,~\eqref{eq:caseD(4)3}.
\end{proof}

Assuming that $\bar{\widetilde{\nu}}=\widetilde{\bar{\nu}}$, we can show $\bar{\widetilde{J}}=\widetilde{\bar{J}}$ for
this case as follows.

\begin{Proposition}
Suppose that $\ell^{(i-1)}<\ell=\ell^{(i)}=\ell_{(i)}$ and $m^{(i)}_\ell(\nu)=2$.
Then we have $\bar{\widetilde{J}}=\widetilde{\bar{J}}$.
\end{Proposition}
\begin{proof}
From the proof of $\bar{\widetilde{\nu}}=\widetilde{\bar{\nu}}$ after Proposition~\ref{prop:caseD(4)1}, we see that it
is enough to look at the string $(\ell,x_\ell)$ of $(\nu,J)^{(i)}$ which is changed by both $\widetilde{f}_i$ and
$\delta$ in both $\widetilde{f}_i\circ\delta$ and $\delta\circ\widetilde{f}_i$.
Then we have
\begin{gather*}
\begin{array}{@{}lllll}
(\ell,x_\ell)&\overset{\widetilde{f}_i}{\longmapsto}& (\ell+1,x_\ell-1)&\overset{\delta}{\longmapsto}&
\big(\ell,P^{(i)}_\ell(\bar{\widetilde{\nu}})\big),
\\
(\ell,x_\ell)&\overset{\delta}{\longmapsto}&
\big(\ell-1,P^{(i)}_{\ell-1}(\bar{\nu})\big)&\overset{\widetilde{f}_i}{\longmapsto}&\big(\ell,P^{(i)}_{\ell-1}(\bar{\nu})-1\big).
\end{array}
\end{gather*}
Recall that we have
\begin{gather*}
\widetilde{\ell}^{(i-1)}<\widetilde{\ell}^{(i)}=\widetilde{\ell}^{(i+1)}
=\widetilde{\ell}_{(i+1)}=\ell<\widetilde{\ell}_{(i)}=\ell+1\leq\widetilde{\ell}_{(i-1)}
\end{gather*}
by the assumption $\ell^{(i-1)}<\ell=\ell^{(i)}=\ell_{(i)}$ and the arguments given in the proof of
$\bar{\widetilde{\nu}}=\widetilde{\bar{\nu}}$.
Then we have $P^{(i)}_\ell(\bar{\widetilde{\nu}})=P^{(i)}_\ell(\nu)-1=x_\ell-1$ since we have $P^{(i)}_\ell(\nu)=x_\ell$
by~\eqref{eq:caseD(4)4}.
Similarly we have $P^{(i)}_{\ell-1}(\bar{\nu})-1=P^{(i)}_{\ell-1}(\nu)-2=x_\ell-1$ since we have
$P^{(i)}_{\ell-1}(\nu)=x_\ell+1$ by~\eqref{eq:caseD(4)4}.
\end{proof}

\subsubsection[Proof of Proposition~\ref{prop:caseD(4)1} for the case $\ell^{(i-1)}=\ell$]{Proof
of Proposition~\ref{prop:caseD(4)1} for the case $\boldsymbol{\ell^{(i-1)}=\ell}$}

\begin{Proposition}
\label{prop:caseD(4)3}
Suppose that $\ell=\ell^{(i)}=\ell_{(i)}$, $m^{(i)}_\ell(\nu)=2$ and $\ell^{(i-1)}=\ell$.
Then we have
\begin{gather}
\label{eq:caseD(4)6}
P^{(i)}_{\ell-1}(\nu)=P^{(i)}_{\ell}(\nu)=x_\ell,
\qquad
P^{(i)}_{\ell+1}(\nu)=x_\ell+1.
\end{gather}
Moreover, we have
\begin{gather}
\label{eq:caseD(4)7}
m^{(i-1)}_\ell(\nu)=1,
\qquad
m^{(i+1)}_\ell(\nu)=2.
\end{gather}
\end{Proposition}

For the proof of this proposition, we prepare the following lemma.

\begin{Lemma}
Suppose that $\ell=\ell^{(i)}=\ell_{(i)}$, $m^{(i)}_\ell(\nu)=2$ and $\ell^{(i-1)}=\ell$.
Then we have
\begin{gather*}
 (1)\quad P^{(i)}_{\ell-1}(\nu)\geq P^{(i)}_\ell(\nu),\qquad
 (2)\quad P^{(i)}_\ell(\nu)<P^{(i)}_{\ell+1}(\nu).
\end{gather*}
\end{Lemma}
\begin{proof}
Since we are assuming that $\ell^{(i)}=\ell_{(i)}$ and $m^{(i)}_\ell(\nu)=2$, we see that the strings $(\ell,x_\ell)$ of
$(\nu,J)^{(i)}$ are singular.
Then we can use the same arguments of the proof of Lemma~\ref{lem:caseD(3)2}.
\end{proof}

\begin{proof}
[Proof of Proposition~\ref{prop:caseD(4)3}] According to the previous lemma, we can express $P^{(i)}_{\ell-1}(\nu)$ and
$P^{(i)}_{\ell+1}(\nu)$ as follows:
\begin{gather*}
P^{(i)}_{\ell-1}(\nu)=x_\ell+\varepsilon_1,
\qquad
P^{(i)}_{\ell+1}(\nu)=x_\ell+\varepsilon_2,
\end{gather*}
where $\varepsilon_1\geq 0$ and $\varepsilon_2\geq 1$.
Then Lemma~\ref{lem:convexity3} with $l=\ell$ reads
\begin{gather*}
-P^{(i)}_{\ell-1}(\nu)+2P^{(i)}_{\ell}(\nu)-P^{(i)}_{\ell+1}(\nu)=-\varepsilon_1-\varepsilon_2
\\
\qquad
{}\geq m^{(i-1)}_\ell(\nu)-2m^{(i)}_\ell(\nu)+m^{(i+1)}_\ell(\nu)
=m^{(i-1)}_\ell(\nu)-4+m^{(i+1)}_\ell(\nu),
\end{gather*}
that is,
\begin{gather*}
4-\varepsilon_1-\varepsilon_2\geq m^{(i-1)}_\ell(\nu)+m^{(i+1)}_\ell(\nu).
\end{gather*}
Since we are assuming that $\ell^{(i-1)}=\ell$, we have $m^{(i-1)}_\ell(\nu)\geq 1$.
Also, from the assumption $\ell^{(i)}=\ell_{(i)}$ we see that $\ell^{(i+1)}=\ell_{(i+1)}$ by def\/inition of~$\delta$.
Thus we have $m^{(i+1)}_\ell(\nu)\geq 2$.
Then the only possibility that is compatible with the above inequality is that $\varepsilon_1=0$, $\varepsilon_2=1$,
$m^{(i-1)}_\ell(\nu)=1$ and $m^{(i+1)}_\ell(\nu)=2$.
\end{proof}

\begin{proof}
[Proof of Proposition~\ref{prop:caseD(4)1} for the case $\boldsymbol{\ell^{(i-1)}=\ell}$] The relation~\eqref{eq:caseD(4)1} is proved
in~\eqref{eq:caseD(4)6}.
As for~\eqref{eq:caseD(4)2}, since $\ell^{(i-1)}=\ell^{(i)}=\ell$ we have
$P^{(i)}_{\ell-1}(\bar{\nu})=P^{(i)}_{\ell-1}(\nu)=x_\ell$ by~\eqref{eq:caseD(4)6}.
Finally, by~\eqref{eq:caseD(4)7} we have $m^{(i-1)}_\ell(\nu)=1$.
Since we are assuming that $\ell^{(i-1)}=\ell$, this string corresponds to $(\ell^{(i-1)},x_{\ell^{(i-1)}})$ which is
dif\/ferent from the string $(\ell_{(i-1)},x_{\ell_{(i-1)}})$.
Since we have $\ell^{(i-1)}\leq\ell_{(i-1)}$ by the def\/inition of~$\delta$, we must have $\ell<\ell_{(i-1)}$, that
is,~\eqref{eq:caseD(4)3}.
\end{proof}

Assuming that $\bar{\widetilde{\nu}}=\widetilde{\bar{\nu}}$, we can show $\bar{\widetilde{J}}=\widetilde{\bar{J}}$ for
this case as follows.

\begin{Proposition}
Suppose that $\ell=\ell^{(i)}=\ell_{(i)}$, $m^{(i)}_\ell(\nu)=2$ and $\ell^{(i-1)}=\ell$.
Then we have $\bar{\widetilde{J}}=\widetilde{\bar{J}}$.
\end{Proposition}
\begin{proof}
Recall that we can assume that $\widetilde{f}_i$ acts on the same string $(\ell,x_\ell)$ of $(\nu,J)^{(i)}$ in both
$\widetilde{f}_i\circ\delta$ and $\delta\circ\widetilde{f}_i$.
Then we have
\begin{gather*}
\begin{array}{@{}lllll}
(\ell,x_\ell)&\overset{\widetilde{f}_i}{\longmapsto}& (\ell+1,x_\ell-1)&\overset{\delta}{\longmapsto}&
\big(\ell,P^{(i)}_\ell(\bar{\widetilde{\nu}})\big),
\\
(\ell,x_\ell)&\overset{\delta}{\longmapsto}&
\big(\ell-1,P^{(i)}_{\ell-1}(\bar{\nu})\big)&\overset{\widetilde{f}_i}{\longmapsto}&\big(\ell,P^{(i)}_{\ell-1}(\bar{\nu})-1\big).
\end{array}
\end{gather*}
Recall that we have
\begin{gather*}
\widetilde{\ell}^{(i-1)}=\widetilde{\ell}^{(i)}=\widetilde{\ell}^{(i+1)}
=\widetilde{\ell}_{(i+1)}=\ell<\widetilde{\ell}_{(i)}=\ell+1\leq\widetilde{\ell}_{(i-1)}
\end{gather*}
by the assumption $\ell^{(i-1)}=\ell^{(i)}=\ell_{(i)}=\ell$ and the arguments given in the proof of
$\bar{\widetilde{\nu}}=\widetilde{\bar{\nu}}$.
Then we have $P^{(i)}_\ell(\bar{\widetilde{\nu}})=P^{(i)}_\ell(\nu)-1=x_\ell-1$ since we have $P^{(i)}_\ell(\nu)=x_\ell$
by~\eqref{eq:caseD(4)6}.
Similarly we have $P^{(i)}_{\ell-1}(\bar{\nu})-1=P^{(i)}_{\ell-1}(\nu)-1=x_\ell-1$ since we have
$P^{(i)}_{\ell-1}(\nu)=x_\ell$ by~\eqref{eq:caseD(4)6}.
\end{proof}

\begin{Example}
Consider the following rigged conf\/iguration $(\nu,J)$ of type $(B^{1,1})^{\otimes 3}\otimes B^{1,3}\otimes
B^{2,1}\otimes B^{2,2}\otimes B^{3,1}$ of $D^{(1)}_5$
\begin{center}
\unitlength 10pt
{\small
\begin{picture}(42,6)
\put(-0.8,2.1){2}
\multiput(-0.8,3.1)(0,1){2}{1}
\put(-1.6,5.1){$-2$}
\put(0,2){\Yboxdim10pt\yng(6,3,2,1)}
\put(1.2,2.1){0}
\put(2.2,3.1){0}
\put(3.2,4.1){1}
\put(6.2,5.1){$-2$}
\put(10,0){
\multiput(-1.6,0.1)(0,1){5}{$-1$}
\put(-0.8,5.1){1}
\put(0,0){\Yboxdim10pt\yng(6,3,3,2,1,1)}
\put(1.2,0.1){$-1$}
\put(1.2,1.1){$-1$}
\put(2.2,2.1){$-1$}
\put(3.2,3.1){$-1$}
\put(3.2,4.0){$-1$}
\put(6.1,5.0){1}
}
\put(19.5,1){
\put(-0.8,0.1){2}
\multiput(-0.8,1.1)(0,1){2}{1}
\put(-1.6,3.1){$-1$}
\put(-1.6,4.1){$-2$}
\put(0,0){\Yboxdim10pt\yng(6,5,3,3,1)}
\put(1.2,0.1){1}
\put(3.2,1.1){1}
\put(3.2,2.1){1}
\put(5.2,3.1){$-1$}
\put(6.2,4.1){$-2$}
}
\put(29.0,4){
\multiput(-0.8,0.1)(0,1){2}{1}
\put(0,0){\Yboxdim10pt\yng(5,3)}
\put(3.2,0.1){1}
\put(5.2,1.1){1}
}
\put(37.0,3){
\multiput(-1.6,0.1)(0,1){3}{$-1$}
\put(0,0){\Yboxdim10pt\yng(5,3,1)}
\put(1.2,0.1){$-1$}
\put(3.2,1.1){$-1$}
\put(5.2,2.1){$-1$}
}
\end{picture}}
%\{\{6,3,2,1\},\{6,3,3,2,1,1\},\{6,5,3,3,1\},\{5,3\},\{5,3,1\}\},\\
%\{\{-2,1,0,0\},\{1,-1,-1,-1,-1,-1\},\{-2,-1,1,1,1\},\{1,1\},\{-1,-1,-1\}\}
\end{center}
The corresponding tensor product $\Phi^{-1}(\nu,J)$ is
\begin{gather*}
\Yboxdim14pt \Yvcentermath1 \young(\mone)\otimes \young(\mone)\otimes \young(\mone)\otimes \young(114)\otimes
\young(\mthree,\mtwo)\otimes \young(2\mthree,\mone\mone)\otimes \young(1,3,\mfive)
\end{gather*}
Then $(\bar{\nu},\bar{J})$ is
\begin{center}
\unitlength 10pt
{\small
\begin{picture}(42,6)
\put(0,2){
\put(-0.8,0.1){1}
\multiput(-0.8,1.1)(0,1){2}{0}
\put(-1.6,3.1){$-2$}
\Yboxdim10pt\yng(5,2,2,1)}
\put(1.2,2.1){0}
\put(2.2,3.1){0}
\put(2.2,4.1){0}
\put(5.2,5.1){$-2$}
\put(10,0){
\multiput(-1.6,0.1)(0,1){5}{$-1$}
\put(-0.8,5.1){1}
\put(0,0){\Yboxdim10pt\yng(6,2,2,2,1,1)}
\put(1.2,0.1){$-1$}
\put(1.2,1.1){$-1$}
\put(2.2,2.1){$-1$}
\put(2.2,3.1){$-1$}
\put(2.2,4.0){$-1$}
\put(6.1,5.0){1}
}
\put(19.5,1){
\multiput(-0.8,0.1)(0,1){3}{2}
\put(-1.6,3.1){$-1$}
\put(-1.6,4.1){$-2$}
\put(0,0){\Yboxdim10pt\yng(6,5,2,2,1)}
\put(1.2,0.1){1}
\put(2.2,1.1){2}
\put(2.2,2.1){2}
\put(5.2,3.1){$-1$}
\put(6.2,4.1){$-2$}
}
\put(29.0,4){
\multiput(-0.8,0.1)(0,1){2}{1}
\put(0,0){\Yboxdim10pt\yng(5,2)}
\put(2.2,0.1){1}
\put(5.2,1.1){1}
}
\put(37.0,3){
\multiput(-1.6,0.1)(0,1){3}{$-1$}
\put(0,0){\Yboxdim10pt\yng(5,2,1)}
\put(1.2,0.1){$-1$}
\put(2.2,1.1){$-1$}
\put(5.2,2.1){$-1$}
}
\end{picture}}
%\{\{5,2,2,1\},\{6,2,2,2,1,1\},\{6,5,2,2,1\},\{5,2\},\{5,2,1\}\},\\
%\{\{-2,0,0,0\},\{1,-1,-1,-1,-1,-1\},\{-2,-1,2,2,1\},\{1,1\},\{-1,-1,-1\}\}
\end{center}
and $(\widetilde{\nu},\widetilde{J})$ is
\begin{center}
\unitlength 10pt
{\small
\begin{picture}(42,6)
\put(0,2){
\put(-0.8,0.1){2}
\multiput(-0.8,1.1)(0,1){2}{1}
\put(-1.6,3.1){$-1$}
\Yboxdim10pt\yng(6,3,2,1)}
\put(1.2,2.1){0}
\put(2.2,3.1){0}
\put(3.2,4.1){1}
\put(6.2,5.1){$-1$}
\put(10,0){
\multiput(-1.6,0.1)(0,1){4}{$-1$}
\put(-1.6,4.1){$-2$}
\put(-1.6,5.1){$-1$}
\put(0,0){\Yboxdim10pt\yng(6,4,3,2,1,1)}
\put(1.2,0.1){$-1$}
\put(1.2,1.1){$-1$}
\put(2.2,2.1){$-1$}
\put(3.2,3.1){$-1$}
\put(4.2,4.0){$-2$}
\put(6.0,5.1){$-1$}
}
\put(19.5,1){
\put(-0.8,0.1){2}
\multiput(-0.8,1.1)(0,1){2}{1}
\put(-0.8,3.1){0}
\put(-1.6,4.1){$-1$}
\put(0,0){\Yboxdim10pt\yng(6,5,3,3,1)}
\put(1.2,0.1){1}
\put(3.2,1.1){1}
\put(3.2,2.1){1}
\put(5.2,3.1){0}
\put(6.2,4.1){$-1$}
}
\put(29.0,4){
\multiput(-0.8,0.1)(0,1){2}{1}
\put(0,0){\Yboxdim10pt\yng(5,3)}
\put(3.2,0.1){1}
\put(5.2,1.1){1}
}
\put(37.0,3){
\multiput(-1.6,0.1)(0,1){3}{$-1$}
\put(0,0){\Yboxdim10pt\yng(5,3,1)}
\put(1.2,0.1){$-1$}
\put(3.2,1.1){$-1$}
\put(5.2,2.1){$-1$}
}
\end{picture}}
%\{\{6,3,2,1\},\{6,4,3,2,1,1\},\{6,5,3,3,1\},\{5,3\},\{5,3,1\}\},
%\{\{-1,1,0,0\},\{-1,-2,-1,-1,-1,-1\},\{-1,0,1,1,1\},\{1,1\},\{-1,-1,-1\}\}
\end{center}
In this example, we have $\ell=\ell^{(1)}=\ell^{(2)}=\ell_{(2)}=3$ and $m^{(2)}_3(\nu)=2$.
Then we have $x_\ell=-1$, $P^{(2)}_2(\nu)=P^{(2)}_3(\nu)=-1$, $P^{(2)}_4(\nu)=0$, $m^{(1)}_3(\nu)=1$ and
$m^{(3)}_3(\nu)=2$ (see Proposition~\ref{prop:caseD(4)3}).
We also have $P^{(2)}_2(\bar{\nu})=-1$ and $\ell_{(1)}=6\geq \ell+1=4$ (see Proposition~\ref{prop:caseD(4)1}).
Finally $(\bar{\widetilde{\nu}},\bar{\widetilde{J}})=(\widetilde{\bar{\nu}},\widetilde{\bar{J}})$ is
\begin{center}
\unitlength 10pt
{\small
\begin{picture}(42,6)
\put(0,2){
\put(-0.8,0.1){1}
\multiput(-0.8,1.1)(0,1){2}{0}
\put(-1.6,3.1){$-1$}
\Yboxdim10pt\yng(5,2,2,1)}
\put(1.2,2.1){0}
\put(2.2,3.1){0}
\put(2.2,4.1){0}
\put(5.2,5.1){$-1$}
\put(10,0){
\put(0,0){
\multiput(-1.6,0.1)(0,1){4}{$-1$}
\put(-1.6,4.1){$-2$}
\put(-1.6,5.1){$-1$}
\Yboxdim10pt\yng(6,3,2,2,1,1)}
\put(1.2,0.1){$-1$}
\put(1.2,1.1){$-1$}
\put(2.2,2.1){$-1$}
\put(2.2,3.1){$-1$}
\put(3.2,4.0){$-2$}
\put(6.1,5.1){$-1$}
}
\put(19.5,1){
\put(0,0){
\multiput(-0.8,0,1)(0,1){3}{2}
\put(-0.8,3.1){0}
\put(-1.6,4.1){$-1$}
\Yboxdim10pt\yng(6,5,2,2,1)}
\put(1.2,0.1){1}
\put(2.2,1.1){2}
\put(2.2,2.1){2}
\put(5.2,3.1){0}
\put(6.2,4.1){$-1$}
}
\put(29.0,4){
\multiput(-0.8,0.1)(0,1){2}{1}
\put(0,0){\Yboxdim10pt\yng(5,2)}
\put(2.2,0.1){1}
\put(5.2,1.1){1}
}
\put(37.0,3){
\multiput(-1.6,0.1)(0,1){3}{$-1$}
\put(0,0){\Yboxdim10pt\yng(5,2,1)}
\put(1.2,0.1){$-1$}
\put(2.2,1.1){$-1$}
\put(5.2,2.1){$-1$}
}
\end{picture}}
%\{\{5,2,2,1\},\{6,3,2,2,1,1\},\{6,5,2,2,1\},\{5,2\},\{5,2,1\}\},\\
%\{\{-1,0,0,0\},\{-1,-2,-1,-1,-1,-1\},\{-1,0,2,2,1\},\{1,1\},\{-1,-1,-1\}\}
\end{center}
\end{Example}

\subsection{Proof for Case E: preliminary steps}

\subsubsection{Classif\/ication}

Recall that the def\/ining condition for Case E is $\ell^{(i)}<\ell$.
We further classify this case as follows.
\begin{enumerate}\itemsep=0pt
\item[(1)] $\ell^{(i)}<\ell$ and $\ell+1<\ell_{(i+1)}$, \item[(2)] $\ell^{(i)}<\ell$ and
$\ell_{(i+1)}\leq\ell+1<\ell_{(i)}=\infty$, \item[(3)] $\ell^{(i)}<\ell$ and
$\ell_{(i+1)}\leq\ell+1\leq\ell_{(i)}<\infty$, \item[(4)] $\ell^{(i)}<\ell$ and $\ell_{(i)}=\ell$, \item[(5)]
$\ell^{(i)}<\ell$ and $\ell_{(i)}<\ell$.
\end{enumerate}

\subsubsection{A common property for Case E}

The fundamental observation for this case is the following property.

\begin{Proposition}
\label{prop:case_E}
If $\ell^{(i)}<\ell$, we have
\begin{gather*}%\label{eq:caseB}
P^{(i)}_{\ell^{(i)}-1}(\bar{\nu})\geq x_\ell.
\end{gather*}
\end{Proposition}
\begin{proof}
We follow the classif\/ication of Lemma~\ref{lem:vacancy_delta}.

 {\bf Case (I).} Let us consider the case $\ell^{(i-1)}=\ell^{(i)}$.
Then by Lemma~\ref{lem:vacancy_delta}, we have to show $P^{(i)}_{\ell^{(i)}-1}(\nu)\geq x_\ell$.
Suppose if possible that $P^{(i)}_{\ell^{(i)}-1}(\nu)<x_\ell$.
If $\ell^{(i)}=1$, then we have $0=P^{(i)}_{\ell^{(i)}-1}(\nu)<x_\ell$ which is a~contradiction since we have
$x_\ell\leq 0$ by $0<\ell$.
Therefore we assume $\ell^{(i)}\geq 2$ in the sequel.
From the assumption $\ell^{(i)}<\ell$ we have $P^{(i)}_{\ell^{(i)}}(\nu)\geq x_{\ell^{(i)}}\geq x_\ell$.
Let $j$ be the length of the longest string of $\nu^{(i)}$ that are strictly shorter than $\ell^{(i)}$.
If $j=\ell^{(i)}-1$, then we have $x_j\leq P^{(i)}_j(\nu)<x_\ell$, which contradicts the minimality of $x_\ell$.
Thus we assume that $j<\ell^{(i)}-1$.
We apply the convexity relation of $P^{(i)}_k$ for $j\leq k\leq\ell^{(i)}$ with the relation
$P^{(i)}_{\ell^{(i)}-1}(\nu)<P^{(i)}_{\ell^{(i)}}(\nu)$ and obtain $x_j\leq
P^{(i)}_j(\nu)<P^{(i)}_{\ell^{(i)}-1}(\nu)<x_\ell$, which is a~contradiction.
Therefore we see that $\nu^{(i)}$ does not contain strings that are strictly shorter than $\ell^{(i)}$.
Again we can use the convexity relation of $P^{(i)}_k$ for $0\leq k\leq\ell^{(i)}$ with the relation
$P^{(i)}_{\ell^{(i)}-1}(\nu)< P^{(i)}_{\ell^{(i)}}(\nu)$ and obtain the relation
$0=P^{(i)}_0(\nu)<P^{(i)}_{\ell^{(i)}-1}(\nu)<x_\ell$.
This is a~contradiction since we have $x_\ell\leq 0$.
Thus we have proved $P^{(i)}_{\ell^{(i)}-1}(\bar{\nu})\geq x_\ell$ for the case $\ell^{(i-1)}=\ell^{(i)}$.

 {\bf Case (II).} Let us consider the case $\ell^{(i-1)}<\ell^{(i)}$.
Then by Lemma~\ref{lem:vacancy_delta}, we have to show $P^{(i)}_{\ell^{(i)}-1}(\nu)>x_\ell$.
Suppose if possible that $P^{(i)}_{\ell^{(i)}-1}(\nu)\leq x_\ell$.
Let $j$ be the length of the longest row of $\nu^{(i)}$ that are strictly shorter than $\ell^{(i)}$.
Suppose that $j=\ell^{(i)}-1$.
Then we have $x_j\leq P^{(i)}_j(\nu)\leq x_\ell$.
By the minimality of $x_\ell$ this relation implies that $x_j=P^{(i)}_j(\nu)$, that is, the string $(j,x_j)$ is singular
in $(\nu,J)^{(i)}$.
This is in contradiction to the def\/inition of $\ell^{(i)}$ since we have $\ell^{(i-1)}\leq j=\ell^{(i)}-1<\ell^{(i)}$.

Therefore we assume $j<\ell^{(i)}-1$ in the sequel.
If $P^{(i)}_{\ell^{(i)}-1}(\nu)<x_\ell$, from $P^{(i)}_{\ell^{(i)}}(\nu)\geq x_{\ell^{(i)}}\geq x_\ell$ we have
$P^{(i)}_{\ell^{(i)}-1}(\nu)<P^{(i)}_{\ell^{(i)}}(\nu)$ and by the convexity relation we have $x_j\leq
P^{(i)}_j(\nu)<P^{(i)}_{\ell^{(i)}-1}(\nu)$, that is $x_j<x_\ell$.
This is a~contradiction.
Therefore we must have $P^{(i)}_{\ell^{(i)}-1}(\nu)=x_\ell$.
If $P^{(i)}_{\ell^{(i)}}(\nu)>x_\ell$, by the similar argument we obtain $x_j<x_\ell$, which is a~contradiction.
Thus we conclude that $P^{(i)}_{\ell^{(i)}}(\nu)=x_\ell$.
Suppose that $\ell^{(i-1)}\leq j$.
Then from the convexity relation of $P^{(i)}_k(\nu)$ for $j\leq k\leq \ell^{(i)}$, we obtain $P^{(i)}_{j}(\nu)\leq
P^{(i)}_{\ell^{(i)}-1}(\nu)=P^{(i)}_{\ell^{(i)}}(\nu)=x_\ell$ and from the requirement $P^{(i)}_j(\nu)\geq x_j\geq
x_\ell$ we obtain $P^{(i)}_j(\nu)=x_j$, that is, the string $(j,x_j)$ is singular.
But this is a~contradiction since $\ell^{(i-1)}\leq j<\ell^{(i)}$.
Thus we have $j<\ell^{(i-1)}$.
Since $m^{(i-1)}_{\ell^{(i-1)}}(\nu)>0$ the convexity relation for $P^{(i)}_k(\nu)$ between $\ell^{(i-1)}-1\leq
k\leq\ell^{(i-1)}+1$ must be strict.
Thus we obtain $x_\ell\geq P^{(i)}_{\ell^{(i-1)}}(\nu)>P^{(i)}_{j}(\nu)\geq x_j$, which is a~contradiction.

Therefore we conclude that $(\nu,J)^{(i)}$ does not contain strings that are strictly shorter than $\ell^{(i)}$.
Then from the convexity relation of $P^{(i)}_k(\nu)$ between $0\leq k\leq \ell^{(i)}$, we obtain $x_\ell\geq
P^{(i)}_{\ell^{(i)}-1}(\nu)\geq P^{(i)}_{\ell^{(i-1)}}(\nu)>P^{(i)}_0(\nu)=0$, which is in contradiction to the
requirement $x_\ell\leq 0$.
Thus we have proved $P^{(i)}_{\ell^{(i)}-1}(\bar{\nu})\geq x_\ell$ for the case $\ell^{(i-1)}<\ell^{(i)}$.
\end{proof}

\subsection{Proof for Case E (1)}

To begin with we prepare the following property.

\begin{Proposition}
Suppose that $\ell^{(i)}<\ell$ and $\ell+1<\ell_{(i+1)}$.
Then we have
\begin{gather*}
P^{(i)}_{\ell_{(i)}-1}(\bar{\nu})>x_\ell.
\end{gather*}
\end{Proposition}
\begin{proof}
We follow the classif\/ication of Lemma~\ref{lem:vacancy_delta}.
From the assumption, we see that $\ell^{(i)}<\ell_{(i+1)}$.
Then we only need to consider cases (V), (VI) and (VII).
During the proof of the proposition, let $j$ be the largest integer such that $j<\ell_{(i)}$ and $m^{(i)}_j(\nu)>0$.
Since $\ell^{(i)}<\ell<\ell_{(i+1)}\leq \ell_{(i)}$ we have $\ell\leq j$.
Proofs for cases (V) and (VII) are the same as the corresponding parts of the proof of Proposition~\ref{lem:case_A}.

{\sloppy  {\bf Case (VI).} In this case, we have to show $P^{(i)}_{\ell_{(i)}-1}(\nu)>x_\ell$.
Suppose if possible that $P^{(i)}_{\ell_{(i)}-1}(\nu)\leq x_\ell$.
First, let us consider the case $j=\ell_{(i)}-1$.
Then the corresponding string satisf\/ies $x_{\ell_{(i)}-1}\leq P^{(i)}_{\ell_{(i)}-1}(\nu)\leq x_\ell$.
Now recall that we have $\ell<\ell_{(i)}-1$ since we have $\ell+1<\ell_{(i+1)}\leq\ell_{(i)}$ by the assumption.
Thus the minimality of $x_\ell$ requires that $x_\ell<x_{\ell_{(i)}-1}$.
This is a~contradiction.
Similarly, if $\ell<j<\ell_{(i)}-1$, we have $x_\ell\geq P^{(i)}_{\ell_{(i)}-1}(\nu)\geq
\min\big\{P^{(i)}_{\ell_{(i)}}(\nu),P^{(i)}_j(\nu)\big\}$ which implies that $x_\ell\geq x_{\ell_{(i)}}$ or $x_\ell\geq
x_j$.
This is in contradiction to the requirements $x_\ell<x_{\ell_{(i)}}$ and $x_\ell<x_j$ derived from $\ell<\ell_{(i)}$ and
$\ell<j$, respectively.

}

Finally let us consider the case $\ell=j$.
Recall that we have $P^{(i)}_\ell(\nu)\geq x_\ell$, $P^{(i)}_{\ell_{(i)}-1}(\nu)\leq x_\ell$ and
$P^{(i)}_{\ell_{(i)}}(\nu)>x_\ell$.
However these relations are in contradiction to the convexity relation of $P^{(i)}_k(\nu)$ between $\ell\leq
k\leq\ell_{(i)}$ which we can use by $\ell<\ell_{(i)}-1$.
\end{proof}

\begin{Proposition}
Suppose that $\ell^{(i)}<\ell$ and $\ell+1<\ell_{(i+1)}$.
Then we have the following identities:
\begin{gather*}
 (1) \quad \bar{\widetilde{\nu}}=\widetilde{\bar{\nu}}, \qquad (2) \quad \bar{\widetilde{J}}=\widetilde{\bar{J}}.
\end{gather*}
\end{Proposition}
\begin{proof}
{\bf Step 1.} Let us consider the case $\delta\circ\widetilde{f}_i$.
Since $\widetilde{f}_i$ does not change coriggings of the untouched strings, and since $\widetilde{f}_i$ creates the
string whose length satisf\/ies $\ell^{(i)}<\ell+1<\ell_{(i+1)}$, we have $\widetilde{\ell}^{(a)}=\ell^{(a)}$ and
$\widetilde{\ell}_{(a)}=\ell_{(a)}$ for all~$a$.

 {\bf Step 2.} Let us consider the case $\widetilde{f}_i\circ\delta$.
Since we are assuming that $\ell^{(i)}<\ell<\ell_{(i)}$ the string $(\ell,x_\ell)$ remains as it is after~$\delta$.
Recall that we have $P^{(i)}_{\ell^{(i)}-1}(\bar{\nu})\geq x_\ell$ by Proposition~\ref{prop:case_E} and
$P^{(i)}_{\ell_{(i)}-1}(\bar{\nu})>x_\ell$ by the previous proposition.
Thus the string $(\ell,x_\ell)$ has the largest length among the strings of $(\widetilde{\nu},\widetilde{J})^{(i)}$ that
have the smallest rigging.

 {\bf Step 3.} To summarize, we have $\bar{\widetilde{\nu}}^{(a)}=\widetilde{\bar{\nu}}^{(a)}=\bar{\nu}^{(a)}$
for all $a\neq i$ and $\bar{\widetilde{\nu}}^{(i)}=\widetilde{\bar{\nu}}^{(i)}$ is obtained by adding a~box to the
$\ell$-th column of $\bar{\nu}^{(i)}$.
The statement $\bar{\widetilde{J}}=\widetilde{\bar{J}}$ is obtained by the equality
$P^{(i)}_k(\bar{\widetilde{\nu}})=P^{(i)}_k(\widetilde{\bar{\nu}})$ for arbitrary $k$ which is the consequence of
$\bar{\widetilde{\nu}}=\widetilde{\bar{\nu}}$.
\end{proof}

\subsection{Proof for Case E (2)}

Let us consider the case $\ell^{(i)}<\ell$ and $\ell_{(i+1)}\leq\ell+1<\ell_{(i)}=\infty$.
We divide the proof into two cases according to whether $\widetilde{f}_i$ creates a~singular string or not.

\subsubsection[Case 1: $\widetilde{f}_i$ creates a~non-singular string]{Case 1: $\boldsymbol{\widetilde{f}_i}$ creates a~non-singular string}

\begin{Proposition}
Suppose that we have $\ell^{(i)}<\ell$ and $\ell_{(i+1)}\leq\ell+1<\ell_{(i)}=\infty$.
Suppose further that $\widetilde{f}_i$ acting on $(\nu,J)$ creates a~non-singular string.
Then we have the following identities:
\begin{gather*}
 (1) \quad \bar{\widetilde{\nu}}=\widetilde{\bar{\nu}}, \qquad (2) \quad \bar{\widetilde{J}}=\widetilde{\bar{J}}.
\end{gather*}
\end{Proposition}
\begin{proof}
{\bf Step 1.} Let us consider the case $\delta\circ\widetilde{f}_i$.
Since $\widetilde{f}_i$ does not change coriggings of the untouched strings, and since $\widetilde{f}_i$ creates
a~non-singular string, we have $\widetilde{\ell}^{(a)}=\ell^{(a)}$ and $\widetilde{\ell}_{(a)}=\ell_{(a)}$ for all~$a$.

{\bf Step 2.} Let us consider the case $\widetilde{f}_i\circ\delta$.
Since we are assuming that $\ell^{(i)}<\ell$ and $\ell_{(i)}=\infty$ the string $(\ell,x_\ell)$ remains as it is after
$\delta$.
Recall that we have $P^{(i)}_{\ell^{(i)}-1}(\bar{\nu})\geq x_\ell$ by Proposition~\ref{prop:case_E}.
Thus the string $(\ell,x_\ell)$ has the largest length among the strings of $(\widetilde{\nu},\widetilde{J})^{(i)}$ that
have the smallest rigging.

 {\bf Step 3.} To summarize, we have $\bar{\widetilde{\nu}}^{(a)}=\widetilde{\bar{\nu}}^{(a)}=\bar{\nu}^{(a)}$
for all $a\neq i$ and $\bar{\widetilde{\nu}}^{(a)}=\widetilde{\bar{\nu}}^{(a)}$ is obtained by adding a~box to the
$\ell$-th column of $\bar{\nu}^{(i)}$.
The statement $\bar{\widetilde{J}}=\widetilde{\bar{J}}$ follows from $\bar{\widetilde{\nu}}=\widetilde{\bar{\nu}}$.
\end{proof}

\subsubsection[Case 2: $\widetilde{f}_i$ creates a~singular string]{Case 2: $\boldsymbol{\widetilde{f}_i}$ creates a~singular string}

\begin{Proposition}
\label{prop:case_E(2)2}
Suppose that we have $\ell^{(i)}<\ell$ and $\ell_{(i+1)}\leq\ell+1<\ell_{(i)}=\infty$.
Suppose further that $\widetilde{f}_i$ acting on $(\nu,J)$ creates a~singular string.
Then we have
\begin{gather*}
\widetilde{f}_i\circ\delta (\nu,J)=0.
\end{gather*}
\end{Proposition}
\begin{proof}
Since $\ell_{(i)}=\infty$,~$\delta$ changes the string $\big(\ell^{(i)},P^{(i)}_{\ell^{(i)}}(\nu)\big)$ of $(\nu,J)^{(i)}$ only.
In particular, from $\ell^{(i)}<\ell$ the string $(\ell,x_\ell)$ remains as it is after the application of~$\delta$.
Since we know that $P^{(i)}_{\ell^{(i)}-1}(\bar{\nu})\geq x_\ell$ by Proposition~\ref{prop:case_E}, the string
$(\ell,x_\ell)$ is the largest string among the strings with the smallest rigging within $(\bar{\nu},\bar{J})^{(i)}$.
Thus the next $\widetilde{f}_i$ acts on this string and produces $(\ell+1,x_\ell-1)$.

Let us compute $P^{(i)}_{\ell+1}(\widetilde{\bar{\nu}})$.
Since $\widetilde{f}_i$ acting on $(\nu,J)$ creates a~singular string, we have $P^{(i)}_{\ell+1}(\nu)=x_\ell+1$ by
Lemma~\ref{lem:f_singular}.
Then we have $P^{(i)}_{\ell+1}(\bar{\nu})=P^{(i)}_{\ell+1}(\nu)-1=x_\ell$ since we have
$\ell^{(i-1)}\leq\ell^{(i)}\leq\ell^{(i+1)}\leq\ell+1$ and $\ell_{(i+1)}\leq\ell+1<\ell_{(i)}=\infty$ by the
assumptions.
Finally, since $\widetilde{f}_i$ adds a~box to the $(\ell+1)$-th column of $(\bar{\nu},\bar{J})^{(i)}$, we have
$P^{(i)}_{\ell+1}(\widetilde{\bar{\nu}})=P^{(i)}_{\ell+1}(\bar{\nu})-2=x_\ell-2$.
Then we see that the string $(\ell+1,x_\ell-1)$ of $(\widetilde{\bar{\nu}},\widetilde{\bar{J}})^{(i)}$ satisf\/ies
$x_\ell-1>P^{(i)}_{\ell+1}(\widetilde{\bar{\nu}})$.
Therefore we conclude that $\widetilde{f}_i\circ\delta (\nu,J)=0$.
\end{proof}

\begin{Example}
Consider the following rigged conf\/iguration $(\nu,J)$ of type $(B^{1,1})^{\otimes 4}\otimes B^{1,3}\otimes
B^{2,1}\otimes B^{2,2}\otimes B^{3,1}$ of $D^{(1)}_5$
\begin{center}
\unitlength 12pt
\begin{picture}(32,6)
\put(0,1){
\multiput(-0.8,0.1)(0,1){2}{1}
\multiput(-0.8,2.1)(0,1){3}{0}
\put(0,0){\Yboxdim12pt\yng(3,3,2,1,1)}
\put(1.2,0.1){0}
\put(1.2,1.1){0}
\put(2.2,2.1){0}
\put(3.2,3.1){$-1$}
\put(3.2,4.1){0}
}
\put(6.0,0){
\multiput(-0.8,0.1)(0,1){5}{1}
\put(-0.8,5.1){0}
\put(0,0){\Yboxdim12pt\yng(4,3,3,2,1,1)}
\put(1.2,0.1){0}
\put(1.2,1.1){1}
\put(2.2,2.1){1}
\put(3.2,3.1){1}
\put(3.2,4.0){1}
\put(4.2,5.0){0}
}
\put(13,0){
\multiput(-0.8,0.1)(0,1){2}{1}
\multiput(-0.8,2.1)(0,1){3}{0}
\put(-1.53,5.1){$-1$}
\put(0,0){\Yboxdim12pt\yng(5,3,3,3,1,1)}
\put(1.2,0.1){0}
\put(1.2,1.1){0}
\put(3.2,2.1){$-1$}
\put(3.2,3.1){$-1$}
\put(3.2,4.1){0}
\put(5.2,5.1){$-1$}
}
\put(21.0,3){
\multiput(-0.8,0.1)(0,1){3}{0}
\put(0,0){\Yboxdim12pt\yng(3,3,1)}
\put(1.2,0.1){$-1$}
\put(3.2,1.1){0}
\put(3.2,2.1){0}
}
\put(27,3){
\multiput(-0.8,0.1)(0,1){2}{0}
\put(-1.53,2.1){$-2$}
\put(0,0){\Yboxdim12pt\yng(5,3,1)}
\put(1.2,0.1){0}
\put(3.2,1.1){0}
\put(5.2,2.1){$-2$}
}
\end{picture}
%\{\{3,3,2,1,1\},\{4,3,3,2,1,1\},\{5,3,3,3,1,1\},\{3,3,1\},\{5,3,1\}\},\\
%\{\{0,-1,0,0,0\},\{0,1,1,1,1,0\},\{-1,0,-1,-1,0,0\},\{0,0,-1\},\{-2,0,0\}\}
\end{center}
$\widetilde{f}_2$ acts on the string $(\ell,x_\ell)=(4,0)$ of $(\nu,J)^{(2)}$ and makes it into $(5,-1)$ of
$(\widetilde{\nu},\widetilde{J})^{(2)}$.
Since we have $P^{(2)}_5(\widetilde{\nu})=-1$, the latter string is singular.
In particular, we see that $\widetilde{f}_2(\nu,J)\neq 0$.
The corresponding tensor product $\Phi^{-1}(\nu,J)$ is
\begin{gather*}
\Yboxdim14pt \Yvcentermath1 \young(\mthree)\otimes \young(\mone)\otimes \young(\mfive)\otimes \young(3)\otimes
\young(11\mone)\otimes \young(1,\mthree)\otimes \young(1\mfive,3\mone)\otimes \young(5,\mfive,\mthree)
\end{gather*}
Then $(\bar{\nu},\bar{J})$ is
\begin{center}
\unitlength 12pt
\begin{picture}(32,6)
\put(0,1){
\multiput(-0.8,0.1)(0,1){5}{0}
\put(0,0){\Yboxdim12pt\yng(3,3,1,1,1)}
\put(1.2,0.1){0}
\put(1.2,1.1){0}
\put(1.2,2.1){0}
\put(3.2,3.1){$-1$}
\put(3.2,4.1){0}
}
\put(6.0,0){
\multiput(-0.8,0.1)(0,1){5}{1}
\put(-0.8,5.1){0}
\put(0,0){\Yboxdim12pt\yng(4,3,3,1,1,1)}
\put(1.2,0.1){0}
\put(1.2,1.1){1}
\put(1.2,2.1){1}
\put(3.2,3.1){1}
\put(3.2,4.0){1}
\put(4.2,5.0){0}
}
\put(13,0){
\multiput(-0.8,0.1)(0,1){2}{1}
\put(-0.8,2.1){0}
\multiput(-1.53,3.1)(0,1){3}{$-1$}
\put(0,0){\Yboxdim12pt\yng(4,3,3,2,1,1)}
\put(1.2,0.1){0}
\put(1.2,1.1){0}
\put(2.2,2.1){0}
\put(3.2,3.1){$-1$}
\put(3.2,4.1){$-1$}
\put(4.2,5.1){$-1$}
}
\put(21.0,3){
\multiput(-0.8,0.1)(0,1){2}{0}
\put(-0.8,2.1){1}
\put(0,0){\Yboxdim12pt\yng(3,2,1)}
\put(1.2,0.1){$-1$}
\put(2.2,1.1){0}
\put(3.2,2.1){0}
}
\put(27,3){
\multiput(-0.8,0.1)(0,1){2}{0}
\put(-1.53,2.1){$-2$}
\put(0,0){\Yboxdim12pt\yng(5,2,1)}
\put(1.2,0.1){0}
\put(2.2,1.1){0}
\put(5.2,2.1){$-2$}
}
\end{picture}
%\{\{3,3,1,1,1\},\{4,3,3,1,1,1\},\{4,3,3,2,1,1\},\{3,2,1\},\{5,2,1\}\},
%\{\{0,-1,0,0,0\},\{0,1,1,1,1,0\},\{-1,-1,-1,0,0,0\},\{0,0,-1\},\{-2,0,0\}\}
\end{center}
Thus we have $\ell^{(2)}=2<\ell=4$ and $\ell_{(3)}=5=\ell+1<\ell_{(2)}=\infty$ (see Proposition~\ref{prop:case_E(2)2}).
Indeed, we have $P^{(2)}_5(\widetilde{\bar{\nu}})=-2$ and thus $\widetilde{f}_2(\bar{\nu},\bar{J})=0$.
\end{Example}

\subsection{Proof for Case E (3)}

Let us consider the case $\ell^{(i)}<\ell$ and $\ell_{(i+1)}\leq\ell+1\leq\ell_{(i)}<\infty$.
As in Case C, we use the following classif\/ication:
\begin{enumerate}\itemsep=0pt
\item[(1)]
$m^{(i)}_{\ell+1}(\nu)>0$ and $\widetilde{f}_i$ creates a~singular string, or\\
$\ell+1=\ell_{(i)}$ and $\widetilde{f}_i$ creates a~non-singular string,
\item[(2)]
$m^{(i)}_{\ell+1}(\nu)=0$ and $\widetilde{f}_i$ creates a~singular string,
\item[(3)] $\ell+1<\ell_{(i)}$ and $\widetilde{f}_i$ creates a~non-singular string.
\end{enumerate}

\subsubsection{Proof for Case 1}

\begin{Proposition}
Suppose that $\ell^{(i)}<\ell$, $\ell_{(i+1)}\leq \ell+1\leq \ell_{(i)}<\infty$, $m^{(i)}_{\ell+1}(\nu)>0$ and
$\widetilde{f}_i$ creates a~singular string.
Then we have the following identities:
\begin{gather*}
 (1)\quad \bar{\widetilde{\nu}}=\widetilde{\bar{\nu}}, \qquad (2) \quad \bar{\widetilde{J}}=\widetilde{\bar{J}}.
\end{gather*}
\end{Proposition}
\begin{proof}
Since we have $P^{(i)}_{\ell^{(i)}-1}(\bar{\nu})\geq x_\ell$ by Proposition~\ref{prop:case_E}, we can use the same
arguments of the proof of Proposition~\ref{prop:C(1)} if we replace $\ell^{(i)}$ (resp.\
$\ell^{(i-1)}$) there by $\ell_{(i)}$ (resp.\
$\ell_{(i+1)}$) and neglect the arguments concerning the string $\ell_{(i)}$ there.
\end{proof}

\begin{Proposition}
Suppose that $\ell^{(i)}<\ell$, $\ell_{(i+1)}\leq \ell+1\leq \ell_{(i)}<\infty$, $\ell+1=\ell_{(i)}$ and
$\widetilde{f}_i$ creates a~non-singular string.
Then we have the following identities:
\begin{gather*}
 (1) \quad \bar{\widetilde{\nu}}=\widetilde{\bar{\nu}}, \qquad (2) \quad \bar{\widetilde{J}}=\widetilde{\bar{J}}.
\end{gather*}
\end{Proposition}
\begin{proof}
As in Case C we can reduce the proof to the previous proposition.
\end{proof}

\subsubsection{Proof for Case 2}

As in Case C, we show the following property:
\begin{Proposition}
\label{prop:caseE(3)2}
Suppose that $\ell^{(i)}<\ell$, $\ell_{(i+1)}\leq \ell+1<\ell_{(i)}<\infty$ and $m^{(i)}_{\ell+1}(\nu)=0$.
If $\widetilde{f}_i$ acting on $(\nu,J)$ creates a~singular string, then the following two conditions are satisfied:
\begin{gather*}
P^{(i)}_{\ell_{(i)}-1}(\bar{\nu})= x_\ell,
%\label{eq:caseE(3)1}
\qquad
m_k^{(i-1)}(\nu)=0\,\text{ for }\,\ell<k<\ell_{(i)}.
%\label{eq:caseE(3)2}
\end{gather*}
\end{Proposition}
\begin{proof}
As in Case C, the proof depend on the following two cases (see Lemma~\ref{lem:C_1});
\begin{gather}
(\ell,x_\ell)\text{ is singular and }\ell=\ell_{(i+1)}-1,
\label{eq:caseE(3)3}
\\
(\ell,x_\ell)\text{ is non-singular.}
\label{eq:caseE(3)4}
\end{gather}
For the case~\eqref{eq:caseE(3)3} (resp.~\eqref{eq:caseE(3)4}) we can use the same arguments of
Section~\ref{sec:caseC(2)3} (resp.\
Section~\ref{sec:caseC(2)4}) by replacing $\ell^{(i)}$ (resp.~$\ell^{(i-1)}$) there by $\ell_{(i)}$ (resp.~$\ell_{(i+1)}$) to show the assertions.
\end{proof}

\begin{Proposition}
Suppose that $\ell^{(i)}<\ell$, $\ell_{(i+1)}\leq \ell+1<\ell_{(i)}<\infty$ and $m^{(i)}_{\ell+1}(\nu)=0$.
Then we have $\bar{\widetilde{\nu}}=\widetilde{\bar{\nu}}$ and $\bar{\widetilde{J}}=\widetilde{\bar{J}}$.
\end{Proposition}
\begin{proof}
Since we have $P^{(i)}_{\ell^{(i)}-1}(\bar{\nu})\geq x_\ell$ by Proposition~\ref{prop:case_E}, we can use the same
arguments of the proof of Proposition~\ref{prop:caseC(2)1} if we replace $\ell^{(i)}$ (resp.~$\ell^{(i-1)}$) there by $\ell_{(i)}$ (resp.~$\ell_{(i+1)}$), neglect the arguments concerning the string $\ell_{(i)}$ there and use Proposition~\ref{prop:caseE(3)2}
instead of Proposition~\ref{prop:caseC1}.
Thus we obtain $\bar{\widetilde{\nu}}=\widetilde{\bar{\nu}}$.
We can also use the same analysis given for Case~C to show $\bar{\widetilde{J}}=\widetilde{\bar{J}}$.
\end{proof}

\subsubsection{Proof for Case 3}

We can use the same arguments in Section~\ref{sec:caseC(3)} if we replace $\ell^{(i)}$ (resp.~$\ell^{(i-1)}$) there by $\ell_{(i)}$ (resp.~$\ell_{(i+1)}$) and neglect the arguments concerning the string $\ell_{(i)}$ there.

\subsection{Proof for Case E (4)}

\subsubsection{Preliminary}

In this case, we assume that $\ell^{(i)}<\ell$ and $\ell_{(i)}=\ell$.
In view of the classif\/ication of Case D we have to consider the following two cases:
\begin{enumerate}\itemsep=0pt
\item[1)] $\ell^{(i)}<\ell$, $\ell_{(i)}=\ell$ and $m^{(i)}_{\ell}(\nu)>1$,
\item[3)] $\ell^{(i)}<\ell$,
$\ell_{(i)}=\ell$ and $m^{(i)}_{\ell}(\nu)=1$.
\end{enumerate}
During the proof of Case E (4) we follow the above classif\/ication.
In both cases, we have the following common property.

\begin{Proposition}
\label{prop:caseE(4)1}
Assume that $\ell^{(i)}<\ell$ and $\ell_{(i)}=\ell$.
Then the following relation holds:
\begin{gather*}
P^{(i)}_{\ell-1}(\bar{\nu})\geq x_\ell.
\end{gather*}
\end{Proposition}
\begin{proof}
We can use the same arguments in the proof of Proposition~\ref{prop:caseD1} if we replace $\ell^{(i)}$ (resp.~$\ell^{(i-1)}$) there by $\ell_{(i)}$ (resp.~$\ell_{(i+1)}$).
Here in order to evaluate the vacancy numbers in Case~(a) and~(b), we note that if we have $\ell-1<\ell_{(i+1)}$ we have
$P^{(i)}_{\ell-1}(\bar{\nu})\geq P^{(i)}_{\ell-1}(\nu)$ since we have $\ell^{(i)}\leq\ell-1$ by the assumption.
\end{proof}

\subsubsection{Proof for Case 1}

\begin{Proposition}
Assume that $\ell^{(i)}<\ell$, $\ell_{(i)}=\ell$ and $m^{(i)}_{\ell}(\nu)>1$.
Then we have the following identities:
\begin{gather*}
 (1) \quad \bar{\widetilde{\nu}}=\widetilde{\bar{\nu}}, \qquad (2) \quad \bar{\widetilde{J}}=\widetilde{\bar{J}}.
\end{gather*}
\end{Proposition}
\begin{proof}
(1) Let us analyze the action of $\widetilde{f}_i$ before and after the application of~$\delta$.
By the assumption $m^{(i)}_{\ell}(\nu)>1$ we can choose two distinct length $\ell$ strings $(\ell,x_\ell)$ and
$\big(\ell,P^{(i)}_\ell(\nu)\big)$ of $(\nu,J)^{(i)}$ where $\widetilde{f}_i$ will act on the former one and~$\delta$ will act
on the latter one.
Recall that we have $P^{(i)}_{\ell^{(i)}-1}(\bar{\nu})\geq x_\ell$ by Proposition~\ref{prop:case_E} and
$P^{(i)}_{\ell_{(i)}-1}(\bar{\nu})\geq x_\ell$ by Proposition~\ref{prop:caseE(4)1}.
Thus both strings $\big(\ell^{(i)}-1,P^{(i)}_{\ell^{(i)}-1}(\bar{\nu})\big)$ and
$\big(\ell_{(i)}-1,P^{(i)}_{\ell_{(i)}-1}(\bar{\nu})\big)$ of $(\bar{\nu},\bar{J})^{(i)}$ are shorter than $(\ell,x_\ell)$ and
their riggings are larger than or equal to $x_\ell$.
Therefore $\widetilde{f}_i$ will act on the same string $(\ell,x_\ell)$ in both $(\nu,J)$ and $(\bar{\nu},\bar{J})$.

Let us analyze~$\delta$.
Recall that $\widetilde{f}_i$ will not change the coriggings of the strings $\big(\ell^{(i)},P^{(i)}_{\ell^{(i)}}(\nu)\big)$ and
$\big(\ell_{(i)},P^{(i)}_{\ell_{(i)}}(\nu)\big)$ of $(\nu,J)^{(i)}$.
Since $\ell^{(i)}<\ell$ we have $\widetilde{\ell}^{(a)}=\ell^{(a)}$ for all~$a$.
Also, since $\ell_{(i)}<\ell+1$ we have $\widetilde{\ell}_{(a)}=\ell_{(a)}$ for all~$a$.
To summarize~$\delta$ chooses the same strings before and after the application of $\widetilde{f}_i$.
Hence we obtain $\bar{\widetilde{\nu}}=\widetilde{\bar{\nu}}$.

(2) Since each $\widetilde{f}_i$ and~$\delta$ chooses the same strings in $\delta\circ\widetilde{f}_i$ and
$\widetilde{f}_i\circ\delta$, the coincidence $\bar{\widetilde{J}}=\widetilde{\bar{J}}$ is a~consequence of
$P^{(i)}_k(\bar{\widetilde{\nu}})=P^{(i)}_k(\widetilde{\bar{\nu}})$.
\end{proof}

\subsubsection{Proof for Case 3}

\begin{Proposition}
Assume that $\ell^{(i)}<\ell$, $\ell_{(i)}=\ell$ and $m^{(i)}_{\ell}(\nu)=1$.
Then we have the following identities:
\begin{gather}
P^{(i)}_{\ell+1}(\nu)=x_\ell+1,
\\
P^{(i)}_{\ell-1}(\bar{\nu})=x_\ell,
\label{eq:caseE(4)1}
\\
\ell+1\leq\ell_{(i-1)}.
\end{gather}
\end{Proposition}
\begin{proof}
The proof is divided into two cases $\ell_{(i+1)}<\ell$ and $\ell_{(i+1)}=\ell$ and we can use the same arguments in the
proof of Proposition~\ref{prop:caseD2} if we replace $\ell^{(i-1)}$, $\ell^{(i)}$ and $\ell^{(i+1)}$ there by
$\ell_{(i+1)}$, $\ell_{(i)}$ and $\ell_{(i-1)}$, respectively.
\end{proof}

\begin{Proposition}
Assume that $\ell^{(i)}<\ell$, $\ell_{(i)}=\ell$ and $m^{(i)}_{\ell}(\nu)=1$.
Then we have the following identities:
\begin{gather*}
 (1) \quad \bar{\widetilde{\nu}}=\widetilde{\bar{\nu}}, \qquad (2)\quad \bar{\widetilde{J}}=\widetilde{\bar{J}}.
\end{gather*}
\end{Proposition}
\begin{proof}
Recall that we have $P^{(i)}_{\ell^{(i)}-1}(\bar{\nu})\geq x_\ell$ by Proposition~\ref{prop:case_E}.
Since we have $P^{(i)}_{\ell_{(i)}-1}(\bar{\nu})=x_\ell$ by~\eqref{eq:caseE(4)1}, the string
$\big(\ell^{(i)}-1,P^{(i)}_{\ell^{(i)}-1}(\bar{\nu})\big)$ of $(\bar{\nu},\bar{J})^{(i)}$ will not be chosen by the next
$\widetilde{f}_i$.
Then we can use the same arguments in the corresponding parts of Case D (3) if we replace $\ell^{(i)}$ (resp.~$\ell^{(i-1)}$) there by $\ell_{(i)}$ (resp.~$\ell_{(i+1)}$) and neglect all the arguments related to~$\ell_{(a)}$.
\end{proof}

\subsection{Proof for Case E (5)}

In this case, we can show the following property.
\begin{Proposition}
\label{prop:caseE(5)1}
Suppose that $\ell_{(i)}<\ell$.
Then we have
\begin{gather*}
P^{(i)}_{\ell_{(i)}-1}(\bar{\nu})\geq x_\ell.
\end{gather*}
\end{Proposition}
\begin{proof}
We can use the arguments of the proof of Proposition~\ref{prop:case_E} if we replace $\ell^{(i)}$ (resp.~$\ell^{(i-1)}$) there by $\ell_{(i)}$ (resp.~$\ell_{(i+1)}$).
\end{proof}
As the consequence of this proposition, we have the following result.
\begin{Proposition}
Suppose that $\ell^{(i)}<\ell$ and $\ell_{(i)}<\ell$.
Then we have $\bar{\widetilde{\nu}}=\widetilde{\bar{\nu}}$ and $\bar{\widetilde{J}}=\widetilde{\bar{J}}$.
\end{Proposition}
\begin{proof}
Since we have $P^{(i)}_{\ell^{(i)}-1}(\bar{\nu})\geq x_\ell$ by Proposition~\ref{prop:case_E} and
$P^{(i)}_{\ell_{(i)}-1}(\bar{\nu})\geq x_\ell$ by Proposition~\ref{prop:caseE(5)1}, we see that $\widetilde{f}_i$ acts
on the same string before and after~$\delta$.
Also, since $\ell^{(i)}\leq\ell_{(i)}<\ell$ and the fact that $\widetilde{f}_i$ does not change coriggings of untouched
strings, we see that~$\delta$ acts on the same strings before and after $\widetilde{f}_i$.
Thus we have $\bar{\widetilde{\nu}}=\widetilde{\bar{\nu}}$.
$\bar{\widetilde{J}}=\widetilde{\bar{J}}$ follows from the fact
$P^{(i)}_k(\bar{\widetilde{\nu}})=P^{(i)}_k(\widetilde{\bar{\nu}})$ which is the consequence of
$\bar{\widetilde{\nu}}=\widetilde{\bar{\nu}}$.
\end{proof}

\section[Proof of Proposition~\ref{th:core2}: $\widetilde{f}_i$ version]{Proof
of Proposition~\ref{th:core2}: $\boldsymbol{\widetilde{f}_i}$ version}\label{sec:main2}

\subsection{Proof for (1)}

\begin{Proposition}
\label{prop:nonzero_zero1}
Let us consider the rigged configuration $(\nu,J)$ of type $B^{1,1}\otimes\bar{B}$.
Suppose that we have the commutativity of $\widetilde{f}_i$ and~$\Phi$ for $\bar{B}$.
Suppose that we have $\widetilde{f}_i(\nu,J)\neq 0$ and $\widetilde{f}_i(\bar{\nu},\bar{J})=0$.
Let $b=\Phi(\nu,J)$ and $b'=\Phi(\bar{\nu},\bar{J})$.
Then we have one of the following two cases:
\begin{enumerate}\itemsep=0pt
\item[$(1)$] $b=i\otimes b'$, $\widetilde{f}_i(b)$ is defined, $\widetilde{f}_i(b')$ is undefined and
$\Phi\bigl(\widetilde{f}_i(\nu,J)\bigr)=(i+1)\otimes b'$,
\item[$(2)$] $b=\overline{i+1}\otimes b'$, $\widetilde{f}_i(b)$
is defined, $\widetilde{f}_i(b')$ is undefined and $\Phi\bigl(\widetilde{f}_i(\nu,J)\bigr)=\overline{i}\otimes b'$.
\end{enumerate}
\end{Proposition}

For the proof, we start by preparing the following properties.

\begin{Proposition}
\label{lem:nonzero_zero1}
$\widetilde{f}_i(\nu,J)\neq 0$ and $\widetilde{f}_i(\bar{\nu},\bar{J})=0$ if and only if one the following conditions is
satisfied:
\begin{enumerate}\itemsep=0pt
\item[$(1)$] $\ell^{(i-1)}\leq \ell+1<\ell^{(i)}=\infty$ and $\widetilde{f}_i$ acting on $(\nu,J)$ creates a~singular
string,
\item[$(2)$] $\ell^{(i)}<\ell$, $\ell_{(i+1)}\leq\ell+1<\ell_{(i)}=\infty$ and $\widetilde{f}_i$ acting on $(\nu,J)$
creates a~singular string.
\end{enumerate}
\end{Proposition}
\begin{proof}
In Section~\ref{sec:main}, we analyze all possible cases such that $\widetilde{f}_i(\nu,J)\neq 0$.
Then only cases such that $\widetilde{f}_i(\nu,J)$ is def\/ined and $\widetilde{f}_i(\bar{\nu},\bar{J})$ is undef\/ined are
the above two cases as described in Proposition~\ref{f_undefined} and Proposition~\ref{prop:case_E(2)2}.
\end{proof}

\begin{Lemma}{\samepage\quad
\begin{enumerate}\itemsep=0pt
\item[$(1)$] Suppose that $\ell^{(i-1)}\leq \ell+1<\ell^{(i)}=\infty$ and $\widetilde{f}_i$ acting on $(\nu,J)$ creates
a~singular string.
Then we have $\widetilde{\ell}^{(i)}=\ell+1$, $\widetilde{\ell}^{(i+1)}=\infty$ and $P^{(i)}_\ell(\bar{\nu})=x_\ell$.

\item[$(2)$] Suppose that $\ell^{(i)}<\ell$, $\ell_{(i+1)}\leq\ell+1<\ell_{(i)}=\infty$ and $\widetilde{f}_i$ acting on $(\nu,J)$
creates a~singular string.
Then we have $\widetilde{\ell}_{(i)}=\ell+1$, $\widetilde{\ell}_{(i-1)}=\infty$ and $P^{(i)}_\ell(\bar{\nu})=x_\ell$.
\end{enumerate}}
\end{Lemma}
\begin{proof}
(1) To begin with we shall show that there is no string of $\bar{\nu}^{(i)}$ that is longer than $\ell$.
Suppose if possible that such string exist.
Let $j$ be the minimal integer satisfying $j>\ell$ and $m^{(i)}_j(\nu)>0$.
Then there are strings $(\ell,x_\ell)$ and $(j,x_{j})$ of $(\nu,J)^{(i)}$ that satisfy $x_{j}>x_\ell$ by def\/inition of
$\ell$.
Since $\ell^{(i)}=\infty$ we have $(\bar{\nu},\bar{J})^{(i)}=(\nu,J)^{(i)}$.
Thus the string $(j,x_{j})$ still exists in $(\bar{\nu},\bar{J})^{(i)}$, which implies that $P^{(i)}_{j}(\bar{\nu})\geq
x_{j}$.
To summarize, we obtain $P^{(i)}_{j}(\bar{\nu})>x_\ell$.

Consider the case $j=\ell+1$.
Then this result is in contradiction to the equality $P^{(i)}_{\ell+1}(\bar{\nu})=x_\ell$ that appears in the proof of
Proposition~\ref{f_undefined}.
Next consider the general case $j>\ell+1$.
By the convexity relation between $\ell$ and $j$, we have $P^{(i)}_{\ell+1}(\bar{\nu})\geq\min
\big\{P^{(i)}_\ell(\bar{\nu}),P^{(i)}_j(\bar{\nu})\big\}$.
Since we know that $P^{(i)}_{\ell+1}(\bar{\nu})=x_\ell<P^{(i)}_{j}(\bar{\nu})$, we deduce that
$P^{(i)}_{\ell+1}(\bar{\nu})\geq P^{(i)}_{\ell}(\bar{\nu})$.
From $P^{(i)}_{\ell+1}(\bar{\nu})=x_\ell$ and $P^{(i)}_{\ell}(\bar{\nu})\geq x_\ell$ (since $\ell^{(i)}=\infty$), we
obtain $P^{(i)}_{\ell}(\bar{\nu})=P^{(i)}_{\ell+1}(\bar{\nu})=x_\ell$.
Then the convexity relation between $\ell$ and $j$ gives $x_\ell=P^{(i)}_{\ell}(\bar{\nu})\geq P^{(i)}_{j}(\bar{\nu})$
which is in contradiction to the relation $P^{(i)}_{j}(\bar{\nu})>x_\ell$.
Hence there is no string of $\bar{\nu}^{(i)}$ that is longer than $\ell$.

The relations $P^{(i)}_{\ell}(\bar{\nu})\geq x_\ell$, $P^{(i)}_{\ell+1}(\bar{\nu})=x_\ell$,
$P^{(i)}_\infty(\bar{\nu})>-\infty$ and the convexity relation bet\-ween~$\ell$ and~$\infty$ gives
$P^{(i)}_{\ell}(\bar{\nu})=P^{(i)}_{\ell+1}(\bar{\nu})= P^{(i)}_{\ell+2}(\bar{\nu})=\cdots$.
This relation forces $m^{(i+1)}_k(\bar{\nu})=0$ for all $k\geq\ell+1$, since otherwise the relation would be strictly
convex.
Since $\ell^{(i)}=\infty$, we have $\bar{\nu}^{(i+1)}=\nu^{(i+1)}$ and thus $\bar{\nu}^{(i+1)}=\widetilde{\nu}^{(i+1)}$.
Therefore we obtain $m^{(i+1)}_k(\widetilde{\nu})=0$ for all $k\geq\ell+1$.

Recall that the def\/inition of $\widetilde{f}_i$ on the rigged conf\/igurations implies that
$\widetilde{\ell}^{(a)}=\ell^{(a)}$ for all $a\leq i-1$.
Then by the assumption we have $\widetilde{\ell}^{(i-1)}\leq\ell +1$.
Again from the assumption there is a~length $\ell+1$ singular string in $\widetilde{\nu}^{(i)}$.
Therefore we have $\widetilde{\ell}^{(i)}=\ell+1$.
Finally from $m^{(i+1)}_k(\widetilde{\nu})=0$ for all $k\geq\ell+1$ we have $\widetilde{\ell}^{(i+1)}=\infty$.

(2) During the proof of Proposition~\ref{prop:case_E(2)2} we show $P^{(i)}_{\ell+1}(\bar{\nu})=x_\ell$.
Then we can use the same arguments of (1) if we replace $\ell^{(i-1)}$, $\ell^{(i)}$ and $\ell^{(i+1)}$ etc.
there by $\ell_{(i+1)}$, $\ell_{(i)}$ and $\ell_{(i-1)}$ etc., respectively.
Note that in this case we have to work under the condition $\ell^{(i)}<\ell$ and $\ell_{(i)}=\infty$.
Then we can use the fact that all strings of $(\nu,J)^{(i-1)}$ that are strictly longer than $\ell^{(i-1)}(<\ell+1)$ do
not change after~$\delta$.
In particular, we can show $m^{(i-1)}_k(\bar{\nu})=m^{(i-1)}_k(\widetilde{\nu})=0$ for all $k\geq\ell+1$.
\end{proof}

\begin{proof}
[Proof of Proposition~\ref{prop:nonzero_zero1}] According to the previous observations, we have to consider two distinct
cases.

(1) Suppose that $\ell^{(i-1)}\leq \ell+1<\ell^{(i)}=\infty$ and $\widetilde{f}_i$ acting on $(\nu,J)$ creates
a~singular string.
From $\ell^{(i-1)}<\infty$ and $\ell^{(i)}=\infty$ we have $\Phi(\nu,J)=i\otimes b'$.
From $\widetilde{\ell}^{(i)}<\infty$ and $\widetilde{\ell}^{(i+1)}=\infty$ we have $\Phi(\widetilde{\nu},\widetilde{J})=
\Phi\circ\widetilde{f}_i(\nu,J)=(i+1)\otimes b''$.
Let us show the coincidence $b''=b'$ by proving $(\bar{\nu},\bar{J})=(\bar{\widetilde{\nu}},\bar{\widetilde{J}})$.
Recall that the dif\/ference between $(\nu,J)$ and $(\widetilde{\nu},\widetilde{J})$ is one box at $(\ell+1)$-th column of
$(\widetilde{\nu},\widetilde{J})^{(i)}$.
Then from $\widetilde{\ell}^{(a)}=\ell^{(a)}$ for all $a\leq i-1$, $\widetilde{\ell}^{(i)}=\ell+1$ and
$\widetilde{\ell}^{(i+1)}=\ell^{(i)}=\infty$, we have $\bar{\nu}=\bar{\widetilde{\nu}}$.
In order to show $\bar{J}=\bar{\widetilde{J}}$ recall that~$\delta$ acting on $(\widetilde{\nu},\widetilde{J})$ creates
the length $\ell$ singular string of $(\bar{\widetilde{\nu}},\bar{\widetilde{J}})^{(i)}$.
On the other hand, since $P^{(i)}_\ell(\bar{\nu})=x_\ell$ all length $\ell$ strings of $(\bar{\nu},\bar{J})^{(i)}$ are
singular.
From $P^{(i)}_\ell(\bar{\nu})=P^{(i)}_\ell(\bar{\widetilde{\nu}})$ we see the coincidence of the corresponding riggings
and hence we deduce that $\bar{J}=\bar{\widetilde{J}}$.
Thus we have $b''=b'$.

As for $b'=\Phi(\bar{\nu},\bar{J})$, we know that $\varphi_i(b')=\varphi_i(\bar{\nu},\bar{J})=0$ by the assumption.
Thus $\widetilde{f}_i\circ\Phi(\nu,J)=\widetilde{f}_i(i)\otimes b'=(i+1)\otimes b'$.
Thus~$\Phi$ and $\widetilde{f}_i$ commutes in this case.

(2) Suppose that $\ell^{(i)}<\ell$, $\ell_{(i+1)}\leq\ell+1<\ell_{(i)}=\infty$ and $\widetilde{f}_i$ acting on $(\nu,J)$
creates a~singular string.
From $\ell_{(i+1)}<\infty$ and $\ell_{(i)}=\infty$ we have $\Phi(\nu,J)=\overline{i+1}\otimes b'$.
From $\widetilde{\ell}_{(i)}<\infty$ and $\widetilde{\ell}_{(i-1)}=\infty$ we have
$\Phi(\widetilde{\nu},\widetilde{J})=\overline{i}\otimes b''$.
The coincidence $b''=b'$ follows from parallel arguments.
Finally, by $\varphi_i(b')=\varphi_i(\bar{\nu},\bar{J})=0$ we have
$\widetilde{f}_i\circ\Phi(\nu,J)=\widetilde{f}_i(\overline{i+1})\otimes b'=\overline{i}\otimes b'$.
Thus~$\Phi$ and $\widetilde{f}_i$ commutes in this case also.
\end{proof}

\subsection{Proof for (2)}

To begin with, let us show the following property.
\begin{Lemma}
\label{lem:sbars}
Suppose that $\ell_{(i+1)}<\infty$, $\ell_{(i)}=\infty$, $\widetilde{\ell}_{(i)}<\infty$ and
$\widetilde{\ell}_{(i-1)}=\infty$.
Let $s$ $($resp.~$\bar{s})$ be the smallest rigging of $(\nu,J)^{(i)}$ $($resp.~$(\bar{\nu},\bar{J})^{(i)})$.
Then we have $s\leq\bar{s}$.
\end{Lemma}
\begin{proof}
Note that if we show that the rigging for the shortened string in $(\bar{\nu},\bar{J})^{(i)}$ is larger than or equal to
$s$, then we have $s\leq\bar{s}$.
In particular, we can assume that $\ell^{(i)}>1$, since otherwise the corresponding rigging will be erased by~$\delta$
and will not contribute to $\bar{s}$, thus $s=\bar{s}$ by $\ell_{(i)}=\infty$.
Since $\ell_{(i)}=\infty$, it is enough to check
\begin{gather*}
P^{(i)}_{\ell^{(i)}-1}(\bar{\nu})\geq s.
\end{gather*}
Then we have the two possibilities:
\begin{gather*}
 \text{(I)}\quad \ell^{(i-1)}=\ell^{(i)} \qquad \text{or} \qquad \text{(II)} \quad \ell^{(i-1)}<\ell^{(i)}.
\end{gather*}
During the proof, let $j$ be the largest integer such that $j<\ell^{(i)}$ and $m^{(i)}_j(\nu)>0$.
If there is no such $j$, we set $j=0$.

{\bf Case (I).} In this case we have to show that $P^{(i)}_{\ell^{(i)}-1}(\nu)\geq s$ under the assumption
$\ell^{(i-1)}=\ell^{(i)}$.
Let us consider the case $j=\ell^{(i)}-1$.
Then we are done since we have $P^{(i)}_{\ell^{(i)}-1}(\nu)\geq x_{\ell^{(i)}-1}\geq s$ by the minimality of $s$.
Next consider the case $0<j<\ell^{(i)}-1$.
Then by the convexity relation between $j$ and $\ell^{(i)}$ we have
$P^{(i)}_{\ell^{(i)}-1}(\nu)\geq\min\big\{P^{(i)}_j(\nu),P^{(i)}_{\ell^{(i)}}(\nu)\big\}\geq s$.
Suppose that $j=0$.
If $s\leq 0$, we can use the same arguments of the previous case since $P^{(i)}_0(\nu)=0$.

Therefore we have to consider the case $j=0$ and $s>0$.
Since $s>0$, $\widetilde{f}_i$ acting on $(\nu,J)$ adds the length $1$ string to $(\nu,J)^{(i)}$.
Recall that $\widetilde{f}_i$ does not change coriggings of untouched strings.
Since we have $\ell^{(i-1)}=\ell^{(i)}>1$, we see that $\widetilde{\ell}^{(a)}=\ell^{(a)}$ and
$\widetilde{\ell}_{(a)}=\ell_{(a)}$ for all~$a$.
However this is in contradiction to the assumptions $\ell_{(i)}=\infty$ and $\widetilde{\ell}_{(i)}<\infty$.
Hence this case cannot happen.

{\bf Case (II).} In this case we have to show that $P^{(i)}_{\ell^{(i)}-1}(\nu)>s$ under the assumption
$\ell^{(i-1)}<\ell^{(i)}$.
Suppose if possible that $P^{(i)}_{\ell^{(i)}-1}(\nu)\leq s$.

Suppose that $j=\ell^{(i)}-1$.
Then we have $x_{\ell^{(i)}-1}\leq P^{(i)}_{\ell^{(i)}-1}(\nu)\leq s$.
Then by the minimality of $s$, we have $x_{\ell^{(i)}-1}=P^{(i)}_{\ell^{(i)}-1}(\nu)=s$, in particular, the string
$(\ell^{(i)}-1,x_{\ell^{(i)}-1})$ is singular.
However this is a~contradiction since its length satisf\/ies $\ell^{(i-1)}\leq \ell^{(i)}-1<\ell^{(i)}$.

Suppose that $\ell^{(i-1)}\leq j<\ell^{(i)}-1$.
By def\/inition of $s$ we have $P^{(i)}_{\ell^{(i)}}(\nu)\geq s$ and by the assumption we have
$P^{(i)}_{\ell^{(i)}-1}(\nu)\leq s$.
Then by the convexity relation of $P^{(i)}_k(\nu)$ between $j\leq k\leq\ell^{(i)}$, we have
$P^{(i)}_{\ell^{(i)}}(\nu)\geq P^{(i)}_{\ell^{(i)}-1}(\nu) \geq P^{(i)}_{j}(\nu)\geq x_j$.
By the minimality of $s$ and $s\geq P^{(i)}_{\ell^{(i)}-1}(\nu)$, we have $P^{(i)}_{j}(\nu)=x_j$, that is, the string
$(j,x_j)$ is singular.
This contradicts the def\/inition of $\ell^{(i)}$.

Suppose that $0<j<\ell^{(i-1)}$.
By $m^{(i-1)}_{\ell^{(i-1)}}(\nu)>0$, the convexity relation of $P^{(i)}_k(\nu)$ between $\ell^{(i-1)}+1\geq
k\geq\ell^{(i-1)}-1$ has to be strictly convex.
Therefore we have $P^{(i)}_{\ell^{(i)}}(\nu)\geq P^{(i)}_{\ell^{(i)}-1}(\nu) >P^{(i)}_{j}(\nu)\geq x_j$ in this case.
Since we have $s\geq P^{(i)}_{\ell^{(i)}-1}(\nu)$, this contradicts the minimality of $s$.

Suppose that $j=0$.
Then we have $s\geq P^{(i)}_{\ell^{(i)}-1}(\nu)>P^{(i)}_{0}(\nu)=0$.
If $s\leq 0$, this is a~contradiction.
Hence we assume that $s>0$.
Then $\widetilde{f}_i$ adds the length 1 string $(1,-1)$ to $(\nu,J)^{(i)}$.
If $\ell^{(i-1)}>1$, we can use the same arguments in the last part of Case (I) to show that such case cannot happen.
Thus assume that $\ell^{(i-1)}=1$.
If $P^{(i)}_1(\widetilde{\nu})>-1$, then the added string $(1,-1)$ is non-singular so that we can use the arguments in
the last part of Case (I) to show that such case cannot happen.

To summarize, we have shown that we have $s\leq\bar{s}$ except for the case $j=0$, $s>0$, $\ell^{(i-1)}=1$ and
$P^{(i)}_1(\widetilde{\nu})=-1$.
Suppose if possible that this case happen.
Since $s>0$, $\widetilde{f}_i$ adds the length 1 string $(1,-1)$ to $(\nu,J)^{(i)}$.
Thus we have $\bar{s}<0$ and $P^{(i)}_1(\nu)=1$.
Since $m^{(i-1)}_{\ell^{(i-1)}}(\nu)=m^{(i-1)}_1(\nu)>0$ and $\ell^{(i)}>1$, the relation $P^{(i)}_k(\nu)$ for $0\leq
k\leq 2$ must be strictly convex.
In view of $P^{(i)}_0(\nu)=0$, this implies that $1=P^{(i)}_1(\nu)\geq P^{(i)}_2(\nu)$.
Since we are assuming that $s>0$, the rigging of the string $\big(\ell^{(i)},P^{(i)}_{\ell^{(i)}}(\nu)\big)$ must satisfy
$P^{(i)}_{\ell^{(i)}}(\nu)>0$.
Then by the convexity relation of $P^{(i)}_k(\nu)$ between $1\leq k\leq\ell^{(i)}$ the only possibility that is
compatible with $1=P^{(i)}_1(\nu)\geq P^{(i)}_2(\nu)$ is the case
$1=P^{(i)}_1(\nu)=P^{(i)}_2(\nu)=\dots=P^{(i)}_{\ell^{(i)}}(\nu)$.
To summarize, we have the following situation:
\begin{gather}
\label{eq:fvanish(2)_1}
s>0,\, \ell^{(i-1)}=1,\, \ell^{(i)}>1
\qquad
\text{and}
\qquad
P^{(i)}_1(\nu)=P^{(i)}_2(\nu)=\dots=P^{(i)}_{\ell^{(i)}}(\nu)=1.
\end{gather}
Since $j=0$ we have $m^{(i)}_1(\nu)=0$.
Then from the convexity relation of Lemma~\ref{lem:convexity3} with $l=1$ we have
\begin{gather*}
-P^{(i)}_0(\nu)+2P^{(i)}_1(\nu)-P^{(i)}_2(\nu)=1
\\
\qquad
{}\geq m^{(i-1)}_1(\nu)-2m^{(i)}_1(\nu)+m^{(i+1)}_1(\nu)
= m^{(i-1)}_1(\nu)+m^{(i+1)}_1(\nu).
\end{gather*}
We have $m^{(i-1)}_1(\nu)\geq 1$ since $\ell^{(i-1)}=1$ and $m^{(i+1)}_1(\nu)\geq 0$ by def\/inition.
From the above inequality we conclude that $m^{(i-1)}_1(\nu)=1$ and $m^{(i+1)}_1(\nu)=0$.
We also have $m^{(i+1)}_k(\nu)=0$ for $1<k<\ell^{(i)}$ by the relation
$P^{(i)}_1(\nu)=P^{(i)}_2(\nu)=\dots=P^{(i)}_{\ell^{(i)}}(\nu)$.
Since $\widetilde{f}_i$ does not change the coriggings of $(\nu,J)^{(a)}$ for $a<i$, we have
$\widetilde{\ell}^{(i-1)}=\ell^{(i-1)}=1$.
Recall that the string $(1,-1)$ of $(\widetilde{\nu},\widetilde{J})^{(i)}$ is singular.
Thus $\widetilde{\ell}^{(i)}=1$.
From $\widetilde{\nu}^{(i+1)}=\nu^{(i+1)}$, we have $m^{(i+1)}_k(\widetilde{\nu})=0$ for $1\leq k<\ell^{(i)}$.
Since $\widetilde{f}_i$ does not change the coriggings of $(\nu,J)^{(i+1)}$, we have
$\widetilde{\ell}^{(i+1)}=\ell^{(i+1)}$, and thus $\widetilde{\ell}^{(a)}=\ell^{(a)}$ for all $i<a$ and
$\widetilde{\ell}_{(a)}=\ell_{(a)}$ for all~$a$.
In particular we have $\widetilde{\ell}_{(i)}=\ell_{(i)}=\infty$.
However this is a~contradiction since we have $\widetilde{\ell}_{(i)}<\infty$ by the assumption.
Thus the situation~\eqref{eq:fvanish(2)_1} cannot happen.
\end{proof}

\begin{Proposition}
Let $b\in B^{1,1}\otimes\bar{B}$ and suppose that we have shown the commutativity of $\widetilde{f}_i$ and~$\Phi$ for
the elements of $\bar{B}$.
Suppose that we have $\widetilde{f}_i(b)\neq 0$ and $\widetilde{f}_i(b')=0$ where $b'\in\bar{B}$ is the corresponding
part of~$b$.
Then we have $\widetilde{f}_i(\nu,J)\neq 0$, $\widetilde{f}_i(\bar{\nu},\bar{J})=0$ and
$\Phi^{-1}(\widetilde{f}_i(b))=\widetilde{f}_i(\Phi^{-1}(b))$.
\end{Proposition}
\begin{proof}
By the assumption, we have the following two possibilities:
\begin{enumerate}\itemsep=0pt
\item[(a)] $b=i\otimes b'$ and $\widetilde{f}_i(b)=(i+1)\otimes b'$, \item[(b)] $b=\overline{i+1}\otimes b'$ and
$\widetilde{f}_i(b)=\overline{i}\otimes b'$.
\end{enumerate}

{\bf Case (a).} In this case we have $\ell^{(i-1)}<\infty$ and $\ell^{(i)}=\infty$.
Then we have $P^{(i)}_\infty(\nu)=P^{(i)}_\infty(\bar{\nu})+1$.
By the induction hypothesis, $\widetilde{f}_i(b')=0$ implies that $\varphi_i(\bar{\nu},\bar{J})=0$.
Then from Theorem~\ref{prop:phi_RC} we have $P^{(i)}_\infty(\bar{\nu})=\bar{s}$ where $\bar{s}$ is the smallest rigging
of $(\bar{\nu},\bar{J})^{(i)}$.
Let $s$ be the smallest rigging of $(\nu,J)^{(i)}$.
Since $\ell^{(i)}=\infty$, we have $s=\bar{s}$.
Again by Theorem~\ref{prop:phi_RC} we have
$\varphi_i(\nu,J)=P^{(i)}_\infty(\nu)-s=P^{(i)}_\infty(\bar{\nu})+1-s=\bar{s}+1-s=1$.
Therefore $\widetilde{f}_i(\nu,J)\neq 0$ as desired.

In order to check $\Phi^{-1}(\widetilde{f}_i(b))=\widetilde{f}_i(\Phi^{-1}(b))$, it is enough to show that
$\widetilde{f}_i$ acting on $(\nu,J)$ creates a~singular string and it satisf\/ies
$\widetilde{\ell}^{(i-1)}\leq\ell+1=\widetilde{\ell}^{(i)}<\widetilde{\ell}^{(i+1)}=\infty$ under the assumptions
$\widetilde{f}_i(\nu,J)\neq 0$ and $\widetilde{f}_i(\bar{\nu},\bar{J})=0$.
Then this is the consequence of Proposition~\ref{lem:nonzero_zero1}.

{\bf Case (b).} In this case we have $\ell_{(i+1)}<\infty$, $\ell_{(i)}=\infty$, $\widetilde{\ell}_{(i)}<\infty$
and $\widetilde{\ell}_{(i-1)}=\infty$.
By $\ell_{(i+1)}<\infty$ and $\ell_{(i)}=\infty$ we have $P^{(i)}_\infty(\nu)=P^{(i)}_\infty(\bar{\nu})+1$.
Since we have $\varphi_i(\bar{\nu},\bar{J})=0$ by $\widetilde{f}_i(b')=0$, Theorem~\ref{prop:phi_RC} asserts that
$P^{(i)}_\infty(\bar{\nu})=\bar{s}$.
Now let us invoke the relation $s\leq\bar{s}$ of Lemma~\ref{lem:sbars}.
Then we have $\varphi_i(\nu,J)=P^{(i)}_\infty(\nu)-s=P^{(i)}_\infty(\bar{\nu})+1-s =\bar{s}-s+1\geq 1$.
Thus $\widetilde{f}_i(\nu,J)\neq 0$.

Finally we can check $\Phi^{-1}(\widetilde{f}_i(b))=\widetilde{f}_i(\Phi^{-1}(b))$ under the assumptions
$\widetilde{f}_i(\nu,J)\neq 0$ and $\widetilde{f}_i(\bar{\nu},\bar{J})=0$ by using Proposition~\ref{lem:nonzero_zero1}.
\end{proof}

\subsection{Proof for (3)}

For the analysis of the present situation, let us prepare several properties of the rigged conf\/igu\-ra\-tions.
Let $s$ (resp.
$\bar{s}$) be the smallest rigging of $(\nu,J)^{(i)}$ (resp.
$(\bar{\nu},\bar{J})^{(i)}$).

\begin{Lemma}
\label{lem:undef_def}
Suppose that we have $\widetilde{f}_i(\nu,J)=0$ and $\widetilde{f}_i(\bar{\nu},\bar{J})\neq 0$.
Then we have the following two possibilities:
\begin{enumerate}\itemsep=0pt
\item[$(1)$] $\ell^{(i)}<\infty$ and $\ell^{(i+1)}=\infty$, or \item[$(2)$] $\ell_{(i)}<\infty$ and $\ell_{(i-1)}=\infty$.
\end{enumerate}
\end{Lemma}
\begin{proof}
Let us show that $\ell^{(i)}<\infty$.
Suppose if possible that $\ell^{(i)}=\infty$.
From Theorem~\ref{prop:phi_RC} we have $\varphi_i(\nu,J)=P^{(i)}_\infty(\nu)-s=0$.
Since $\ell^{(i)}=\infty$, we have $s=\bar{s}$ and $P^{(i)}_\infty(\bar{\nu})\leq P^{(i)}_\infty(\nu)$.
Thus $\varphi_i(\bar{\nu},\bar{J})=P^{(i)}_\infty(\bar{\nu})-s\leq 0$ contradicting the assumption
$\widetilde{f}_i(\bar{\nu},\bar{J})\neq 0$.
Hence we conclude that $\ell^{(i)}<\infty$.

Suppose that $\ell^{(i+1)}<\infty$.
Then we have the following four possibilities:
\begin{enumerate}\itemsep=0pt
\item[(a)] $\ell_{(i+1)}=\infty$, \item[(b)] $\ell_{(i+1)}<\infty$ and $\ell_{(i)}=\infty$, \item[(c)]
$\ell_{(i)}<\infty$ and $\ell_{(i-1)}=\infty$, \item[(d)] $\ell_{(i-1)}<\infty$.
\end{enumerate}
Let us analyze these possibilities case by case.
Note that by the assumption $\widetilde{f}_i(\nu,J)=0$, we have $\varphi_i(\nu,J)=P^{(i)}_\infty(\nu)-s=0$ from
Theorem~\ref{prop:phi_RC}.

{\bf Case (a).} By the assumptions $\ell^{(i+1)}<\infty$ and $\ell_{(i+1)}=\infty$, we have
$P^{(i)}_\infty(\bar{\nu})=P^{(i)}_\infty(\nu)$.
Then the assumption $\widetilde{f}_i(\bar{\nu},\bar{J})\neq 0$ implies that
\begin{gather*}
\varphi_i(\bar{\nu},\bar{J})=P^{(i)}_\infty(\bar{\nu})-\bar{s} =P^{(i)}_\infty(\nu)-\bar{s}=s-\bar{s}>0.
\end{gather*}
By the assumption $\ell_{(i+1)}=\infty$, we see that $s>\bar{s}$ can happen only if the string
$\big(\ell^{(i)}-1,P^{(i)}_{\ell^{(i)}-1}(\bar{\nu})\big)$ which is created by~$\delta$ has the rigging strictly smaller than
$s$.
Thus the present assumptions imply that $P^{(i)}_{\ell^{(i)}-1}(\bar{\nu})<s$.
To summarize we are left with the following two possibilities:
\begin{enumerate}\itemsep=0pt
\item[(I)] $\ell^{(i-1)}=\ell^{(i)}$, then $P^{(i)}_{\ell^{(i)}-1}(\nu)<s$ by
$P^{(i)}_{\ell^{(i)}-1}(\bar{\nu})=P^{(i)}_{\ell^{(i)}-1}(\nu)$, \item[(II)] $\ell^{(i-1)}<\ell^{(i)}$, then
$P^{(i)}_{\ell^{(i)}-1}(\nu)\leq s$ by $P^{(i)}_{\ell^{(i)}-1}(\bar{\nu})=P^{(i)}_{\ell^{(i)}-1}(\nu)-1$.
\end{enumerate}
Let us show that these cases cannot happen.
Let $j$ be the largest integer such that $j<\ell^{(i)}$ and $m^{(i)}_j(\nu)>0$.
If there is no such $j$, set $j=0$.

{\bf Case (a-I).} Suppose that $j=\ell^{(i)}-1$.
Then we have $x_j\leq P^{(i)}_j(\nu)<s$, which contradicts the minimality of $s$.
Suppose that $0<j<\ell^{(i)}-1$.
Recall that by the minimality of $s$ we have $P^{(i)}_{\ell^{(i)}}(\nu)\geq s$.
Thus we have $P^{(i)}_{\ell^{(i)}}(\nu)>P^{(i)}_{\ell^{(i)}-1}(\nu)$.
Then by the convexity relation of $P^{(i)}_k(\nu)$ between $j\leq k\leq\ell^{(i)}$ we have
$P^{(i)}_{\ell^{(i)}}(\nu)>P^{(i)}_{\ell^{(i)}-1}(\nu) >P^{(i)}_j(\nu)\geq x_j$.
Since we are assuming that $s>P^{(i)}_{\ell^{(i)}-1}(\nu)$, this contradicts the minimality of $s$.

Finally suppose that $j\!=\!0$.
Then the previous inequalities become \mbox{$s\!>\!P^{(i)}_{\ell^{(i)}-1}(\nu)\!>\!P^{(i)}_0(\nu)\!=\!0$}.
Note that this relation is always valid since we have $\ell^{(i)}\geq 1$ by the def\/inition.
Therefore we have $s>0$ which implies that $\widetilde{f}_i$ acting on $(\nu,J)$ adds the length~1 string $(1,-1)$ to
$(\nu,J)^{(i)}$.
Since we are assuming that $\widetilde{f}_i(\nu,J)=0$, the rigging of the string $(1,-1)$ of
$(\widetilde{\nu},\widetilde{J})^{(i)}$ must be larger than the corresponding vacancy number.
Thus we have $P^{(i)}_1(\widetilde{\nu})<-1$, that is, $P^{(i)}_1(\nu)\leq 0$ by
$P^{(i)}_1(\widetilde{\nu})=P^{(i)}_1(\nu)-2$.
Since $P^{(i)}_0(\nu)=0$, the convexity relation of $P^{(i)}_k(\nu)$ between $0\leq k\leq\ell^{(i)}$ implies that
$0=P^{(i)}_0(\nu)\geq P^{(i)}_1(\nu)\geq P^{(i)}_{\ell^{(i)}}(\nu)$.
Therefore the rigging of the string $\big(\ell^{(i)},P^{(i)}_{\ell^{(i)}}(\nu)\big)$ of $(\nu,J)^{(i)}$ satisf\/ies that
$P^{(i)}_{\ell^{(i)}}(\nu)\leq 0$.
This is a~contradiction since we have $s>0$.

{\bf Case (a-II).} Suppose that $j=\ell^{(i)}-1$.
Then we have $x_j\leq P^{(i)}_j(\nu)\leq s$.
By the minimality of $s$ we have $x_j=P^{(i)}_j(\nu)=s$, in particular, the string $(j,x_j)$ of $(\nu,J)^{(i)}$ is
singular.
However this is in contradiction to the def\/inition of $\ell^{(i)}$ since we have $\ell^{(i-1)}\leq j<\ell^{(i)}$.
Suppose that $\ell^{(i-1)}\leq j<\ell^{(i)}-1$.
Then by the convexity relation of $P^{(i)}_k(\nu)$ between $j\leq k\leq\ell^{(i)}$ we have $s\geq
P^{(i)}_{\ell^{(i)}-1}(\nu)\geq P^{(i)}_{j}(\nu)\geq x_j$.
Then by the minimality of $s$ we have $P^{(i)}_j(\nu)=x_j$, in particular, the string $(j,x_j)$ of $(\nu,J)^{(i)}$ is
singular.
Again this is in contradiction to the def\/inition of $\ell^{(i)}$ since $\ell^{(i-1)}\leq j<\ell^{(i)}$.
Suppose that $0<j<\ell^{(i-1)}$.
Since $m^{(i-1)}_{\ell^{(i-1)}}(\nu)>0$, the convexity relation of $P^{(i)}_k(\nu)$ between $\ell^{(i-1)}-1\leq
k\leq\ell^{(i-1)}+1$ is strictly convex.
Therefore we have $s\geq P^{(i)}_{\ell^{(i)}-1}(\nu)>P^{(i)}_{j}(\nu)\geq x_j$.
This contradicts the minimality of $s$.

Finally suppose that $j=0$.
Then the above relation becomes $s\geq P^{(i)}_{\ell^{(i)}-1}(\nu)>P^{(i)}_0(\nu)=0$, that is, $s>0$.
By the similar arguments of the latter part of Case (a-I), we obtain $P^{(i)}_1(\nu)\leq 0$.
Since $P^{(i)}_0(\nu)=0$ and $m^{(i-1)}_{\ell^{(i-1)}}(\nu)>0$, the convexity relation of $P^{(i)}_k(\nu)$ between
$0\leq k\leq\ell^{(i)}$ becomes $0=P^{(i)}_0(\nu)\geq P^{(i)}_{\ell^{(i-1)}}(\nu)> P^{(i)}_{\ell^{(i)}}(\nu)$.
Hence the rigging of the string $\big(\ell^{(i)},P^{(i)}_{\ell^{(i)}}(\nu)\big)$ of $(\nu,J)^{(i)}$ satisf\/ies that
$P^{(i)}_{\ell^{(i)}}(\nu)<0$.
This is a~contradiction since we have $s>0$.

{\bf Case (b).} By the assumptions $\ell_{(i+1)}<\infty$ and $\ell_{(i)}=\infty$, we have
$P^{(i)}_\infty(\bar{\nu})=P^{(i)}_\infty(\nu)-1$.
Then the assumption $\widetilde{f}_i(\bar{\nu},\bar{J})\neq 0$ implies that
\begin{gather*}
\varphi_i(\bar{\nu},\bar{J})=P^{(i)}_\infty(\bar{\nu})-\bar{s} =P^{(i)}_\infty(\nu)-1-\bar{s}=s-\bar{s}-1>0.
\end{gather*}
By the assumption $\ell_{(i)}=\infty$, the situation $s-1>\bar{s}$ can happen only if the rigging of the string
$\big(\ell^{(i)}-1,P^{(i)}_{\ell^{(i)}-1}(\bar{\nu})\big)$ which is created by~$\delta$ satisf\/ies
$P^{(i)}_{\ell^{(i)}-1}(\bar{\nu})<s-1$.
Then we can use the same arguments of the above Case (a) to show that this case cannot happen.

{\bf Case (c).} By the assumptions $\ell_{(i)}<\infty$ and $\ell_{(i-1)}=\infty$, we have
$P^{(i)}_\infty(\bar{\nu})=P^{(i)}_\infty(\nu)+1$.
Then the assumption $\widetilde{f}_i(\bar{\nu},\bar{J})\neq 0$ implies that
\begin{gather*}
\varphi_i(\bar{\nu},\bar{J})=P^{(i)}_\infty(\bar{\nu})-\bar{s} =P^{(i)}_\infty(\nu)+1-\bar{s}=s-\bar{s}+1>0
\qquad
\Longleftrightarrow
\qquad
s\geq\bar{s}.
\end{gather*}
Thus, basically we have the following two cases:
\begin{gather*}
 (1) \quad s>\bar{s}, \qquad (2) \quad s=\bar{s}.
\end{gather*}

{\bf Case (c-1)} To begin with let us consider the case $s>\bar{s}$.
This situation can happen if at least one of the riggings for the strings
$\big(\ell^{(i)}-1,P^{(i)}_{\ell^{(i)}-1}(\bar{\nu})\big)$ or $\big(\ell_{(i)}-1,P^{(i)}_{\ell_{(i)}-1}(\bar{\nu})\big)$ is smaller than~$s$.
If $P^{(i)}_{\ell^{(i)}-1}(\bar{\nu})<s$, we can use the same arguments of Case~(a) to show that this case cannot
happen.
Thus $P^{(i)}_{\ell^{(i)}-1}(\bar{\nu})\geq s>\bar{s}$.

Therefore we suppose that $P^{(i)}_{\ell^{(i)}-1}(\bar{\nu})>\bar{s}$ and $P^{(i)}_{\ell_{(i)}-1}(\bar{\nu})=\bar{s}$ in
the sequel.
Let us follow the classif\/ication (III) to (VII) of Lemma~\ref{lem:vacancy_delta}.
We observe that the Cases (III) and (IV) cannot happen in this setting, since in these settings we have
$\ell^{(i)}=\ell_{(i)}$.
This is a~contradiction since we are assuming that
$P^{(i)}_{\ell^{(i)}-1}(\bar{\nu})>P^{(i)}_{\ell_{(i)}-1}(\bar{\nu})$.
Let us consider the remaining cases under the assumption $\ell^{(i)}<\ell_{(i)}$.

{\bf Case (c-1-V).} In this case we have $P^{(i)}_{\ell_{(i)}-1}(\nu)=\bar{s}-1$.
Let $j$ be the largest integer such that $j<\ell_{(i)}$ and $m^{(i)}_j(\nu)>0$.
Since $\ell^{(i)}<\ell_{(i)}$ we have $\ell^{(i)}\leq j$.
If $j=\ell_{(i)}-1$, we obtain the immediate contradiction $s>\bar{s}-1=P^{(i)}_{\ell_{(i)}-1}(\nu)\geq
x_{\ell_{(i)}-1}$.
Thus suppose that $\ell^{(i)}\leq j<\ell_{(i)}-1$.
Since $s$ is the minimal rigging of $(\nu,J)^{(i)}$ we see that $P^{(i)}_{\ell^{(i)}}(\nu)\geq s>\bar{s}$.
Therefore we have $P^{(i)}_{\ell^{(i)}}(\nu)>P^{(i)}_{\ell^{(i)}-1}(\nu)$.
Then by the convexity relation of $P^{(i)}_k(\nu)$ between $j\leq k\leq\ell^{(i)}$ we obtain
$P^{(i)}_{\ell^{(i)}-1}(\nu)>P^{(i)}_j(\nu)\geq x_j$, in particular, $s>x_j$.
This contradicts the minimality of $s$.

{\bf Case (c-1-VI).} In this case we have $P^{(i)}_{\ell_{(i)}-1}(\nu)=\bar{s}$.
Then we can use the same arguments of Case (c-1-V) to show that this case cannot happen.

{\bf Case (c-1-VII).} In this case we have $P^{(i)}_{\ell_{(i)}-1}(\nu)=\bar{s}+1$ under the assumption
$\ell_{(i+1)}<\ell_{(i)}$.
If $\bar{s}+1<s$ we can use the same arguments of Case (c-1-V) to show that such case cannot happen.
Thus suppose that $\bar{s}+1=s$.
Let $j$ be the largest integer such that $j<\ell_{(i)}$ and $m^{(i)}_j(\nu)>0$.
Since $\ell^{(i)}<\ell_{(i)}$ we have $\ell^{(i)}\leq j$.

Suppose that $j=\ell_{(i)}-1$.
Then we see that the string $(\ell_{(i)}-1,x_{\ell_{(i)}-1})$ is singular by the minimality of $s$, since we have
$s=P^{(i)}_{\ell_{(i)}-1}(\nu)\geq x_{\ell_{(i)}-1}$.
However this is in contradiction to the def\/inition of $\ell_{(i)}$ since its length satisf\/ies that
$\ell_{(i+1)}\leq\ell_{(i)}-1<\ell_{(i)}$.

Suppose that $\ell_{(i+1)}\leq j<\ell_{(i)}-1$.
Since $s$ is the minimal rigging of $(\nu,J)^{(i)}$, we have $P^{(i)}_{\ell_{(i)}}(\nu)\geq s$ and $P^{(i)}_j(\nu)\geq
s$.
Since we have $P^{(i)}_{\ell_{(i)}-1}(\nu)=s$, the only possibility that is compatible with the convexity relation of
$P^{(i)}_k(\nu)$ between $j\leq k\leq\ell_{(i)}$ is the case
$s=P^{(i)}_{\ell_{(i)}}(\nu)=P^{(i)}_{\ell_{(i)}-1}(\nu)=\dots =P^{(i)}_j(\nu)\geq x_j$.
By the minimality of $s$ we have $P^{(i)}_j(\nu)=x_j$, in particular, the string $(j,x_j)$ is singular.
However this is in contradiction to the def\/inition of $\ell_{(i)}$ since we have $\ell_{(i+1)}\leq j<\ell_{(i)}$.

Suppose that $\ell^{(i)}\leq j<\ell_{(i+1)}$.
As in the previous paragraph, we obtain $s=P^{(i)}_{\ell_{(i)}}(\nu)=P^{(i)}_{\ell_{(i)}-1}(\nu)=\dots
=P^{(i)}_j(\nu)\geq x_j$.
However, this is a~contradiction since we have $m^{(i+1)}_{\ell_{(i+1)}}(\nu)>0$, which implies that the vacancy numbers
$P^{(i)}_k(\nu)$ are strictly convex function of $k$.

To summarize, we have shown that Case (c-1) cannot happen.

{\bf Case (c-2).} In this case we have $s=\bar{s}$.
This case is indeed possible as the example at the end of this subsection shows.

{\bf Case (d).} By the assumption $\ell_{(i-1)}<\infty$, we have $P^{(i)}_\infty(\bar{\nu})=P^{(i)}_\infty(\nu)$.
Then the assumptions $\widetilde{f}_i(\bar{\nu},\bar{J})\neq 0$ and $\widetilde{f}_i(\nu,J)=0$ imply that
\begin{gather*}
\varphi_i(\bar{\nu},\bar{J})=P^{(i)}_\infty(\bar{\nu})-\bar{s} =P^{(i)}_\infty(\nu)-\bar{s}=s-\bar{s}>0
\qquad
\Longleftrightarrow
\qquad
s>\bar{s}.
\end{gather*}
We can use the same arguments of Case (c-1) to show that this case cannot happen.
\end{proof}

\begin{Lemma}
\label{cor:undef_def}
Suppose that we have $\widetilde{f}_i(\nu,J)=0$ and $\widetilde{f}_i(\bar{\nu},\bar{J})\neq 0$.
Then we have the following two possibilities:{\samepage
\begin{enumerate}\itemsep=0pt
\item[$(1)$] $\ell^{(i)}<\infty$, $\ell^{(i+1)}=\infty$ and $s=\bar{s}$,
\item[$(2)$] $\ell_{(i)}<\infty$, $\ell_{(i-1)}=\infty$ and $s=\bar{s}$.
\end{enumerate}}
\end{Lemma}

\begin{proof}
(1) {\bf Step 1.} Let us show $s\leq\bar{s}$.
It is enough to show that the new string of $(\bar{\nu},\bar{J})^{(i)}$ has the rigging larger than or equal to $s$,
that is, $P^{(i)}_{\ell^{(i)}-1}(\bar{\nu})\geq s$.
Suppose that $P^{(i)}_{\ell^{(i)}-1}(\bar{\nu})<s$.
Then we can use the same arguments of Case (a) of the proof of Lemma~\ref{lem:undef_def} to show that this case cannot
happen.
Thus $s\leq\bar{s}$.

{\bf Step 2.} Next, let us show $s\geq\bar{s}$.
Then it is enough to show $P^{(i)}_{\ell^{(i)}}(\nu)\geq\bar{s}$ since the only rigging that is changed by $\delta^{-1}$
is the one associated with the lengthened string (see Remark~\ref{rem:delta_inverse} for the explanation of
$\delta^{-1}$).
By the assumption $\ell^{(i+1)}=\infty$ we have $P^{(i)}_{\ell^{(i)}}(\nu)=P^{(i)}_{\ell^{(i)}}(\bar{\nu})-1$.
Thus the above relation is equivalent to $P^{(i)}_{\ell^{(i)}}(\bar{\nu})>\bar{s}$.

Suppose if possible that $P^{(i)}_{\ell^{(i)}}(\bar{\nu})\leq\bar{s}$.
Suppose if possible that there are strings of $(\bar{\nu},\bar{J})^{(i)}$ that are longer than or equal to $\ell^{(i)}$.
Let $j$ be the length of the shortest such string.
By the minimality of $\bar{s}$ we have $P^{(i)}_{\ell^{(i)}-1}(\bar{\nu})\geq\bar{s}$.
Then we have $P^{(i)}_{\ell^{(i)}-1}(\bar{\nu})\geq P^{(i)}_{\ell^{(i)}}(\bar{\nu})\geq P^{(i)}_{j}(\bar{\nu})\geq x_j$
by the convexity relation between $\ell^{(i)}-1$ and $j$ (if $j=\ell^{(i)}$ we can directly show the inequality without
the second term).
By the minimality of $\bar{s}$, we have $\bar{s}\leq x_j$.
Therefore we obtain $\bar{s}=x_j$ which implies that the string $(j,x_j)$ is singular with length larger than
$\ell^{(i)}-1$.
However, we know that $\delta^{-1}$ will add a~box to the length $\ell^{(i)}-1$ string of $(\bar{\nu},\bar{J})^{(i)}$.
This is a~contradiction since $\ell^{(i+1)}=\infty$ implies that~$\delta^{-1}$ will add a~box to the longest possible
string.
Therefore we conclude that the longest string of $(\bar{\nu},\bar{J})^{(i)}$ has length $\ell^{(i)}-1$.
Since we have $P^{(i)}_{\ell^{(i)}-1}(\bar{\nu})\geq P^{(i)}_{\ell^{(i)}}(\bar{\nu})$, the convexity relation between
$\ell^{(i)}-1$ and $\infty$ gives $P^{(i)}_{\ell^{(i)}-1}(\bar{\nu})\geq P^{(i)}_{\ell^{(i)}}(\bar{\nu}) \geq\dots\geq
P^{(i)}_{\infty}(\bar{\nu})$.
However this is a~contradiction since we already know that $P^{(i)}_{\ell^{(i)}}(\bar{\nu})\leq\bar{s}$ and
$P^{(i)}_\infty(\bar{\nu})>\bar{s}$ by $\varphi_i(\bar{\nu},\bar{J})=P^{(i)}_\infty(\bar{\nu})-\bar{s}>0$.
Thus we have shown that the relation $P^{(i)}_{\ell^{(i)}}(\bar{\nu})>\bar{s}$ is always satisf\/ied.
Hence we obtain $s\geq\bar{s}$.

By combining the two inequalities, we conclude that $s=\bar{s}$.

(2) During the proof of the previous lemma, we have shown that the only possible case under the assumption is
$s=\bar{s}$.
\end{proof}

\begin{Example}
Consider the following rigged conf\/iguration $(\nu,J)$ of type $(B^{1,1})^{\otimes 4}\otimes B^{1,3}\otimes
B^{2,1}\otimes B^{2,2}\otimes B^{3,1}$ of $D^{(1)}_5$
\begin{center}
\unitlength 12pt
\begin{picture}(36,6)
\multiput(-0.8,1.1)(0,1){5}{1}
\put(0,1){\Yboxdim12pt\yng(3,3,2,1,1)}
\put(1.2,1.1){0}
\put(1.2,2.1){0}
\put(2.2,3.1){1}
\put(3.2,4.1){1}
\put(3.2,5.1){1}
\put(6.0,0){
\put(-0.8,0.1){1}
\multiput(-0.8,1.1)(0,1){4}{0}
\put(-1.53,5.1){$-1$}
\put(0,0){\Yboxdim12pt\yng(5,3,3,2,2,1)}
\put(1.2,0.1){1}
\put(2.2,1.1){0}
\put(2.2,2.1){0}
\put(3.2,3.1){0}
\put(3.2,4.0){0}
\put(5.1,5.1){$-1$}
}
\put(14.8,0){
\put(-0.8,0.1){1}
\put(-0.8,1.1){0}
\multiput(-1.53,2.1)(0,1){4}{$-2$}
\put(0,0){\Yboxdim12pt\yng(5,4,3,3,2,1)}
\put(1.2,0.1){1}
\put(2.2,1.1){0}
\put(3.2,2.1){$-2$}
\put(3.2,3.1){$-2$}
\put(4.2,4.1){$-2$}
\put(5.2,5.1){$-3$}
}
\put(23.0,3){
\put(-0.8,0.1){0}
\put(-0.8,1.1){1}
\put(-0.8,2.1){2}
\put(0,0){\Yboxdim12pt\yng(5,2,1)}
\put(1.2,0.1){0}
\put(2.2,1.1){1}
\put(5.2,2.1){2}
}
\put(31,3){
\put(-0.8,0.1){0}
\multiput(-1.53,1.1)(0,1){2}{$-1$}
\put(0,0){\Yboxdim12pt\yng(4,4,1)}
\put(1.2,0.1){0}
\put(4.2,1.1){$-1$}
\put(4.2,2.1){$-1$}
}
\end{picture}
%\{\{3,3,2,1,1\},\{5,3,3,2,2,1\},\{5,4,3,3,2,1\},\{5,2,1\},\{4,4,1\}\},\\
%\{\{1,1,1,0,0\},\{-1,0,0,0,0,1\},\{-3,-2,-2,-2,0,1\},\{2,1,0\},\{-1,-1,0\}\}
\end{center}
This is the example for Case (c-2) of the proof of Lemma~\ref{lem:undef_def}.
The corresponding image $\Phi^{-1}(\nu,J)$~is
\begin{gather*}
\Yboxdim14pt \Yvcentermath1 \young(\mtwo)\otimes \young(\mone)\otimes \young(1)\otimes \young(\mfour)\otimes
\young(5\mthree\mone)\otimes \young(\mfive,5)\otimes \young(14,4\mfive)\otimes \young(1,\mfive,\mtwo)
\end{gather*}
We have $\widetilde{f}_2(\nu,J)=0$ since we have $P^{(2)}_6(\widetilde{\nu})=-3$.
$(\bar{\nu},\bar{J})$ is given as follows:
\begin{center}
\unitlength 12pt
\begin{picture}(36,6)
\multiput(-0.8,1.1)(0,1){3}{0}
\multiput(-0.8,4.1)(0,1){2}{1}
\put(0,1){\Yboxdim12pt\yng(3,3,1,1,1)}
\put(1.2,1.1){0}
\put(1.2,2.1){0}
\put(1.2,3.1){0}
\put(3.2,4.1){1}
\put(3.2,5.1){1}
\put(6,0){
\multiput(-0.8,0.1)(0,1){2}{1}
\multiput(-0.8,2.1)(0,1){3}{0}
\put(-1.53,5.1){$-1$}
\put(0,0){\Yboxdim12pt\yng(4,3,3,2,1,1)}
\put(1.2,0.1){1}
\put(1.2,1.1){1}
\put(2.2,2.1){0}
\put(3.2,3.1){0}
\put(3.2,4.0){0}
\put(4.1,5.0){$-1$}
}
\put(14.8,0){
\multiput(-0.8,0.1)(0,1){2}{1}
\multiput(-1.53,2.1)(0,1){4}{$-2$}
\put(0,0){\Yboxdim12pt\yng(5,3,3,3,1,1)}
\put(1.2,0.1){1}
\put(1.2,1.1){1}
\put(3.2,2.1){$-2$}
\put(3.2,3.1){$-2$}
\put(3.2,4.1){$-2$}
\put(5.2,5.1){$-3$}
}
\put(23.0,3){
\multiput(-0.8,0.1)(0,1){2}{0}
\put(-0.8,2.1){2}
\put(0,0){\Yboxdim12pt\yng(5,1,1)}
\put(1.2,0.1){0}
\put(1.2,1.1){0}
\put(5.2,2.1){2}
}
\put(31.0,3){
\multiput(-0.8,0.1)(0,1){2}{0}
\put(-1.53,2.1){$-1$}
\put(0,0){\Yboxdim12pt\yng(4,3,1)}
\put(1.2,0.1){0}
\put(3.2,1.1){0}
\put(4.2,2.1){$-1$}
}
\end{picture}
%\{\{3,3,1,1,1\},\{4,3,3,2,1,1\},\{5,3,3,3,1,1\},\{5,1,1\},\{4,3,1\}\},\\
%\{\{1,1,0,0,0\},\{-1,0,0,0,1,1\},\{-3,-2,-2,-2,1,1\},\{2,0,0\},\{-1,0,0\}\}
\end{center}
Note that the minimal rigging $-1$ of $(\nu,J)^{(2)}$ is preserved before and after~$\delta$ in a~nontrivial way.
We have $\widetilde{f}_2(\bar{\nu},\bar{J})\neq 0$ since we have $P^{(2)}_5(\widetilde{\bar{\nu}})=-2$.
\end{Example}

\begin{Proposition}%\label{prop:zero_nonzero}
Let us consider the rigged configuration $(\nu,J)$ of type $B^{1,1}\otimes\bar{B}$.
Suppose that we have the commutativity of $\widetilde{f}_i$ and~$\Phi$ for $\bar{B}$.
Suppose that $\widetilde{f}_i(\nu,J)=0$ and $\widetilde{f}_i(\bar{\nu},\bar{J})\neq 0$.
Let $b=\Phi(\nu,J)$ and $b'=\Phi(\bar{\nu},\bar{J})$.
Then we have $\widetilde{f}_i(b)=0$, $\widetilde{f}_i(b')\neq 0$ and
$\Phi(\widetilde{f}_i(\nu,J))=\widetilde{f}_i(\Phi(\nu,J))$.
\end{Proposition}
\begin{proof}
According to Lemma~\ref{lem:undef_def}, we see that there are only two possibilities.
The f\/irst case is $\ell^{(i)}<\infty$ and $\ell^{(i+1)}=\infty$.
In this case we have $b=(i+1)\otimes b'$ and $P^{(i)}_\infty(\bar{\nu})=P^{(i)}_\infty(\nu)+1$.
The second case is $\ell_{(i)}<\infty$ and $\ell_{(i-1)}=\infty$.
In this case we have $\overline{i}\otimes b'$ and $P^{(i)}_\infty(\bar{\nu})=P^{(i)}_\infty(\nu)+1$.

Recall that we have $\varphi_i(b')=\varphi_i(\bar{\nu},\bar{J})>0$ by the induction hypothesis.
Let us show that in fact we have $\varphi_i(\bar{\nu},\bar{J})=1$.
From Theorem~\ref{prop:phi_RC}, we compute
\begin{gather*}
\varphi_i(\bar{\nu},\bar{J})=P^{(i)}_\infty(\bar{\nu})-\bar{s}
= P^{(i)}_\infty(\nu)-s+1=1.
\end{gather*}
Here we have used the relation $s\!=\!\bar{s}$ from Lemma~\ref{cor:undef_def} and the assumption
$\varphi(\nu,J)\!=\!P^{(i)}_\infty(\nu)-s\!=\!0$.
Thus $\varphi_i(\bar{\nu},\bar{J})=1$.

Then by the induction hypothesis we have $\varphi_i(b')=\varphi_i(\bar{\nu},\bar{J})=1$.
On the other hand, recall that we have $\varepsilon_i(i+1)=\varepsilon_i(\minusi)=1$.
Then from the tensor product rule at Section~\ref{se:crystal}, we conclude that $\widetilde{f}_i(b)=0$.
Hence we have $\Phi(\widetilde{f}_i(\nu,J))=\widetilde{f}_i(\Phi(\nu,J))$.
\end{proof}

\subsection{Proof for (4)}

\begin{Proposition}
Let $b\in B^{1,1}\otimes\bar{B}$ and suppose that we have shown the commutativity of $\widetilde{f}_i$ and~$\Phi$ for
the elements of $\bar{B}$.
Suppose that we have $\widetilde{f}_i(b)=0$ and $\widetilde{f}_i(b')\neq 0$ where $b'\in\bar{B}$ is the corresponding
part of~$b$.
Then we have $\widetilde{f}_i(\nu,J)=0$, $\widetilde{f}_i(\bar{\nu},\bar{J})\neq 0$ and
$\Phi^{-1}(\widetilde{f}_i(b))=\widetilde{f}_i(\Phi^{-1}(b))$.
\end{Proposition}

\begin{proof}
By the assumption, we have $b=(i+1)\otimes b'$ or $b=\overline{i}\otimes b'$ and $\varphi_i(b')=1$.
In the f\/irst case we have $\ell^{(i)}<\infty$ and $\ell^{(i+1)}=\infty$ and in the second case we have
$\ell_{(i)}<\infty$ and $\ell_{(i-1)}=\infty$.
In both cases we have $P^{(i)}_\infty(\nu)=P^{(i)}_\infty(\bar{\nu})-1$.
From the condition $\varphi_i(b')=1$ together with the induction hypothesis we have
$\varphi_i(b')=\varphi_i(\bar{\nu},\bar{J})=1$, which implies $P^{(i)}_\infty(\bar{\nu})=\bar{s}+1$ by
Theorem~\ref{prop:phi_RC}.
Here $\bar{s}$ is the minimal rigging of $(\bar{\nu},\bar{J})^{(i)}$.
Then we have
\begin{gather*}
\varphi_i(\nu,J)=P^{(i)}_\infty(\nu)-s =P^{(i)}_\infty(\bar{\nu})-1-s=\bar{s}-s.
\end{gather*}
If we show $s\geq\bar{s}$ we have $\varphi_i(\nu,J)=0$ as desired.
Note that this result leads to $\Phi^{-1}(\widetilde{f}_i(b))=\widetilde{f}_i(\Phi^{-1}(b))$ also.

Let us show $s\geq\bar{s}$.
If $\ell^{(i)}<\infty$ and $\ell^{(i+1)}=\infty$ we can use the same arguments of Step~2 of the proof of
Lemma~\ref{cor:undef_def}(1) since we do not use the assumption $\varphi_i(\nu,J)=0$ there.
If $\ell_{(i)}<\infty$ and $\ell_{(i-1)}=\infty$, again we can use the same arguments of Step~2 of the proof of
Lemma~\ref{cor:undef_def}(1) by replacing $\ell^{(i)}$ and $\ell^{(i+1)}$ there by $\ell_{(i)}$ and $\ell_{(i-1)}$
respectively.
\end{proof}

\section[Proof of Proposition~\ref{th:core2}: $\widetilde{e}_i$ version]{Proof
of Proposition~\ref{th:core2}: $\boldsymbol{\widetilde{e}_i}$ version}\label{sec:main2_e}

In this section we use a~parallel notation that are introduced in Section~\ref{subsec:notation} except for changing the
role of $\widetilde{f}_i$ there by $\widetilde{e}_i$.
For example, $\widetilde{e}_i$ acts on the string $(\ell,x_\ell)$ of $(\nu,J)^{(i)}$.

\subsection{Proof for (5)}

\begin{Proposition}\label{prop:nonzero(5)_1}
Let us consider the rigged configuration $(\nu,J)$ of type $B^{1,1}\otimes\bar{B}$.
Suppose that we have the commutativity of $\widetilde{e}_i$ and $\widetilde{f}_i$ with~$\Phi$ for $\bar{B}$.
Suppose that $\widetilde{e}_i(\nu,J)\neq 0$ and $\widetilde{e}_i(\bar{\nu},\bar{J})=0$.
Let $b=\Phi(\nu,J)$ and $b'=\Phi(\bar{\nu},\bar{J})$.
Then we have one of the following two cases:
\begin{enumerate}\itemsep=0pt
\item[$(1)$] $b=(i+1)\otimes b'$, $\widetilde{e}_i(b)$ is defined, $\widetilde{e}_i(b')$ is undefined and
$\Phi(\widetilde{e}_i(\nu,J))=i\otimes b'$,
\item[$(2)$] $b=\overline{i}\otimes b'$, $\widetilde{e}_i(b)$ is
defined, $\widetilde{e}_i(b')$ is undefined and $\Phi(\widetilde{e}_i(\nu,J))=\overline{i+1}\otimes b'$.
\end{enumerate}
\end{Proposition}

For the proof, we show the following properties.

\begin{Proposition}
\label{prop:nonzero(5)_2}
$\widetilde{e}_i(\nu,J)$ is defined and $\widetilde{e}_i(\bar{\nu},\bar{J})$ is undefined if and only if one of the
following conditions is satisfied:
\begin{enumerate}\itemsep=0pt
\item[$(1)$] $\ell^{(i)}<\infty$, $\ell^{(i+1)}=\infty$ and all the riggings of $(\nu,J)^{(i)}$ are non-negative except
for the longest string whose rigging is $-1$.
Moreover we have
\begin{gather*}
P^{(i)}_{\ell^{(i)}-1}(\bar{\nu})=P^{(i)}_{\ell^{(i)}}(\bar{\nu})= \dots=P^{(i)}_\infty(\bar{\nu})=0.
\end{gather*}
\item[$(2)$] $\ell_{(i)}<\infty$, $\ell_{(i-1)}=\infty$ and all the riggings of $(\nu,J)^{(i)}$ are non-negative except
for the longest string$($s$)$ whose rigging is~$-1$.
Moreover we have
\begin{gather*}
P^{(i)}_{\ell_{(i)}-1}(\bar{\nu})=P^{(i)}_{\ell_{(i)}}(\bar{\nu})= \dots=P^{(i)}_\infty(\bar{\nu})=0.
\end{gather*}
\end{enumerate}
\end{Proposition}
\begin{proof}
By the assumptions $\widetilde{e}_i(\nu,J)\neq 0$ and $\widetilde{e}_i(\bar{\nu},\bar{J})=0$, we see that~$\delta$
changes $(\nu,J)^{(i)}$.
Thus we have $\ell^{(i)}<\infty$.
Also, from $\widetilde{e}_i(\bar{\nu},\bar{J})=0$, we see that all the riggings of $(\bar{\nu},\bar{J})^{(i)}$ are
non-negative.

(1) Suppose that we have $\ell^{(i+1)}=\infty$.
Suppose if possible that we have $m^{(i)}_j(\bar{\nu})>0$ for some $\ell^{(i)}\leq j$.
Let us choose the minimal such $j$.
Note that by the condition $\ell^{(i)}\leq j$, the string $(j,x_j)$ is dif\/ferent from the string
$\big(\ell^{(i)}-1,P^{(i)}_{\ell^{(i)}-1}(\bar{\nu})\big)$ created by~$\delta$.
By the def\/inition of $\delta^{-1}$, $\ell^{(i)}<\infty$ and $\ell^{(i+1)}=\infty$ imply that the string
$\big(\ell^{(i)}-1,P^{(i)}_{\ell^{(i)}-1}(\bar{\nu})\big)$ is the longest singular string of $(\bar{\nu},\bar{J})^{(i)}$.
In particular, the string $(j,x_j)$ is not singular so that we have $P^{(i)}_j(\bar{\nu})>x_j\geq 0$.
On the other hand, we have $P^{(i)}_{\ell^{(i)}-1}(\bar{\nu})\geq 0$ since $\widetilde{e}_i(\bar{\nu},\bar{J})=0$.
Let us show $P^{(i)}_{\ell^{(i)}}(\bar{\nu})>0$.
If $\ell^{(i)}=j$ it is already proved.
Thus suppose that $\ell^{(i)}<j$.
Since we have $m^{(i)}_k(\bar{\nu})=0$ for all $\ell^{(i)}<k<j$, the only possibility that is compatible with
$P^{(i)}_{\ell^{(i)}-1}(\bar{\nu})\geq 0$, $P^{(i)}_j(\bar{\nu})>0$ and the convexity relation of $P^{(i)}_k(\bar{\nu})$
for $\ell^{(i)}-1\leq k\leq j$ is $P^{(i)}_{\ell^{(i)}}(\bar{\nu})>0$.
In conclusion, we have $P^{(i)}_{\ell^{(i)}}(\nu)\geq 0$ by
$P^{(i)}_{\ell^{(i)}}(\nu)=P^{(i)}_{\ell^{(i)}}(\bar{\nu})-1$ which follows from the relations $\ell^{(i)}<\infty$ and
$\ell^{(i+1)}=\infty$.
Since all the riggings of $(\bar{\nu},\bar{J})^{(i)}$ are non-negative, the relation $P^{(i)}_{\ell^{(i)}}(\nu)\geq 0$
implies that all the riggings of $(\nu,J)^{(i)}$ are non-negative.
This contradicts the assumption $\widetilde{e}_i(\nu,J)\neq 0$.
Hence the string $\big(\ell^{(i)}-1,P^{(i)}_{\ell^{(i)}-1}(\bar{\nu})\big)$ is the longest string of
$(\bar{\nu},\bar{J})^{(i)}$.

Let us show that $P^{(i)}_{\ell^{(i)}}(\bar{\nu})=0$.
Suppose if possible that we have $P^{(i)}_{\ell^{(i)}}(\bar{\nu})>0$.
Then we have $P^{(i)}_{\ell^{(i)}}(\nu)\geq 0$ by $P^{(i)}_{\ell^{(i)}}(\nu)=P^{(i)}_{\ell^{(i)}}(\bar{\nu})-1$.
This contradicts the assumption $\widetilde{e}_i(\nu,J)\neq 0$.
Therefore we have $P^{(i)}_{\ell^{(i)}}(\bar{\nu})=0$.
In particular, the string $\big(\ell^{(i)},P^{(i)}_{\ell^{(i)}}(\nu)\big)=(\ell^{(i)},-1)$ is the only string of $(\nu,J)^{(i)}$
which has the negative rigging.

Finally let us show that $P^{(i)}_{\ell^{(i)}-1}(\bar{\nu})=\dots=P^{(i)}_\infty(\bar{\nu})=0$.
By the convexity relation of $P^{(i)}_k(\bar{\nu})$ between $\ell^{(i)}-1\leq k<\infty$ and the condition
$P^{(i)}_\infty(\bar{\nu})>-\infty$, we have $P^{(i)}_{\ell^{(i)}-1}(\bar{\nu})\leq P^{(i)}_{\ell^{(i)}}(\bar{\nu}) \leq
\dots \leq P^{(i)}_{\infty}(\bar{\nu})$.
Recall that we have $P^{(i)}_{\ell^{(i)}-1}(\bar{\nu})\geq 0$ by $\widetilde{e}_i(\bar{\nu},\bar{J})=0$ and
$P^{(i)}_{\ell^{(i)}}(\bar{\nu})=0$ by the previous paragraph.
Then the only possibility is $P^{(i)}_{\ell^{(i)}-1}(\bar{\nu})=P^{(i)}_{\ell^{(i)}}(\bar{\nu}) =\dots
=P^{(i)}_{\infty}(\bar{\nu})=0$.

(2) In this case we assume that $\ell^{(i+1)}<\infty$.

{\bf Step 1.} We shall show that $\ell_{(i+1)}<\infty$.
Suppose if possible that $\ell_{(i+1)}=\infty$.
Consider the case $\ell^{(i)}=\ell^{(i+1)}$.
Then we have $P^{(i)}_{\ell^{(i)}}(\nu)=P^{(i)}_{\ell^{(i)}}(\bar{\nu})$.
This implies that all the riggings of $(\nu,J)^{(i)}$ are non-negative.
This contradicts the assumption $\widetilde{e}_i(\nu,J)\neq 0$.

Thus we are left with the case $\ell^{(i)}<\ell^{(i+1)}$.
Then we have $P^{(i)}_{\ell^{(i)}}(\nu)=P^{(i)}_{\ell^{(i)}}(\bar{\nu})-1$.
In order to have $\widetilde{e}_i(\nu,J)\neq 0$, we must have $P^{(i)}_{\ell^{(i)}}(\bar{\nu})=0$.
Suppose if possible that there are strings of $(\bar{\nu},\bar{J})^{(i)}$ which are longer than $\ell^{(i)}-1$.
Let $j$ be the minimal integer such that $m^{(i)}_j(\bar{\nu})>0$ and $\ell^{(i)}\leq j$.
Then the string $(j,x_j)$ of $(\bar{\nu},\bar{J})^{(i)}$ must satisfy $P^{(i)}_j(\bar{\nu})\geq x_j\geq 0$.
Similarly, the string $\big(\ell^{(i)}-1,P^{(i)}_{\ell^{(i)}-1}(\bar{\nu})\big)$ must satisfy
$P^{(i)}_{\ell^{(i)}-1}(\bar{\nu})\geq 0$.
Then the only possibility that is compatible with the convexity relation of $P^{(i)}_k(\bar{\nu})$ between
$\ell^{(i)}-1\leq k\leq j$ is the relation
\begin{gather}
\label{eq:proof_mainth2_(5)_1}
P^{(i)}_{\ell^{(i)}-1}(\bar{\nu})=P^{(i)}_{\ell^{(i)}}(\bar{\nu})=\dots =P^{(i)}_j(\bar{\nu})=0.
\end{gather}
From~\eqref{eq:proof_mainth2_(5)_1} we have $m^{(i+1)}_k(\bar{\nu})=0$ for all $\ell^{(i)}\leq k<j$.
Recall that by def\/inition of $\ell^{(i+1)}$ there is the string $(\ell^{(i+1)}-1,P^{(i+1)}_{\ell^{(i+1)}-1}(\bar{\nu}))$
of $(\bar{\nu},\bar{J})^{(i+1)}$.
Then we see that its length satisf\/ies $j\leq\ell^{(i+1)}-1$ by $\ell^{(i)}\leq\ell^{(i+1)}-1$.
Also from~\eqref{eq:proof_mainth2_(5)_1} and the assumption $\widetilde{e}_i(\bar{\nu},\bar{J})=0$, we see that the
string $(j,x_j)$ is singular since we have $0=P^{(i)}_j(\bar{\nu})\geq x_j\geq 0$.
By the relation $j\leq\ell^{(i+1)}-1$, $\delta^{-1}$ will add a~box to the string $(j,x_j)$ of
$(\bar{\nu},\bar{J})^{(i)}$ instead of the string $\big(\ell^{(i)}-1,P^{(i)}_{\ell^{(i)}-1}(\bar{\nu})\big)$.
This is a~contradiction since we have $\ell^{(i)}-1<j$.
In conclusion we see that the string $\big(\ell^{(i)}-1,P^{(i)}_{\ell^{(i)}-1}(\bar{\nu})\big)$ is the longest string of
$(\bar{\nu},\bar{J})^{(i)}$.

We know that $P^{(i)}_{\ell^{(i)}}(\bar{\nu})=0$.
Since $\widetilde{e}_i(\bar{\nu},\bar{J})=0$, we have $P^{(i)}_{\ell^{(i)}-1}(\bar{\nu})\geq 0$ and, in particular,
$P^{(i)}_{\ell^{(i)}-1}(\bar{\nu})\geq P^{(i)}_{\ell^{(i)}}(\bar{\nu})$.
Since there is no string of $(\bar{\nu},\bar{J})^{(i)}$ which is longer than $\ell^{(i)}-1$, we can use the convexity
relation of $P^{(i)}_k(\bar{\nu})$ between $\ell^{(i)}-1\leq k<\infty$ with $P^{(i)}_\infty(\bar{\nu})>-\infty$ to
conclude that $P^{(i)}_{\ell^{(i)}-1}(\bar{\nu})\leq P^{(i)}_{\ell^{(i)}}(\bar{\nu})\leq\dots\leq
P^{(i)}_{\infty}(\bar{\nu})$.
By the relation $P^{(i)}_{\ell^{(i)}-1}(\bar{\nu})\geq P^{(i)}_{\ell^{(i)}}(\bar{\nu})$ we conclude that
$P^{(i)}_{\ell^{(i)}-1}(\bar{\nu})=P^{(i)}_{\ell^{(i)}}(\bar{\nu})=\dots =P^{(i)}_{\infty}(\bar{\nu})$ is the only
possibility that is compatible with the convexity relation.
Then we have $m^{(i+1)}_k(\bar{\nu})=0$ for all $\ell^{(i)}\leq k$.
This is a~contradiction since we are assuming that $\ell^{(i)}<\ell^{(i+1)}$.
Hence we have $\ell_{(i+1)}<\infty$.

{\bf Step 2.} We shall show that $\ell_{(i)}<\infty$.
Suppose if possible that $\ell_{(i)}=\infty$.
In order to have $\widetilde{e}_i(\nu,J)\neq 0$, we must have
$P^{(i)}_{\ell^{(i)}}(\nu)=P^{(i)}_{\ell^{(i)}}(\bar{\nu})-1$.
Under the condition $\ell_{(i)}=\infty$, this is only possible if we have $\ell^{(i)}<\ell^{(i+1)}$.
Then we can use the same argument of the corresponding part of Step 1 to deduce a~contradiction.
Thus we have $\ell_{(i)}<\infty$.

{\bf Step 3.} We shall show that $\ell_{(i-1)}=\infty$.
Suppose if possible that $\ell_{(i-1)}<\infty$.
Then we can use the same arguments of Step 1 if we replace $\ell^{(i)}$ and $\ell^{(i+1)}$ there by $\ell_{(i)}$ and
$\ell_{(i-1)}$ respectively to deduce a~contradiction.
Thus we have $\ell_{(i-1)}=\infty$.

{\bf Step 4.} We shall show that the length $\ell^{(i)}$ string of $(\nu,J)^{(i)}$ created by $\delta^{-1}$ has
non-negative rigging if $\ell^{(i)}<\ell_{(i)}$.
For this we have to show that $P^{(i)}_{\ell^{(i)}}(\nu)\geq 0$.
Recall that we have $P^{(i)}_{\ell^{(i)}-1}(\bar{\nu})\geq 0$ by the assumption $\widetilde{e}_i(\bar{\nu},\bar{J})=0$.
On the other hand, let $j$ be the smallest integer such that $m^{(i)}_j(\bar{\nu})>0$ and $\ell^{(i)}\leq j$.
By $\ell^{(i)}<\ell_{(i)}$ we have $j\leq\ell_{(i)}$.
Then by the assumption $\widetilde{e}_i(\bar{\nu},\bar{J})=0$ we have $P^{(i)}_{j}(\bar{\nu})\geq 0$.

If $\ell^{(i)}=\ell^{(i+1)}$, we have $P^{(i)}_{\ell^{(i)}}(\nu)\geq P^{(i)}_{\ell^{(i)}}(\bar{\nu})$ by
$\ell^{(i)}<\ell_{(i)}$.
Then by the convexity relation of~$P^{(i)}_k(\bar{\nu})$ between $\ell^{(i)}-1\leq k\leq j$ we have
$P^{(i)}_{\ell^{(i)}}(\nu)\geq P^{(i)}_{\ell^{(i)}}(\bar{\nu})\geq
\min\big\{P^{(i)}_{\ell^{(i)}-1}(\bar{\nu}),P^{(i)}_{j}(\bar{\nu})\big\}\geq 0$.
Next let us consider the case $\ell^{(i)}<\ell^{(i+1)}$.
Then we have $P^{(i)}_{\ell^{(i)}}(\nu)=P^{(i)}_{\ell^{(i)}}(\bar{\nu})-1$.
Then one can attain $P^{(i)}_{\ell^{(i)}}(\nu)=-1$ only if
$P^{(i)}_{\ell^{(i)}-1}(\bar{\nu})=\dots=P^{(i)}_{j}(\bar{\nu})=0$.
This relation implies that $m^{(i+1)}_k(\bar{\nu})=0$ for all $\ell^{(i)}-1<k<j$.
Then we have $j\leq\ell^{(i+1)}-1$ since we have $\ell^{(i)}-1<\ell^{(i+1)}-1$ by the assumption.
Let us consider the string $(j,x_j)$ of $(\bar{\nu},\bar{J})^{(i)}$.
Then we have $P^{(i)}_{j}(\bar{\nu})\geq x_j\geq 0$ by the assumption $\widetilde{e}_i(\bar{\nu},\bar{J})=0$.
Since $P^{(i)}_{j}(\bar{\nu})=0$ we see that the string $(j,x_j)$ is singular and its length satisf\/ies
$\ell^{(i)}-1<j\leq\ell^{(i+1)}-1$.
Thus $\delta^{-1}$ will add a~box to the string $(j,x_j)$ instead of the length $\ell^{(i)}-1$ string.
This is a~contradiction.
To summarize, we have shown that the length~$\ell^{(i)}$ string of $(\nu,J)^{(i)}$ created by $\delta^{-1}$ has
non-negative rigging if $\ell^{(i)}<\ell_{(i)}$.

Then we can use a~parallel arguments of the previous case (1) to deduce that the string
$\big(\ell_{(i)}-1,P^{(i)}_{\ell_{(i)}-1}(\bar{\nu})\big)$ is the longest string of $(\bar{\nu},\bar{J})^{(i)}$.

{\bf Step 5.} Finally let us show that $P^{(i)}_{\ell_{(i)}-1}(\bar{\nu})=\dots=P^{(i)}_\infty(\bar{\nu})=0$.
Since the length $\ell_{(i)}-1$ string is the longest string of $(\bar{\nu},\bar{J})^{(i)}$, we have
$P^{(i)}_{\ell_{(i)}-1}(\bar{\nu})\leq P^{(i)}_{\ell_{(i)}}(\bar{\nu})\leq\dots\leq P^{(i)}_\infty(\bar{\nu})$ by the
convexity relation of $P^{(i)}_k(\bar{\nu})$ between $\ell_{(i)}-1\leq k<\infty$ and
$P^{(i)}_\infty(\bar{\nu})>-\infty$.
Since we have $P^{(i)}_{\ell_{(i)}}(\nu)=-1$, we have $P^{(i)}_{\ell_{(i)}}(\bar{\nu})=0$ by
$P^{(i)}_{\ell_{(i)}}(\nu)=P^{(i)}_{\ell_{(i)}}(\bar{\nu})-1$.
Since we have $P^{(i)}_{\ell_{(i)}-1}(\bar{\nu})\geq 0$ by $\widetilde{e}_i(\bar{\nu},\bar{J})=0$, we conclude that
$P^{(i)}_{\ell_{(i)}-1}(\bar{\nu})=P^{(i)}_{\ell_{(i)}}(\bar{\nu})=\dots= P^{(i)}_\infty(\bar{\nu})=0$.
\end{proof}

\begin{proof}
[Proof of Proposition~\ref{prop:nonzero(5)_1}] We divide the proof into two cases following
Proposition~\ref{prop:nonzero(5)_2}.

(1) Let us consider the case $\ell^{(i)}<\infty$ and $\ell^{(i+1)}=\infty$.
In this case we have $b=(i+1)\otimes b'$.
By the assumption we have $\varepsilon_i(b')=\varepsilon_i(\bar{\nu},\bar{J})=0$.
From Proposition~\ref{prop:nonzero(5)_2} we know that $P^{(i)}_\infty(\bar{\nu})=0$ and the minimal rigging of
$(\bar{\nu},\bar{J})^{(i)}$ is 0.
Therefore we have $\varphi_i(\bar{\nu},\bar{J})=0$ by Theorem~\ref{prop:phi_RC}.
Then by the assumption we have $\varphi_i(b')=\varphi_i(\bar{\nu},\bar{J})=0$.
Thus $\varepsilon_i(b)=1$, that is, $\widetilde{e}_i(b)\neq 0$.

Finally let us show $\Phi(\widetilde{e}_i(\nu,J))=i\otimes b'$.
Recall that the longest string of $(\bar{\nu},\bar{J})^{(i)}$ has length $\ell^{(i)}-1$.
Therefore we have $m^{(i)}_{\ell^{(i)}}(\nu)=1$.
Thus if $\ell^{(i-1)}=\ell^{(i)}$ we see that we cannot remove a~box from $(\widetilde{e}_i(\nu,J))^{(i)}$ by~$\delta$
since the length of the longest string of $(\widetilde{e}_i(\nu,J))^{(i)}$ is $\ell^{(i)}-1$ and thus all the strings of
$(\widetilde{e}_i(\nu,J))^{(i)}$ are strictly shorter than $\ell^{(i-1)}$.
In this case we have $\Phi(\widetilde{e}_i(\nu,J))=i\otimes b'$.
Therefore let us suppose that $\ell^{(i-1)}<\ell^{(i)}$ in the following.
Recall that we have $P^{(i)}_{\ell^{(i)}-1}(\bar{\nu})=0$.
Then we have $P^{(i)}_{\ell^{(i)}-1}(\nu)=1$ by $\ell^{(i-1)}<\ell^{(i)}$.
We know that $\widetilde{e}_i$ changes the string $(\ell^{(i)},-1)$ of $(\nu,J)^{(i)}$ into $(\ell^{(i)}-1,0)$.
In particular, the latter string is non-singular by
$P^{(i)}_{\ell^{(i)}-1}(\widetilde{\nu})=P^{(i)}_{\ell^{(i)}-1}(\nu)=1$.
Thus $\widetilde{\ell}^{(i)}=\infty$ also.
In conclusion we have $\Phi(\widetilde{e}_i(\nu,J))=i\otimes b'$.

(2) The proof for this case is almost identical to the previous case if we replace $\ell^{(i-1)}$ and $\ell^{(i)}$ by
$\ell_{(i)}$ and $\ell_{(i+1)}$ respectively.
\end{proof}

\begin{Example}
Consider the following rigged conf\/iguration $(\nu,J)$ of type $(B^{1,1})^{\otimes 4}\otimes B^{1,3}\otimes
B^{2,1}\otimes B^{2,2}\otimes B^{3,1}$ of $D^{(1)}_5$
\begin{center}
\unitlength 12pt
\begin{picture}(35,4)(0,2)
\put(-0.8,3.1){1}
\put(-0.8,4.1){0}
\put(-0.8,5.1){0}
\put(0,3){\Yboxdim12pt\yng(4,3,2)}
\put(2.2,3.1){1}
\put(3.2,4.1){$-2$}
\put(4.2,5.1){0}
\put(6.7,0){
\put(-0.8,2.1){1}
\put(-0.8,3.1){2}
\put(-0.8,4.1){2}
\put(-1.5,5.1){$-1$}
\put(0,2){\Yboxdim12pt\yng(6,4,2,1)}
\put(1.2,2.1){0}
\put(2.2,3.1){2}
\put(4.2,4.0){2}
\put(6.1,5.1){$-1$}
}
\put(15.3,0){
\put(-0.8,2.1){1}
\put(-0.8,3.1){0}
\put(-0.8,4.1){0}
\put(-0.8,5.1){0}
\put(0,2){\Yboxdim12pt\yng(5,4,3,1)}
\put(1.2,2.1){0}
\put(3.2,3.1){$-1$}
\put(4.2,4.1){0}
\put(5.2,5.1){0}
}
\put(23.0,3){
\put(-1.53,1.1){$-1$}
\put(-1.53,2.1){$-1$}
\put(0,1){\Yboxdim12pt\yng(5,2)}
\put(2.2,1.1){$-1$}
\put(5.2,2.1){$-1$}
}
\put(31.0,3){
\put(-1.53,1.1){$-1$}
\put(-0.8,2.1){0}
\put(0,1){\Yboxdim12pt\yng(4,2)}
\put(2.2,1.1){$-1$}
\put(4.2,2.1){0}
}
\end{picture}
%\{\{4,3,2\},\{6,4,2,1\},\{5,4,3,1\},\{5,2\},\{4,2\}\},\\
%\{\{0,-2,1\},\{-1,2,2,0\},\{0,0,-1,0\},\{-1,-1\},\{0,-1\}\}
\end{center}
The corresponding image $\Phi^{-1}(\nu,J)$ is
\begin{gather*}
\Yboxdim14pt \Yvcentermath1 \young(\mtwo)\otimes \young(\mone)\otimes \young(5)\otimes \young(\mtwo)\otimes
\young(11\mtwo)\otimes \young(2,3)\otimes \young(22,34)\otimes \young(1,\mfour,\mthree)
\end{gather*}
We have $\widetilde{e}_2(\nu,J)\neq 0$ since $(\nu,J)^{(2)}$ has a~negative rigging.
Here $\widetilde{e}_2$ changes the string $(6,-1)$ of $(\nu,J)^{(2)}$ into $(5,0)$.
Since $P^{(i)}_5(\widetilde{\nu})=1$, the latter string is non-singular.
Thus the leftmost letter of $\Phi(\widetilde{\nu},\widetilde{J})$ is $\mthree$ and the rest part is the same as the
corresponding part of the above path.
$(\bar{\nu},\bar{J})$ is given as follows:
\begin{center}
\unitlength 12pt
\begin{picture}(34,4)(0,2)
\put(-0.8,3.1){2}
\put(-0.8,4.1){0}
\put(-0.8,5.1){0}
\put(0,3){\Yboxdim12pt\yng(4,3,1)}
\put(1.2,3.1){2}
\put(3.2,4.1){$-2$}
\put(4.2,5.1){0}
\put(7.0,0){
\put(-0.8,2.1){1}
\put(-0.8,3.1){1}
\put(-0.8,4.1){2}
\put(-0.8,5.1){0}
\put(0,2){\Yboxdim12pt\yng(5,4,1,1)}
\put(1.2,2.1){0}
\put(1.2,3.1){1}
\put(4.2,4.0){2}
\put(5.1,5.1){0}
}
\put(15.3,0){
\put(-0.8,2.1){1}
\put(-1.53,3.1){$-1$}
\put(-1.53,4.1){$-1$}
\put(-0.8,5.1){0}
\put(0,2){\Yboxdim12pt\yng(4,3,3,1)}
\put(1.2,2.1){0}
\put(3.2,3.1){$-1$}
\put(3.2,4.1){$-1$}
\put(4.2,5.1){0}
}
\put(23.0,3){
\put(-1.53,1.1){$-1$}
\put(-1.53,2.1){$-1$}
\put(0,1){\Yboxdim12pt\yng(4,2)}
\put(2.2,1.1){$-1$}
\put(4.2,2.1){$-1$}
}
\put(30.7,3){
\put(-1.53,1.1){$-1$}
\put(-0.8,2.1){0}
\put(0,1){\Yboxdim12pt\yng(3,2)}
\put(2.2,1.1){$-1$}
\put(3.2,2.1){0}
}
\end{picture}
%\{\{4,3,1\},\{5,4,1,1\},\{4,3,3,1\},\{4,2\},\{3,2\}\}\\
%\{\{0,-2,2\},\{0,2,1,0\},\{0,-1,-1,0\},\{-1,-1\},\{0,-1\}\}
\end{center}
Since all the riggings of $(\bar{\nu},\bar{J})^{(2)}$ are non-negative, we have $\widetilde{e}_2(\bar{\nu},\bar{J})=0$.
\end{Example}

\subsection{Proof for (6)}

\begin{Proposition}
Let $b\in B^{1,1}\otimes\bar{B}$ and suppose that we have shown the commutativity of $\widetilde{e}_i$ and
$\widetilde{f}_i$ with~$\Phi$ for the elements of $\bar{B}$.
Suppose that $\widetilde{e}_i(b)\neq 0$ and $\widetilde{e}_i(b')=0$ where $b'\in\bar{B}$ is the corresponding part of~$b$.
Then $\widetilde{e}_i(\nu,J)\neq 0$, $\widetilde{e}_i(\bar{\nu},\bar{J})=0$ and
$\Phi^{-1}(\widetilde{e}_i(b))=\widetilde{e}_i(\Phi^{-1}(b))$.
\end{Proposition}
\begin{proof}
By the assumption, we have the following two possibilities:
\begin{enumerate}\itemsep=0pt
\item[(a)] $b=(i+1)\otimes b'$ and $\widetilde{e}_i(b)=i\otimes b'$, \item[(b)] $b=\overline{i}\otimes b'$ and
$\widetilde{e}_i(b)=\overline{i+1}\otimes b'$.
\end{enumerate}
In both cases we have $\varphi_i(b')=0$.
Note that we have $\widetilde{e}_i(\bar{\nu},\bar{J})=0$ by the assumption.

{\bf Case (a).} By the assumption we have $\varepsilon_i(\bar{\nu},\bar{J})=0$ and
$\varphi_i(\bar{\nu},\bar{J})=0$.
Then from Theorem~\ref{prop:phi_RC} we have
$P^{(i)}_\infty(\bar{\nu})=\varphi_i(\bar{\nu},\bar{J})-\varepsilon_i(\bar{\nu},\bar{J})=0$.
Let $l:=\bar{\nu}^{(i)}_1$ be the length of the longest string of $(\bar{\nu},\bar{J})^{(i)}$.
Then we can use the convexity relation of $P^{(i)}_k(\bar{\nu})$ between $l\leq k\leq\infty$ to deduce that
$P^{(i)}_l(\bar{\nu})=P^{(i)}_{l+1}(\bar{\nu})=\dots =P^{(i)}_\infty(\bar{\nu})=0$.
Consider the string $(l,x_l)$ of $(\bar{\nu},\bar{J})^{(i)}$.
Then from $\widetilde{e}_i(\bar{\nu},\bar{J})=0$ we have $P^{(i)}_l(\bar{\nu})\geq x_l\geq 0$.
Thus we have $x_l=0$ and, in particular, the string $(l,x_l)=(l,0)$ is singular.
Since we have $b=(i+1)\otimes b'$, we see that $\delta^{-1}$ adds a~box to the string $(l,0)$ of
$(\bar{\nu},\bar{J})^{(i)}$, that is, $l+1=\ell^{(i)}$.
Since $\ell^{(i+1)}=\infty$, we see that $P^{(i)}_{\ell^{(i)}}(\nu)=P^{(i)}_{\ell^{(i)}}(\bar{\nu})-1=-1$.
To summarize, the longest string of $(\nu,J)^{(i)}$ is $(\ell^{(i)},-1)$ and therefore we have
$\widetilde{e}_i(\nu,J)\neq 0$.

Since $\delta^{-1}$ does not change the riggings of the strings except for the longest one, we see that the unique
string of $(\nu,J)^{(i)}$ which has negative rigging is $(\ell^{(i)},-1)$.
Thus $\widetilde{e}_i$ acts on the string $(\ell^{(i)},-1)$ and changes it to $(\ell^{(i)}-1,0)$.
If $\ell^{(i-1)}=\ell^{(i)}$ we see that~$\delta$ cannot remove a~box from $(\widetilde{\nu},\widetilde{J})^{(i)}$ since
all the strings are strictly shorter than $\ell^{(i-1)}$.
Thus $\Phi^{-1}(\widetilde{e}_i(b))=i\otimes b'$.
So suppose that $\ell^{(i-1)}<\ell^{(i)}$.
Then we have $P^{(i)}_{\ell^{(i)}-1}(\widetilde{\nu})=P^{(i)}_{\ell^{(i)}-1}(\nu)
=P^{(i)}_{\ell^{(i)}-1}(\bar{\nu})+1=1$.
Thus the string $(\ell^{(i)}-1,0)$ created by $\widetilde{e}_i$ is non-singular.
Therefore we have $\widetilde{\ell}^{(i)}=\infty$ and hence $\Phi^{-1}(\widetilde{e}_i(b))=i\otimes b'$.

{\bf Case (b).} Proof of this case is almost identical to the previous Case (a).
\end{proof}

\subsection{Proof for (7)}

\begin{Proposition}
The situation $\widetilde{e}_i(\nu,J)=0$ and $\widetilde{e}_i(\bar{\nu},\bar{J})\neq 0$ cannot happen.
\end{Proposition}
\begin{proof}
Suppose if possible that we have $\widetilde{e}_i(\nu,J)=0$ and $\widetilde{e}_i(\bar{\nu},\bar{J})\neq 0$.
Then~$\delta$ must change some strings of $(\nu,J)^{(i)}$.
Thus we have $\ell^{(i)}<\infty$.

{\bf Step 1.} Let us consider the case $\ell_{(i)}=\infty$.
Then the string $\big(\ell^{(i)}-1,P^{(i)}_{\ell^{(i)}-1}(\bar{\nu})\big)$ is the unique string of $(\bar{\nu},\bar{J})^{(i)}$
which has a~negative rigging.
In particular, we have $P^{(i)}_{\ell^{(i)}-1}(\bar{\nu})<0$.
Let $j$ be the largest integer such that $m^{(i)}_j(\bar{\nu})>0$ and $j<\ell^{(i)}-1$.
If there is no such string, set $j=0$.
Then we must have $P^{(i)}_j(\bar{\nu})\geq 0$.
From the convexity relation of $P^{(i)}_k(\bar{\nu})$ between $j\leq k\leq \ell^{(i)}-1$, we see that
$P^{(i)}_{\ell^{(i)}-2}(\bar{\nu})>P^{(i)}_{\ell^{(i)}-1}(\bar{\nu})$.

On the other hand, from $\widetilde{e}_i(\nu,J)=0$ the string $\big(\ell^{(i)},P^{(i)}_{\ell^{(i)}}(\nu)\big)$ must have
a~non-negative rigging so that we have $P^{(i)}_{\ell^{(i)}}(\nu)\geq 0$.

1.~Let us consider the case $\ell^{(i)}<\ell^{(i+1)}$.
Then we have $P^{(i)}_{\ell^{(i)}}(\nu)=P^{(i)}_{\ell^{(i)}}(\bar{\nu})-1$.
Therefore we have $P^{(i)}_{\ell^{(i)}}(\bar{\nu})\geq 1$.
Thus
\begin{gather*}
-P^{(i)}_{\ell^{(i)}-2}(\bar{\nu}) +2P^{(i)}_{\ell^{(i)}-1}(\bar{\nu}) -P^{(i)}_{\ell^{(i)}}(\bar{\nu})\leq -3.
\end{gather*}
On the other hand, since we cannot have more than one strings of $(\bar{\nu},\bar{J})^{(i)}$ which have negative
riggings, we have $m^{(i)}_{\ell^{(i)}-1}(\bar{\nu})=1$.
Thus
\begin{gather*}
m^{(i-1)}_{\ell^{(i)}-1}(\bar{\nu}) -2m^{(i)}_{\ell^{(i)}-1}(\bar{\nu}) +m^{(i+1)}_{\ell^{(i)}-1}(\bar{\nu})\geq -2.
\end{gather*}
This contradicts the convexity relation of Lemma~\ref{lem:convexity3}.

2.~Next, let us consider the case $\ell^{(i)}=\ell^{(i+1)}<\ell_{(i+1)}$.
Then we have $P^{(i)}_{\ell^{(i)}}(\nu)=P^{(i)}_{\ell^{(i)}}(\bar{\nu})$ so that $P^{(i)}_{\ell^{(i)}}(\bar{\nu})\geq
0$.
Thus
\begin{gather*}
-P^{(i)}_{\ell^{(i)}-2}(\bar{\nu}) +2P^{(i)}_{\ell^{(i)}-1}(\bar{\nu}) -P^{(i)}_{\ell^{(i)}}(\bar{\nu})\leq -2.
\end{gather*}
On the other hand, we have $m^{(i)}_{\ell^{(i)}-1}(\bar{\nu})=1$ and $m^{(i+1)}_{\ell^{(i)}-1}(\bar{\nu})\geq 1$.
Thus
\begin{gather*}
m^{(i-1)}_{\ell^{(i)}-1}(\bar{\nu}) -2m^{(i)}_{\ell^{(i)}-1}(\bar{\nu}) +m^{(i+1)}_{\ell^{(i)}-1}(\bar{\nu})\geq -1.
\end{gather*}
Again this contradicts the convexity relation of Lemma~\ref{lem:convexity3}.

3.~Finally let us consider the case $\ell^{(i)}=\ell^{(i+1)}=\ell_{(i+1)}$.
Then we have $P^{(i)}_{\ell^{(i)}}(\nu)=P^{(i)}_{\ell^{(i)}}(\bar{\nu})+1$ so that $P^{(i)}_{\ell^{(i)}}(\bar{\nu})\geq
-1$.
Thus
\begin{gather*}
-P^{(i)}_{\ell^{(i)}-2}(\bar{\nu}) +2P^{(i)}_{\ell^{(i)}-1}(\bar{\nu}) -P^{(i)}_{\ell^{(i)}}(\bar{\nu})\leq -1.
\end{gather*}
In this case, we have $m^{(i)}_{\ell^{(i)}-1}(\bar{\nu})=1$ and $m^{(i+1)}_{\ell^{(i)}-1}(\bar{\nu})\geq 2$.
Thus
\begin{gather*}
m^{(i-1)}_{\ell^{(i)}-1}(\bar{\nu}) -2m^{(i)}_{\ell^{(i)}-1}(\bar{\nu}) +m^{(i+1)}_{\ell^{(i)}-1}(\bar{\nu})\geq 0.
\end{gather*}
This contradicts the convexity relation of Lemma~\ref{lem:convexity3}.

{\bf Step 2.} Let us consider the case $\ell_{(i)}<\infty$.
If $\ell^{(i)}=\ell_{(i)}$ we can use the previous arguments to deduce a~contradiction.
Note that we have $\ell^{(i+1)}=\ell_{(i+1)}$ in this case.
Thus we assume that $\ell^{(i)}<\ell_{(i)}$.
Suppose that we have $\ell^{(i)}=\ell_{(i)}-1$.
Suppose that we also have $\ell_{(i)}-1>\ell_{(i+1)}-1$.
The full condition for this case is as follows
\begin{gather*}
\ell^{(i)}-1=\ell^{(i+1)}-1=\ell_{(i+1)}-1=(\ell_{(i)}-1)-1.
\end{gather*}
Thus we have $P^{(i)}_{\ell_{(i)}-1}(\nu)=P^{(i)}_{\ell_{(i)}-1}(\bar{\nu})+1$.
Since the string $\big(\ell^{(i)},P^{(i)}_{\ell^{(i)}}(\nu)\big)$ of $(\nu,J)^{(i)}$ has a~non-negative rigging, we have
$P^{(i)}_{\ell_{(i)}-1}(\nu)\geq 0$.
To summarize, we have $P^{(i)}_{\ell_{(i)}-1}(\bar{\nu})\geq -1$ in this case.

1.~Consider the case $P^{(i)}_{\ell_{(i)}-1}(\bar{\nu})\geq 0$.
Since $\widetilde{e}_i(\bar{\nu},\bar{J})\neq 0$ there is a~string of $(\bar{\nu},\bar{J})^{(i)}$ with negative rigging.
Therefore we must have $P^{(i)}_{\ell^{(i)}-1}(\bar{\nu})<0$ since $P^{(i)}_{\ell_{(i)}-1}(\bar{\nu})\geq 0$.
Let $j$ be the largest integer such that $m^{(i)}_j(\bar{\nu})>0$ and $j<\ell^{(i)}-1$.
If there is no such $j$, set $j=0$.
Then we have $P^{(i)}_j(\bar{\nu})\geq 0$.
By the convexity relation of $P^{(i)}_k(\bar{\nu})$ between $j\leq k\leq\ell^{(i)}-1$, we have
$P^{(i)}_{\ell^{(i)}-2}(\bar{\nu})>P^{(i)}_{\ell^{(i)}-1}(\bar{\nu})$.
Since $\ell^{(i)}=\ell_{(i)}-1$, we have
\begin{gather*}
-P^{(i)}_{\ell^{(i)}-2}(\bar{\nu}) +2P^{(i)}_{\ell^{(i)}-1}(\bar{\nu}) -P^{(i)}_{\ell^{(i)}}(\bar{\nu})\leq -2.
\end{gather*}
On the other hand, since $\ell^{(i)}-1=\ell^{(i+1)}-1=\ell_{(i+1)}-1$ we have $m^{(i+1)}_{\ell^{(i)}-1}(\bar{\nu})\geq
2$.
Recall that the string $\big(\ell^{(i)}-1,P^{(i)}_{\ell^{(i)}-1}(\bar{\nu})\big)$ is the unique string of
$(\bar{\nu},\bar{J})^{(i)}$ with the negative rigging.
Thus we have $m^{(i)}_{\ell^{(i)}-1}(\bar{\nu})=1$.
Therefore we have
\begin{gather*}
m^{(i-1)}_{\ell^{(i)}-1}(\bar{\nu}) -2m^{(i)}_{\ell^{(i)}-1}(\bar{\nu}) +m^{(i+1)}_{\ell^{(i)}-1}(\bar{\nu})\geq 0.
\end{gather*}
This contradicts the convexity relation of Lemma~\ref{lem:convexity3}.

2.~Consider the case $P^{(i)}_{\ell_{(i)}-1}(\bar{\nu})=-1$ and $\ell_{(i-1)}>\ell_{(i)}$.
We will show that this case cannot happen in Lemma~\ref{lem:e_vanish(7)_1}.

3.~Consider the case $P^{(i)}_{\ell_{(i)}-1}(\bar{\nu})=-1$ and $\ell_{(i-1)}=\ell_{(i)}$.
We will show that this case cannot happen in Lemma~\ref{lem:e_vanish(7)_2}.

On the other hand, together with the present assumption $\ell^{(i)}=\ell_{(i)}-1$, suppose that we also have
$\ell_{(i)}-1=\ell_{(i+1)}-1$.
We will show that this case cannot happen in Lemma~\ref{lem:e_vanish(7)_3}.

Finally let us assume that $\ell^{(i)}<\ell_{(i)}-1$.
Recall that the only strings of $(\bar{\nu},\bar{J})^{(i)}$ which may have negative riggings are
$\big(\ell^{(i)}-1,P^{(i)}_{\ell^{(i)}-1}(\bar{\nu})\big)$ and $\big(\ell_{(i)}-1,P^{(i)}_{\ell_{(i)}-1}(\bar{\nu})\big)$.
Suppose that we have $P^{(i)}_{\ell^{(i)}-1}(\bar{\nu})<0$.
In this case, as in Step 1, we have $P^{(i)}_{\ell^{(i)}-2}(\bar{\nu})>P^{(i)}_{\ell^{(i)}-1}(\bar{\nu})$ and
$P^{(i)}_{\ell^{(i)}}(\nu)\geq 0$.
Then by the same arguments of Step 1 we can deduce a~contradiction.
Thus suppose that we have $P^{(i)}_{\ell^{(i)}-1}(\bar{\nu})\geq 0$ in the following.
Then by the assumption $\widetilde{e}_i(\bar{\nu},\bar{J})\neq 0$, we have $P^{(i)}_{\ell_{(i)}-1}(\bar{\nu})<0$.
Let $j$ be the largest integer such that $j<\ell_{(i)}-1$ and $m^{(i)}_j(\bar{\nu})>0$.
Then we have $P^{(i)}_j(\bar{\nu})\geq 0$ and by the convexity relation of $P^{(i)}_k(\bar{\nu})$ between $j\leq
k\leq\ell_{(i)}-1$, we have $P^{(i)}_{\ell_{(i)}-2}(\bar{\nu})>P^{(i)}_{\ell_{(i)}-1}(\bar{\nu})$.
On the other hand, from the assumption $\widetilde{e}_i(\nu,J)=0$ we see that the string
$\big(\ell_{(i)},P^{(i)}_{\ell_{(i)}}(\nu)\big)$ of $(\nu,J)^{(i)}$ must have a~non-negative rigging.
Thus we have $P^{(i)}_{\ell_{(i)}}(\nu)\geq 0$.
Then we can use the same arguments of Step 1 if we replace $\ell^{(i)}$ and $\ell^{(i+1)}$ there by $\ell_{(i)}$ and
$\ell_{(i-1)}$ respectively to deduce a~contradiction.
Thus we have conf\/irmed that the f\/inal case $\ell^{(i)}<\ell_{(i)}-1$ cannot happen.
\end{proof}

\begin{Lemma}
\label{lem:e_vanish(7)_1}
The following situation cannot happen.
$\widetilde{e}_i(\nu,J)=0$, $\widetilde{e}_i(\bar{\nu},\bar{J})\neq 0$,
\begin{gather*}
\ell^{(i)}-1=\ell^{(i+1)}-1=\ell_{(i+1)}-1=(\ell_{(i)}-1)-1,
\end{gather*}
$\ell_{(i-1)}>\ell_{(i)}$ and $P^{(i)}_{\ell_{(i)}-1}(\bar{\nu})=-1$.
\end{Lemma}
\begin{proof}
Since $P^{(i)}_{\ell_{(i)}}(\nu)=P^{(i)}_{\ell_{(i)}}(\bar{\nu})-1$, we have $P^{(i)}_{\ell_{(i)}}(\bar{\nu})\geq 1$ by
$\widetilde{e}_i(\nu,J)=0$.
Let us show that $m^{(i)}_{\ell_{(i)}-1}(\bar{\nu})=1$.
Since we know that $\ell^{(i)}-1<\ell_{(i)}-1$, we see that $\delta^{-1}$ will add boxes to the strings
$\big(\ell^{(i)}-1,P^{(i)}_{\ell^{(i)}-1}(\bar{\nu})\big)$ and $\big(\ell_{(i)}-1,P^{(i)}_{\ell_{(i)}-1}(\bar{\nu})\big)$.
Thus if we have $m^{(i)}_{\ell_{(i)}-1}(\bar{\nu})>1$, there will be a~remaining negative rigging in $(\nu,J)^{(i)}$ by
$P^{(i)}_{\ell_{(i)}-1}(\bar{\nu})<0$, which is in contradiction to the assumption $\widetilde{e}_i(\nu,J)=0$.
Therefore we have $m^{(i)}_{\ell_{(i)}-1}(\bar{\nu})=1$.
Then $P^{(i)}_{\ell_{(i)}-2}(\bar{\nu})$, $P^{(i)}_{\ell_{(i)}-1}(\bar{\nu})$ and $(P^{(i)}_{\ell_{(i)}}(\bar{\nu})-2)$
must satisfy the convexity relation by Lemma~\ref{lem:convexity2}.
Since $P^{(i)}_{\ell_{(i)}-1}(\bar{\nu})=-1$ and $P^{(i)}_{\ell_{(i)}}(\bar{\nu})\geq 1$, we see that
\begin{gather}
\label{eq:e_vanish_0}
P^{(i)}_{\ell_{(i)}-2}(\bar{\nu})=P^{(i)}_{\ell^{(i)}-1}(\bar{\nu})\leq -1.
\end{gather}
From this, we also have $m^{(i)}_{\ell^{(i)}-1}(\bar{\nu})=1$ as otherwise we would have a~remaining negative rigging in
$(\nu,J)^{(i)}$.
Let $j$ be the largest integer such that $m^{(i)}_j(\bar{\nu})>0$ and $j<\ell^{(i)}-1$.
If there is no such~$j$, set $j=0$.
Then we have $P^{(i)}_j(\bar{\nu})\geq 0$ since the only strings of $(\bar{\nu},\bar{J})^{(i)}$ with negative riggings
are $\big(\ell^{(i)}-1,P^{(i)}_{\ell^{(i)}-1}(\bar{\nu})\big)$ and $\big(\ell_{(i)}-1,P^{(i)}_{\ell_{(i)}-1}(\bar{\nu})\big)$.
Then from the convexity relation of $P^{(i)}_k(\bar{\nu})$ between $j\leq k\leq \ell^{(i)}-1$ and~\eqref{eq:e_vanish_0},
we have
\begin{gather*}%\label{eq:e_vanish_4}
P^{(i)}_{\ell^{(i)}-2}(\bar{\nu})>P^{(i)}_{\ell^{(i)}-1}(\bar{\nu}).
\end{gather*}

Then the left hand side of the convexity relation of Lemma~\ref{lem:convexity3} is
\begin{gather*}%\label{eq:e_vanish_1}
-P^{(i)}_{\ell^{(i)}-2}(\bar{\nu})+2P^{(i)}_{\ell^{(i)}-1}(\bar{\nu})-P^{(i)}_{\ell^{(i)}}(\bar{\nu}) \leq -1
\end{gather*}
by $\ell_{(i)}-1=\ell^{(i)}$ and the right hand side of it is
\begin{gather*}%\label{eq:e_vanish_2}
m^{(i-1)}_{\ell^{(i)}-1}(\bar{\nu})-2m^{(i)}_{\ell^{(i)}-1}(\bar{\nu})+m^{(i+1)}_{\ell^{(i)}-1}(\bar{\nu}) \geq 0.
\end{gather*}
Here we use the fact $m^{(i+1)}_{\ell^{(i)}-1}(\bar{\nu})\geq 2$ which follows from
$\ell^{(i)}-1=\ell^{(i+1)}-1=\ell_{(i+1)}-1$.
This is a~contradiction.
\end{proof}

\begin{Lemma}
\label{lem:e_vanish(7)_2}
The following situation cannot happen.
$\widetilde{e}_i(\nu,J)=0$, $\widetilde{e}_i(\bar{\nu},\bar{J})\neq 0$,
\begin{gather*}
\ell^{(i)}-1=\ell^{(i+1)}-1=\ell_{(i+1)}-1=(\ell_{(i)}-1)-1=(\ell_{(i-1)}-1)-1
\end{gather*}
and $P^{(i)}_{\ell_{(i)}-1}(\bar{\nu})=-1$.
\end{Lemma}
\begin{proof}
Since $P^{(i)}_{\ell_{(i)}}(\nu)=P^{(i)}_{\ell_{(i)}}(\bar{\nu})$, we have
\begin{gather*}
P^{(i)}_{\ell_{(i)}}(\bar{\nu})\geq 0
\end{gather*}
by $\widetilde{e}_i(\nu,J)=0$.
As in Lemma~\ref{lem:e_vanish(7)_1}, we have $m^{(i)}_{\ell_{(i)}-1}(\bar{\nu})=1$ by $\ell^{(i)}-1<\ell_{(i)}-1$.
If we use $P^{(i)}_{\ell_{(i)}-1}(\bar{\nu})=-1$, the left hand side of the convexity relation of
Lemma~\ref{lem:convexity3} with $l=\ell_{(i)}-1$ is
\begin{gather*}
-P^{(i)}_{\ell_{(i)}-2}(\bar{\nu})+2P^{(i)}_{\ell_{(i)}-1}(\bar{\nu})-P^{(i)}_{\ell_{(i)}}(\bar{\nu}) \leq
-P^{(i)}_{\ell_{(i)}-2}(\bar{\nu})-2
\end{gather*}
and the right hand side of it is
\begin{gather*}
m^{(i-1)}_{\ell_{(i)}-1}(\bar{\nu})-2m^{(i)}_{\ell_{(i)}-1}(\bar{\nu})+m^{(i+1)}_{\ell_{(i)}-1}(\bar{\nu}) \geq
-1+m^{(i+1)}_{\ell_{(i)}-1}(\bar{\nu})\geq -1.
\end{gather*}
Here we have used $m^{(i-1)}_{\ell_{(i)}-1}(\bar{\nu})\geq 1$ by $\ell_{(i-1)}-1=\ell_{(i)}-1$.
Therefore we have
\begin{gather}
\label{eq:e_vanish_7}
P^{(i)}_{\ell_{(i)}-2}(\bar{\nu})=P^{(i)}_{\ell^{(i)}-1}(\bar{\nu})\leq -1.
\end{gather}
From this and $P^{(i)}_{\ell_{(i)}-1}(\bar{\nu})=-1$ we deduce that $m^{(i)}_{\ell^{(i)}-1}(\bar{\nu})=1$ as otherwise
we must have $\widetilde{e}_i(\nu,J)\neq 0$.
Let $j$ be the largest integer such that $m^{(i)}_j(\bar{\nu})>0$ and $j<\ell^{(i)}-1$.
If there is no such $j$, set $j=0$.
Then we have $P^{(i)}_j(\bar{\nu})\geq 0$ since the only strings of $(\bar{\nu},\bar{J})^{(i)}$ which may have negative
riggings are $\big(\ell^{(i)}-1,P^{(i)}_{\ell^{(i)}-1}(\bar{\nu})\big)$ and $\big(\ell_{(i)}-1,P^{(i)}_{\ell_{(i)}-1}(\bar{\nu})\big)$.
Then by the convexity relation of $P^{(i)}_k(\bar{\nu})$ between $j\leq k\leq \ell^{(i)}-1$ we have
\begin{gather*}
P^{(i)}_{\ell^{(i)}-2}(\bar{\nu})> P^{(i)}_{\ell^{(i)}-1}(\bar{\nu}).
\end{gather*}

Recall that the left hand side of the convexity relation of Lemma~\ref{lem:convexity3} with $l=\ell^{(i)}-1$ is
\begin{gather*}%\label{eq:e_vanish_5}
-P^{(i)}_{\ell^{(i)}-2}(\bar{\nu})+2P^{(i)}_{\ell^{(i)}-1}(\bar{\nu})-P^{(i)}_{\ell^{(i)}}(\bar{\nu}) \leq
P^{(i)}_{\ell^{(i)}-1}(\bar{\nu})
\end{gather*}
by $P^{(i)}_{\ell^{(i)}-2}(\bar{\nu})> P^{(i)}_{\ell^{(i)}-1}(\bar{\nu})$ and $P^{(i)}_{\ell^{(i)}}(\bar{\nu})=-1$ by
$\ell^{(i)}=\ell_{(i)}-1$.
On the other hand, the right hand side of the convexity relation of Lemma~\ref{lem:convexity3} is
\begin{gather*}%\label{eq:e_vanish_6}
m^{(i-1)}_{\ell^{(i)}-1}(\bar{\nu})-2m^{(i)}_{\ell^{(i)}-1}(\bar{\nu})+m^{(i+1)}_{\ell^{(i)}-1}(\bar{\nu})
=m^{(i-1)}_{\ell^{(i)}-1}(\bar{\nu})-2+m^{(i+1)}_{\ell^{(i)}-1}(\bar{\nu}) \geq m^{(i-1)}_{\ell^{(i)}-1}(\bar{\nu}).
\end{gather*}
Here we use the fact $m^{(i+1)}_{\ell^{(i)}-1}(\bar{\nu})\geq 2$ which follows from
$\ell^{(i)}-1=\ell^{(i+1)}-1=\ell_{(i+1)}-1$.
This is a~contradiction since we have~\eqref{eq:e_vanish_7}.
\end{proof}

\begin{Lemma}
\label{lem:e_vanish(7)_3}
The following situation cannot happen.
$\widetilde{e}_i(\nu,J)=0$, $\widetilde{e}_i(\bar{\nu},\bar{J})\neq 0$ and
\begin{gather*}
\ell^{(i)}-1=(\ell_{(i+1)}-1)-1=(\ell_{(i)}-1)-1.
\end{gather*}
\end{Lemma}
\begin{proof}\sloppy
In this case we have $P^{(i)}_{\ell_{(i)}-1}(\nu)\leq P^{(i)}_{\ell_{(i)}-1}(\bar{\nu})$.
Suppose if possible that we have $P^{(i)}_{\ell_{(i)}-1}(\bar{\nu})<0$.
Then the string $(\ell_{(i)}-1,P^{(i)}_{\ell_{(i)}-1}(\nu))$ of $(\nu,J)^{(i)}$ has a~negative rigging.
This is a~contradiction since we have $\widetilde{e}_i(\nu,J)=0$.
Thus we can assume that $P^{(i)}_{\ell_{(i)}-1}(\bar{\nu})=P^{(i)}_{\ell^{(i)}}(\bar{\nu})\geq 0$ by
$\ell^{(i)}=\ell_{(i)}-1$.
Then we see that we must have
\begin{gather*}
P^{(i)}_{\ell^{(i)}-1}(\bar{\nu})<0
\end{gather*}
since we have $\widetilde{e}_i(\bar{\nu},\bar{J})\neq 0$.

Let $j$ be the largest integer such that $j<\ell^{(i)}-1$ and $m^{(i)}_j(\bar{\nu})>0$.
If there is no such string, set $j=0$.
Then we have $P^{(i)}_j(\bar{\nu})\geq 0$.
From the convexity relation of $P^{(i)}_k(\bar{\nu})$ between $j\leq k\leq\ell^{(i)}-1$, we have
\begin{gather*}
P^{(i)}_{\ell^{(i)}-2}(\bar{\nu})>P^{(i)}_{\ell^{(i)}-1}(\bar{\nu})
\end{gather*}
since $P^{(i)}_{\ell^{(i)}-1}(\bar{\nu})<0$.

Let us show that $m^{(i)}_{\ell^{(i)}-1}(\bar{\nu})=1$.
Suppose if possible that we have $m^{(i)}_{\ell^{(i)}-1}(\bar{\nu})>1$.
Since we have $P^{(i)}_{\ell^{(i)}-1}(\bar{\nu})<0$ and $\ell^{(i)}-1<\ell_{(i)}-1$, we see that there is the remaining
negative rigging in $(\nu,J)^{(i)}$.
This is a~contradiction since we have $\widetilde{e}_i(\nu,J)=0$.
Thus we have $m^{(i)}_{\ell^{(i)}-1}(\bar{\nu})=1$.

{\bf Case 1.} Suppose that we have $\ell^{(i)}-1=(\ell^{(i+1)}-1)-1$.
Then we have $P^{(i)}_{\ell^{(i)}}(\nu)=P^{(i)}_{\ell^{(i)}}(\bar{\nu})-1$.
From $\widetilde{e}_i(\nu,J)=0$ we have $P^{(i)}_{\ell^{(i)}}(\bar{\nu})\geq 1$.
To summarize, we have
\begin{gather*}
P^{(i)}_{\ell^{(i)}-2}(\bar{\nu})>P^{(i)}_{\ell^{(i)}-1}(\bar{\nu})<0,
\qquad
P^{(i)}_{\ell^{(i)}}(\bar{\nu})\geq 1.
\end{gather*}
Let $n^{(a)}_l$ be the length of the $l$th column of $\bar{\nu}^{(a)}$
\begin{gather*}
n^{(a)}_l:=
\sum\limits_{l\leq k}
m^{(a)}_k(\bar{\nu}).
\end{gather*}
Then the above relations imply that
\begin{gather}
\label{eq:e_vanish_12}
n^{(i-1)}_{\ell^{(i)}-1}-2n^{(i)}_{\ell^{(i)}-1}+n^{(i+1)}_{\ell^{(i)}-1}\leq -1,
\qquad
n^{(i-1)}_{\ell^{(i)}}-2n^{(i)}_{\ell^{(i)}}+n^{(i+1)}_{\ell^{(i)}}\geq 2.
\end{gather}
Since we have $m^{(i)}_{\ell^{(i)}-1}(\bar{\nu})=1$, we have $n^{(i)}_{\ell^{(i)}}=n^{(i)}_{\ell^{(i)}-1}-1$ here.
Then the second relation of~\eqref{eq:e_vanish_12} becomes
$n^{(i-1)}_{\ell^{(i)}}-2n^{(i)}_{\ell^{(i)}-1}+n^{(i+1)}_{\ell^{(i)}}\geq 0$.
Combining this with the f\/irst relation of~\eqref{eq:e_vanish_12}, we obtain
\begin{gather*}
0>\left(n^{(i-1)}_{\ell^{(i)}-1}-n^{(i-1)}_{\ell^{(i)}}\right)+
\left(n^{(i+1)}_{\ell^{(i)}-1}-n^{(i+1)}_{\ell^{(i)}}\right).
\end{gather*}
This is a~contradiction since we have $n^{(a)}_l\geq n^{(a)}_{l+1}$ by def\/inition of $n^{(a)}_l$.

{\bf Case 2.} Suppose that we have $\ell^{(i)}-1=\ell^{(i+1)}-1$.
In particular, we have $m^{(i+1)}_{\ell^{(i)}-1}(\bar{\nu})>0$.
As in the previous case, we have
\begin{gather*}
P^{(i)}_{\ell^{(i)}-2}(\bar{\nu})>P^{(i)}_{\ell^{(i)}-1}(\bar{\nu})<0,
\qquad
P^{(i)}_{\ell^{(i)}}(\bar{\nu})\geq 0.
\end{gather*}
Therefore we have
\begin{gather}
\label{eq:e_vanish_13}
n^{(i-1)}_{\ell^{(i)}-1}-2n^{(i)}_{\ell^{(i)}-1}+n^{(i+1)}_{\ell^{(i)}-1}\leq -1,
\qquad
n^{(i-1)}_{\ell^{(i)}}-2n^{(i)}_{\ell^{(i)}}+n^{(i+1)}_{\ell^{(i)}}\geq 1.
\end{gather}
From $n^{(i)}_{\ell^{(i)}}=n^{(i)}_{\ell^{(i)}-1}-1$, the second relation of~\eqref{eq:e_vanish_13} gives
$n^{(i-1)}_{\ell^{(i)}}-2n^{(i)}_{\ell^{(i)}-1}+n^{(i+1)}_{\ell^{(i)}}\geq -1$.
Then from the f\/irst relation of~\eqref{eq:e_vanish_13}, we obtain
\begin{gather*}
0\geq\left(n^{(i-1)}_{\ell^{(i)}-1}-n^{(i-1)}_{\ell^{(i)}}\right)+
\left(n^{(i+1)}_{\ell^{(i)}-1}-n^{(i+1)}_{\ell^{(i)}}\right).
\end{gather*}
Therefore we have $n^{(i-1)}_{\ell^{(i)}-1}=n^{(i-1)}_{\ell^{(i)}}$ and
$n^{(i+1)}_{\ell^{(i)}-1}=n^{(i+1)}_{\ell^{(i)}}$.
In particular, we have $m^{(i+1)}_{\ell^{(i)}-1}(\bar{\nu})=0$, which is a~contradiction.
Hence this case cannot happen.
\end{proof}

\subsection{Proof for (8)}

\begin{Proposition}
Let $b=b'\otimes\bar{b}\in B^{1,1}\otimes\bar{B}$ and $i\in I_0$.
Then the situation $\widetilde{e}_i(b)=0$ and $\widetilde{e}_i(\bar{b})\neq 0$ cannot happen.
\end{Proposition}
\begin{proof}
The proof of this assertion follows from the tensor product rule at Section~\ref{se:crystal}.
\end{proof}

\section{Comments on Deka--Schilling's work}\label{sec:DekaSchilling}

In Appendix C of the original preprint version of~\cite{DS:2006}, Deka and Schilling proved the compatibility of the
classical Kashiwara operators and the rigged conf\/iguration bijection of type $A^{(1)}_n$.
I am pleased to admit that they proved some portion of the entire proof and also I would like to sincerely admire that
they embarked on such a~dif\/f\/icult problem.
However, due to the following crucial problems, the author assumes that they did not f\/inish their proof.
\begin{enumerate}\itemsep=0pt
\item[(a)] In~\cite{DS:2006}, the proofs for $\widetilde{e}_i$ case are completely omitted since they claim at Lemmas~C.2 and~C.3 (i.e., not in the main proof) that the proofs are ``similar'' to $\widetilde{f}_i$ case.
Actually, as we have shown in Sections~\ref{sec:main2} and~\ref{sec:main2_e} of the present paper, the logical structure
of two proofs are entirely dif\/ferent.
In particular, the proofs for $\widetilde{e}_i$ case require the result of $\widetilde{f}_i$ case which was proved
beforehand.
It is natural as the statements of the necessary propositions are already dif\/ferent.
\item[(b)] In~\cite{DS:2006}, the authors give entirely no words about the riggings.
As we have shown in the present paper, the coincidence of the riggings is non-trivial in many cases since the Kashiwara
operators on the rigged conf\/igurations modify riggings.
Thus they strongly deserve individual and careful analysis.
\item[(c)] In~\cite{DS:2006}, the authors give no words to Cases A and E of the present paper.
They simply claim that the only non-trivial cases are the Cases B, C and D.~However, as we have shown in the present
paper, each proof of Case A or E breaks into subcases and requires somewhat elaborated arguments even if we restrict to
the type $A^{(1)}_n$.
Thus it is necessary to mention about these cases.
\end{enumerate}

\subsection*{Acknowledgements}

The author would like to thank Professor Masato Okado for valuable discussion throughout the project.
In particular, the project began when both of us worked together to understand the details of Appendix C of the preprint
version of~\cite{DS:2006}.
The author is also grateful to Professor Anne Schilling for valuable discussion on her results.
For both of them, the author is very grateful for the fruitful collaboration on closely related
projects~\cite{OSS:2012,OSSS:2012} which provides a~motivation and an important application of the present work.
I would like to thank anonymous referees for valuable suggestions which greatly help the author to improve the original
manuscript.
This work is partially supported by Grants-in-Aid for Scientif\/ic Research No.~21740114 from JSPS.

%\pdfbookmark[1]{References}{ref}
\LastPageEnding


\addcontentsline{toc}{section}{References}
\begin{thebibliography}{99}
\footnotesize\itemsep=0pt

\bibitem{AFLT}
Alba V.A., Fateev V.A., Litvinov A.V., Tarnopolskiy G.M., On combinatorial
  expansion of the conformal blocks arising from {AGT} conjecture,
  \href{http://dx.doi.org/10.1007/s11005-011-0503-z}{\textit{Lett. Math. Phys.}} \textbf{98} (2011), 33--64, \href{http://arxiv.org/abs/1012.1312}{arXiv:1012.1312}.

\bibitem{BB}
Belavin A., Belavin V., A{GT} conjecture and integrable structure of conformal
  f\/ield theory for {$c=1$}, \href{http://dx.doi.org/10.1016/j.nuclphysb.2011.04.014}{\textit{Nuclear Phys.~B}} \textbf{850} (2011),
  199--213, \href{http://arxiv.org/abs/1102.0343}{arXiv:1102.0343}.

\bibitem{CJ}
Cai W., Jing N., Applications of a {L}aplace--{B}eltrami operator for {J}ack
  polynomials, \href{http://dx.doi.org/10.1016/j.ejc.2011.11.003}{\textit{European~J. Combin.}} \textbf{33} (2012), 556--571,
  \href{http://arxiv.org/abs/1101.5544}{arXiv:1101.5544}.

\bibitem{DS:2006}
Deka L., Schilling A., New fermionic formula for unrestricted {K}ostka
  polynomials, \href{http://dx.doi.org/10.1016/j.jcta.2006.01.003}{\textit{J.~Combin. Theory Ser.~A}} \textbf{113} (2006),
  1435--1461, \href{http://arxiv.org/abs/math.CO/0509194}{math.CO/0509194}.

\bibitem{DLM}
Desrosiers P., Lapointe L., Mathieu P., Superconformal f\/ield theory and {J}ack
  superpolynomials, \href{http://dx.doi.org/10.1007/JHEP09(2012)037}{\textit{J.~High Energy Phys.}} \textbf{2012} (2012), no.~9,
  037, 42~pages, \href{http://arxiv.org/abs/1205.0784}{arXiv:1205.0784}.

\bibitem{EPSS}
Estienne B., Pasquier V., Santachiara R., Serban D., Conformal blocks in
  {V}irasoro and {W} theories: duality and the {C}alogero--{S}utherland model,
  \href{http://dx.doi.org/10.1016/j.nuclphysb.2012.03.007}{\textit{Nuclear Phys.~B}} \textbf{860} (2012), 377--420, \href{http://arxiv.org/abs/1110.1101}{arXiv:1110.1101}.

\bibitem{FeFu}
Feigin B.L., Fuchs D.B., Representations of the {V}irasoro algebra, in
  Representation of {L}ie Groups and Related Topics, \textit{Adv. Stud.
  Contemp. Math.}, Vol.~7, Gordon and Breach, New York, 1990, 465--554.

\bibitem{FHV}
Feigin M.V., Halln{\"a}s M.A., Veselov A.P., Baker--{A}khiezer functions and
  generalised {M}acdonald--{M}ehta integrals, \href{http://dx.doi.org/10.1063/1.4804615}{\textit{J.~Math. Phys.}}
  \textbf{54} (2013), 052106, 22~pages, \href{http://arxiv.org/abs/1210.5270}{arXiv:1210.5270}.

\bibitem{FOY}
Fukuda K., Yamada Y., Okado M., Energy functions in box ball systems,
  \href{http://dx.doi.org/10.1142/S0217751X00000616}{\textit{Internat.~J. Modern Phys.~A}} \textbf{15} (2000), 1379--1392,
  \href{http://arxiv.org/abs/math.QA/9908116}{math.QA/9908116}.

\bibitem{HHIKTT}
Hatayama G., Hikami K., Inoue R., Kuniba A., Takagi T., Tokihiro T., The
  {$A^{(1)}_M$} automata related to crystals of symmetric tensors,
  \href{http://dx.doi.org/10.1063/1.1322077}{\textit{J.~Math. Phys.}} \textbf{42} (2001), 274--308,
  \href{http://arxiv.org/abs/math.QA/9912209}{math.QA/9912209}.

\bibitem{HKT1}
Hatayama G., Kuniba A., Takagi T., Soliton cellular automata associated with
  crystal bases, \href{http://dx.doi.org/10.1016/S0550-3213(00)00105-X}{\textit{Nuclear Phys.~B}} \textbf{577} (2000), 619--645,
  \href{http://arxiv.org/abs/solv-int/9907020}{solv-int/9907020}.

\bibitem{HKT2}
Hatayama G., Kuniba A., Takagi T., Simple algorithm for factorized dynamics of
  the {${\mathfrak g}_n$}-automaton, \href{http://dx.doi.org/10.1088/0305-4470/34/48/331}{\textit{J.~Phys.~A: Math. Gen.}}
  \textbf{34} (2001), 10697--10705, \href{http://arxiv.org/abs/nlin.CG/0103022}{nlin.CG/0103022}.

\bibitem{Kac}
Kac V.G., Inf\/inite-dimensional {L}ie algebras, 3rd ed., \href{http://dx.doi.org/10.1017/CBO9780511626234}{Cambridge University
  Press}, Cambridge, 1990.

\bibitem{Kashiwara:1991}
Kashiwara M., On crystal bases of the {$Q$}-analogue of universal enveloping
  algebras, \href{http://dx.doi.org/10.1215/S0012-7094-91-06321-0}{\textit{Duke Math.~J.}} \textbf{63} (1991), 465--516.

\bibitem{Kashiwara:book}
Kashiwara M., Bases cristallines des groupes quantiques, \textit{Cours
  Sp\'ecialis\'es}, Vol.~9, Soci\'et\'e Math\'ematique de France, Paris, 2002.

\bibitem{KN:1994}
Kashiwara M., Nakashima T., Crystal graphs for representations of the
  {$q$}-analogue of classical {L}ie algebras, \href{http://dx.doi.org/10.1006/jabr.1994.1114}{\textit{J.~Algebra}} \textbf{165}
  (1994), 295--345.

\bibitem{KKR}
Kerov S.V., Kirillov A.N., Reshetikhin N.Y., Combinatorics, the {B}ethe ansatz
  and representations of the symmetric group, \href{http://dx.doi.org/10.1007/BF01247087}{\textit{J.~Soviet Math.}}
  \textbf{41} (1988), 916--924.

\bibitem{KR}
Kirillov A.N., Reshetikhin N.Y., The {B}ethe ansatz and the combinatorics of
  {Y}oung tableaux, \href{http://dx.doi.org/10.1007/BF01247088}{\textit{J.~Soviet Math.}} \textbf{41} (1988), 925--955.

\bibitem{KSS:2002}
Kirillov A.N., Schilling A., Shimozono M., A bijection between
  {L}ittlewood--{R}ichardson tableaux and rigged conf\/igurations,
  \href{http://dx.doi.org/10.1007/s00029-002-8102-6}{\textit{Selecta Math.~(N.S.)}} \textbf{8} (2002), 67--135,
  \href{http://arxiv.org/abs/math.CO/9901037}{math.CO/9901037}.

\bibitem{KOSTY}
Kuniba A., Okado M., Sakamoto R., Takagi T., Yamada Y., Crystal interpretation
  of {K}erov--{K}irillov--{R}eshetikhin bijection, \href{http://dx.doi.org/10.1016/j.nuclphysb.2006.02.005}{\textit{Nuclear Phys.~B}}
  \textbf{740} (2006), 299--327, \href{http://arxiv.org/abs/math.QA/0601630}{math.QA/0601630}.

\bibitem{KSY2}
Kuniba A., Sakamoto R., Yamada Y., Generalized energies and integrable
  {$D_n^{(1)}$} cellular automaton, in \href{http://dx.doi.org/10.1142/9789814324373_0012}{New Trends in Quantum Integrable
  Systems}, World Sci. Publ., Hackensack, NJ, 2011, 221--242,
  \href{http://arxiv.org/abs/1001.1813}{arXiv:1001.1813}.

\bibitem{MY}
Mimachi K., Yamada Y., Singular vectors of the {V}irasoro algebra in terms of
  {J}ack symmetric polynomials, \href{http://dx.doi.org/10.1007/BF02099610}{\textit{Comm. Math. Phys.}} \textbf{174} (1995),
  447--455.

\bibitem{OS:2011}
Okado M., Sakamoto R., Stable rigged conf\/igurations for quantum af\/f\/ine algebras
  of nonexceptional types, \href{http://dx.doi.org/10.1016/j.aim.2011.06.012}{\textit{Adv. Math.}} \textbf{228} (2011), 1262--1293,
  \href{http://arxiv.org/abs/1008.0460}{arXiv:1008.0460}.


\bibitem{OSS:2012}
Okado M., Sakamoto R., Schilling A., Af\/f\/ine crystal structure on rigged
  conf\/igurations of type {$D_n^{(1)}$}, \href{http://dx.doi.org/10.1007/s10801-012-0383-z}{\textit{J.~Algebraic Combin.}}
  \textbf{37} (2013), 571--599, \href{http://arxiv.org/abs/1109.3523}{arXiv:1109.3523}.

\bibitem{OSSS:2012}
Okado M., Sakamoto R., Schilling A., Kerov--Kirillov--Reshetikhin type
  bijection for $D^{(1)}_n$, {i}n preparation.

\bibitem{OSano}
Okado M., Sano N., K{KR} type bijection for the exceptional af\/f\/ine algebra
  {$E_6^{(1)}$}, in Algebraic Groups and Quantum Groups, \href{http://dx.doi.org/10.1090/conm/565/11181}{\textit{Contemp.
  Math.}}, Vol.~565, Amer. Math. Soc., Providence, RI, 2012, 227--242,
  \href{http://arxiv.org/abs/1105.1636}{arXiv:1105.1636}.

\bibitem{OSS:2003}
Okado M., Schilling A., Shimozono M., A crystal to rigged conf\/iguration
  bijection for nonexceptional af\/f\/ine algebras, in \href{http://dx.doi.org/10.1142/9789812775405_0005}{Algebraic Combinatorics and
  Quantum Groups}, World Sci. Publ., River Edge, NJ, 2003, 85--124,
  \href{http://arxiv.org/abs/math.QA/0203163}{math.QA/0203163}.

\bibitem{Sak2}
Sakamoto R., Kirillov--{S}chilling--{S}himozono bijection as energy functions
  of crystals, \href{http://dx.doi.org/10.1093/imrn/rnn140}{\textit{Int. Math. Res. Not.}} \textbf{2009} (2009), 579--614,
  \href{http://arxiv.org/abs/0711.4185}{arXiv:0711.4185}.

\bibitem{S:review}
Sakamoto R., Ultradiscrete soliton systems and combinatorial representation
  theory, \href{http://arxiv.org/abs/1212.2774}{arXiv:1212.2774}.

\bibitem{SSAFR}
Sakamoto R., Shiraishi J., Arnaudon D., Frappat L., Ragoucy E., Correspondence
  between conformal f\/ield theory and {C}alogero--{S}utherland model,
  \href{http://dx.doi.org/10.1016/j.nuclphysb.2004.10.005}{\textit{Nuclear Phys.~B}} \textbf{704} (2005), 490--509,
  \href{http://arxiv.org/abs/hep-th/0407267}{hep-th/0407267}.

\bibitem{S:2005}
Schilling A., A bijection between type {$D^{(1)}_n$} crystals and rigged
  conf\/igurations, \href{http://dx.doi.org/10.1016/j.jalgebra.2004.12.010}{\textit{J.~Algebra}} \textbf{285} (2005), 292--334,
  \href{http://arxiv.org/abs/math.QA/0406248}{math.QA/0406248}.

\bibitem{Sch:2006}
Schilling A., Crystal structure on rigged conf\/igurations, \href{http://dx.doi.org/10.1155/IMRN/2006/97376}{\textit{Int. Math.
  Res. Not.}} \textbf{2006} (2006), Art.~ID 97376, 27~pages,
  \href{http://arxiv.org/abs/math.QA/0508107}{math.QA/0508107}.

\bibitem{S:2008}
Schilling A., Combinatorial structure of {K}irillov--{R}eshetikhin crystals of
  type $D^{(1)}_n$, $B^{(1)}_n$, $A^{(2)}_{2n-1}$, \href{http://dx.doi.org/10.1016/j.jalgebra.2007.10.020}{\textit{J.~Algebra}}
  \textbf{319} (2008), 2938--2962, \href{http://arxiv.org/abs/0704.2046}{arXiv:0704.2046}.

\bibitem{SS}
Schilling A., Shimozono M., {$X=M$} for symmetric powers, \href{http://dx.doi.org/10.1016/j.jalgebra.2005.04.023}{\textit{J.~Algebra}}
  \textbf{295} (2006), 562--610, \href{http://arxiv.org/abs/math.QA/0412376}{math.QA/0412376}.

\bibitem{Sutherland1}
Sutherland B., Exact results for a quantum many-body problem in one dimension,
  \href{http://dx.doi.org/10.1103/PhysRevA.4.2019}{\textit{Phys. Rev.~A}} \textbf{4} (1971), 2019--2021.

\bibitem{Sutherland2}
Sutherland B., Exact results for a quantum many-body problem in one
  dimension.~II, \href{http://dx.doi.org/10.1103/PhysRevA.5.1372}{\textit{Phys. Rev.~A}} \textbf{5} (1972), 1372--1376.

\bibitem{TS}
Takahashi D., Satsuma J., A soliton cellular automaton, \href{http://dx.doi.org/10.1143/JPSJ.59.3514}{\textit{J.~Phys. Soc.
  Japan}} \textbf{59} (1990), 3514--3519.

\end{thebibliography}
\end{document}